\documentclass[reqno]{amsart}
\usepackage{xcolor}
\usepackage{dirtytalk}
\usepackage{amsmath}
\usepackage{amssymb}
\usepackage{amsthm}
\usepackage{mathtools}

\renewcommand{\emph}[1]{\texttt{#1}}
\usepackage[
backend=bibtex,
style=numeric,
giveninits=true,
maxbibnames=99,
autolang=hyphen,
doi=false,isbn=false,url=false,
sortcites
]{biblatex}
\renewbibmacro{in:}{}
\usepackage[colorlinks=true,hypertexnames=false]{hyperref}
\usepackage[stretch=15]{microtype} 
\usepackage{todonotes}
\usepackage{bm}
\usepackage{amsmath}
\usepackage{amssymb}
\usepackage{amsthm}
\usepackage{mathtools}
\usepackage{enumitem}

\usepackage[nameinlink]{cleveref}
\usepackage{thmtools} 

\usepackage[bb=dsserif]{mathalpha}

\usepackage[OT1]{eulervm}

\usepackage[only,llbracket,rrbracket]{stmaryrd}

\numberwithin{equation}{section}
\theoremstyle{plain}
\newtheorem{theorem}{Theorem}
\numberwithin{theorem}{section}

\newtheorem{lemma}[theorem]{Lemma}
\newtheorem{example}[theorem]{Example}
\newtheorem{proposition}[theorem]{Proposition}
\newtheorem{corollary}[theorem]{Corollary}
\newtheorem{remark}[theorem]{Remark}
\newtheorem*{remark*}{Remark}

\newtheorem{lemmaInThm}{Sublemma}
\Crefname{lemmaInThm}{Sublemma}{Sublemmas}
\newtheorem{assumptions}{Assumptions}
\Crefname{assumptions}{Assumptions}{Assumptions}

\theoremstyle{definition}
\newtheorem{definition}[theorem]{Definition}
\newtheorem{notation}[theorem]{Notation}

\RenewCommandCopy{\theHtheorem}{\thetheorem}
\RenewCommandCopy{\theHlemma}{\thelemma}
\RenewCommandCopy{\theHexample}{\theexample}
\RenewCommandCopy{\theHproposition}{\theproposition}
\RenewCommandCopy{\theHcorollary}{\thecorollary}
\RenewCommandCopy{\theHremark}{\theremark}
\RenewCommandCopy{\theHdefinition}{\thedefinition}
\RenewCommandCopy{\theHnotation}{\thenotation}

\newlist{enumerateInLemma}{enumerate}{1}
\setlist[enumerateInLemma]{label=(\alph*),ref=(\alph*)}
\crefname{enumerateInLemmai}{part}{parts}
\Crefname{enumerateInLemmai}{Lemma \thelemma\spaceskip1sp}{Lemmas \thelemma\spaceskip1sp}

\newlist{enumerateInExample}{enumerate}{1}
\setlist[enumerateInExample]{label=(\roman*)}
\crefname{enumerateInExamplei}{part}{parts}
\Crefname{enumerateInExamplei}{}{}

\newcommand{\eps}{\varepsilon}

\renewcommand{\phi}{\varphi}
\renewcommand{\theta}{\vartheta} 

\newcommand{\bR}{\mathbb{R}}
\newcommand{\bRExt}{\widehat{\mathbb{R}}}

\newcommand{\bN}{\mathbb{N}}

\newcommand{\bD}{\mathbb{D}}

\newcommand{\bS}{\mathbb{S}}

\newcommand{\mF}{\mathfrak{F}}
\newcommand{\mf}[1]{\mathfrak{#1}}
\newcommand{\mc}[1]{\mathcal{#1}}

\newcommand{\ms}[1]{\mathsf{#1}}

\newcommand{\bx}{\mathbb{x}}
\newcommand{\by}{\mathbb{y}}
\newcommand{\bz}{\mathbb{z}}
\newcommand{\ba}{\mathbb{a}}
\newcommand{\bb}{\mathbb{b}}
\newcommand{\bc}{\mathbb{c}}
\newcommand{\bd}{\mathbb{d}}
\newcommand{\be}{\mathbb{e}}
\newcommand{\bv}{\mathbb{v}}

\newcommand{\cP}{\mathcal{P}}

\newcommand{\cR}{\mathcal{R}}

\newcommand{\partitions}[2][]{\mathfrak{P}^{#1}_{#2}}

\newcommand{\lebesguePoints}[1]{\mathrm{Leb} \del*{#1}}
\newcommand{\ball}[2][]{\mathsf{B}_{#1} \del*{#2}}

\newcommand{\cl}[1]{\overline{#1}}
\newcommand{\clball}[2][]{\cl{\mathsf{B}}_{#1}\del*{#2}}

\newcommand{\sobolev}[1]{W^{ #1 }}

\newcommand{\locSobolev}[1]{W_{\mathrm{loc}}^{ #1 }}

\newcommand{\integrable}[1]{L^{ #1 }}

\newcommand{\locIntegrable}[1]{L_{\mathrm{loc}}^{ #1 }}

\newcommand{\continuousAndBounded}[1][]{C_{\mathrm{b}}^{#1}}
\newcommand{\smoothCompactSupport}{C_{\mathrm{c}}^{ \infty }}

\newcommand{\holder}[1]{C^{#1}}
\newcommand{\holderBounded}[1]{C_{\mathrm{b}}^{#1}}

\newcommand{\hajlasz}[1]{M^{ #1 }}
\newcommand{\locHajlasz}[1]{M_{\mathrm{loc}}^{ #1 }}

\newcommand{\multipoint}[2][]{\mc{M}_{#1}^{#2}}

\DeclarePairedDelimiterX\setc[2]{\{}{\}}{\,#1 \;\delimsize\vert\; #2\,}
\DeclarePairedDelimiterX\genc[2]{\langle}{\rangle}{#1 \delimsize\vert #2}

\DeclareMathOperator*{\diff}{ \text{$\Delta$} }

\newcommand{\poly}[2][]{\mathsf{P}_{#1}^{#2} }
\newcommand{\polyset}[2][]{\poly[#1]{#2}}
\newcommand{\polyfam}[2][]{\poly[#1]{#2}}
\newcommand{\polynum}[2][]{\poly[#1]{#2}}
\newcommand{\polygen}[2][]{\poly[#1]{#2}}

\renewcommand{\t}{\cdot}
\newcommand{\blank}[1][\cdot]{{\mspace{1mu}#1\mspace{1mu}}}

\newcommand{\wave}[1]{\widetilde{ #1}}

\newcommand{\indicator}[1]{\mathbb{1}_{#1}}

\newcommand{\restMaxFun}[3][]{\mathsf{M}^{#1}_{#2} \del*{#3}}
\newcommand{\maxFun}[2][]{\mathsf{M}^{#1} \del*{#2}}
\newcommand{\ballAv}[2]{\mathsf{B}_{#1} \del*{ #2 } }

\newcommand{\integral}[4][]{  \int^{#1}_{ #2} #3 \, \mathrm{d}#4}

 \usepackage{ esint }
\newcommand{\sintegral}[4][]{\fint^{#1}_{#2} #3 \, \mathrm{d}#4}

\DeclareMathOperator{\dist}{\mathrm{dist}}
\DeclareMathOperator{\sign}{\mathrm{sgn}}
\DeclareMathOperator{\supp}{\mathrm{supp}}

\newcommand{\metric}[2][]{\mathsf{d}_{#1}\del*{ #2 }} 
\newcommand{\metricAlone}[1][]{\mathsf{d}_{#1}}

\newcommand{\normAlone}[1][]{\| \, \blank \, \|_{#1}}
\newcommand{\absAlone}[1][]{| \, \blank \, |}
\newcommand{\seminormAlone}[1][]{\seminorm{ \, \blank \,  }_{#1}}

\newcommand{\measure}[2][]{\mu_{#1}\del*{ #2 }}
\newcommand{\measureAlone}[1][]{\mu_{#1}}

\let\oldth\th
\renewcommand{\th}{\mathrm{th}}

\DeclareFontFamily{U}{MnSymbolA}{}
\DeclareFontShape{U}{MnSymbolA}{m}{n}{
  <-6> MnSymbolA5
  <6-7> MnSymbolA6
  <7-8> MnSymbolA7
  <8-9> MnSymbolA8
  <9-10> MnSymbolA9
  <10-12> MnSymbolA10
  <12-> MnSymbolA12}{}
\DeclareFontShape{U}{MnSymbolA}{b}{n}{
  <-6> MnSymbolA-Bold5
  <6-7> MnSymbolA-Bold6
  <7-8> MnSymbolA-Bold7
  <8-9> MnSymbolA-Bold8
  <9-10> MnSymbolA-Bold9
  <10-12> MnSymbolA-Bold10
  <12-> MnSymbolA-Bold12}{}

\DeclareSymbolFont{MnSyA}{U}{MnSymbolA}{m}{n}

\DeclareMathSymbol{\upmapsto}{\mathrel}{MnSyA}{41}

\DeclarePairedDelimiter{\abs}{|}{|}
\DeclarePairedDelimiter{\del}{(}{)}
\DeclarePairedDelimiter{\norm}{\|}{\|}
\DeclarePairedDelimiter{\set}{\{}{\}}
\DeclarePairedDelimiter{\sbr}{[}{]}
\DeclarePairedDelimiter{\card}{|}{|}
\DeclarePairedDelimiter{\intoc}{(}{]}
\DeclarePairedDelimiter{\intcc}{[}{]}
\DeclarePairedDelimiter{\intco}{[}{)}
\DeclarePairedDelimiter{\intoo}{(}{)}
\DeclarePairedDelimiter{\inner}{\langle}{\rangle}
\DeclarePairedDelimiter{\gen}{\langle}{\rangle}

\DeclarePairedDelimiter{\seminorm}{[}{]}

\DeclarePairedDelimiter{\dsbr}{\llbracket}{\rrbracket}

\usepackage[margin=2.5cm, top=2.5cm, bottom=2.5cm]{geometry}

\title[Multipoint Characterization of Higher-Order Sobolev Spaces]{Multipoint Characterization\\of Higher-Order Sobolev Spaces}
\author{Kacper Kazimierz Kurowski}
\thanks{e-mail: kacper.kurowski.dokt@pw.edu.pl}
\date{\today}


\begin{document}

\begin{abstract}
	In the paper, we prove a rather general characterization of higher-order Sobolev spaces.
	We show that the $k\th$-order regularity, where $k \in \bN$, is captured via inequalities involving $2^k$-tuples of points. 
	In~fact, in full generality, the obtained results characterize higher-order Sobolev spaces based on Banach function spaces.
	Moreover, we show an analogous characterization of higher-order H\"{o}lder spaces.
	Finally, we propose a way to use the obtained results to define higher-order Sobolev and H\"{o}lder spaces on metric measure spaces.
\end{abstract}

\maketitle

\section{Introduction}
\label{sec::introduction}
In his paper from 1996, Hajlasz \cite{Hajlasz} showed that,
if $p \in \intoc{1,\infty}$
and 
$ \Omega \subseteq \bR^n$ is open and such that there exists a bounded extension operator
$
	\mathsf{E} \colon \sobolev{1,p}\del*{ \Omega } \to \sobolev{1,p}\del*{ \bR^n },
$
then the following conditions are equivalent:
\begin{enumerate}[label=(\alph*)]
\item
	$ f \in \sobolev{1,p}\del*{ \Omega }$;
\item 
	$ f \in \integrable{p}\del*{ \Omega }$ 
	and there exists 
	$ g \in \integrable{p}\del*{ \Omega }$ 
	that is a Hajlasz gradient of $f$, i.e., a nonnegative
	measurable function for which there exists a set 
	$ F \subseteq \Omega $
	of full measure such that
	\begin{equation*}
		\forall x, y \in F
	\qquad 
		\abs*{
			f(x) - f(y)
		}
		\le 
		\norm*{ x- y }
		\del*{
			g(x) + g(y)
		}.
	\end{equation*}
\end{enumerate}
Using this result, he then introduced first-order Sobolev spaces on metric measure spaces.
The theory of~thusly introduced Hajłasz--Sobolev
spaces
$ \hajlasz{1,p} $
is already rather abundant
\cite{Hajlasz, sobolev_met_poincare, holder_quasicontinuity, kalamajska_zwarte, equivalence_density_of_measure, dachun_fractional, gorka_looking_for_compactness, joint_work}.

Within the paper, we will generalize the 
Hajłasz's result to all orders $k \in \bN$
by proving the following theorem.
\begin{theorem}
\label{thm::sobolev_characterization}
	Let $n, k \in\bN$, 
	$p \in \intoc{1, \infty}$,
	and $\Omega \subseteq \bR^n$
	be open and such that there exists a bounded extension operator
	$ 
		\mathsf{E} 
		\colon 
		\sobolev{k,p}(\Omega)
		\to
		\sobolev{k,p}(\bR^n).
	$
	Then the following conditions are equivalent:
	\begin{enumerate}[label=(\alph*)]
	\item 
	\label{thmpart::thm::sobolev_characterization::sobolev_condition}
		$f \in \sobolev{k,p}(\Omega)$;
	\item 
	\label{thmpart::thm::sobolev_characterization::hypergradient_condition}
		$f \in \integrable{p}(\Omega)$
		and there exists 
		$ G \in \integrable{p}(\Omega)$
		such that
		$ G \in \bD^{k}_{\lambda}(f)$,
		i.e., $G$ is a nonnegative measurable function for which there exists a set $F \subseteq \Omega$ of full measure such that the~%
		following condition is satisfied:
		\begin{equation*}
			%
		\forall 
			\bx = \set{ x_I }_{ I \subseteq [k] }
			\subseteq 
			F
		\qquad 
			\abs*{
				\diff_{ I = \emptyset }^{ [k] }
					f(x_I)
			}
			\le 
			\polygen{ \gen{k} }(\bx)
			\sum_{I = \emptyset }^{ [k ] }
				G(x_I),
		\end{equation*}
		where, letting 
		$ \partitions{[k]} $
		denote the family of partitions of 
		$[k] \coloneq \set*{1, \, \ldots, \, k }$,
		\begin{equation*}
			\diff_{ I = \emptyset }^{ [k] }
				f(x_I)
			\coloneq
			\sum_{ I = \emptyset }^{ [k] }
				(-1)^{ \card*{ [k] \setminus I } }
				f\del*{ x_I }
		\quad \text{and} \quad 
			\polygen{ \gen{k} }(\bx)
			\coloneq
			\sum_{ \cP \in \partitions{[k]} }
				\prod_{ S \in \cP }
					\sum_{ A = \emptyset }^{ [k] \setminus S }
						\norm*{
							\sum_{ I= \emptyset }^{ S }
								(-1)^{ \card*{ S \setminus I } }
								x_{I \cup A }
						}.
		\end{equation*}
	\end{enumerate}
	Moreover, there are constants
	$ C_1,C_2 > 0$ 
	such that
	\begin{equation*}
	\forall 
		f \in \sobolev{k,p}( \Omega)
	\qquad 
		C_1 \norm{ f }_{ \sobolev{k,p}( \Omega) }
		\le 
		\norm{f}_{ \integrable{p}(\Omega) }
		+
		\inf_{ G \in \bD_{\lambda}^{k}(f) }
			\norm{G}_{ \integrable{p}(\Omega) }
		\le 
		C_2 \norm{ f }_{ \sobolev{k,p}( \Omega) }.
	\end{equation*}
\end{theorem}
In fact, the above theorem will be a quick consequence of a more general characterization of higher-order Sobolev spaces based on Banach function spaces.
Furthermore, we will show an analogous characterization of higher-order H\"{o}lder spaces.

These characterizations introduce a potentially novel approach to defining function spaces.
Indeed, there are many function spaces that are defined using conditions involving two points of~the~underlying space, for example,
H\"{o}lder spaces, 
Hajłasz--Sobolev spaces, Newtonian spaces \cite{Heinonen_Koskela_Shanmugalingam_Tyson_2015, Newtonian}, or even, arguably, Sobolev--Slobodeckij spaces \cite{fractional_sobolev}.
However, the possibility of introducing conditions that involve more than two points does not seem to be well-explored yet.
As considering more numerous tuples of points is precisely what will allow us to~characterize higher-order Sobolev and H\"{o}lder spaces, 
this approach could be a fruitful ground for~new research.
It is also the reason why we~have decided to call our characterization a ``multipoint'' one.

The introduction of first-order Sobolev spaces
on metric measure spaces using the Hajłasz's result is~fairly straightforward.
However, it is not immediately clear if we can use
\Cref{thm::sobolev_characterization}
to do the same for higher-order Sobolev spaces.
Indeed, the main issue is the presence of the
$
	\norm*{
		\sum_{ I= \emptyset }^{ S }
			(-1)^{ \card*{ S \setminus I } }
			x_{I \cup A }
	}
$
term
in~the~definition of~%
$
	\polygen{ \gen{k} }\del*{ \bx }.
$
Nevertheless, as we will see in \cref{sec::metric_spaces},
this issue can be overcome,
and the resulting function spaces 
$ \multipoint{s,p}\del*{X}$
work rather well with the Hajłasz--Sobolev spaces $\hajlasz{1,p}\del*{X}$.
For example, 
we have
$
	\multipoint{1,p}\del*{X}
	\cong 
	\hajlasz{1,p}\del*{X}
$
and
$
	\multipoint{s,p}\del*{X}
	\hookrightarrow
	\hajlasz{1,p}\del*{X}
$
if $s > 1$.

Let us note that
the introduction of higher-order Sobolev spaces on metric measure spaces has already been done in the past.
Indeed, for example,
in \cite{Wheeden},
the authors introduced such spaces on 
metric spaces of homogeneous type in the sense of Coifman and Weiss.
In this approach, the spaces are defined with respect to an appropriate family of functions, therein called
\say{polynomial functions,} that has to be provided. 
For the second-order spaces, a different approach was developed in 
\cite{Gigli}.
There, the authors, treating the theory of weak Hessians on Riemannian manifolds as a guideline, 
develop the theory of second-order Sobolev spaces on $\mathsf{RCD}\del*{ K, \infty} $ spaces.

Other examples of such results are present in 
\cite{Verdera}
and
\cite{AveragesOnBalls},
where the authors provide various characterizations of higher-order Sobolev spaces on $\bR^n$ and use them to define such spaces on metric measure spaces.
Furthermore, the authors also characterize fractional order Sobolev spaces.
However,  it is not clear whether one can obtain similar characterizations of higher-order Sobolev spaces when the underlying space $\Omega$ is a proper subset of $\bR^n$.%
\footnote{
	It should be noted that in
	\cite{Sobolev_on_sphere},  
	the authors showed that higher-order Sobolev spaces defined on the unit sphere $\bS^{n-1}$ could be characterized in a similar way as the one present
	in \cite{Verdera}.
}
As such, especially for integer orders, the results that we obtain here could be considered more general than the ones present in the two papers. 
Moreover, in contrast to the other approaches, we define higher-order Sobolev spaces on completely arbitrary metric spaces with a measure.

The rest of the paper is structured as follows.

In \cref{sec::preliminaries}, we recall some basic theory regarding locally integrable functions, Hardy--Littlewood's maximal function, and weak derivatives.
\Cref{sec::combinatorial_preludium} is devoted to many combinatorial lemmas that will be used in~the~next sections.
In \cref{sec::characterizations}, we prove the main theorems of the paper. 
In particular, we prove our characterization of higher-order Sobolev spaces based on Banach function spaces and give several examples of Banach function spaces for which it can be used.
In this section, we also prove the characterization of higher-order H\"{o}lder spaces.
Finally, in 
\cref{sec::metric_spaces},
we explain how the obtained results can be used to define higher-order Sobolev and 
H\"{o}lder spaces on metric measure spaces.

\section{Preliminaries}
\label{sec::preliminaries}
Let $\Omega \subseteq \bR^n$ be open,
$(V, \normAlone)$ be a finite-dimensional normed space,
$A \subseteq \Omega$ be such that $0 < \abs{A} < \infty$,
and let $f \colon \Omega \to V$
be measurable. 
We define the average of $f$ over $A$ by 
\begin{equation*}
    \sintegral{A}{f}{x}
    \coloneq
    \frac{1}{\abs{A} }
    \integral{A}{f}{x}
\end{equation*}
as long as the expression on the right hand side is well-defined.
We will say that $x \in \Omega$ is a \emph{Lebesgue point} of $f$ if 
\begin{equation*}
    \lim_{ r \to 0^+ }
    \sintegral{ \ball{x, r} }{
        \norm{ f(x) - f(y) }
    }{y}
    =
    0.
\end{equation*}
We will denote the set of all Lebesgue points of $f$ by 
$\lebesguePoints{f}$.
We will say that  
$ f $ is \emph{locally integrable}
(which we will denote by writing $f \in \locIntegrable{1}(\Omega;V)$)
if for every $y \in \Omega$
there exists an open neighborhood
$U$ of $y$ such that
$
    \integral{U}{\norm{f}}{x} < \infty.
$
Note that this is equivalent to stating that
$
    \integral{K}{\norm{f}}{x} < \infty
$
for all compact $K \subseteq \Omega$.
It is worth noting that 
if $f \in \locIntegrable{1}(\Omega;V)$, then
$\lebesguePoints{f}$ is of full measure in $\Omega$
by the Lebesgue differentiation theorem.

Now, assume that $\Omega = \bR^n$.
For $r > 0$ we define the operator $\mathsf{B}_r$ by
\begin{equation*}
    \forall x \in \bR^n
    \qquad 
    \ballAv{r}{f}(x)
    \coloneq
    \sintegral{ \ball{x, r} }{ f }{ x }.
\end{equation*}
We define the \emph{Hardy--Littlewood maximal function} of $f$ by 
\begin{equation*}
    \forall x \in \bR^n 
    \qquad 
    \maxFun{f}(x)
    \coloneq
    \sup_{ r > 0 }
    \sintegral{\ball{x, r}}{ \norm{f} }{ x }.
\end{equation*}
For $R > 0$,
we also define the \emph{$R$-restricted Hardy--Littlewood maximal function} of $f$ by 
\begin{equation*}
    \forall x \in \bR^n 
    \qquad 
    \restMaxFun{R}{f}(x)
    \coloneq
    \sup_{ r \in \intoc{0, R} }
    \sintegral{\ball{x, r}}{ \norm{f} }{ x }.
\end{equation*}
We also define
\begin{equation*}
    \forall x \in \bR^n 
    \qquad 
    \restMaxFun{0}{f}(x)
    \coloneq
    \limsup_{ r \to 0^+ }
    \sintegral{\ball{x, r}}{ \norm{f} }{ x }.
\end{equation*}
In this paper, 
we consider vector-valued functions primarily to make expressions like 
$
    \maxFun{ \nabla f }
$
well-defined.
Note that while 
$ \ballAv{r}{f} $ 
might not be well-defined for given $r > 0$ and 
measurable function $f \colon \bR^n \to V$, 
$\mathsf{M}(f)$ and $\mathsf{M}_R(f)$, where $R \ge 0$,
are always well-defined.
However, if $f \in \locIntegrable{1}(\bR^n; V)$,
then $ \ballAv{r}{f} $ is well-defined.

Let us now recall some of the properties of the three operators mentioned above:
\begin{itemize}
    \item 
    For any $r > 0$,
    operator $\mathsf{B}_r$ is linear, 
    whereas the operators: $\mathsf{M}$ and $\mathsf{M}_R$, where $R \ge 0$, are sublinear,
    \item
    For all $r > 0$ and $R \ge r$,
    we have
    $
        \ballAv{r}{f}
        \le 
        \restMaxFun{R}{f}
        \le
        \maxFun{f}
    $
    everywhere.
    \item 
    If $x \in \lebesguePoints{f}$, then for any $R \ge 0$,
    we have
    $
        f(x)
        \le 
        \restMaxFun{R}{f}.
    $
    \item 
    If 
    $ \norm{f} \le \norm{g} $ 
    almost everywhere, 
    then for all $R \ge 0$,
    we have
    $
        \restMaxFun{R}{f}
        \le 
        \restMaxFun{R}{g}
    $
    and 
    $
        \maxFun{f} \le \maxFun{g}
    $
    everywhere.    
\end{itemize}
\begin{lemma}\label{lem::composing_maximal_functions}
    Let $n \in \bN$,
    $\del*{ V, \normAlone }$ be a finite-dimensional normed space,
    and $f \colon \bR^n \to V$ be~measurable.
    Then for all $m \in \bN$ and $R \ge 0$
    we have
    $
        \restMaxFun[m]{R}{f}
        \le 
        \restMaxFun[m+1]{R}{f}
    $
    and 
    $
        \maxFun[m]{f} \le \maxFun[m+1]{f}
    $
    everywhere.
\end{lemma}
\begin{proof}
    Fix $R \ge 0$.
    Denote 
    \begin{equation*}
        \Omega 
        \coloneq 
        \setc*{
            x \in \bR^n 
        }{
            \integral{ \ball{x,r} }{ 
            	\norm{f} 
        	}{x} < \infty 
            \text{ for some }
            r > 0
        }.
    \end{equation*}
    It is clear that $\Omega$ is open and that 
    $ 
        f \rvert_{\Omega} 
        \in 
        \locIntegrable{1}\del*{ \Omega; V },
    $
    hence 
    $ 
        \norm{f}
        \le 
        \restMaxFun{0}{f} 
        \le 
        \restMaxFun{R}{f} 
        \le 
        \maxFun{f} 
    $
    almost everywhere in~$\Omega$,
    where the first inequality follows by the Lebesgue differentiation theorem.
    Next, suppose that $x \notin \Omega$.
    Then 
    $
    	\ballAv{r}{ 
    		\norm{f}  
    	}(x) = \infty
    $ 
    for all $r > 0$, hence
    $
    	\restMaxFun{R}{f}(x) = \infty 
    $ 
    and 
    $
    	\maxFun{f}(x) = \infty.
    $
    Thus, 
    $
        \norm{f(x)} 
        \le 
        \restMaxFun{R}{f}(x) 
        \le 
        \maxFun{f}(x).
    $
    We have showed that both
    $
    	\norm{f} 
    	\le 
    	\restMaxFun{R}{f}
    $
    and 
    $
    	\norm{f}
    	\le 
    	\maxFun{f}
    $
    almost everywhere in $\Omega$ and everywhere in 
    $\bR^n \setminus \Omega$, 
    hence almost everywhere in $\bR^n$.
    It follows that 
    $ 
    	\restMaxFun{R}{f} 
    	\le 
    	\restMaxFun[2]{R}{f}
    $
    and $\maxFun{f} \le \maxFun[2]{f}$ 
    everywhere in~$\bR^n$.

    Now, suppose that for some $m \in \bN$ and $R \ge 0$
    we have 
    $ \restMaxFun[m]{R}{f} \le \restMaxFun[m+1]{R}{f}$
    and $\maxFun[m]{f} \le \maxFun[m+1]{f}$ everywhere.
    Then also 
    $ \restMaxFun[m+1]{R}{f} \le \restMaxFun[m+2]{R}{f}$
    and $\maxFun[m+1]{f} \le \maxFun[m+2]{f}$
    everywhere.
    %
\end{proof}
\begin{lemma}\label{lem::shifting_operators}
    Let $n \in \bN$, $r > 0$, $R \ge 0$,
    and $\del*{ V, \normAlone }$ 
    be a finite-dimensional 
    normed space.
    Let $f \colon \bR^n \to V$ be 
    measurable.
    Then for all $x,y \in \bR^n$,
    \begin{gather*}
        \ballAv{r}{
        	f\del*{ \blank + y } 
        }(x)
        =
        \ballAv{r}{ f }(x+y),
        \quad
        \restMaxFun{R}{
        	f\del*{ \blank + y} 
        }(x)
        =
        \restMaxFun{R}{f}(x+y),
    \\
    \text{and} \qquad 
        \maxFun{
        	f\del*{ \blank + y } 
        }(x)
        =
        \maxFun{f}(x+y),
    \end{gather*}
    where the first of the above equalities is understood in the sense that if one side is well-defined, then so is the other, and the equality is satisfied;
    in the latter ones, both sides are always well-defined.
\end{lemma}
\begin{proof}
    For the first equality, let us first assume that $f \ge 0$.
    Since $\abs{ \ball{x, r} } = \abs{ \ball{x+y, r } }$,
    \begin{multline*}
        \ballAv{r}{
        	f\del*{ \blank + y }
        }(x)
        =
        \sintegral{ \ball{x, r} }{ f(z + y) }{z}
        =
        \dfrac{1}{ \abs{\ball{x, r}} }
        \integral{ \ball{x, r} }{ f(z + y) }{z}
    \\
        =
        \dfrac{1}{ \abs{\ball{x + y, r}} }
        \integral{ \ball{x + y, r} }{f(z)}{z}
        =
        \sintegral{ \ball{x + y, r} }{f(z)}{z}
        =
        \ballAv{r}{ f }(x+y).
    \end{multline*}
    In the general case, we have that
    $ f = f^+ - f^- $,
    where 
    $
    	f^+ \coloneq \max\del*{ f, 0 }
    $
    and
    $
    	f^- \coloneq \max\del*{ -f, 0 }.
    $
    Hence, 
    \begin{align*}
        \ballAv{r}{
        	f\del*{ \blank + y }
        }(x)
        =
        \ballAv{r}{
        	f^+\del*{ \blank + y }
        }(x)	
        -
        \ballAv{r}{
        	f^-\del*{ \blank + y }
        }(x)
    &
    \\
        =
        \ballAv{r}{ f^+ }(x+y)
        -
        \ballAv{r}{ f^- }(x+y)
    &=
        \ballAv{r}{ f }(x+y),
    \end{align*}
    where the equalities are valid as long as
    $
    	\ballAv{r}{
        	f\del*{ \blank + y }
		}(x)
    $
    or
    $
    	\ballAv{r}{ f }(x+y)
    $
    is well-defined.

    The second equality follows from the fact that,
    for $R > 0$,
    we have
    \begin{equation*}
        \restMaxFun{R}{
        	f\del*{ \blank + y } 
        }(x)
        =
        \sup_{r \in \intoc{0,R}}
            \ballAv{r}{ 
            	\norm{  
            		f\del*{ \blank + y} 
            	}
            }(x)
        =
        \sup_{r \in \intoc{0,R}}
            \ballAv{r}{ 
            	\norm{ f}
            }(x+y)
        =
        \restMaxFun{R}{f}(x+y),
    \end{equation*}
    while for $R = 0$ we have
    \begin{equation*}
        \restMaxFun{0}{
        	f\del*{ \blank + y } 
        }(x)
        =
        \limsup_{r \to 0^+}
            \ballAv{r}{ 
            	\norm{ 
            		f\del*{ \blank + y }
            	}
            }(x)
        =
        \limsup_{r \to 0^+}
            \ballAv{r}{ 
            	\norm{f}
            }(x+y)
        =
        \restMaxFun{0}{f}(x+y).
    \end{equation*}
    Finally, we have the last equality since
    \begin{equation*}
        \maxFun{
        	f\del*{ \blank + y }
        }(x)
        =
        \sup_{r > 0}
            \ballAv{r}{ 
            	\norm{  
            		f\del*{ \blank + y }
            	}
            }(x)
        =
        \sup_{r > 0 }
            \ballAv{r}{ 
            	\norm{ f }
            }(x+y)
        =
        \maxFun{f}(x+y). \qedhere
    \end{equation*}
\end{proof}
\begin{lemma}\label{lem::centered_max_fun_bounds_shifted_average}
    Let $n \in \bN$, $x \in \bR^n$, $r > 0$,
    and $\del*{ V, \normAlone }$ 
    be a finite-dimensional normed space. 
    Let $f \in \locIntegrable{1}\del*{ \bR^n; V}$.
    Then for every 
    $R \ge r$ and $y \in \clball{x, r}$ we have
    $
        \ballAv{r}{ \norm{f} }(y) 
        \le 
        2^n
        \restMaxFun{2R}{f}(x).
    $
\end{lemma}
\begin{proof}
    Let us notice that 
    $
        \abs{ \ball{y, r} }
        =
        \abs{ \ball{x, r} },
    $
    $
        \frac{ \abs{ \ball{x, 2r} } }{ \abs{ \ball{x, r} } }= 2^n,
    $
    and 
    $
        \ball{y, r} \subseteq \ball{x, 2r}.
    $
    Therefore,
    \begin{align*}
        \ballAv{r}{ \norm{f} }(y) 
    &=
        \frac{1}{ \abs{ \ball{y, r} } }
        \integral{ \ball{y, r} }{ \norm{f} }{z}
    \\
    &\le 
        \frac{1}{ \abs{ \ball{y, r} } }
        \integral{ \ball{x, 2r} }{ \norm{f} }{z}
        =
        \frac{ \abs{ \ball{x, 2r} } }{ \abs{ \ball{x, r} } }
        \sintegral{ \ball{x, 2r} }{ \norm{f} }{z}
        \le
        2^n \restMaxFun{2R}{f}(x).
        \tag*{\qedhere}
    \end{align*}
\end{proof}

Let $n \in \bN$ and $\Omega \subseteq \bR^n$ be open.
We will write 
$ \smoothCompactSupport(\Omega)$
to denote the family of smooth functions 
$\phi \colon \Omega \to \bR$ 
whose support is a compact subset of $\Omega$.
We will say that $h  \in \locIntegrable{1}(\Omega)$ 
is an $\alpha\th$ weak derivative of 
$ f \in \locIntegrable{1}(\Omega)$,
where 
$ \alpha \in \bN_0^n  $,
if 
\begin{equation*}
	\forall \phi \in \smoothCompactSupport(\Omega)
\qquad  
	\integral{\Omega}{ h \phi }{ x }
	=
	(-1)^{ \abs{\alpha} }
	\integral{\Omega}{ f  \partial^{\alpha}\phi }{ x }.
\end{equation*}
If the above condition is satisfied, we will write $h = \partial^{\alpha} f$.
Note that, in particular, $\partial^0  f =f $.

For a fixed $k \in \bN_0$, we then denote
\begin{equation*}
	%
	\locSobolev{k,1}(\Omega)
	\coloneq 
	\setc*{
		f \in \locIntegrable{1}(\Omega)
	}{
		\partial^{\alpha} f 
		\in 
		\locIntegrable{1}\del*{ \Omega }
		\text{ for all }
		\alpha \in \bN_0^n
		\text{ with }
		\abs{ \alpha } \le k
	}.
\end{equation*}
In particular, we have
$ \locSobolev{0,1}(\Omega)=\locIntegrable{1}(\Omega)$.
For $f \in \locSobolev{k,1}(\Omega)$ and $j \in [k]$,
we will often write 
$ \nabla^j f$ to denote the $j\th$ order (weak) gradient of $f$.

\subsection{Conventions and Notation}
We adopt the following conventions and notation:
\begin{itemize}
    %
    %
    %
    %
    \item 
    We will often not distinguish between the measurable functions 
    and their equivalence classes in the relation of equality $\measureAlone$-almost everywhere.
    \item 
    $ 0 \cdot \infty = 0$,
    \item
    Suppose we have an expression of the form 
    $
        \abs*{ \sum_i v_i}.
    $
    If for any of the $v_i$'s we have 
    $\abs{v_i} = \infty$, 
    we assign the value of $\infty$ to this expression.
    \item 
    The sum over an empty set of indices equals $0$,
    \item 
    The product over an empty set of indices equals $1$,
    \item 
    For $k \in \bN$ we denote 
    $[k] \coloneq \set{1, \, \ldots, \, k }$. 
    We also put
    $[0] \coloneq \emptyset$.
    \item 
    We will write $\delta_a^b$ to denote the Kronecker's delta, 
    i.e., a function such that $\delta_a^b = 1$ if $a = b$, and 
    $\delta_a^b = 0$ otherwise.
    \item 
    We denote the symmetric difference of sets by $\div$,
    i.e.
    $
    	A \div B 
    	\coloneq 
    	(A \setminus B) \cup (B \setminus A).
    $
    \item 
    We will usually use $\absAlone$ to denote the cardinality of a set,
    Lebesgue measure of a (Lebesgue) measurable set, as well as the order of a multiindex.
    However, if there would be a need to explicitly distinguish between these notions, we will use $\#$ for the cardinality and $\lambda$ for the Lebesgue measure.
    %
    \item 
    We will denote the indicator function of a set $A$ by $\indicator{A}$.
    \item 
    We will use both the parentheses and curly braces to refer to a tuple. 
    For example, we will consider expressions of the form
    $
        \del{ x_i }_{ i \in I }
    $
    and 
    $
        \set{ x_i }_{ i \in I }
    $
    equivalent.
    However, we will always use the latter form when writing expressions of the form
    $
        \set{ x_i }_{ i \in I } \subseteq A
    $
    to indicate that every element of the tuple is an element of $A$.
    Also, when we will index the elements of a tuple by subsets of some set $B$, we will write 
    $
        \set{ x_I }_{ I \subseteq B }
    $
    instead of 
    $
        \set{ x_I }_{ I \in 2^B },
    $
    \item 
    If $\del{X, \Sigma, \measureAlone}$ is a measure space,
    then we will write 
    $ \measureAlone \forall x \in X$
    to indicate that what follows is satisfied
    for $\measureAlone$-almost every $x \in X$.
    Also, by writing, for example,
    $ 
    	\measureAlone 
    	\forall 
    	\bx = 
    	\set*{ x_I }_{
    		I\subseteq [k]
    	}
    	\subseteq X
    $,
    we mean that there exists a~set 
    $ F \subseteq X$ of full measure such that whatever follows is satisfied for all tuples
    $     	
    	\bx = 
    	\set*{ x_I }_{
    		I\subseteq [k]
    	}
    	\subseteq F.
    $
    \item 
    Let $\del{\Omega, \Sigma, \measureAlone}$
    be a measure space and 
    $ 
    	\mc{F}(\Omega)
    $
    be a Banach function space over it.
	Then we might write 
	$ \nabla^j f \in \mc{F}(\Omega)$
	without explicitly writing the codomain
	to mean that we have
	$ \partial^{\alpha} f \in \mc{F}(\Omega) $
	for all multiindices 
	$ \alpha $
	such that 
	$ \abs{\alpha} = j$.
	Moreover, we put
	$
		\norm*{
			\nabla^j f
		}_{ \mc{F}\del*{ \Omega} }
		\coloneq
		\norm*{
			\,
			\norm*{
				\nabla^j f
			}
			\,
		}_{ \mc{F}\del*{ \Omega} },
	$  
	where 
	$
		\norm*{
			\nabla^j f
		}
	$
	is the Euclidean norm applied pointwise to 
	$
		\nabla^j f.
	$
\end{itemize}

\section{Combinatorial Preludium} 
\label{sec::combinatorial_preludium}
In this section, we prove many combinatorial lemmas that will be needed for the characterization of higher-order Sobolev spaces.
Since this part will be fairly sizable, we have decided to separate it into its own section.

\begin{definition}[Partition of a finite set]
\label{def::partition}
    Let $S$ be a finite nonempty set. 
    We will say that $\cP \subseteq 2^{ S }$
    is a \emph{partition} of $ S $ if the elements of $\cP$  
        are nonempty and pairwise disjoint and
        $\bigcup_{ P \in \cP} P = S$.
    We will write 
    $ \partitions{S } $
    to denote the family of partitions of $S $. 
    We also define the family
    \begin{equation*}
        \partitions[j]{ S }
        \coloneq 
        \setc*{
            \cP \in \partitions{ S }
        }{
            \card{ \cP } = j
        }.
    \end{equation*}
\end{definition}
\begin{remark}
    Let us remark that for any finite nonempty set $S$,
    we have
    $
        \card{ \partitions{S} }
        =
        \partitions{ \card{S} },
    $
    where $\partitions{n}$ 
    denotes the $n\th$ Bell number.
    Also, for such $S$ we have
    $ \partitions[1]{S} = \set*{ \set*{ S } } $
    and 
    $ 
        \partitions[ \card{S} ]{ S }
        =
        \set*{
            \setc*{
                \set*{ j } 
            }{
                j \in S
            }
        }.
    $
\end{remark}
\begin{notation}
    We will frequently use subsets of 
    $ [k] $, where $k \in \bN$,
    as indices over which we will iterate.
    To make it easier to differentiate between
    the set-indices and the number-indices, we will use the uppercase letters for the former and lowercase letters for the latter.
    Furthermore, when writing $\sum_{I = A}^B$,
    we sum over
    all $I$ such that $A \subseteq I \subseteq B$.
    For example, for every $f \colon \bN_0 \to \bR$ 
    and all $A \subseteq B \subseteq [k]$,
    we have 
    \begin{equation} \label{eqnot::example_of_summing_with_sets}
        \sum_{ I = A }^{ B } f( \card{I} )
        =
        \sum_{i = \card{A}}^{ \card{B} } 
            \binom{ \card{ B \setminus A } }{ i - \card{A} }
            f(i)
        =
        \sum_{i = 0}^{ \card{B \setminus A } } 
            \binom{ \card{ B \setminus A } }{ i }
            f( i + \card{A} ),
    \end{equation}
    where we use the fact that for all 
    $i \in \intcc{ \, \card{A}, \card{B} \, }$ 
    there are exactly 
    $ \binom{ \card{ B \setminus A } }{ i - \card{A} } $
    sets $I \subseteq [k] $ such that $A \subseteq I \subseteq B$ 
    and 
    $ \card{ I } = i $.
    (Note that $\binom{0}{0} = 1 $.)
\end{notation}
\begin{lemma}\label{lem::alternating_sum_over_cardinalities}
    Let $k \in \bN$ and $A \subseteq B \subseteq [k]$. Then 
    $
        \sum_{I = A}^{B}
            (-1)^{ \card{I} }
        =
        (-1)^{\card{A}}
        \delta_{A  }^{B }.
    $
\end{lemma}
\begin{proof}
    If $A = B$, then 
    $
        \sum_{I = A}^{B}
            (-1)^{ \card{I} }
        =
        (-1)^{ \card{A} }.
    $
    Suppose now that $A \ne B$. 
    Using \eqref{eqnot::example_of_summing_with_sets} 
    and the Binomial Theorem, we get
    \begin{equation*}
        \sum_{I = A}^{B}
            (-1)^{ \card{I} }
        =
        \sum_{i = 0}^{ \card{B \setminus A } } 
            \binom{ \card{ B \setminus A } }{ i }
            (-1)^{ i + \card{A} }
        =
        (-1)^{\card{A}}
        (1-1)^{\card{ B \setminus A }}
        =
        0.
        \tag*{\qedhere}
    \end{equation*}
\end{proof}
\begin{lemma}
\label{lem::between_A_and_B_cardinalities}
    Let $k \in \bN$ and $A \subseteq B \subseteq [k]$ with 
    $B \setminus A \ne \emptyset$. Then
    \begin{equation*}
        \card{
            \setc*{
                I \subseteq B 
            }{
                A \subseteq I \text{ and }
                \card{I} \text{ is odd}
            }
        }
        =
        \card{
            \setc*{
                I \subseteq B 
            }{
                A \subseteq I \text{ and }
                \card{I} \text{ is even}
            }
        }
        =
        2^{ \card{B} - \card{A} - 1 }.
    \end{equation*}
\end{lemma}
\begin{proof}
    First of all, let us notice that the set
    $
        \setc*{
            I \subseteq B
        }{
            A \subseteq I
        }
    $
    has $2^{ \card{B} - \card{A} }$ elements.
    Indeed, this follows from the fact that by 
    \eqref{eqnot::example_of_summing_with_sets} 
    we have
    \begin{equation*}
	    \card*{
	    	\setc*{
				I \subseteq B
			}{
				A \subseteq I
			}
	    }
	    =
	    \sum_{ I = A }^{ B }
	    	1
	    =
	    \sum_{i = 0}^{ \card{ B \setminus A } }
	    	\binom{ \card{ B \setminus A } }{ i }
	    =
	    2^{ \card{ B \setminus A } }.
    \end{equation*}
    Therefore, once we prove that the set
    $
        \setc*{
            I \subseteq B
        }{
            A \subseteq I
        }
    $ 
    has the same number of 
    elements with odd cardinality and even cardinality, the claim will follow.
    This, however, follows from the fact that by 
    \Cref{lem::alternating_sum_over_cardinalities},
    \begin{equation*}
        \sum_{ \substack{ I = A \\ \card{I} \text{ is even}} }^{
            B
        }
            1
        -
        \sum_{ \substack{ I = A \\ \card{I} \text{ is odd}} }^{
            B
        }
            1
        =
        \sum_{I = A}^{B}
            (-1)^{ \card{I} }
        =
        0.
        \tag*{\qedhere}
    \end{equation*}
\end{proof}
\begin{lemma} \label{lem::sum_over_x_IuS_gives_2^|S|_before_the_sum}
    Let $(G, +)$ be an Abelian group, $k \in \bN$, and $S \subseteq [k]$. 
    Fix 
    $
        \bx = \set*{x_I}_{ I \subseteq [k] } 
        \subseteq G.
    $
    Then 
    \begin{equation*}
        \sum_{ I = \emptyset }^{ [k] }
            x_{ I \cup S} 
        =
        2^{ \card{S} }
        \sum_{ I = S}^{ [k] }
            x_I.
    \end{equation*}
\end{lemma}
\begin{proof}
    Let us notice that $(J, L) \mapsto J \cup L$
    is a bijection from $2^{S} \times 2^{ [k] \setminus S}$ 
    to $2^{ [k] }$.
    Also, 
    $ L \mapsto L \cup S$ 
    is a bijection from 
    $2^{ [k] \setminus S}$ 
    to 
    $ \setc*{ I \subseteq [k] }{ S \subseteq I}$.
    In consequence,
    \begin{equation*}
        \sum_{ I = \emptyset }^{ [k] }
            x_{ I \cup S} 
        =
        \sum_{ J = \emptyset }^{ S }
            \sum_{ L = \emptyset }^{ [k] \setminus S }
                x_{ J \cup L \cup S }
        =
        \sum_{ J = \emptyset }^{ S }
            \sum_{ L = \emptyset }^{ [k] \setminus S }
                x_{ L \cup S }
        =
        2^{ \card{S} }
        \sum_{ L = \emptyset }^{ [k] \setminus S }
            x_{ L \cup S }
        =
        2^{ \card{S} }
        \sum_{ I = S }^{ [k]  }
            x_{ I }. \tag*{\qedhere}
    \end{equation*}
\end{proof}
\begin{lemma}
\label{lem::folding_hypercubes_with_a_sum}
    Let $(G, \! +)$ be an Abelian group, $k \in \bN,$ $m \in [k]$, and $i \in [m]$.
    Fix
    $
    	{
        \bx = \set*{ x_I }_{ I \subseteq [m] }
        \subseteq
        G.
        }
    $
    Define $\by = \set*{ y_J }_{ J \subseteq [k] }$ by 
    \begin{equation*}
        \forall I \subseteq [m] 
        \quad \forall L \subseteq [k] \setminus [m]
        \qquad
        y_{ I \cup L} 
        \coloneq 
        \begin{cases}
            x_I & \text{if } \card{L} \text{ is even}, \\
            x_{I \div \set*{i} } & \text{otherwise}.
        \end{cases}
    \end{equation*}
    Then
    \begin{equation*}
        \sum_{ J = \emptyset }^{ [k] }
            y_J
        =
        2^{k-m}
        \sum_{ I = \emptyset }^{ [m] }
            x_I.
    \end{equation*}
\end{lemma}
\begin{proof}
    First, let us notice that 
    $(I, L) \mapsto I \cup L$ is a bijection from 
    $2^{[m]} \times 2^{ [k] \setminus [m]}$ to 
    $2^{ [k] }$. 
    Also, 
    $
        I \mapsto I \div \set*{i}
    $
    is a bijection from $2^{ [m] }$ to itself.
    Therefore, 
    \begin{align*}
        \sum_{ J = \emptyset }^{ [k] }
            y_J
        =
        \sum_{ L = \emptyset }^{ [k] \setminus [m] }
            \sum_{ I = \emptyset }^{ [m] }
                y_{ I \cup L }
        &=
        \sum_{ \substack{   
            L = \emptyset 
            \\ \card{L} \text{ is even}
        }}^{ [k] \setminus [m] }
            \sum_{ I = \emptyset }^{ [m] }
                y_{ I \cup L }
        +
        \sum_{ \substack{   
            L = \emptyset 
            \\ \card{L} \text{ is odd}
        }}^{ [k] \setminus [m] }
            \sum_{ I = \emptyset }^{ [m] }
                y_{ I \cup L }
    \\
        &=
        \sum_{ \substack{   
            L = \emptyset 
            \\ \card{L} \text{ is even}
        }}^{ [k] \setminus [m] }
            \sum_{ I = \emptyset }^{ [m] }
                x_I
        +
        \sum_{ \substack{   
            L = \emptyset 
            \\ \card{L} \text{ is odd}
        }}^{ [k] \setminus [m] }
            \sum_{ I = \emptyset }^{ [m] }
                x_{ I \div \set*{i} }
    \\
        &=
        \sum_{ \substack{   
            L = \emptyset 
            \\ \card{L} \text{ is even}
        }}^{ [k] \setminus [m] }
            \sum_{ I = \emptyset }^{ [m] }
                x_I
        +
        \sum_{ \substack{   
            L = \emptyset 
            \\ \card{L} \text{ is odd}
        }}^{ [k] \setminus [m] }
            \sum_{ I = \emptyset }^{ [m] }
                x_{ I }
        =
        \sum_{ L = \emptyset }^{ [k] \setminus [m] }
            \sum_{ I = \emptyset }^{ [m] }
                x_I
        =
        2^{k-m}
        \sum_{ I = \emptyset }^{ [m] }
            x_I,
    \end{align*}
    as claimed.
\end{proof}
\begin{definition}
\label{def::diff}
    Let $(G, +)$ be an Abelian group, $k \in \bN$, and let 
    $
        \set*{ x_I }_{ I \subseteq [k] }
        \subseteq 
        G.
    $
    For given sets 
    $A \subseteq B \subseteq [k]$, 
    we define
    \begin{equation*} 
        \diff_{ I = A }^{ B } 
            x_I
        \coloneq 
        \sum_{ I = A }^{ B }
            (-1)^{ \card{ B \setminus I} }
            x_I.
    \end{equation*}
    For disjoint $A, S \subseteq [k]$,
    we also introduce the notation
    \begin{equation*}
        %
        \diff_{ I \upmapsto A }^{ S }
            x_I
        \coloneqq 
        \diff_{ I = A}^{ S \cup A }
            x_I.
    \end{equation*}
\end{definition}
We will frequently use the following remark.
\begin{remark}
\label{rem::diff_with_an_arrow}
	Let $\del*{G, +}$ be an Abelian group, 
	$ k \in \bN $, 
	and 
	$
		\set*{ x_I }_{ I \subseteq [k] }
		\subseteq 
		G.
	$
	Then for all disjoint
	$ A, S \subseteq [k]$,
    since 
    $ A \cap S = \emptyset $,
    the mapping 
    $ I \mapsto I \cup A $
    is a bijection from 
    $ 2^S $ to 
    $ \setc*{ I \subseteq A \cup S }{ A \subseteq I } $.
    Hence,
    \begin{equation*} 
        %
        \diff_{ I \upmapsto A }^{ S }
            x_I
        =
        \diff_{ I = A}^{ S \cup A }
            x_I
        =
        \sum_{I = A }^{S \cup A}
            (-1)^{ \card{ (S \cup A) \setminus I} }
            x_{I}
        =
        \sum_{I = \emptyset}^{S}
            (-1)^{ \card{ S \setminus I} }
            x_{I \cup A}
        =
        \diff_{ I = \emptyset}^{S}
            x_{I \cup A}.
    \end{equation*}
\end{remark}
\begin{lemma} \label{lem::diff_constant_term}
    Let $(G, +)$ be an Abelian group, $g \in G$, $k \in \bN$,
    and 
	$
		\set*{ x_I }_{ I \subseteq [k] }
		\subseteq 
		G.
	$
    Let $A, S \subseteq [k]$ be disjoint. 
    If $x_I = g$ for all $I \subseteq [k]$ 
    satisfying $A \subseteq I \subseteq A \cup S$,
    then
    \begin{equation*}
        \diff_{ I \upmapsto A}^{ S } x_I 
        =
        \begin{cases}
            g 
        & 
        	\text{if } S = \emptyset, 
        \\ 
            0 
        &
        	\text{otherwise}.
        \end{cases}
    \end{equation*}
\end{lemma}
\begin{proof}
    Since 
    $I \mapsto S \setminus I$ is a bijection from 
    $2^S$ to itself, by 
    \Cref{lem::alternating_sum_over_cardinalities}
    we have
    \begin{equation*}
        \diff_{ I \upmapsto A}^{ S } x_I 
        =
        \sum_{ I = \emptyset }^{ S }
            (-1)^{ \card{ S \setminus I } } 
            x_{ I \cup A }
        =
        \sum_{ I = \emptyset }^{ S }
            (-1)^{ \card{ S \setminus I } } 
            g
        =
        \del*{
            \sum_{ I = \emptyset }^{ S }
                (-1)^{ \card{ S \setminus I } } 
        }
        g
        =
        \del*{
            \sum_{ I = \emptyset }^{ S }
                (-1)^{ \card{ I } } 
        }
        g
        =
        \delta_{ \emptyset }^{ S }g.
        \tag*{\qedhere}
    \end{equation*}
\end{proof}
\begin{lemma} \label{lem::diff_B_diff_C_to_diff_BuC}
    Let $(G, +)$ be an Abelian group, $k \in \bN$, and 
    $
        \set*{
            x_{I}
        }_{
            I \subseteq [k]
        }
        \subseteq 
        G.
    $
    Then for every pairwise disjoint
    $
        A,B,C \subseteq [k]
    $
    we have
    \begin{equation*}
        \diff_{ I \upmapsto A }^{ B }
            \diff_{ J \upmapsto I }^{ C}
                x_J
        =
        \diff_{ I \upmapsto A }^{ B \cup C }
            x_I.
    \end{equation*}
    Moreover, 
    for all $i \in [k] \setminus (A \cup C)$,
    \begin{equation*}
        \diff_{ I \upmapsto A }^{ \set{i} \cup C }
            x_I
        =
        \diff_{ I \upmapsto A }^{ \set*{i} }
            \diff_{ J \upmapsto I }^{ C}
                x_J
        =
        \diff_{ J \upmapsto A \cup \set*{i} }^{ C}
                x_J
        - \diff_{ J \upmapsto A  }^{ C}
                x_J
        =
        \diff_{ J \upmapsto A }^{ C}
            \del{
                x_{J \cup \set*{i} } - x_{J  }
            }.
    \end{equation*}
    In particular, when $A = \emptyset$, the above equalities become the following ones:
    \begin{equation*}
		\diff_{ I = \emptyset }^{ \set{i} \cup C }
			 x_I
	        =
        \diff_{ I = \emptyset }^{ \set*{i} }
            \diff_{ J \upmapsto I }^{ C}
                x_J
        =
        \diff_{ J = \set*{i} }^{ C}
                x_J
        - \diff_{ J = \emptyset  }^{ C}
                x_J
        =
        \diff_{ J = \emptyset }^{ C}
            \del{
                x_{J \cup \set*{i} } - x_{J  }
            }.	    
    \end{equation*}
\end{lemma}
\begin{proof}
    By direct calculation,
    \begin{align*}
        \diff_{ I \upmapsto A }^{ B }
            \diff_{ J \upmapsto I }^{ C}
                x_J
    &=
        \sum_{ I = \emptyset }^{  B }
            (-1)^{ \card{B \setminus I} }
            \diff_{ J \upmapsto I \cup A }^{ C}
                x_J
    \\
	&=
        \sum_{ I = \emptyset }^{  B }
            (-1)^{ \card{ B \setminus I} }
            \sum_{ J = \emptyset }^{ C }
                (-1)^{ \card{ C \setminus J } }
                x_{J \cup I \cup A}
    	=
        \sum_{ 
            \substack{ 
                I \subseteq B \\ 
                J \subseteq C 
            } 
        } 
            (-1)^{ \card{ (B \cup C) \setminus (I \cup J) } }
            x_{ J \cup I \cup A }
    \\
    &\phantom{
        \sum_{ I = \emptyset }^{  B }
            (-1)^{ \card{ B \setminus I} }
            \sum_{ J = \emptyset }^{ C }
                (-1)^{ \card{ C \setminus J } }
                x_{J \cup I \cup A}
    	=    \,
    }
        =
        \sum_{ L = \emptyset }^{ B \cup C }
            (-1)^{ \card{ (B \cup C) \setminus L} }
            x_{ L \cup A }
        =
        \diff_{ L \upmapsto A }^{ B \cup C }
            x_L,
    \end{align*}
    where in the fourth equality we used the fact that
    since $B$ and $C$ are disjoint,
    the mapping
    $
        (I, J) \mapsto I \cup J
    $
    is~a~bijection from $2^B \times 2^C$ to $2^{B \cup C}$.
    The resulting equality also implies the first equality in the \say{moreover} part of the lemma. 
    The other equalities follow from the fact that,
    for all 
    $ 
    	i \in [k] \setminus (A \cup C),
    $ 
    we have
    \begin{equation*}
        \diff_{ J \upmapsto A \cup \set*{i} }^{ C }
            x_J 
        =
        \sum_{ J = \emptyset }^{ C }
            (-1)^{ \abs{ C \setminus J } }
            x_{ J \cup A \cup \set*{i} }
        =
        \diff_{ J \upmapsto A  }^{ C }
            x_{J \cup \set*{i} }. 
    \end{equation*}
    Finally, the equalities in the 
    ``in particular'' part
    follow from substituting $A = \emptyset$
    in the ones from the ``moreover'' part
    and using 
    \Cref{rem::diff_with_an_arrow}.
\end{proof}
%
%
%
%
\begin{lemma} \label{lem::partitioned_diff_gives_0}
    Let $k \in \bN$ and $\set*{C, D}$ be a partition of $[k]$. 
    Let $(G,+)$ be an Abelian group, $ g \in G$, and 
    $
        \set*{ v_{L} }_{ L \subseteq C}, 
        \set*{w_{J}}_{J \subseteq D} 
        \subseteq G
    $. 
    For $I \subseteq [k]$ let
    $
        x_{I}
        \coloneq 
        g
        + v_{ I \cap C}
        + w_{ I \cap D}.   
    $
    Suppose that $B \subseteq [k]$ 
    is~such that
    $B \cap C \ne \emptyset$ and 
    $B \cap D \ne \emptyset$.
    Then 
    \begin{equation*}
        \forall
        A \subseteq [k] \setminus B
        \qquad
        \diff_{ I \upmapsto A }^{ B }
            x_I
        =
        0.
    \end{equation*}
\end{lemma}
\begin{proof}
    By direct calculation,
    \begin{align*}
        \diff_{ I \upmapsto A }^{ B }
            x_I
        &=
        \diff_{ I \upmapsto A }^{ B }
            \del{
                g
                + v_{ I \cap C}
                + w_{ I \cap D}
            }
    \\
        &=
        \diff_{ I \upmapsto A }^{ B }
            g
        +
        \diff_{ I \upmapsto A }^{ B }
            \del{ v_{ I \cap C}}
        +
        \diff_{ I \upmapsto A }^{ B }
            \del{ v_{ I \cap D} }
    \\
        &=
        \diff_{ I \upmapsto A }^{ B }
            g
        +
        \diff_{ I \upmapsto A }^{ B \cap C }
            \del*{ 
                \diff_{ J \upmapsto I }^{ B \cap D }
                    v_{ J \cap C}
            }
        +
        \diff_{ I \upmapsto A }^{ B \cap D}
            \del*{
                \diff_{ J \upmapsto I }^{ B \cap C }
                    v_{ J \cap D}
            }
    \\
        &=
        0
        +
        \diff_{ I \upmapsto A }^{ B \cap C }
            0
        +
        \diff_{ I \upmapsto A }^{ B \cap D }
            0
    \\
        &=
        0,
    \end{align*}
    where in the third line we have used 
    \Cref{lem::diff_B_diff_C_to_diff_BuC}
    and in the fourth we have used \Cref{lem::diff_constant_term}.
\end{proof}
\begin{lemma}
\label{lem::difference_to_S_and_difference_from_S_into_single_one}
    Let $(G, +)$ be an Abelian group, $k \in \bN$, $S \subseteq [k]$, and 
    $
        \set*{
            x_{I,J}
        }_{
            I, J \subseteq [k]
        }
        \subseteq 
        G.
    $
    Define 
    $ \set*{ y_I }_{ I \subseteq [k] }$
    by~the~formula 
    $ y_I = x_{I \cap S, I \cup S}$
    for $I \subseteq [k]$.
    Then
    \begin{equation*}
        \diff_{ I = \emptyset }^{ S }
            \diff_{ J = S }^{ [k]  }
                x_{I, J}
        =
        \diff_{ I = \emptyset }^{ [k] }
            y_I.
    \end{equation*}
\end{lemma}
\begin{proof}
    Let us notice that 
    $L \mapsto (L \cap S, L \cup S)$ is a bijection from 
    $2^{[k] }$ to 
    $
    	2^{S} 
    	\times 
    	\setc*{ 
    		J \subseteq [k] 
    	}{ 
    		S \subseteq J 
	    }.
    $
    Also, for all $L \subseteq [k]$, we have
    \begin{equation}
    \label{proofeq::lem::difference_to_S_and_difference_from_S_into_single_one::product_of_minus_one_to_power}
        (-1)^{ \card{ S \setminus (L \cap S) } }
        (-1)^{ \card{ [k] \setminus (L \cup S) } }
        =
        (-1)^{ \card{ S \setminus L } + \card{ 
            ([k] \setminus S) \setminus L
        }}
        =
        (-1)^{ \card{ [k] \setminus L } }
    \end{equation}
    Therefore, we have
    \begin{align*}
        \diff_{ I = \emptyset }^{ S }
            \diff_{ J = S }^{ [k]  }
                x_{I, J}
    &=
        \sum_{ I = \emptyset }^{ S }
            (-1)^{ \card{ S \setminus I } }
            \sum_{ J = S }^{ [k] }
                (-1)^{ \card{ [k] \setminus J } }
                    x_{I,J} 
    \\
    &=
        \sum_{ I = \emptyset }^{ S }
        \sum_{ J = S }^{ [k] }
            (-1)^{ \card{ S \setminus I } }
            (-1)^{ \card{ [k] \setminus J } }
                    x_{I,J} 
    \\
    &=
        \sum_{ L = \emptyset }^{ [k] }
            (-1)^{ \card{ S \setminus (L \cap S) } }
            (-1)^{ \card{ [k] \setminus (L \cup S) } }
                    x_{L \cap S, L \cup S}
        ,
    \intertext{
    	which, by 
    	\eqref{proofeq::lem::difference_to_S_and_difference_from_S_into_single_one::product_of_minus_one_to_power},
    }
    &=
        \sum_{ L = \emptyset }^{ [k] }
            (-1)^{ \card{ [k] \setminus L } }
                    x_{L \cap S, L \cup S}
    	=
        \sum_{ L = \emptyset }^{ [k] }
            (-1)^{ \card{ [k] \setminus L } }
                y_L
        =
        \diff_{ L = \emptyset }^{ [k] }
            y_L,
    \end{align*}
    as claimed. 
\end{proof}
\begin{lemma} \label{lem::diag_diff_into_row_and_col_diffs}
    Let $(G, +)$ be an Abelian group, $k \in \bN$, and 
    $
        \set*{
            x_{I,J}
        }_{
            I, J \subseteq [k]
        }
        \subseteq 
        G.
    $
    Then
    \begin{equation*}
        \diff_{ I = \emptyset }^{ [k] }
            x_{I,I}
        =
        \sum_{
            S \subseteq [k]
        }
            \diff_{ I = \emptyset }^{ S }
            \diff_{ J = S }^{ [k]  }
                x_{I, J}.
    \end{equation*}
\end{lemma}
\begin{proof}
    First of all, let us notice that for all 
    $I \subseteq J \subseteq [k]$,
    the function 
    $S \mapsto  J \setminus S$ 
    is a bijection 
    from
    $
        \setc*{
            S \subseteq J
        }{
            I \subseteq S
        }
    $
    to
    $2^{J \setminus I}$ .
    Therefore, by \Cref{lem::alternating_sum_over_cardinalities},
    \begin{equation*}
        \sum_{ S = I }^{ J }
            (-1)^{ \card{ J \setminus S } }
        =
        \sum_{ S = \emptyset }^{ J \setminus I }
            (-1)^{ \card{ S } }
        =
        \delta_{ \emptyset }^{ J \setminus I  }
        =
        \delta_{I}^{J}.
    \end{equation*}
    Also, for all $I, S, J \subseteq [k]$ such that
    $I \subseteq S \subseteq J$, we have
    \begin{equation*}
        (-1)^{ \card{ S \setminus I } } 
        (-1)^{ \card{ [k] \setminus J } } 
        =
        (-1)^{ 
            \card{S} - \card{I}
            + \card{ [k] } - \card{J} 
        }
        =
        (-1)^{ 
            \card{[k]} - \card{I}
            -( \card{J} - \card{S} )
        }
        =
        (-1)^{ \card{ [k] \setminus I } }
        (-1)^{ \card{ J \setminus S} }.
    \end{equation*}
    In consequence,
    \begin{align*}
        \sum_{S \subseteq [k] }
            \diff_{I = \emptyset}^{ S }
                \diff_{J = S}^{ [k] }
                    x_{ I, J}
        &=
        \sum_{S \subseteq [k] }
            \sum_{ I = \emptyset }^{ S }
                (-1)^{ \card{ S \setminus I } } 
                \sum_{ J = S }^{ [k] }
                    (-1)^{ \card{ [k] \setminus J } } 
                        x_{ I, J }
    \\
        &=
        \sum_{ \substack{ (I, S, J) \colon \\ I \subseteq S \subseteq J \subseteq [k]} }
                (-1)^{ \card{ S \setminus I } } 
                (-1)^{ \card{ [k] \setminus J } } 
                        x_{ I, J }
    \\
        &=
        \sum_{ \substack{ (I, S, J) \colon \\ I \subseteq S \subseteq J \subseteq [k]} }
                (-1)^{ \card{ [k] \setminus I } } 
                (-1)^{ \card{ J \setminus S } } 
                        x_{ I, J }
    \\
        &=
        \sum_{ I = \emptyset }^{ [k] }
            (-1)^{ \card{ [k] \setminus I } } 
            \sum_{ J = I }^{ [k] }
                x_{I, J} 
                    \sum_{ S = I }^{ J }
                        (-1)^{ \card{ J \setminus S } }
    \\
        &=
        \sum_{ I = \emptyset }^{ [k] }
            (-1)^{ \card{ [k] \setminus I } } 
            \sum_{ J = I }^{ [k] }
                x_{I, J} 
                \delta_{I}^{J}
        =
        \sum_{ I = \emptyset }^{ [k] }
            (-1)^{ \card{ [k] \setminus I } } 
            x_{I, I} 
        =
        \diff_{I = \emptyset}^{ [k] } x_{I, I},
    \end{align*}
    as claimed.
\end{proof}
%
%
\begin{lemma}
\label{lem::characterization_of_diff_0_with_indicators}
    Let $(G, +)$ be an Abelian group, $k \in \bN$, 
    and 
    $
        \bx = \set*{ x_I }_{ I \subseteq [k] }
        \subseteq 
        G.
    $
    Suppose that for every 
    $
        x \in 
        \setc*{ x_I }{ I \subseteq [k] }
    $
    we have
    %
	$
        \diff_{ I = \emptyset }^{ [k] }
            \indicator{ \set*{x} }( x_I )
        =
        0.
    $
    %
    Then 
    %
	$
        \diff_{ I = \emptyset }^{ [k] }
            x_I
        =
        0.
    $
        %
    %
\end{lemma}
\begin{proof}
    Denote 
    $
        X \coloneq \setc*{ x_I }{ I \subseteq [k] }.
    $
    Then
    \begin{align*}
        \diff_{ I = \emptyset }^{ [k] }
            x_I
    &=
        \diff_{ I = \emptyset }^{ [k] }
            \indicator{X}(x_I)
            x_I
    \\
    &=
        \diff_{ I = \emptyset }^{ [k] }
            \sum_{ x \in X }
                \indicator{ \set*{x} }( x_I)
            x_I
    	=
        \diff_{ I = \emptyset }^{ [k] }
            \sum_{ x \in X }
                \del*{
                    \indicator{ \set*{x} }( x_I)
                    x
                }
    \\
    &
    \phantom{
    	=
		\diff_{ I = \emptyset }^{ [k] }
		\sum_{ x \in X }
		   \indicator{ \set*{x} }( x_I)
		x_I
		\:
    }
        =
        \sum_{ x \in X }
            \del*{
                \diff_{ I = \emptyset }^{ [k] }
                    \indicator{ \set*{x} }( x_I)
            }
            x
        =
        \sum_{ x \in X }
            0x
        =
        0,
    \end{align*}
    as claimed.
\end{proof}
%
%
\begin{lemma}
\label{lem::diff_[k]_to_diff_[m]_diff_[k-m]}
    Let $(G, +)$ be an Abelian group, $k \in \bN$,
    and 
    $
        \set*{ x_J }_{ J \subseteq [k] }
        \subseteq
        G.
    $
    Then for every $m \in [k]$ we have
    \begin{equation*}
        \diff_{ J = \emptyset }^{ [k] }
            x_J
        =
        \diff_{ L = \emptyset }^{ [m] }
            \diff_{ I = \emptyset }^{ [k-m] }
                x_{ L \cup (m+I) }
        =
        \diff_{ I = \emptyset }^{ [k-m] }
            \diff_{ L = \emptyset }^{ [m] }
                x_{ L \cup (m+I) }.
    \end{equation*}
\end{lemma}
\begin{proof}
    First of all, let us notice that,
    for all $L \subseteq [m]$
    and $I \subseteq [k-m]$,
    we have
    \begin{align}
        (-1)^{ \card*{ 
            [k] \setminus \del*{ L \cup (m+I) }
        }}
        =
        (-1)^{ \card*{ 
            \del*{ [m] \cup \del*{m+[k-m]}} 
            \setminus \del*{ L \cup (m+I) }
        }}
   	&
    \notag
    \\
        =
        (-1)^{ \card*{ [m] \setminus L }
            + \card*{ \del{ m + [k-m]} 
                \setminus \del*{ m+ I } }
        }
    &=
        (-1)^{ \card{ [m] \setminus L } }
        (-1)^{ \card{ [k-m] \setminus I } }.
        \label{proofeq::lem::diff_[k]_to_diff_[m]_diff_[k-m]::(-1)_to_powers}
    \end{align}
    Therefore, since    
    $
        (L, I) \mapsto L \cup (m+I)
    $
    is a bijection from
    $
        2^{ [m] } \times 2^{ [k-m] }
    $
    to 
    $ 2^{ [k] } $,
    \begin{align*}
        \diff_{ J = \emptyset }^{ [k] }
            x_J
        =
        \sum_{ J = \emptyset }^{ [k] }
            (-1)^{ \card{ [k] \setminus J } }
            x_J
        &=
        \sum_{ L = \emptyset}^{ [m] }
            \sum_{ I = \emptyset }^{ [k-m] }
                (-1)^{ \card{ 
                    [k] \setminus \del{ L \cup (m+I)} } }
                x_{ L \cup (m+I) }, 
    \intertext{
        which, by \eqref{proofeq::lem::diff_[k]_to_diff_[m]_diff_[k-m]::(-1)_to_powers},
    }
        &\hspace{-1cm}=
        \sum_{ L = \emptyset}^{ [m] }
            \sum_{ I = \emptyset }^{ [k-m] }
                (-1)^{ \card{ [m] \setminus L } }
                (-1)^{ \card{ [k-m] \setminus I } }
                x_{ L \cup (m+I) }
        =
        \diff_{ L = \emptyset }^{ [m] }
            \diff_{ I = \emptyset }^{ [k-m] }
                x_{ L \cup (m+I) },
    \end{align*}
    proving the first equality.
    The second one follows from the fact that
    \begin{align*}
        \diff_{ L = \emptyset }^{ [m] }
            \diff_{ I = \emptyset }^{ [k-m] }
                x_{ L \cup (m+I) }
        &=
        \sum_{ L = \emptyset}^{ [m] }
            \sum_{ I = \emptyset }^{ [k-m] }
                (-1)^{ \card{ [m] \setminus L } }
                (-1)^{ \card{ [k-m] \setminus I } }
                x_{ L \cup (m+I) }
    \\
        &=
        \sum_{ I = \emptyset }^{ [k-m] }
            \sum_{ L = \emptyset}^{ [m] }
                (-1)^{ \card{ [m] \setminus L } }
                (-1)^{ \card{ [k-m] \setminus I } }
                x_{ L \cup (m+I) }
        =
        \diff_{ I = \emptyset }^{ [k-m] }
            \diff_{ L = \emptyset }^{ [m] }
                x_{ L \cup (m+I) },
    \end{align*}    
    as claimed.
\end{proof}
%
%
%
\begin{lemma}
\label{lem::folding_hypercubes_in_groups}
    Let $( G, + )$ be an Abelian group, $k \in \bN$, $m \in [k]$, 
    and $i \in [m]$.
    Fix
    $
        \set*{ x_I }_{ I \subseteq [m] }
        \subseteq 
        G
    $
    and define 
    $
        \set*{ y_J }_{ J \subseteq [k] }
    $
    by 
    \begin{equation*}
        \forall I \subseteq [m]
        \quad 
        \forall L \subseteq [k] \setminus [m]
        \qquad 
        y_{ I \cup L }
        \coloneqq 
        \begin{cases}
            x_{I} &\text{if } \abs{ L } \text{ is even}, \\
            x_{ I \div \set*{i} } &\text{otherwise}.
        \end{cases}
    \end{equation*}
    Then for all disjoint $A, S \subseteq [k]$,
    \begin{equation*}
        \diff_{ J \upmapsto A  }^{ S}
            y_J
        =
        \eta_{A}^{S }
        2^{
            \card{
                S \cap \phi^{-1}\sbr{ \set*{i} } 
            }
            - \card{
                \phi \sbr{ S } \cap \set*{i} 
            }
        }
        \diff_{
            I \upmapsto 
            \Psi( A)
            \setminus \phi \sbr{ S }
        }^{
            \phi \sbr{ S  } 
        }
            x_I,
    \end{equation*}
    where $ \abs*{ \eta_{A}^{S } } = 1$, 
    $\phi \colon [k] \to [m]$ is defined by
    \begin{align*}
        \forall j \in [k]
        \qquad 
        \phi(j)
        &\coloneq
        \begin{cases}
            j &\text{if } j \in [m], \\ 
            i &\text{otherwise},
        \end{cases} 
    \intertext{
        and 
        $\Psi \colon 2^{[k]} \to 2^{[m]}$
        is defined by
    }
        \forall I \subseteq [m] 
        \quad \forall L \subseteq [k] \setminus [m]
        \qquad
        \Psi(I \cup L)
        &\coloneq
        \begin{cases}
            I &\text{if } \abs{L} \text{ is even,} \\ 
            I \div \set*{i} &\text{otherwise}.
        \end{cases} 
    \end{align*}
    In particular,
    \begin{equation*}
        \diff_{J = \emptyset}^{ [k] }
            y_J
        =
        (-1)^{ k-m }
        2^{k-m}
        \diff_{I = \emptyset}^{ [m] }
            x_I.
    \end{equation*}
\end{lemma}
\begin{proof}
    Fix $S \subseteq [k]$ and $A \subseteq [k] \setminus S$.
    Throughout the proof, we will use the following notation:
    \begin{equation*}
        S_m \coloneq S \cap [m],
        \quad  
        A_m \coloneq A \cap [m],
        \quad 
        S_k \coloneq S \setminus [m],
        \quad \text{and} \quad
        A_k \coloneq A \setminus [m]
    \end{equation*}
    We will prove the statement by considering several cases.
    \begin{enumerate}[label={Case \Roman*:}, ref={Case \Roman*}, leftmargin=1em, itemindent=2.5em]
        \item
        \label{proofitem::lem::folding_hypercubes_in_groups::Case::S_k_is_empty}
        $ S_k = \emptyset $.

        In this case we have 
        $
            \phi[ S ]
            =
            \phi[ S_m ]
            =
            S_m.
        $
        Moreover, if $i \in S$, then
        \begin{equation*}
            \card{
                S \cap \phi^{-1}[ \set*{i} ]
            }
            - \card{
                \phi[S] \cap \set*{i} 
            }
            =
            \card{
                \set*{i}
            }
            - \card{
                \set*{i} 
            }
            =
            0,
        \end{equation*}
        and if $i \notin S$, then also
        \begin{equation*}
            \card{
                S \cap \phi^{-1}[ \set*{i} ]
            }
            - \card{
                \phi[S] \cap \set*{i} 
            }
            =
            \card{
                \emptyset
            }
            - \card{
                \emptyset
            }
            =
            0.
        \end{equation*}
        %
        \begin{enumerate}[label={Subcase \arabic*:}, ref={Subcase \arabic*}, leftmargin=1em, itemindent=3.5em]
            \item
            \label{proofitem::lem::folding_hypercubes_in_groups::Case::S_k_is_empty::Subcase::A_k_is_even}
            $ \card{ A_k }$ is even.

            In this subcase, we have
            \begin{equation*}
                \diff_{ J \upmapsto A }^{ S }
                    y_J
                {=}
                \diff_{ J \upmapsto A_m \cup A_k }^{ S_m }
                    y_J
                {=}
                \sum_{ J = \emptyset }^{ S_m }
                    (-1)^{ \card{ S_m \setminus J } }
                    y_{ J \cup A_m \cup A_k }
                {=}
                \sum_{ J = \emptyset }^{ S_m }
                    (-1)^{ \card{ S_m \setminus J } }
                    x_{ J \cup A_m}
                {=}
                \diff_{ J \upmapsto A_m }^{ S_m }
                    x_J,
            \end{equation*}
            which matches the desired form as
            $
                \Psi ( A )
                =
                \Psi( A_m \cup A_k )
                =
                A_m
            $
            and $A_m$ is disjoint with 
            $\phi[S] = S_m$.
            \item 
            $ \card{ A_k }$ is odd and $i \notin S$.

            This time, we have 
            $
                y_{ J \cup A_m \cup A_k }
                =
                x_{ (J \cup A_m) \div \set*{i} }
                =
                x_{ J \cup (A_m \div \set*{i}) }
            $
            for all $J \subseteq S_m$, hence
            \begin{align*}
                \diff_{ J \upmapsto A }^{ S }
                    y_J
                =
                \diff_{ J \upmapsto A_m \cup A_k }^{ S_m }
                    y_J
            &=
                \sum_{ J = \emptyset }^{ S_m }
                    (-1)^{ \card{ S_m \setminus J } }
                    y_{ J \cup A_m \cup A_k }
            \\
            &=
                \sum_{ J = \emptyset }^{ S_m }
                    (-1)^{ \card{ S_m \setminus J } }
                    x_{ J \cup (A_m \div \set*{i}) }
                =
                \diff_{ J \upmapsto A_m \div \set*{i} }^{ S_m }
                    y_J.
            \end{align*}
            This matches our form since 
            $
                \Psi(A_m \cup A_k)
                =
                A_m \div \set*{i}
            $
            and 
            $ A_m \div \set*{i} $
            is disjoint with
            $
                \phi[S] = S_m.
            $
            \item 
            $ \card{ A_k }$ is odd and $i \in S$.
            
            If $i \in S (= S_m)$, then for all 
            $J \subseteq S_m \setminus \set*{i}$,
            \begin{equation*}
                \diff_{I \upmapsto J \cup A }^{\set*{i}}
                    y_I
                =
                y_{ J \cup A_m \cup A_k \cup  \set*{i}}
                - y_{ J \cup A_m \cup A_k }
                =
                x_{ J \cup A_m}
                - x_{J \cup A_m \cup \set*{i}}
                =
                - \diff_{I \upmapsto J \cup A_m}^{ \set*{i} }
                    x_I.
            \end{equation*}
            Hence, by 
            \Cref{lem::diff_B_diff_C_to_diff_BuC},
            \begin{align*}
                \diff_{ J \upmapsto A }^{ S }
                    y_J
                =
                \diff_{ J \upmapsto A }^{ S_m }
                    y_J
                =
                \diff_{ J \upmapsto A }^{ S_m \setminus \set*{i} }
                    \diff_{ I \upmapsto J}^{\set*{i} }
                        y_I
            &=
                \diff_{ J = \emptyset }^{ S_m \setminus \set*{i} }
                    \diff_{ I \upmapsto J \cup A}^{\set*{i} }
                        y_I
            \\
            &=
                \diff_{ J = \emptyset }^{ S_m \setminus \set*{i} }
                    \del*{ 
                        -
                        \diff_{ I \upmapsto J \cup A_m}^{\set*{i} }
                            x_I
                    }
                =
                -
                \diff_{ J \upmapsto A_m }^{ S_m }
                    x_J.
            \end{align*}
            Since $i \in S_m$, and $A_m$ and $S_m$ are disjoint,
            \begin{equation*}
                \Psi(A) \setminus \phi[S]
                =
                (A_m \div \set*{i}) \setminus S_m
                =
                (A_m \cup \set*{i}) \setminus \set*{i}
                =
                A_m,
            \end{equation*}
            so we can get the desired form by setting 
            $ \eta_A^S = -1$.
        \end{enumerate}
        \item \label{proofitem::lem::folding_hypercubes_in_groups::Case::S_m_is_empty}
        $S_m = \emptyset$ and $S_k \ne \emptyset$.

        In this case we have 
        $
            \phi[ S ]
            =
            \phi[ S_k ]
            =
            \set*{i},
        $ 
        \begin{equation*}
            \card{
                S \cap \phi^{-1}[ \set*{i} ]
            }
            - \card{
                \phi[S] \cap \set*{i} 
            }
            =
            \card{ S_k } - \card{ \set*{i} }
            =
            \card{ S_k } - 1,
        \end{equation*}
        and since 
        $\Psi(A) = A_m$ or $\Psi(A) = A_m \div \set*{i}$,
        we also have
        \begin{equation*}
            \Psi(A) \setminus \phi[S]
            =
            \Psi(A) \setminus \set*{i}
            =
            A_m \setminus \set*{i}.
        \end{equation*}
        %
		Therefore,
        \begin{align*}
            \diff_{J \upmapsto A}^{S}
                y_J
        &=
            \sum_{ J = \emptyset }^{ S_k }
                (-1)^{ \card{S_k \setminus J } }
                y_{J \cup A_m \cup A_k}
        \\
        &=
            (-1)^{\card{S_k \cup A_k}}
            \sum_{ J = \emptyset }^{ S_k }
                (-1)^{ \card{A_k \cup J } }
                y_{J \cup A_m \cup A_k}
        \\
       	&=
        	(-1)^{\card{S_k \cup A_k}}
        	\sum_{ J = \emptyset }^{ S_k }
        	(-1)^{ \card{A_k \cup J } }
        	y_{A_m \cup J \cup  A_k}
       	\\
        &=
            (-1)^{\card{S_k \cup A_k}}
            \del*{
	                \sum_{ 
	                    \substack{
	                        J = \emptyset 
	                    \\ 
	                        \card{A_k \cup J} \text{ is even}
	                    }
	                }^{ S_k }
	                    y_{A_m \cup J \cup  A_k}
                -
	                \sum_{ 
	                    \substack{
	                        J = \emptyset 
	                    \\ 
	                        \card{A_k \cup J} \text{ is odd}
	                    }
	                }^{ S_k }
	                    y_{A_m \cup J \cup  A_k}
            },
        \intertext{
        	which, since $S_k \ne \emptyset$, by \Cref{lem::between_A_and_B_cardinalities},      	
        }
        &=
            (-1)^{\card{S_k \cup A_k}}
            2^{ \card{S_k} - 1}
            \del*{
                x_{A_m}
                - x_{A_m \div \set*{i}}
            }.
        \end{align*}
        Next, notice that if $i \notin A_m$, then
        \begin{equation*}
            x_{A_m} - x_{A_m \div \set*{i}}
            =
            x_{A_m} - x_{A_m \cup \set*{i}}
            =
            - \diff_{I \upmapsto A_m}^{ \set*{i} }
                x_I
            =
            (-1)^{ \card{ \set*{i} \setminus A_m } }
            \diff_{I \upmapsto A_m \setminus \set*{i} }^{ \set*{i} }
                x_I,
        \end{equation*}
        and if $i \in A_m$, then also
        \begin{equation*}
            x_{A_m} - x_{A_m \div \set*{i}}
            =
            x_{A_m} - x_{A_m \setminus \set*{i}}
            =
            \diff_{I \upmapsto A_m \setminus \set*{i}}^{ \set*{i }}
                x_I
            =
            (-1)^{ \card{ \set*{i} \setminus A_m } }
            \diff_{I \upmapsto A_m \setminus \set*{i} }^{ \set*{i} }
                x_I.
        \end{equation*}
        In consequence, 
        \begin{align*}
            \diff_{J \upmapsto A}^{S}
                y_J
        &=
            (-1)^{\card{S_k \cup A_k}}
            2^{ \card{S_k} - 1}
            \del*{
                x_{A_m}
                - x_{A_m \div \set*{i}}
            }
        \\
        &=
            (-1)^{
                \card{S_k \cup A_k} 
                + \card{ \set*{i} \setminus A_m }
            }
            2^{ \card{S_k} - 1}
            \diff_{I \upmapsto A_m \setminus \set*{i} }^{ \set*{i} }
                x_I.
        \end{align*}
        It remains to note that we get the desired form by 
        setting 
        $
            \eta_A^S
            =
            (-1)^{
                \card{S_k \cup A_k} 
                + \card{ \set*{i} \setminus A_m }
            }.
        $
        %
        \item \label{proofitem::lem::folding_hypercubes_in_groups::Case::both_S_m_and_S_k_are_nonempty}
        $S_m \ne \emptyset$ and $S_k \ne \emptyset$.

        This time, using \Cref{lem::diff_B_diff_C_to_diff_BuC}, 
        we have 
        \begin{equation*}
            \diff_{ J \upmapsto A}^{S}
                y_J
            =
            \diff_{ J \upmapsto A}^{S_m}
                \diff_{ I \upmapsto J}^{S_k}
                    y_I
            =
            \diff_{ J = \emptyset}^{S_m}
                \diff_{ I \upmapsto J \cup A}^{S_k}
                    y_I
            =
            \diff_{ J = \emptyset}^{S_m}
                \del*{
                    \eta_{J \cup A}^{S \setminus S_m}
                    2^{ \card{S_k} - 1}
                    \diff_{I \upmapsto (J \cup A_m) \setminus \set*{i} }^{ \set*{i} }
                        x_I
                },
        \end{equation*}
        where the final equality follows by
        \ref{proofitem::lem::folding_hypercubes_in_groups::Case::S_m_is_empty}
        with 
        $ S_k $ and $J \cup A $ playing the roles of
        $ S $ and $A$, respectively.
        Hence, we have
        $
            \eta_{J\cup A}^{S \setminus S_m}
            =
            (-1)^{
                \card{S_k \cup A_k} 
                + \card{ \set*{i} \setminus (J \cup A_m) }
            }
        $
        for all $J \subseteq S_m$.
        Since $S_k \ne \emptyset$, in this case we have
        $
            \phi[S] = S_m \cup \set*{i}.
        $
        Thus,
        we have
        $
            \Psi(A) \setminus \phi[S]
            =
            A_m \setminus \set*{i}.
        $
        \begin{enumerate}[label={Subcase \arabic*:}, ref={Subcase \arabic*}, leftmargin=1em, itemindent=3.5em]
            \item 
            \label{proofitem::lem::folding_hypercubes_in_groups::Case::both_S_m_and_S_k_are_nonempty::subcase::i_not_in_Sm}
            $i \notin S_m$.
            
            If $i \notin S_m$, then 
            $
                \eta_{J\cup A}^{S \setminus S_m}
                =
                \eta_{A}^{S \setminus S_m}.
            $
            Hence,
            by \Cref{lem::diff_B_diff_C_to_diff_BuC},
            \begin{align*}
                \diff_{ J = \emptyset}^{S_m}
                    \del*{
                        \eta_{J \cup A}^{S \setminus S_m}
                        2^{ \card{S_k} - 1}
                        \diff_{I \upmapsto (J \cup A_m) \setminus \set*{i} }^{ \set*{i} }
                            x_I
                    }
            &=
                \diff_{ J = \emptyset}^{S_m}
                    \del*{
                        \eta_{A}^{S \setminus S_m}
                        2^{ \card{S_k} - 1}
                        \diff_{I \upmapsto J \cup (A_m \setminus \set*{i}) }^{ \set*{i} }
                            x_I
                    }
            \\
            &=
                \eta_{A}^{S \setminus S_m}
                2^{ \card{S_k} - 1}
                \diff_{ J = \emptyset}^{S_m}
                    \del*{
                        \diff_{I \upmapsto J \cup (A_m \setminus \set*{i}) }^{ \set*{i} }
                           x_I
                    }
            \\
            &=
                \eta_{A}^{S \setminus S_m}
                2^{ \card{S_k} - 1}
                \diff_{ J \upmapsto A_m \setminus \set*{i} }^{S_m}
                    \diff_{I \upmapsto J }^{ \set*{i} }
                        x_I
            \\
            &=
                \eta_{A}^{S \setminus S_m}
                2^{ \card{S_k} - 1}
                \diff_{ J \upmapsto A_m \setminus \set*{i} }^{S_m \cup \set*{i}}
                    x_J.
            \end{align*}
            This gives us the desired form by setting 
            $
                \eta_A^S = \eta_{A}^{S \setminus S_m}.
            $
            Indeed, since $i \notin S$ and $S_k \ne \emptyset$, we have
            \begin{equation*}
                \card{
                    S \cap \phi^{-1}[ \set*{i} ]
                }
                - \card{
                    \phi[S] \cap \set*{i} 
                }
                =
                \card{ S_k } - \card{ \set*{i} }
                =
                \card{ S_k } - 1,
            \end{equation*}
            as needed.
            \item 
            \label{proofitem::lem::folding_hypercubes_in_groups::Case::both_S_m_and_S_k_are_nonempty::subcase::i_in_Sm}
            $i \in S_m$.

            If $i \in S_m$, then $i \notin A_m$. 
            Hence,
            \begin{equation*}
                \eta_{\set*{i} \cup A}^{S \setminus S_m}
                =
                (-1)^{
                    \card{S_k \cup A_k} 
                    + \card{ \set*{i} \setminus (\set*{i}  \cup A_m) }
                }
                =
                (-1)^{
                    \card{S_k \cup A_k} 
                    + 0
                }
                =
                (-1)^{
                    \card{S_k \cup A_k} 
                }
            \end{equation*}
            and 
            \begin{equation*}
                \eta_{\emptyset \cup A}^{S \setminus S_m}
                =
                (-1)^{
                    \card{S_k \cup A_k} 
                    + \card{ \set*{i} \setminus (\emptyset  \cup A_m) }
                }
                =
                (-1)^{
                    \card{S_k \cup A_k} 
                    + 1
                }
                =
                -
                (-1)^{
                    \card{S_k \cup A_k} 
                }.
            \end{equation*}
            Moreover, we have
            \begin{equation*}
                \del{ \set*{i} \cup A_m } \setminus \set*{i}
                =
                A_m
                \quad\text{and}\quad 
                \del{ 
                    \emptyset \cup A_m 
                } \setminus \set*{i}
                =
                A_m
            \end{equation*}
            Therefore, 
            by the \say{moreover} part of \Cref{lem::diff_B_diff_C_to_diff_BuC},
            \begin{align*}
                \diff_{ J \upmapsto A}^{S}
                y_J
            &=
                \diff_{ J \upmapsto A \cup \set*{i}}^{S \setminus \set*{i}}
                    y_J
                -
                \diff_{ J \upmapsto A  }^{S \setminus \set*{i}}
                    y_J,
            \intertext{
                which, if $S_m = \set*{i}$, by \ref{proofitem::lem::folding_hypercubes_in_groups::Case::S_m_is_empty},
                and if $ S_m \ne \set*{i}$, by 
                \ref{proofitem::lem::folding_hypercubes_in_groups::Case::both_S_m_and_S_k_are_nonempty::subcase::i_not_in_Sm}
                of \ref{proofitem::lem::folding_hypercubes_in_groups::Case::both_S_m_and_S_k_are_nonempty},
            }
                &=
                \eta_{\set*{i} \cup A}^{S \setminus S_m}
                2^{ \card{S_k} - 1}
                \diff_{ J \upmapsto (\set*{i} \cup A_m) \setminus \set*{i} }^{S_m \cup \set*{i}}
                    x_J
                -
                \eta_{ \emptyset \cup A }^{S \setminus S_m}
                2^{ \card{S_k} - 1}
                \diff_{ J \upmapsto (\emptyset \cup A_m)  \setminus \set*{i} }^{S_m \cup \set*{i}}
                    x_J
            \\
                &=
                (-1)^{
                    \card{S_k \cup A_k} 
                }
                2^{ \card{S_k} - 1}
                \diff_{ J \upmapsto A_m }^{S_m \cup \set*{i}}
                    x_J
                +
                (-1)^{
                    \card{S_k \cup A_k} 
                }
                2^{ \card{S_k} - 1}
                \diff_{ J \upmapsto A_m }^{S_m \cup \set*{i}}
                    x_J
            \\
                &=
                (-1)^{
                    \card{S_k \cup A_k} 
                }
                2^{ \card{S_k}}
                \diff_{ J \upmapsto A_m }^{S_m \cup \set*{i}}
                    x_J.
            \end{align*}
            This matches the desired form by setting 
            $ 
                \eta_A^S 
                = (-1)^{
                    \card{S_k \cup A_k} 
                }.
            $            
            Indeed, since $i \notin A_m$, we have
            $
                \Psi(A) \setminus \phi[S]
                =
                A_m \setminus \set*{i}
                =
                A_m.
            $
            Also, 
            \begin{equation*}
                \card{
                    S \cap \phi^{-1}[ \set*{i} ]
                }
                - \card{
                    \phi[S] \cap \set*{i} 
                }
                =
                \card{ S_k \cup \set*{i} } - \card{ \set*{i} }
                =
                \card{ S_k } + 1 - 1
                =
                \card{ S_k },
            \end{equation*}
            as needed.
        \end{enumerate}
    \end{enumerate}
    It remains to prove the ``in particular'' part.
    Let $A = \emptyset$ and $S = [k] $.
    Then 
    $ A_k = \emptyset $, 
    $A_m = \emptyset$,
    $S_k = [k] \setminus [m]$, 
    and 
    $S_m = [m]$.
    We have two possibilities:
    \begin{itemize}[leftmargin=0.5cm]
        \item
        If $k = m$, then 
        $ A_k = \emptyset $,
        $S_k = [k] \setminus [m] = \emptyset$,
        and
        $
        	S_m = [m] = [k] = S.
        $
		Hence, by 
        \ref{proofitem::lem::folding_hypercubes_in_groups::Case::S_k_is_empty::Subcase::A_k_is_even}
        of
        \ref{proofitem::lem::folding_hypercubes_in_groups::Case::S_k_is_empty},
        \begin{equation*}
            \diff_{J = \emptyset}^{ [k] }
                y_J
            =
            \diff_{J \upmapsto \emptyset}^{ [k] }
                y_J 
    		=
			\diff_{J \upmapsto A }^{ S }
				y_J     
    		=
			\diff_{J \upmapsto A_m }^{ S_m }
				x_J   
            =
            \diff_{J \upmapsto \emptyset}^{ [m] }
                x_J
            =
            \diff_{J = \emptyset}^{ [m] }
                x_J
            =
            (-1)^{k-m}
            2^{k-m}
            \diff_{I = \emptyset}^{ [m] }
                x_I.
        \end{equation*}
        \item
        If $k > m$, then 
        $ S_k = [k] \setminus [m] \ne \emptyset$
        and
        $ i \in [m] = S_m $.
		Hence, by
        \ref{proofitem::lem::folding_hypercubes_in_groups::Case::both_S_m_and_S_k_are_nonempty::subcase::i_in_Sm}
        of
        \ref{proofitem::lem::folding_hypercubes_in_groups::Case::both_S_m_and_S_k_are_nonempty},
        \begin{align*}
            \diff_{J = \emptyset}^{ [k] }
                y_J
            =
            \diff_{J \upmapsto \emptyset}^{ [k] }
                y_J 
    		=
			\diff_{J \upmapsto A }^{ S }
				y_J
    	&=
            (-1)^{ \card{ S_k  \cup A_k } }
            2^{ \card{ S_k } }
			\diff_{J \upmapsto A_m }^{ S_m }
				x_J        
        \\
        &=
            (-1)^{ \card{ ([k] \setminus [m])  \cup \emptyset } }
            2^{ \card{ [k] \setminus [m] } }
            \diff_{I = \emptyset}^{ [m] }
                x_I
            =
            (-1)^{k-m}
            2^{k-m}
            \diff_{I = \emptyset}^{ [m] }
                x_I.
        \end{align*}
    \end{itemize}
    In both cases, the resulting expression matches the desired form.
\end{proof}
The following proposition can be thought of as a motivation for the next definition.
Function $ c_{\bx}$ that will be introduced therein will also be used in 
\Cref{prop::higher_order_fundamental_theorem_of_calculus_inequality},
which could be considered a motivation for our \say{multipoint} characterization of higher-order Sobolev spaces.
\begin{proposition}
\label{prop::interpolation_estimate}
	Let $\del*{ V, \normAlone }$ be a normed space,
	$k \in \bN$,
	and 
	$ 
		\bx = \set*{ x_I }_{ I \subseteq [k] }
		\subseteq 
		V.
	$
	Let
	$c_{\bx} \colon [0,1]^k \to V$
	be a
	$k$-affine function from $[0,1]^k$ to $V$ such that
	\begin{equation*}
		\forall I \subseteq [k]
	\qquad 
		c_{\bx}\del*{ \indicator{I}(1), \, \ldots, \, \indicator{I}(k) }
		= 
		x_I.
	\end{equation*}
	(%
		We will call such a function a 
		\emph{$k$-affine interpolant of $\bx$}.%
	)
	
	\noindent Then, the following formulas hold:%
	\footnote{Note that the formulas imply the uniqueness of such $c_{\bx}$.}
	\begin{equation}
	\label{eqprop::prop::interpolation_estimate::deifnitions_of_c_x}
		\forall t \in [0,1]^k
		\qquad 
		c_{\bx}(t) 
		=
		\sum_{ I = \emptyset}^{ [k] }
			\alpha_I(t)\beta_{ [k] \setminus I }(t)x_I
		=
		\sum_{ I = \emptyset}^{ [k] }
			\alpha_I(t)
			\diff_{ A = \emptyset }^{ I }
				x_A,
	\end{equation}
	where, for all 
	$I \subseteq \bN$ 
	and 
	$ m \in \bN $ 
	such that $I \subseteq [m]$, 
	functions 
	$ \alpha_I, \beta_I \colon [0,1]^m \to \bR$
	are defined by the formulas
	\begin{equation*}
		\forall t = \del*{ t_1, \, \ldots, \, t_m } \in [0,1]^m
		\qquad 
		\alpha_I(t)
		\coloneq 
		\prod_{ i \in I }
		t_i
		\quad \text{and} \quad 
		\beta_I(t)
		\coloneq 
		\prod_{ i \in I }
		\del*{1-t_i}.
	\end{equation*}
	Moreover, for all 
	$S \subseteq [k]$ 
	and 
	$t \in \intoo{0,1}^k$, 
	we have
	\begin{equation*}
		\norm*{ \partial^S c_{\bx}(t) }
		\le 
		\sum_{ A = \emptyset }^{ [k] \setminus S }
			\norm*{
				\diff_{ I \upmapsto A }^{ S }
					x_I
			},
	\quad \text{where we use the notation} \quad
		\partial^S
		\coloneq 
		\frac{ 
			\partial^{ \card{S} } 
		}{
			\prod_{ i \in S }
				\partial_{ x_i } 	
		}.
	\end{equation*}
\end{proposition}
\begin{proof}
	We will begin the proof by showing that the two formulas in 
	\eqref{eqprop::prop::interpolation_estimate::deifnitions_of_c_x}
	hold using induction over $k \in \bN$.
	They are correct for $k=1$ since in this case
	\begin{equation*}
		c_{\bx}(t)
		=
		(1-t)x_{\emptyset} + t x_{\set*{1} }
		=
		\alpha_{\emptyset}(t)\beta_{[1] \setminus \emptyset}(t) x_{\emptyset}
		+ \alpha_{ \set*{1} }(t) \beta_{ [1] \setminus \set*{1} }(t) x_{\set*{1}}
		=
		\sum_{I=\emptyset}^{ [1] }
			\alpha_{I}(t) \beta_{ [k] \setminus I }(t) x_I
	\end{equation*}
	and 
	\begin{equation*}
		c_{\bx}(t)
		=
		(1-t)x_{\emptyset} + t x_{\set*{1} }
		=
		x_{\emptyset} + t\del*{ x_{\set*{1} } - x_{ \emptyset } }
		=
		\alpha_{\emptyset}(t) x_{\emptyset}
		+ \alpha_{ \set*{1} }(t) 
			\diff_{ A = \emptyset }^{ \set*{1} }
				x_A
		=
		\sum_{I=\emptyset}^{ [1] }
			\alpha_{I}(t) 
			\diff_{ A = \emptyset }^{ I }
				x_A.
	\end{equation*}
	Now, fix $k \in \bN$ such that $k \ge 2$ and suppose that the formulas are correct for $k-1$.
	Let us note that, writing 
	$ t' = \del*{ t_1, \, \ldots, \, t_{k-1} }$,
	the functions
	\begin{equation*}
		[0,1]^{k-1} 
		\ni 
		t'
		\mapsto
		c_{\bx}\del*{ t', 0}
	\quad \text{and} \quad 
		[0,1]^{k-1} 
		\ni 
		t'
		\mapsto
		c_{\bx}\del*{ t',1 }		
	\end{equation*}
	%
%
	are the $(k-1)$-affine interpolants of
	$
		\bx' 
		\coloneq 
		\set*{ 
			x_I 
		}_{ 
			I \subseteq [k-1 ] 
		}
	$
	and 
	$
		\bx'' 
		\coloneq 
		\set*{ 
			x_{I \cup \set*{k} } 
		}_{ 
			I \subseteq [k-1 ] 
		},
	$
	respectively. 
	Therefore,
	\begin{align*}
		c_{\bx}(t)
	&=
		(1-t_k) c_{\bx}\del*{t', 0 }
		+ t_k c_{\bx}\del*{t', 1 }
	\\
	&=
		(1-t_k)c_{\bx'} \del*{t'}
		+ t_k c_{\bx''} \del*{t'}
	\\
	&=
		(1-t_k) 
		\sum_{I = \emptyset}^{ [k-1] }
			\alpha_I\del*{t'} \beta_{[k-1] \setminus I}\del*{t'}x_{I}
		+
		t_k
		\sum_{I = \emptyset}^{ [k-1] }
			\alpha_I\del*{t'} \beta_{[k-1] \setminus I}\del*{t'}x_{I \cup \set*{k}}
	\\
	&=
		\sum_{I = \emptyset}^{ [k-1] }
			\alpha_I\del*{t} \beta_{[k] \setminus I}\del*{t}x_{I}
		+
		\sum_{I = \emptyset}^{ [k-1] }
			\alpha_{I\cup \set*{k}}\del*{t} 
			\beta_{[k] \setminus \del*{I \cup \set*{k} }}\del*{t}
				x_{I \cup \set*{k}},
	\intertext{
		which, since 
		$ I \mapsto I \cup \set*{k} $
		is a bijection from 
		$ 2^{ [k-1] } $
		to
		$ 
			\setc*{ 
				J \subseteq [k] 
			}{
				\set*{k} \subseteq J
			},
		$
	}
	&=
		\sum_{I = \emptyset}^{ [k-1] }
			\alpha_I\del*{t} \beta_{[k] \setminus I}\del*{t}x_{I}
		+
		\sum_{I = \set*{k}}^{ [k] }
			\alpha_I\del*{t} \beta_{[k] \setminus I}\del*{t}x_{I}
		=
		\sum_{I = \emptyset}^{ [k] }
			\alpha_I\del*{t} \beta_{[k] \setminus I}\del*{t}x_{I}.
	\end{align*}
	Similarly,
	\begin{align*}
		c_{\bx}(t)
	&=
		(1-t_k)c_{\bx'} \del*{t'}
		+ t_k c_{\bx''} \del*{t'}
	\\
	&=
		(1-t_k) 
		\sum_{I = \emptyset}^{ [k-1] }
			\alpha_I\del*{t'} 
			\diff_{ A = \emptyset }^{ I } 
				x_A
		+
		t_k
		\sum_{I = \emptyset}^{ [k-1] }
			\alpha_I\del*{t'} 
			\diff_{ A = \emptyset }^{ I } 
				x_{A \cup \set*{k} }
	\\
	&=
		\sum_{I = \emptyset}^{ [k-1] }
			\alpha_I\del*{t} 
			\diff_{ A = \emptyset }^{ I } 
				x_A
		+
		\sum_{I = \emptyset}^{ [k-1] }
			\alpha_{I \cup \set*{k} }\del*{t}
			\diff_{ A = \emptyset }^{ I }
				\del*{ 
					x_{ A \cup \set*{k} }
					-
					x_{ A }
				},
	\intertext{
		which, by the 
		\say{in particular} part of
		\Cref{lem::diff_B_diff_C_to_diff_BuC},
	}
	&=
		\sum_{I = \emptyset}^{ [k-1] }
			\alpha_I\del*{t} 
			\diff_{ A = \emptyset }^{ I } 
				x_A
		+
		\sum_{I = \emptyset}^{ [k-1] }
			\alpha_{I \cup \set*{k} }\del*{t}
			\diff_{ A = \emptyset }^{ I \cup \set*{k} }
				x_A
	\\
	&=
		\sum_{I = \emptyset}^{ [k-1] }
			\alpha_I\del*{t} 
			\diff_{ A = \emptyset }^{ I } 
				x_A
		+
		\sum_{I = \set*{k}}^{ [k] }
			\alpha_{I }\del*{t}
			\diff_{ A = \emptyset }^{ I }
				x_A
		=
		\sum_{I = \emptyset}^{ [k] }
			\alpha_I\del*{t} 
			\diff_{ A = \emptyset }^{ I } 
				x_A.
	\end{align*}
	This proves the inductive step. 
	Hence, by induction, our formulas are correct for all $k \in \bN$.

	Next, let us prove the estimate for $\partial^S c_{\bx}$.
	Fix $ S \subseteq [k]$ and define
	$ 
		\by = \set*{ y_A }_{ A \subseteq [k] }
	$
	by the formula
	\begin{equation*}
		\forall A \subseteq [k]
	\qquad 
		y_A \coloneq 
		\begin{cases}
			\diff_{ J \upmapsto A }^{ S } 
				x_J
			&\text{if } A \subseteq [k] \setminus S, \\ 
			0 
			&\text{otherwise}.
		\end{cases}
	\end{equation*}
	Let us note that 
	\begin{equation*}
		\forall I \subseteq [k] 
	\quad 
		\forall t \in \intoo{0,1}^k
	\qquad 
		\partial^S \alpha_I(t)
		=
		\begin{cases}
			\alpha_{I \setminus S}(t) &\text{if } S \subseteq I, \\ 
			0 &\text{otherwise}.
		\end{cases}
	\end{equation*}
	Also, define the function
	$ \pi_{[k] \setminus S} \colon [0,1]^k \to [0,1]^k$
	by the formula
	\begin{equation*}
		\forall t = \del*{ t_1, \, \ldots, \, t_k} \in [0,1]^k
	\qquad 
		\pi_{[k] \setminus S}(t) 
		\coloneq
		\del*{ t_1 \indicator{{[k] \setminus S}}(1), \, \ldots, \, t_k \indicator{{[k] \setminus S}}(k) }.
	\end{equation*}
	%
	Then, for all $t \in [0,1]^k$ and
	$ I \subseteq [k]$, 
	we have
	\begin{align*}
		\alpha_I\del*{ \pi_{ [k] \setminus S }(t) }
	&=
		\prod_{ i \in I }
			t_i
			\indicator{[k] \setminus S }(i)
	\\
	&=
		\del*{
			\prod_{ 
				i \in I 
				\cap \del*{
					[k] \setminus S
				}
			}
				t_i \t 1	
		}
		\del*{
			\prod_{ 
				i \in I 
				\setminus \del*{
					[k] \setminus S
				}}
				t_i \t 0		
		}	
		=
		\begin{cases}
			\alpha_I(t) &\text{if } I \subseteq [k] \setminus S, 
		\\ 
			0 &\text{otherwise}, 
		\end{cases}		
	\end{align*}
	and,
	\begin{align*}
		\beta_I\del*{ \pi_{ [k] \setminus S }(t) }
	&=
		\prod_{ i \in I }
			\del*{
				1 - t_i \indicator{ [k] \setminus S }(i)
			}
	\\
	&=
		\del*{
			\prod_{ 
				i \in I 
				\cap 
				\del*{
				[k] \setminus S
			} }
				\del*{
					1 - t_i \t 1
				}
		}
		\del*{
			\prod_{ 
				i \in I 
				\setminus 
				\del*{
				[k] \setminus S
			} }
				\del*{
					1 - t_i \t 0
				}
		}	
		=
		\beta_{I \setminus S}(t),			
	\end{align*}
	where the final equality follows from the fact that
	$
		I \cap 
		\del*{
			[k] \setminus S
		} 
		=
		I \setminus S.
	$
	Hence, for all
	$
		t 
		\in [0,1]^k,
	$
	\begin{equation}
	\label{proofeq::prop::interpolation_estimate::relationship_between_alpha_and_pi}
		\forall I \subseteq [k]
	\quad \, \,
		\alpha_I\del*{ \pi_{ [k] \setminus S }(t) }
		=
		\begin{cases}
			\alpha_I(t) &\text{if } I \subseteq [k] \setminus S, \\ 
			0 &\text{otherwise}, 
		\end{cases}
	\,\,\, \text{and} \,\,\, 
		\beta_I\del*{ \pi_{ [k] \setminus S }(t) }
		=
		\beta_{I \setminus S}(t).
	\end{equation}
	%
	
	In consequence, for all $t \in \intoo{0,1}^k$,
	\begin{align*}
		\partial^S c_{\bx}(t)
	&=
		\partial^S \del*{ 
			\sum_{I = \emptyset}^{ [k] } 
				\alpha_I(t)
				\diff_{ A = \emptyset}^{ I }
					x_A
		}
	\\
	&=
		\sum_{ I = S }^{ [k] }
			\alpha_{I \setminus S}(t)
			\diff_{ A = \emptyset}^{ I }
				x_A
		=
		\sum_{I = \emptyset}^{ [k] \setminus S }
			\alpha_I(t)
			\diff_{A = \emptyset}^{ I \cup S }
				x_A,
	\intertext{
		which, by \Cref{lem::diff_B_diff_C_to_diff_BuC},	
	}
	&=
		\sum_{I = \emptyset}^{ [k] \setminus S }
			\alpha_I(t)
			\diff_{A = \emptyset}^{ I }
				\diff_{J \upmapsto A }^{ S }
					x_J
		=
		\sum_{I = \emptyset}^{ [k] \setminus S }
			\alpha_I(t)
			\diff_{A = \emptyset}^{ I }
				y_A,
	\intertext{
		which, by 
		\eqref{proofeq::prop::interpolation_estimate::relationship_between_alpha_and_pi},
	}
	&=
		\sum_{I = \emptyset}^{ [k]}
			\alpha_I\del*{ \pi_{[k] \setminus S}(t) }
			\diff_{A = \emptyset}^{ I }
				y_A,
	\intertext{
		which, by 
		\eqref{eqprop::prop::interpolation_estimate::deifnitions_of_c_x}
		used for the tuple $\by$,
	}
	&=
		\sum_{ A = \emptyset }^{ [k] }
			\alpha_A\del*{ \pi_{[k] \setminus S}(t) }
			\beta_{[k] \setminus A}\del*{ \pi_{[k] \setminus S}(t) }
			y_A
	\intertext{
		which, by 
		\eqref{proofeq::prop::interpolation_estimate::relationship_between_alpha_and_pi}
		and the definition of $y_A$,
	}
	&=
		\sum_{ A = \emptyset }^{ [k] \setminus S }
			\alpha_A(t)
			\beta_{[k] \setminus (A\cup S)}(t)
				\diff_{ I \upmapsto A }^S
					x_I.
	\end{align*}
	Finally, note that for all $A \subseteq [k]$ we have
	$
		\abs*{ \alpha_A } \le 1
	$
	and 
	$
		\abs*{ \beta_A } \le 1
	$
	everywhere in $[0,1]^k$.
	Thus, for all $t \in \intoo{0,1}^k $,
	\begin{equation*}
		\norm*{
			\partial^S c_{\bx}(t)
		}
		=
		\norm*{
			\sum_{ A = \emptyset }^{ [k] \setminus S }
				\alpha_A(t)
				\beta_{[k] \setminus \del*{ A \cup S } }(t)
				\diff_{ I \upmapsto A }^S
					x_I
		}
		\le 
		\sum_{ A = \emptyset }^{ [k] \setminus S }
			\norm*{
				\diff_{ I \upmapsto A }^{ S }
					x_I
			},
	\end{equation*}
	as claimed.
\end{proof}
\begin{definition}
\label{def::poly_normed_space}
    Let $(V, \| \cdot \|)$ be a normed space, 
    $k \in \bN$,
    $s \in \intoc{k-1, k}$,
    and 
    $
        \bx = \set*{ x_I }_{ I \subseteq [k] } \subseteq V.
    $   
    \begin{itemize}[itemsep=0em,parsep=0em]
        \item 
        For all 
        $S \subseteq [k]$
        and
        $A \subseteq [k] \setminus S$, 
        we define
        \begin{equation}
        \label{eqdef::def::poly_normed_space::poly_A^S}
            \polyset[ A ]{S}( \bx)
            \coloneq 
            \norm*{
                \diff_{ I \upmapsto A }^{ S }
                        x_{I}
                }.
        \end{equation}
        \item 
        For all 
        $S \subseteq [k]$, 
        we define
        \begin{equation}
       	\label{eqdef::def::poly_normed_space::poly^S}
            \polyset{S}( \bx)
            \coloneq 
            \sum_{ A = \emptyset }^{ [k] \setminus S }
                \polyset[ A ]{S}( \bx ).
        \end{equation}
        \item 
        For all 
        $
            \cP \in 
            \partitions{[k]},
        $
        we define
        \begin{equation*}
            \polyfam{\cP}( \bx)
            \coloneq
            \prod_{ S \in \cP }
                \polyset{S}( \bx)
            .
        \end{equation*}
        \item 
        For all 
        $
            j \in [k-1],
        $
        we denote 
        \begin{equation*}
            \polynum{ j }( \bx)
            \coloneq
            \sum_{ \cP \in \partitions[j]{[k]} }
                \polyfam{ \cP }( \bx).
        \end{equation*}
        We also define\footnote{Here, for $k=1$, we use the convention that the value of the empty product is $1$. In consequence, 
        $
            \prod_{ j=1 }^{ 1-1 }
                \polyset{ \set*{j} }( \bx ) 
            =1.
        $}
        \begin{equation*}
            \polynum{ s }( \bx)
            \coloneq
            \polyset{ \set*{k} }(\bx)^{s-k+1}
            \prod_{ j=1}^{k-1}
                \polyset{ \set*{j} }( \bx)
            .
        \end{equation*}
        \item 
        Finally, we define 
        \begin{equation*}
            \polygen{ \gen{s} }( \bx)
            \coloneq
            \polynum{ s }( \bx)
            +
            \sum_{ j = 1 }^{ k-1 }
                \polynum{ j }( \bx)  
        \end{equation*}
        and 
        \begin{equation*}
            \ell( \bx )
            \coloneq 
            \sum_{ j = 1 }^{ k }
                \polyset{ \set*{j} }( \bx).
        \end{equation*}
    \end{itemize}
    Moreover, 
    recalling that $[0] = \emptyset$,
    we extend the definition of 
    $ \polygen{ \gen{s} } $
    to $ s = 0 $
    by the formula    
    \begin{equation*}
        \forall 
        \bx = \set*{ x_I }_{ I \subseteq [0] }
        \subseteq 
        V
    \qquad 
        \polygen{ \gen{0} }( \bx )
        \coloneq
        1.
    \end{equation*}
\end{definition}
\begin{example}\label{ex::special_values_of_P_s_k}
    Let $(V, \| \cdot \|)$ be a normed space, 
    $ k \in \bN $,
    $s \in \intoc{k-1, k}$,
    and 
    $
        \bx = \set*{ x_I }_{ I \subseteq [k] } \subseteq V.
    $   
    \begin{enumerateInExample}[ref=\Cref{ex::special_values_of_P_s_k}(\roman*)]
        \item \label{ex::special_values_of_P_s_k::item::k=1}
        If $s \ge 1$,
        then
        $
            \polynum{ 1 }(\bx)
            =
            \norm*{
                \diff_{I = \emptyset}^{[k]}
                    x_I
            }.
        $
        Indeed, this follows from the fact that $\set*{ [k] }$ is the unique partition of~$[k]$ that has only one element. In consequence,
        \begin{equation*}
            \polynum{ 1 }(\bx)
            =
            \polyfam{ \set*{  [k] } }\del*{ \bx }
            =
            \polyset{ [k] }( \bx )
            =
            \polyset[\emptyset]{ [k] }( \bx )
            =
            \norm*{
                \diff_{I \upmapsto \emptyset}^{[k]}
                    x_I
            }
            =
            \norm*{
                \diff_{I = \emptyset}^{[k]}
                    x_I
            }.
        \end{equation*}
        In particular, if $s=1$,
        then
        $
            \polynum{ 1 }(\bx)
            =
            \norm{
                x_{ \set*{1}}
                - x_{\emptyset}
            }.
        $
        \item 
        When $s \in \intoc{0,1}$, then
        $
            \polygen{ \gen{s} }(\bx)
            =
            \norm{
                x_{\set*{1}} - x_{\emptyset}
            }^s.
        $
        \item \label{ex::special_values_of_P_s_k::item::k>=s}
        $
            \polynum{ s }(\bx)
            =
            \del*{
                \sum_{
                        I = \emptyset
                    }^{ [k-1] }
                        \norm{
                            x_{ I \cup \set*{k} }
                            - x_I
                        }
            }^{s-k+1}
            \prod_{j=1}^{k-1}
                \sum_{
                    I = \emptyset
                }^{ [k] \setminus \set*{j} }
                    \norm{
                            x_{ I \cup \set*{j} }
                            - x_I
                        }
            .
        $
    \item 
    \label{ex::special_values_of_P_s_k::item::polyen_of_k}
        $
        	\polygen{ \gen{k}}\del*{\bx }
        	=
	        \sum_{ \cP \in \partitions{ [k] } }
	            \polyfam{ \cP }( \bx ).
        $
    \end{enumerateInExample}
\end{example}
\begin{lemma}
\label{cor::If_P^set(j)(x)=0_for_any_j_then_diff_f(x)=0}
    Let $( V, \norm{ \cdot })$ be a normed space, 
    $ X \subseteq V $,
    $k \in \bN$, and 
    $
        \bx = \set*{ x_I }_{ I \subseteq [k] }
        \subseteq 
        X.
    $
    Let $ (G, +) $ be~an~Abelian group and
    $
        f \colon X \to G.
    $
    If there exists $j \in [k]$ such that
    $
        \polyset{ \set*{j} }( \bx ) = 0,
    $
    then
    $
        \diff_{ I = \emptyset }^{ [k] }
            f(x_I)
        = 
        0.
    $    
\end{lemma}
\begin{proof}
    Let us
    suppose that $\polyset{ \set*{j} }( \bx ) = 0$
    for some $j \in [k]$. 
    Then, for all $A \subseteq [k] \setminus \set*{j}$,
    we have
    $
        x_{A \cup \set*{j} } = x_A
    $
    since
    $
        0
        =
        \polyset[A]{ \set*{j} }( \bx ) 
        =
        \norm{ x_{A \cup \set*{j} } - x_A }.
    $
    Therefore, 
    for all
    $
        y \in \setc*{ f(x_I) }{ I \subseteq [k] },
    $
    by the 
    \say{in particular}
    part of 
    \Cref{lem::diff_B_diff_C_to_diff_BuC},
    we~have
    \begin{align*}
        \diff_{ I = \emptyset }^{ [k] }
            \indicator{ \set*{y} }\del{ f(x_I) }
        =
        \diff_{ I = \emptyset }^{ [k] \setminus \set*{j} }
            \del*{
                \indicator{ \set*{y} }\del*{ f\del*{x_{I \cup \set*{ j } } } } 
                -
                \indicator{ \set*{y} }\del*{ f(x_I) }
            }&
    \\
        =
        \diff_{ I = \emptyset }^{ [k] \setminus \set*{j} }
            \del*{
                \indicator{ \set*{y} }\del*{ f(x_I) }  
                -
                \indicator{ \set*{y} }\del*{ f(x_I) }
            }&
        =
        \diff_{ I = \emptyset }^{ [k] \setminus \set*{j} }
            0
        =
        0.
    \end{align*}
    Therefore, 
    $
        \diff_{ I = \emptyset }^{ [k] }
            f(x_I)
        = 
        0
    $  
    by
    \Cref{lem::characterization_of_diff_0_with_indicators}.
\end{proof}
\begin{lemma}
\label{lem::invariant_transformations_of_poly}
    Let $(V, \normAlone)$ be a normed space, 
    $k \in \bN$,
    $s \in \intoc{k-1, k}$,
    and 
    $
        \bx = \set*{ x_I }_{ I \subseteq [k] } \subseteq V.
    $   
    \begin{itemize}
        \item
        Let $T \colon V \to V$ be an isometry.\footnote{i.e., a bijective, distance-preserving map.}
        Define $\ba = \set*{ a_I }_{ I \subseteq [k] }$
        by $a_I \coloneq T(x_I)$
        for all $ I \subseteq [k]$.
        \item
        Fix $i \in [k]$ and define 
        $\bb = \set*{ b_I }_{ I \subseteq [k] }$
        by $b_I \coloneq x_{ I \div \set*{i} }$
        for all $ I \subseteq [k]$.
        \item
        Let 
        $ \sigma $
        be a permutation of 
        $ [k] $.
        Define 
        $\bc = \set*{ c_I }_{ I \subseteq [k] }$
        by $c_I \coloneq x_{ \sigma\sbr{ I } }$
        for all $ I \subseteq [k]$.
    \end{itemize}
    %
    %
    Then~%
    \begin{center}
        \begin{tabular}{r r|l|l|l l}
            &
            $ ( \cdot )  = $
            & \multicolumn{1}{|c|}{ $( \ba ) $ }
            & \multicolumn{1}{|c|}{ $( \bb ) $ }
            & \multicolumn{1}{|c}{ $( \bc ) $ }
            &
        \\
            \hline \hline
            For all nonempty $S\subseteq [k],$ 
            & 
            $\polyset{S}( \cdot ) = $ 
            &  
            $\polyset{S}(\bx)$
            &  
            $\polyset{S}(\bx)$
            &  
            $\polyset{ \sigma \sbr{S} }(\bx)$
            &
        \\
            \hline
            For all $ \cP \in \partitions{ [k] },$
            &
            $\polyfam{\cP}( \cdot ) = $ 
            &  
            $\polyfam{\cP}(\bx)$
            &  
            $\polyfam{\cP}(\bx)$
            &  
            $\polyfam{ \sigma \dsbr{\cP} }(\bx)$ 
            &
        \\
            \hline
            For all $ j \in [k]$, 
            &
            $\polynum{ j }( \cdot ) = $ 
            &  
            $\polynum{ j }(\bx)$
            &  
            $\polynum{ j }(\bx)$
            &  
            $\polynum{ j }(\bx)$ 
            &
        \\
            \hline
            &
            $\polynum{ s }( \cdot ) = $ 
            &  
            $\polynum{ s }(\bx)$
            &  
            $\polynum{ s }(\bx)$
            &  
            \multicolumn{2}{l}{
                $\polynum{ s }(\bx)$ 
                if $\sigma(k) = k$
            } 
        \\
            \hline
            &
            $\polygen{ \gen{s} }( \cdot ) = $ 
            &  
            $\polygen{ \gen{s} }(\bx)$
            &  
            $\polygen{ \gen{s} }(\bx)$
            &  
            \multicolumn{2}{l}{
                $\polygen{ \gen{s} }(\bx)$ 
                if $\sigma(k) = k$
            } 
        \\  
            \hline
            &
            $\ell( \cdot )=$
            & 
            $\ell( \bx )$
            &
            $\ell( \bx )$
            &
            $\ell( \bx )$
            &
        \end{tabular}
    \end{center}
\end{lemma}
\begin{proof}
    We will prove the statements for each column separately.
    \begin{enumerate}
        \item[$(\ba)$]
        It is sufficient to only show that 
        $
            \polyset[A]{S}(\ba)
            = \polyset[A]{S}(\bx)
        $ for all $S \subseteq [k]$,
        $S \ne \emptyset$, and 
        $A \subseteq [k] \setminus S$.
        Fix such $S$ and $A$.
        By the Mazur--Ulam theorem \cite{quick_Mazur-Ulam_proof},
        $T$ is an affine transformation. 
        In consequence, 
        there are $v \in V$ and~a~linear isometry 
        $L \colon V\to V$
        such that 
        $T(x) = v + L(x)$ 
        for all 
        $x \in V$.
        Therefore, 
        \begin{align*}
            \polyset[A]{S}(\ba) 
            =
            \norm*{
                \diff_{ I \upmapsto A }^{ S }
                    a_I
            }
        &=
            \norm*{
                \diff_{ I \upmapsto A }^{ S }
                    \del{ v + L(x_I) }
            }
        \\
        &=
            \norm*{
                \diff_{ I \upmapsto A }^{ S }
                    v 
                + 
                \diff_{ I \upmapsto A }^{ S }
                    L(x_I)
            }
            =
            \norm*{
                \diff_{ I \upmapsto A }^{ S }
                    L(x_I)
            }
        \\
        &\;
        \phantom{
        	=
        	\norm*{
        	                \diff_{ I \upmapsto A }^{ S }
        	                    v 
        	                + 
        	                \diff_{ I \upmapsto A }^{ S }
        	                    L(x_I)
        	            }
        }
            =
            \norm*{
                L \del*{
                    \diff_{ I \upmapsto A }^{ S }
                        x_I
                }
            }
            =
            \norm*{
                \diff_{ I \upmapsto A }^{ S }
                        x_I
            }
            =
            \polyset[A]{S}( \bx),
        \end{align*}
        where 
        $
            \diff_{ I \upmapsto A }^{ S }
                    v 
            =
            0
        $
        by \Cref{lem::diff_constant_term}.
        \item[$(\bb)$]
        It is sufficient to only show that 
        $\polyset{ S }( \bb) = \polyset{ S }(\bx)$
        for all $S \subseteq [k]$.
        Fix disjoint $A, S \subseteq [k]$. 
        If~$i \in S$, then $i \notin A$, hence,
        by the \say{moreover} part of
        \Cref{lem::diff_B_diff_C_to_diff_BuC},
        \begin{align*}
            \polyset[A]{S}( \bb) 
            =
            \norm*{
                \diff_{ I \upmapsto A }^{ S }
                    b_I
            }
        &=
            \norm*{
                \diff_{ I \upmapsto A }^{ 
                    S \setminus \set*{i} 
                }
                    \del*{
                        b_{ I \cup \set*{i} } 
                        - b_I
                    }
            }
        \\
        &=
            \norm*{
                \diff_{ I \upmapsto A }^{ 
                    S \setminus \set*{i} 
                }
                    \del*{
                        x_{ (I \cup \set*{i}) \div \set*{i} } 
                        - x_{I \div \set*{i} } 
                    }
            }
            =
            \norm*{
                \diff_{ I \upmapsto A }^{ 
                    S \setminus \set*{i} 
                }
                    \del*{
                        x_{ I }
                        - x_{I \cup \set*{i} } 
                    }
            }
        \\
        &\;
        \phantom{
        	=
				  \norm*{
				      \diff_{ I \upmapsto A }^{ 
				          S \setminus \set*{i} 
				      }
				          \del*{
				            x_{ (I \cup \set*{i}) \div \set*{i} } 
				              - x_{I \div \set*{i} } 
				          }
  }       
        }
        	=
            \norm*{
                - \diff_{ I \upmapsto A}^{S}
                    x_I
            }
            =
            \polyset[A]{S}( \bx).
        \end{align*}
        Summing the resulting equality over all $A \subseteq [k] \setminus S$,
        we get 
        $
            \polyset{ S }(\bb)
            = \polyset{ S }( \bx).
        $
        If $i \notin S$, then
        \begin{align*}
            \polyset[ A ]{ S }( \bb) 
            =
            \norm*{
                \diff_{ I \upmapsto A }^{ S }
                    b_I
            }
        	=
            \norm*{
                \diff_{ I = \emptyset }^{ S }
                    b_{ I \cup A }
            }
        &=
            \norm*{
                \diff_{ I = \emptyset }^{ S }
                    x_{ (I \cup A) \div \set*{i } }
            }
        \\
        &=
            \norm*{
                \diff_{ I = \emptyset }^{ S }
                    x_{ I \cup ( A \div \set*{i }) }
            }
            =
            \norm*{
                \diff_{ I \upmapsto A \div \set*{i} }^{ S }
                    x_I
            }
            =
            \polyset[A \div \set*{i}]{ S }( \bx).
        \end{align*}
        Since $i \notin S$, 
        the map 
        $ A \mapsto A \div \set*{i}$ is a bijection
        from $2^{ [k] \setminus S }$ to itself.
        Hence, 
        \begin{equation*}
            \polyset{ S }(\bb)
            =
            \sum_{ A = \emptyset }^{ [k] \setminus S}
                \polyset[ A ]{ S }( \bb)
            =
            \sum_{ A = \emptyset }^{ [k] \setminus S}
                \polyset[A \div \set*{i}]{S}( \bx)
            =
            \sum_{ A = \emptyset }^{ [k] \setminus S}
                \polyset[ A ]{ S }( \bx)
            =
            \polyset{ S }( \bx).
        \end{equation*}
        Thus, we have 
        $\polyset{ S }( \bb) = \polyset{ S }(\bx)$
        for all $S \subseteq [k]$.
        \item[$(\bc)$]
        Fix disjoint $A, S \subseteq [k]$. 
        Since $\sigma \colon [k] \to [k]$ is a bijection, 
        for all $I \subseteq S$ we have that
        $
            \card{ S \setminus I}
            =
            \card{ \sigma[S] \setminus \sigma[I] }.
        $
        Hence,
        \begin{multline*}
            \polyset[A]{S}(\bc) 
            =
            \norm*{
                \diff_{ I \upmapsto A }^{ S }
                    c_I
            }
    	=
            \norm*{ 
                \sum_{ I = \emptyset }^{ S }
                    (-1)^{ \card{ S \setminus I } }
                    c_{ I \cup A }
            }
        \\
        =
            \norm*{ 
                \sum_{ I = \emptyset }^{ S }
                    (-1)^{ \card{ S \setminus I } }
                    x_{ \sigma \sbr{ I \cup A } }
            }
            =
            \norm*{ 
                \sum_{ I = \emptyset }^{ S }
                    (-1)^{ 
                        \card*{ \sigma[S] \setminus \sigma[I] } 
                    } 
                    x_{ 
                        \sigma \sbr{ I } 
                        \cup \sigma \sbr{ A } 
                    }
            }
        \\
        	=
            \norm*{ 
                \sum_{ I = \emptyset }^{ 
                    \sigma \sbr{S} 
                }
                    (-1)^{ \card{ \sigma[S] \setminus I } }
                    x_{ 
                        I \cup \sigma \sbr{ A } 
                    }
            }
            =
            \norm*{
                \diff_{ I \upmapsto \sigma \sbr{A} }^{
                    \sigma \sbr{S} 
                }
                    x_I
            }
            =
            \polyset[ \sigma \sbr{A} ]{ \sigma\sbr{S} }
                (\bx),
        \end{multline*}
        where the first equality of the third line follows from the fact that
        $
            I \mapsto \sigma \sbr{I}
        $
        is~a~bijection from $2^S$
        to $2^{ \sigma\sbr{S} }$.
        Next, notice that 
        $A \mapsto \sigma [A]$ 
        is a bijection from 
        $2^{ [k] \setminus S}$ to 
        $2^{ [k] \setminus \sigma\sbr{S} }$.
        Thus, 
        \begin{equation*}
            \polyset{ S }( \bc)
            =
            \sum_{ A = \emptyset }^{ [k] \setminus S}
                \polyset[ A ]{ S }(\bc)
            =
            \sum_{ A = \emptyset }^{ [k] \setminus S}
                \polyset[ \sigma\sbr{A} ]{ \sigma \sbr{S} }(\bx)
            =
            \sum_{ A = \emptyset }^{ 
                [k] \setminus \sigma\sbr{S} 
            }
                \polyset[ A ]{ \sigma\sbr{S} }( \bx)
            =
            \polyset{ \sigma\sbr{S} }( \bx).
        \end{equation*}
        Hence,
        \begin{equation*}
            \ell( \bc )
            =
            \sum_{ j = 1 }^{ k }
                \polyset{ \set*{j} }( \bc )
            =
            \sum_{ j = 1 }^{ k }
                \polyset{ \sigma \sbr{ \set*{j} } }( \bx )
            =
            \sum_{ j = 1 }^{ k }
                \polyset{ \set*{ \sigma(j) } }( \bx )
            =
            \sum_{ j = 1 }^{ k }
                \polyset{ \set*{ j } }( \bx )
            =
            \ell( \bx ).
        \end{equation*}
        Also, for all 
        $\cP \in \partitions{[k]}$ we have
        \begin{equation*}
            \polyfam{ \cP }( \bc)
            =
            \prod_{ S \in \cP }
                \polyset{ S }( \bc)
            =
            \prod_{ S \in \cP }
                \polyset{ \sigma\sbr{S} }( \bx)
            =
            \prod_{ S \in \sigma \dsbr{\cP} }
                \polyset{ S }( \bx)
            =
            \polyfam{ \sigma \dsbr{ \cP } } ( \bx ),
            %
        \end{equation*}
        where
        $
            \sigma \dsbr{ \cP }
            \coloneq 
            \setc*{
                \sigma [S]
            }{
                S \in \cP
            }.
        $
        Now, notice that for all 
        $j \in [k]$, 
        the mapping
        $ \cP \mapsto \sigma \dsbr{ \cP }$
        is~a~bijection from $\partitions[j]{[k]}$
        to itself.
        Hence, for all such $j$,
        \begin{equation*}
            \polynum{ j }(\bc)
            =
            \sum_{ \cP \in \partitions[j]{[k]} }
                \polyfam{ \cP }(\bc)
            =
            \sum_{ \cP \in \partitions[j]{[k]} }
                \polyfam{ \sigma \dsbr{ \cP} }(\bx)
            =
            \sum_{ \cP \in \partitions[j]{[k]} }
                \polyfam{ \cP }(\bx)
            =
            \polynum{ j }(\bx).
        \end{equation*}
        Finally, assume that $\sigma(k) = k$. Then
        $ j \mapsto \sigma(j) $ is a bijection from 
        $[k-1]$ to itself.
        Hence,
        \begin{multline*}
            \polynum{ s }(\bc)
            =
            \polyset{ \set*{k} }(\bc)^{s-k+1}
            \prod_{j=1}^{k-1}
                \polyset{ \set*{j} }(\bc)
            =
            \polyset{ \sigma\sbr{\set*{k} } }(\bx)^{s-k+1}
            \prod_{j=1}^{k-1}
                \polyset{ \sigma\sbr{\set*{j}} }(\bx)
            \\
            =
            \polyset{ \set*{ \sigma(k) }  }(\bx)^{s-k+1}
            \prod_{j=1}^{k-1}
                \polyset{ \set*{\sigma(j)} }(\bx)
            =
            \polyset{ \set*{k}  }(\bx)^{s-k+1}
            \prod_{j=1}^{k-1}
                \polyset{ \set*{j} }(\bx)
            =
            \polynum{ s }(\bx)
        \end{multline*}
        and 
        \begin{equation*}
            \polygen{ \gen{s} }(\bc)
            =
            \polynum{ s }(\bc)
            +
            \sum_{j=1}^{k-1}
                \polynum{ j }(\bc)
            =
            \polynum{ s }(\bx)
            +
            \sum_{j=1}^{k-1}
                \polynum{ j }(\bx)
            =
            \polygen{ \gen{s} }(\bx).
            \tag*{\qedhere} 
        \end{equation*}
    \end{enumerate}
\end{proof}
\begin{lemma} \label{lem::poly_partitioned_gives_0}
    Let $\del*{V, \normAlone}$ be a normed space, 
    $k \in \bN$, 
    and $\set*{C, D}$ be a partition of $[k]$.
    Let~%
    $g \in V$
    and
    $
        \set*{ v_{L} }_{ L \subseteq C}, 
        \set*{w_{J}}_{J \subseteq D} 
        \subseteq V.
    $
    For $I \subseteq [k]$, let 
    $
        x_{I}
        \coloneq 
        g
        + v_{ I \cap C}
        + w_{ I \cap D}.   
    $
    Suppose that
    $\cP \in \partitions{[k]}$
    has~an~element $B \in \cP$ such that
    $B \cap C \ne \emptyset$
    and $B \cap D \ne \emptyset$.
    Then
    $
        \polyfam{ \cP }(\bx)
        =
        0.
    $
\end{lemma}
\begin{proof}
    Let us fix such $\cP$ and let $B \in \cP$
    be such that
    $B \cap C \ne \emptyset$
    and $B \cap D \ne \emptyset$.
    By \Cref{lem::partitioned_diff_gives_0}
    for all 
    $A \subseteq [k] \setminus B$
    we have
    $\diff_{I \upmapsto A}^{B} x_I = 0.$
    In consequence,
    \begin{equation*}
        \polyset{ B }( \bx)
        =
        \sum_{ A = \emptyset }^{ [k] \setminus B }
            \norm*{
                \diff_{ I \upmapsto A }^{ B }
                    x_{I}
            }
        =
        0.
    \end{equation*}
    Therefore,
    \begin{equation*}
        \polyfam{ \cP }( \bx)
        =
        \prod_{ S \in \cP }
            \polyset{ S }( \bx)
        =
        0
    \end{equation*}
    since at least one of the factors is $0$.
\end{proof}
\begin{corollary}
\label{cor::value_of_Poly_when_tuple_is_a_hyperparallelogram}
    Let $(V, \| \cdot \|)$ be a normed space, 
    $k \in \bN$, and
    $s \in \intoc{k-1, k}$.
    Fix $g \in V$ and 
    $
        \set*{ y_i }_{ i \in [k] }
        \subseteq
        V.
    $
    For $I \subseteq [k]$,
    let 
    $
        x_I
        \coloneq 
        g 
        + \sum_{ i \in I } 
            y_i.
    $
    Then 
    \begin{equation*}
        \polygen{ \gen{s} }(\bx) 
        = 
        \polynum{ s }(\bx)
        =
        2^{s(k-1)}
        \norm{ y_k }^{ s-k + 1}
        \prod_{ j = 1 }^{ k - 1 }
            \norm{ y_j }.
    \end{equation*} 
\end{corollary}
\begin{proof}
    To prove the first equality, 
    it is sufficient to show that for 
    all $j \in [k-1 ]$ we have
    $\polynum{ j }( \bx ) =0$, 
    or, equivalently, that for all 
    $j \in [k-1]$ and $\cP \in \partitions[j]{[k]}$
    we have
    $
        \polyfam{ \cP }(\bx) = 0.
    $
    Fix such $j$ and $\cP$. 
    Then there is $B \in \cP$ such that
    $\card{ B } > 1.$
    Fix $c \in B$ and let 
    $C \coloneq \set*{c}$
    and 
    $D \coloneq [k] \setminus C$. 
    Then
    $ \set*{C, D} $
    is a partition of $[k]$.
    For $L \subseteq C$, 
    let 
    $v_L \coloneq x_L - g$,
    and for $J \subseteq D$, 
    let 
    $w_J \coloneq x_J - g$. 
    Then for all $I \subseteq [k]$ we have
    \begin{equation*}
        x_I
        =
        g
        + \sum_{ i \in I }
            y_i
        =
        g
        + \sum_{ i \in I \cap C}
            y_i
        +
        \sum_{ i \in I \cap D}
            y_i
        =
        x_{I \cap C}
        + x_{I \cap D}
        - g
        =
        g
        + v_{I \cap C}
        + w_{I \cap D}.
    \end{equation*}
    Since $B \cap C \ne \emptyset$ 
    and $B \cap D \ne \emptyset$,
    by \Cref{lem::poly_partitioned_gives_0}
    we have 
    $
        \polyfam{ \cP }(\bx) = 0,
    $
    as needed.

    To prove the second equality, notice that for all
    $j \in [k]$ and all 
    $A \subseteq [k] \setminus \set*{j}$ we have
    \begin{equation*}
        \polyset[ A ]{ \set*{j} }(\bx)
        =
        \norm*{
            \diff_{ I \upmapsto A}^{ \set*{j } }
                x_I
        }
        =
        \norm{
            x_{ A \cup \set*{j} }
            - x_{A}
        }
        =
        \norm*{
            \del*{ 
                g
                + \sum_{ i \in A \cup \set*{j} } 
                y_i
            }
            -
            \del*{
                g
                + \sum_{ i \in A } 
                    y_i
            }
        }
        =
        \norm{ y_j },
    \end{equation*}
    so
    \begin{equation*}
        \polyset{ \set*{j} }(\bx)
        =
        \sum_{ A = \emptyset }^{ [k] \setminus \set*{j} } 
            \polyset[ A ]{ \set*{j} }( \bx)
        =
        \sum_{ A = \emptyset }^{ [k] \setminus \set*{j} } 
            \norm{ y_j }
        =
        2^{k-1} \norm{ y_j }.
    \end{equation*}
    Therefore,
    \begin{multline*}
        \polynum{ s }( \bx)
        =
        \polyset{ \set*{k} }(\bx)^{s-k+1}
        \prod_{ j=1}^{k-1}
            \polyset{ \set*{j} }( \bx)
        =
        \del*{ 2^{k-1} \norm{ y_k} }^{s-k+1}
        \prod_{ j=1}^{k-1}
            \del*{ 2^{k-1} \norm{ y_j} }
    \\
        =
        2^{(k-1)(s-k+1) + (k-1)^2}
        \norm{ y_k}^{s-k+1}
        \prod_{ j=1}^{k-1}
            \norm{ y_j}
        =
        2^{s(k-1)}
        \norm{ y_k}^{s-k+1} 
        \prod_{ j=1}^{k-1}
            \norm{ y_j},
    \end{multline*}
    as claimed.
\end{proof}
%
%
\begin{corollary}
\label{cor::tuple_part_tuple_part_vectors}
    Let $\del*{V, \normAlone}$ be a normed space,
    $ k  \in \bN$,
    $ m \in [k-1] $,
    and 
    $
        s \in \intoc{k-1,k}.
    $
    Let~%
    $
        \bx = \set*{ x_I }_{ I \subseteq [k-m] }
        \subseteq 
        V
    $
    and 
    $
        \set*{ v_i }_{ j = 1}^{ m }
        \subseteq 
        V.
    $
    Define 
    $   
        \bv = \set*{ v_L }_{ L \subseteq [m] } 
    $
    by the formula
    $
        v_L
        \coloneq 
        \sum_{ \ell \in L }
            v_{\ell}
    $
    for $L \subseteq [m]$, 
    and 
    $
        \by = \set*{ y_J }_{ J \subseteq [k] }
    $
    by~the~formula
    \begin{equation*}
        \forall L \subseteq [m]
        \quad 
        \forall I \subseteq [k-m]
        \qquad 
        y_{ L \cup (m+I)}
        \coloneq 
        v_L + x_I.
    \end{equation*}
    Then 
    \begin{equation*}
        \polygen{ \gen{ s} }( \by )
        \le
        2^{ (2k-1)m }
        \polygen{ \gen{ s - m} }(\bx)
        \prod_{ \ell = 1 }^{ m }
            \norm{ v_{\ell} }.
    \end{equation*}
\end{corollary}
%
%
\begin{proof}
    %
    We will start by showing that 
    \begin{align}
        \label{proofeq::cor::tuple_part_tuple_part_vectors::P^l(y)_for_l_in_[m]}
        \forall \ell \in [m]
        \qquad 
        &\polyset{ \set*{ \ell } }(\by) 
        =
        2^{k-1} \norm{ v_{ \ell } }
    \intertext{and}
        \label{proofeq::cor::tuple_part_tuple_part_vectors::P^m+S(y)_for_S_in_[k-m]}
        \forall S \subseteq [k-m]
        \qquad 
        S \ne \emptyset 
        \implies
        &\polyset{ m + S }( \by )
        =
        2^m \polyset{ S }( \bx ).
    \end{align}
    Fix 
    $ \ell' \in [m]$,
    $
        A \subseteq [m] \setminus \set*{ \ell' },
    $
    and 
    $
        B \subseteq [k-m].
    $
    Then 
    $
        A \cup \set*{ \ell' } \subseteq [m],
    $
    hence
    \begin{align*}
        \polyset[A \cup (m+B)]{\set*{\ell'}}( \by )
        =
        \norm*{
            \diff_{ J \upmapsto A \cup (m+B) }^{ \set*{\ell'} }
                y_J
        }
        &=
        \norm*{
            y_{ A \cup (m+B) \cup \set*{\ell'} }
            - y_{ A \cup (m+B) }
        }
    \\
        &=
        \norm*{
            \del*{
                \sum_{ \ell \in A \cup \set*{ \ell' } }
                    v_{\ell }
                +
                x_B
            }
            - 
            \del*{
                \sum_{ \ell \in A }
                    v_{\ell }
                +
                x_B
            }
        }
        =
        \norm{ v_{\ell'} }.
    \end{align*}
    In consequence,
    since 
    $ (A, B) \mapsto A \cup (m+B)$
    is a bijection from
    $ 2^{[m] \setminus \set*{\ell'} } \times 2^{[k-m]}$
    to 
    $ 2^{ [k] \setminus \set*{\ell'} }$,
    \begin{align*}
        \polyset{ \set*{\ell' } }( \by)
        =
        \sum_{ C = \emptyset }^{ 
            [k] \setminus \set*{ \ell' }
        }
            \polyset[C]{ \set*{\ell' } }( \by)
    	=
        \sum_{ A = \emptyset }^{ 
            [m] \setminus \set*{ \ell'}  
        }
            \sum_{ B = \emptyset }^{ [k-m] }
                \polyset[ A \cup (m+B) ]{ \set*{\ell' } }( \by)
    &
    \\
    	=
        \sum_{ A = \emptyset }^{ 
            [m] \setminus \set*{ \ell'}  
        }
            \sum_{ B = \emptyset }^{ [k-m] }
                \norm{ v_{\ell'} }
    &=
        2^{ m-1 } 2^{ k-m } \norm{ v_{\ell'} }
        =
        2^{k-1}
        \norm{ v_{\ell'} },   
    \end{align*}
    which proves 
    \eqref{proofeq::cor::tuple_part_tuple_part_vectors::P^l(y)_for_l_in_[m]}.
    
    Next, for the proof of 
    \eqref{proofeq::cor::tuple_part_tuple_part_vectors::P^m+S(y)_for_S_in_[k-m]},
    fix 
    $S \subseteq [k-m]$ that is nonempty,
    $ A \subseteq [m]$,
    and 
    $ B \subseteq [k-m] \setminus S $.
    Since
    $ J \mapsto m+J $ 
    is a bijection from $2^S$ to $2^{ m+S }$,
    \begin{align*}
        \diff_{ J \upmapsto m + B }^{ m + S }
            x_{ J -m }
        &=
        \sum_{ J = \emptyset }^{ m + S }
            (-1)^{ \card{ (m+S) \setminus J } }
            x_{ \del*{ J \cup (m + B)} - m }
    \\
        &=
        \sum_{ J = \emptyset }^{ S }
            (-1)^{ \card{ (m+S) \setminus (m+J) } }
            x_{ \del*{ (m+J) \cup (m+B) } - m }
        =
        \sum_{ J = \emptyset }^{ S }
            (-1)^{ \card{ S \setminus J } }
            x_{ J \cup B }
        =
        \diff_{ J \upmapsto B }^{ S }
            x_J.
    \end{align*}
    Therefore, because 
    $A \cap (m+S) = \emptyset$
    and 
    $ m + S \ne \emptyset$,
    by \Cref{lem::diff_constant_term} we have
    \begin{align*}
        \polyset[ A \cup (m+B) ]{ m + S}( \by )
        =
        \norm*{
            \diff_{ J \upmapsto A \cup (m+B) }^{ m + S }
                y_J
        }
        =
        \norm*{
            \diff_{ J \upmapsto m+B }^{ m + S }
                y_{A \cup J}
        }
        &=
        \norm*{
            \diff_{ J \upmapsto m+B }^{ m + S }
                \del*{
                    v_{A} + x_{J - m}
                }
        }
    \\
    &=
        \norm*{
            \diff_{ J \upmapsto m+B }^{ m + S }
                x_{J - m}
        }
        =
        \norm*{
            \diff_{ J \upmapsto B }^{ S }
                x_{J}
        }
        =
        \polyset[ B ]{ S }( \bx ).
    \end{align*}
    Therefore, since 
    $ (A, B) \mapsto A \cup (m+B)$
    is a bijection from 
    $ 2^{ [m] } \times 2^{ [k-m] \setminus S }$
    to 
    $ 2^{ [k] \setminus S }$,
    \begin{equation*}
        \polyset{m+S}( \by )
        =
        \sum_{ C = \emptyset }^{ [k] \setminus S }
            \polyset[ C ]{ m+S }( \by )
        =
        \sum_{ A = \emptyset }^{ [m] }
            \sum_{ B = \emptyset }^{ [k-m] \setminus S }
                \polyset[ A \cup (m+B) ]{ m + S }( \by )
        =
        \sum_{ A = \emptyset }^{ [m] }
            \sum_{ B = \emptyset }^{ [k-m] \setminus S }
                \polyset[ B ]{ S }( \bx )
        =
        2^m
        \polyset{ S }( \bx ),
    \end{equation*}
    proving 
    \eqref{proofeq::cor::tuple_part_tuple_part_vectors::P^m+S(y)_for_S_in_[k-m]}.

    Next, we will show that if we have 
    $
        \polyfam{ \cP }( \by ) \ne 0
    $
    for some $\cP \in \partitions{ [k] }$,
    then 
    \begin{equation*}
        \cP 
        = 
        \cP'
        \cup 
        \setc*{
            \set*{\ell}
        }{
            \ell \in [m]
        }
    \end{equation*}
    for some 
    $ \cP' \in \partitions{ m + [k-m] }$.
    Fix 
    $\cP \in \partitions{ [k] }$.
    It is sufficient to show that 
    $ \set*{ \ell } \in \cP $
    for all 
    $\ell \in [m]$.
    Suppose this is not the case and there exists 
    $\ell' \in [m] $ 
    such that
    $ \set*{ \ell' } \notin \cP $.
    Since 
    $\ell' \in [m] \subseteq [k]$,
    there exists 
    $B \in \cP$ 
    such that 
    $\ell' \in B$.
    Denote 
    $C \coloneq \set*{ \ell' }$
    and 
    $D \coloneq [k] \setminus \set*{ \ell' }$.
    Then 
    $B \cap C \ne \emptyset$
    and 
    $B \cap D \ne \emptyset$.
    Let us now define 
    $
        \set*{ \wave{v}_{L} }_{ L \subseteq C }
    $
    by~the~formula
    \begin{equation*}
        \forall L \subseteq C
        \qquad 
        \wave{v}_{L}
        \coloneq
        \sum_{ \ell \in L }
            v_{ \ell }.
    \end{equation*}
    Next, noting that $m + [k-m] \subseteq D$,
    define
    $
        \set*{ \wave{w}_{J} }_{ J \subseteq D }
    $
    by the formula
    \begin{equation*}
        \forall L \subseteq D \cap [m]
        \quad 
        \forall I \subseteq [k-m]
        \qquad 
        \wave{w}_{L \cup (m+I)}
        \coloneq
        \sum_{ \ell \in  L}
            v_{\ell}
        +
        x_I.
    \end{equation*}
    Then, for all $L \subseteq [m]$ and all $I \subseteq [k-m]$, we have
    \begin{align*}
        y_{L \cup (m+I)}
    &=
        v_L 
        + x_I
    \\
    &=
        \sum_{ \ell \in L }
            v_{\ell}
        + x_I
    	=
        \sum_{ \ell \in L \cap C }
            v_{\ell}
        +
        \sum_{ \ell \in L \cap D }
            v_{\ell}
        + x_I
    \\
    &
    \phantom{
    	\sum_{ \ell \in L }
    	            v_{\ell}
    	        + x_I
    	===\,\,
    }
    =
        \wave{v}_{ L \cap C }
        + \wave{w}_{ L \cap D \cup (m+I)}
    	=
        \wave{v}_{ \del*{ L \cup (m+I) } \cap C }
        +
        \wave{w}_{ \del*{ L \cup (m+I) } \cap D },
    \end{align*}
    or, equivalently, 
    $
        y_J
        =
        \wave{v}_{ J \cap C}
        + \wave{w}_{ J \cap D}
    $
    for all $J \subseteq [k]$.
    Therefore, by 
    \Cref{lem::poly_partitioned_gives_0}
    we have
    $
        \polyfam{\cP}( \by ) = 0.
    $
    We obtained a contradiction, 
    so if 
    $ \polyfam{\cP}( \by ) \ne 0, $
    then necessarily 
    $ \set*{ \ell } \in \cP $
    for all 
    $ \ell \in [m] $.

    Since 
    $
        \wave{ \cP } 
        \mapsto 
        \setc*{
            m + S 
        }{
            S \in \wave{ \cP }
        }
    $
    is a bijection from 
    $ \partitions{ [k-m] } $
    to 
    $ \partitions{ m + [k-m] } $,
    we see that for every 
    $
        \cP  \in  \partitions{ [k] }
    $
    such that 
    $
        \polyfam{ \cP }( \by ) \ne 0,
    $
    there exists a unique
    $
        \wave{ \cP } \in \partitions{ [k-m] }
    $
    such that
    \begin{equation}
    \label{proofeq::cor::tuple_part_tuple_part_vectors::mapping_P_to_P_wave}
        \cP 
        =
        \setc*{ m + S }{ S \in \wave{\cP} }
        \cup
        \setc*{ \set*{\ell} }{ \ell \in [m] }
        .
    \end{equation}
    For such $\cP$ and $\wave{ \cP }$ we also have
    $
        \card*{ \cP } = \big| \wave{ \cP } \big| + m.
    $
    We can draw two conclusions from this:
    \begin{enumerate}[label=Conclusion (\alph*):, ref=Conclusion (\alph*),leftmargin=*]
        \item
        \label{proofitem::cor::tuple_part_tuple_part_vectors::conclusion_P^j(y)=0}
        $ \polynum{ j }( \by ) = 0$ for all $j \in [m]$,
        \item 
        \label{proofitem::cor::tuple_part_tuple_part_vectors::existence_of_a_mapping}
        For every $j \in [k-m]$, the mapping $\cP \mapsto \wave{ \cP }$ resulting from \eqref{proofeq::cor::tuple_part_tuple_part_vectors::mapping_P_to_P_wave}
        is an injection from
        $
            \setc*{ \cP \in \partitions[ m+j ]{ [k] } }{
                \polyfam{\cP}( \by ) \ne 0
            }
        $
        to 
        $ \partitions[j]{ [k-m] }$.         
    \end{enumerate}
    Now, fix $\cP \in \partitions{ [k] }$ such that
    $
        \polyfam{ \cP }( \by ) \ne 0
    $
    and let $\wave{\cP} \in \partitions{ [m] }$
    be such that we have \eqref{proofeq::cor::tuple_part_tuple_part_vectors::mapping_P_to_P_wave}.
    Denote 
    $
        \cP' 
        \coloneq 
        \setc*{ m + S }{ S \in \wave{ \cP } }.
    $
    Then, thanks to 
    \eqref{proofeq::cor::tuple_part_tuple_part_vectors::P^l(y)_for_l_in_[m]}
    and 
    \eqref{proofeq::cor::tuple_part_tuple_part_vectors::P^m+S(y)_for_S_in_[k-m]},
    \begin{multline*}
        \polyfam{ \cP }(\by)
        =
        \prod_{S \in \cP}
            \polyset{S}(\by)
        =
        \del*{
            \prod_{S \in \cP'}
                \polyset{S}(\by)
        }
        \prod_{ \ell = 1 }^{ m }
            \polyset{ \set*{\ell} }( \by)  
    \\
        =
        \del*{
            \prod_{S \in \wave{ \cP } }
                \polyset{ m + S }(\by)  
        }
        \prod_{ \ell = 1 }^{ m }
            \polyset{ \set*{\ell} }( \by)
        =
        \del*{
            \prod_{S \in \wave{ \cP } }
                2^m \polyset{ S }(\bx)
        }
        \prod_{ \ell = 1 }^{ m }
            2^{ k - 1 } \norm{ v_{ \ell } }
    \\
        =
        2^{m \card{ \wave{\cP}} + (k-1)m}
        \polyfam{ \wave{\cP } }( \bx )
        \prod_{ \ell = 1 }^{ m }
            \norm{ v_{ \ell } }
        \le 
        2^{ (2k-1)m }
        \polyfam{ \wave{\cP } }( \bx )
        \prod_{ \ell = 1 }^{ m }
            \norm{ v_{ \ell } }.
    \end{multline*}
    %
    Therefore, thanks to \ref{proofitem::cor::tuple_part_tuple_part_vectors::existence_of_a_mapping}, for all $j \in [k-m]$,
    \begin{align}
    \label{proofeq::cor::tuple_part_tuple_part_vectors::bound_for_P^(m+j)(y)}
        \polynum{ m+j }( \by )
        &\le
        \sum_{ \cP \in \partitions[m+j]{ [k] } }
            2^{ (2k-1)m }
            \polyfam{ \wave{\cP } }( \bx )
            \prod_{ \ell = 1 }^{ m }
                \norm{ v_{ \ell } }
        \notag
    \\
        &\le 
        \sum_{ \wave{\cP} \in \partitions[j]{ [k-m] } }
            2^{ (2k-1)m }
            \polyfam{ \wave{\cP } }( \bx )
            \prod_{ \ell = 1 }^{ m }
                \norm{ v_{ \ell } }
        =
        2^{ (2k-1)m }
        \polynum{ j }( \bx )
        \prod_{ \ell = 1 }^{ m }
            \norm{ v_{ \ell } }
        .
    \end{align}
    Moreover, 
    since $m \le k-1$,
    we have
    %
    \begin{align*}
        \polynum{ s }( \by )
        &=
        \polyset{ \set*{k} }( \by)^{s-k+1}
        \prod_{ \ell = 1 }^{ k-1}
            \polyset{ \set*{ \ell } }(\by)
    \\[-0.2em]
        &=
        \del*{
            \polyset{ \set*{k} }( \by)^{s-k+1}
            \prod_{ \ell = m+1 }^{ k-1}
                \polyset{ \set*{ \ell } }(\by)
        }
        \prod_{ \ell = 1 }^{ m}
            \polyset{ \set*{\ell} }(\by)
    \\
        &=
        \del*{
            \polyset{ m + \set*{k-m} }( \by)^{s-k+1}
            \prod_{ \ell = 1 }^{ k-m-1}
                \polyset{ m + \set*{ \ell } }(\by)
        }
        \prod_{ \ell = 1 }^{ m}
            \polyset{ \set*{\ell} }(\by),
    \intertext{
        which, thanks to 
        \eqref{proofeq::cor::tuple_part_tuple_part_vectors::P^l(y)_for_l_in_[m]}
    and 
    \eqref{proofeq::cor::tuple_part_tuple_part_vectors::P^m+S(y)_for_S_in_[k-m]},
    }
        &=
        \del*{
            \del*{
                2^m \polyset{\set*{k - m}}( \bx)
            }^{s-k+1}
            \prod_{ \ell = 1 }^{ k-m-1}
                \del*{
                    2^m \polyset{\set*{ \ell }}( \bx)
                }
        }
        \prod_{ \ell = 1 }^{ m}
            \del*{
                 2^{k-1} \norm{ v_{ \ell } }
            }
    \\
        &=
        2^{ m(s+k-m-1) }
        \del*{
            \polyset{\set*{k - m}}( \bx)^{ (s-m) - (k-m) +1}
            \prod_{ \ell = 1 }^{ k-m-1 }
                \polyset{\set*{ \ell }}( \bx)
        }
        \prod_{ \ell = 1 }^{ m }
            \norm{ v_{ \ell } }
    \\
        &=
        2^{ m(s+k-m-1) }
        \polynum{ s-m }( \bx )
        \prod_{ \ell = 1 }^{ m}
            \norm{ v_{ \ell } }
    \\
        &\le 
        2^{ m(2k-1) }
        \polynum{ s-m }( \bx )
        \prod_{ \ell = 1 }^{ m}
            \norm{ v_{ \ell } }
        .
        %
    \end{align*}
    %
    In consequence, thanks to \ref{proofitem::cor::tuple_part_tuple_part_vectors::conclusion_P^j(y)=0} and \eqref{proofeq::cor::tuple_part_tuple_part_vectors::bound_for_P^(m+j)(y)},
    \begin{align*}
        \polygen{ \gen{s} }(\by)
        =
        \polynum{ s }(\by)
        + \sum_{  j =1 }^{ k-1 }
            \polynum{ j }(\by)
        &=
        \polynum{ s }(\by)
        + \sum_{  j = m + 1 }^{ k-1 }
            \polynum{ j }(\by)
    \\
        &=
        \polynum{ s }(\by)
        + \sum_{  j = 1 }^{ k-m-1 }
            \polynum{ m + j }(\by)
    \\[-0.2em]
        &\le 
        2^{ m(2k-1) }
        \polynum{ s-m }( \bx )
        \prod_{ \ell = 1 }^{  m}
            \norm{ v_{ \ell } }
        + 
        \sum_{  j = 1 }^{ k-m-1 }
                2^{ (2k-1)m }
                \polynum{ j }( \bx )
                \prod_{ \ell = 1 }^{ m }
                    \norm{ v_{ \ell } }
    \\
        &=
        2^{ (2k-1)m }
        \del*{
            \polynum{ s-m }( \bx )
            +
            \sum_{  j = 1 }^{ k-m-1 }
                \polynum{ j }( \bx )
        }
        \prod_{ \ell = 1 }^{  m }
                    \norm{ v_{ \ell } }
    \\[-0.2em]
        &=
        2^{ (2k-1)m }
        \polygen{ \gen{ s-m } }( \bx )
        \prod_{ \ell = 1 }^{  m }
            \norm{ v_{ \ell } },
    \end{align*}
    giving us the desired inequality.
    %
\end{proof}
\begin{lemma}\label{lem::adding_k+1_direction_to_poly_j_gives_poly_j+1}
    Let $(V, \| \cdot \|)$ be a normed space, 
    $k \in \bN$, 
    $s \in \intoc{k, k+1}$. 
    Let 
    $
        \bx 
        = 
        \set*{x_I}_{ I \subseteq [k+1]}
        \subseteq 
        V
    $
    and denote
    $
        \bx'
        \coloneq 
        \set*{ x_I }_{ I \subseteq [k] }
    $
    and 
    $
        \bx''
        \coloneq 
        \set*{ x_{I \cup \set*{k+1} } }_{ I \subseteq [k] }
        .
    $
    Then 
    \begin{enumerateInLemma}[label=(\roman*)]
    \item
    \label{itemlem::lem::adding_k+1_direction_to_poly_j_gives_poly_j+1::item::ell}
    	$
	        \ell( \bx' ) + \ell( \bx'' ) = \ell( \bx ) - \polyset{ \set*{k+1} }(\bx).
		$ 
    \item 
    \label{itemlem::lem::adding_k+1_direction_to_poly_j_gives_poly_j+1::item::poly_j}
        For all $j \in [k]$ we have
		$
            \polyset{ \set*{k+1} }( \bx) 
            \del*{
                \polynum{ j }( \bx' )
                + \polynum{ j }( \bx'' )
            }
            \le 
            \polynum{ j+1  }( \bx).
        $
	\item
    \label{itemlem::lem::adding_k+1_direction_to_poly_j_gives_poly_j+1::item::poly_s}
	    $
	        \polyset{ \set*{k+1}  }( \bx)^{s-k} 
	        \del*{
	            \polynum{ k }( \bx' )
	            + \polynum{ k }( \bx'' )
	        }
	        \le 
	        \polynum{ s }( \bx).
	    $
    \end{enumerateInLemma}
\end{lemma}
\begin{proof}
    First, let us note that,
    for all disjoint
    $A, S \subseteq [k]$, 
    we have
    $
        \polyset[A]{S}( \bx')
        =
        \norm*{
            \diff_{ I \upmapsto A }^{ S }
                x_I
        }
        =
        \polyset[A]{S}( \bx)
    $
    and 
    \begin{equation*}
        \polyset[A]{S}( \bx'')
        =
        \norm*{
            \diff_{ I \upmapsto A }^{ S }
                x_{I \cup \set*{k+1} }
        }
        =
        \norm*{
            \diff_{ I \upmapsto A\cup \set*{k+1} }^{ S }
                x_{I  }
        }
        =
        \polyset[A \cup \set*{k+1}]{S}( \bx).
    \end{equation*}
    Hence, for all $ S \subseteq [k]$, we have
    \begin{align*}
        \polyset{S}( \bx')
        + \polyset{S}(\bx'')
    &=
        \sum_{ A = \emptyset}^{ [k] \setminus S }
            \del*{
                \polyset[A]{S}(\bx')
                + \polyset[A]{S}(\bx'')
            }
    \\
    &=
        \sum_{ A = \emptyset}^{ [k] \setminus S }
            \del*{
                \polyset[A]{S}(\bx)
                + \polyset[A \cup \set*{k+1}]{S}(\bx)
            }
        =
        \sum_{ A = \emptyset}^{ [k+1] \setminus S }
            \polyset[A]{S}(\bx)
        =
        \polyset{S}(\bx).
    \end{align*}
    In consequence,
    \begin{equation*}
        \ell( \bx' ) 
        + \ell( \bx '' )
        =
        \sum_{ j = 1 }^{ k }
            \del*{
                \polyset{ \set*{j}} \del*{ \bx' }
                + \polyset{ \set*{j} }\del*{  \bx'' }
            }
        =
        \sum_{ j = 1 }^{ k }
            \polyset{ \set*{j}} (\bx )
        =
        \ell( \bx )
        - \polyset{ \set*{k+1} }(\bx),
    \end{equation*}
    proving \labelcref{itemlem::lem::adding_k+1_direction_to_poly_j_gives_poly_j+1::item::ell}.
    
    Next, fix $j \in [k]$ 
    and note that, 
    for all 
    $\cP \in \partitions[j]{[k]}$, 
    we have
    \begin{align*}
        \polyset{ \set*{k+1} }(\bx)
        \del*{
            \polyfam{ \cP }( \bx') 
            + \polyfam{ \cP }( \bx'')
        }
        &=
        \polyset{ \set*{k+1} }(\bx)
        \del*{
            \prod_{S \in \cP }
                \polyset{ S }(\bx')    
            +
            \prod_{S \in \cP }
                \polyset{ S }(\bx'')
        }
    \\
        &\le 
        \polyset{ \set*{k+1} }(\bx)
        \prod_{S \in \cP }
            \del*{
                \polyset{ S }(\bx')
                + \polyset{ S }(\bx'')
            }
    \\
        &= 
        \polyset{ \set*{k+1} }(\bx)
        \prod_{S \in \cP }
            \polyset{ S }(\bx)        
        =
        \prod_{ S \in \cP \cup \set*{ \set*{k+1} } }
        \! \!
            \polyset{ S }(\bx)
        =
        \polyfam{ \cP \cup \set*{ \set*{k+1} } }( \bx).
    \end{align*}
    Therefore, since
    $\cP \mapsto \cP \cup \set*{ \set*{k+1} }$
    is an injection from 
    $\partitions[j]{[k]}$ to $\partitions[j+1]{[k]}$,
    \begin{align*}
        \polyset{ \set*{k+1} }(\bx)
        \del*{
            \polynum{ j }( \bx') 
            + \polynum{ j }( \bx'') 
        }
        &=
        \sum_{ \cP \in \partitions[j]{[k]}}
            \polyset{ \set*{k+1} }(\bx)
            \del*{
                \polyfam{ \cP }( \bx') 
                + \polyfam{ \cP }( \bx'')
            }
    \\
        &\le
        \sum_{ \cP \in \partitions[j]{[k]}}
            \polyfam{ \cP \cup \set*{ \set*{k+1} } }( \bx)
        \le 
        \sum_{ \cP \in \partitions[j+1]{[k]}}
            \polyfam{ \cP }( \bx)
        =
        \polynum{ j+1 }(\bx),
    \end{align*}
    proving 
    \labelcref{itemlem::lem::adding_k+1_direction_to_poly_j_gives_poly_j+1::item::poly_j}.
    
    Finally, for the proof of  
    \labelcref{itemlem::lem::adding_k+1_direction_to_poly_j_gives_poly_j+1::item::poly_s},
    let us notice that
    \begin{align*}
        \polyset{ \set*{k+1} }( \bx)^{s-k} 
        \del*{
            \polynum{ k }( \bx' )
            + \polynum{ k }( \bx'') 
        }
		=
        \polyset{ \set*{k+1} }( \bx)^{s-k}
        \del*{
            \prod_{j=1}^{k}
                \polyset{ \set*{j} }(\bx')        
            +
            \prod_{j=1}^{k}
                \polyset{ \set*{j} }(\bx'')
        }
    &
    \\
    	\le
        \polyset{ \set*{k+1} }( \bx)^{s-k} 
        \prod_{j=1}^{k}
            \del*{
                \polyset{ \set*{j} }( \bx' )
                + \polyset{ \set*{j} }( \bx'' )
            }
    &
    \\
    	=
        \polyset{ \set*{k+1} }( \bx)^{s-k} 
        \prod_{j=1}^{k}
            \polyset{ \set*{j} }( \bx ) 
    &=
        \polynum{ s }(\bx),
    \end{align*}
    as claimed.
\end{proof}
\begin{lemma}
\label{lem::normed_if_ScS'_then_PS'<=PS}
    Let $\del*{ V, \normAlone } $ 
    be a normed space, $k \in \bN$,
    and 
    $
    	\bx = \set*{ x_I }_{ I \subseteq [k] }
    	\subseteq 
    	V.
    $
    Then for every $S, T \subseteq [k]$ 
    such that $S \subseteq T$, we have
    $
        \polyset{T}( \bx )
        \le 
        \polyset{S}( \bx ).
    $
\end{lemma}
\begin{proof}
    For every $ A \subseteq [k] \setminus T$, 
    by \Cref{lem::diff_B_diff_C_to_diff_BuC}
    we have
    \begin{align*}
        \polyset[A]{T}( \bx )
        =
        \norm*{
            \diff_{ I \upmapsto A }^{ T }
                x_I
        }
        =
        \norm*{
            \diff_{ I \upmapsto A }^{ T \setminus S}
                \diff_{ J \upmapsto I }^{ S }
                    x_J
        }
        =
        \norm*{
            \sum_{ I = \emptyset  }^{ T \setminus S}
                (-1)^{ \card{ (T \setminus S) \setminus I} }
                \diff_{ J \upmapsto I \cup A }^{ S }
                    x_J
        }
    &
    \\
        \le 
        \sum_{ I = \emptyset  }^{ T \setminus S}
        \norm*{
            \diff_{ J \upmapsto I \cup A }^{ S }
                    x_J
        }
    &=
        \sum_{ I = \emptyset  }^{ T \setminus S}
            \polyset[ I \cup A ]{ S }( \bx ).
    \end{align*}
	Hence,
    \begin{equation*}
        \polyset{T}( \bx )
        =
        \sum_{ A = \emptyset }^{ [k] \setminus T }
            \polyset[A]{T}( \bx )
        \le 
        \sum_{ A = \emptyset }^{ [k] \setminus T }
            \sum_{ I = \emptyset  }^{ T \setminus S}
                \polyset[ I \cup A ]{ S }( \bx )
        =
        \sum_{ A = \emptyset }^{ [k] \setminus S }
            \polyset[A]{ S }( \bx )
        =
        \polyset{S}( \bx ). 
        \tag*{\qedhere}
    \end{equation*}
\end{proof}
\begin{lemma}
\label{lem::folding_hypercubes_in_normed_spaces}
    Let 
    $\del*{ V, \normAlone }$ 
    be a normed space.
    Let $k \in \bN$ be such that $k \ge 2$, let $m \in [k-1]$, $i \in [m]$, 
    and $s \in \intoc{ k-1, k }$. 
    Fix
    $
        \bx = \set*{ x_I }_{ I \subseteq [m] }
        \subseteq 
        V
    $
    and define 
    $
        \by = \set*{ y_J }_{ J \subseteq [k] }
    $
    by 
    \begin{equation*}
        \forall I \subseteq [m]
        \quad 
        \forall L \subseteq [k] \setminus [m]
        \qquad 
        y_{ I \cup L }
        \coloneqq 
        \begin{cases}
            x_{I} &\text{if } \abs{ L } \text{ is even}, \\
            x_{ I \div \set*{i} } &\text{otherwise}.
        \end{cases}
    \end{equation*}
    Then 
    \begin{equation*}
        \polygen{ \gen{s} }( \by )
        \le 
        2^{k (k-m)}
        \partitions{k} \polygen{ \gen{m} }( \bx )
        \del{
            1 + \polyset{ \set*{ i} }( \bx )^{ k-m }
        }.
    \end{equation*}
\end{lemma}
\begin{proof}
    By \Cref{lem::folding_hypercubes_in_groups}, 
    for all 
    $S \subseteq [m]$, 
    $A \subseteq [m] \setminus S$,
    $S' \subseteq [k] \setminus [m]$,
    and 
    $A' \subseteq [k] \setminus ([m] \cup S')$ 
    such that $S \cup S'$ is nonempty,
    \begin{align*}
        \polyset[A \cup A']{S \cup S'}( \by )
        =
        \norm*{
            \diff_{ J \upmapsto A \cup A' }^{ S \cup S' }
                y_J
        }
    &=
        \norm*{
            2^{
                \card{
                    (S \cup S') \cap \phi^{-1}\sbr{ \set*{i} } 
                }
                - \card{
                    \phi \sbr{ S \cup S' } \cap \set*{i} 
                }
            }
            \diff_{
                I \upmapsto 
                \Psi( A \cup A' )
                \setminus \phi \sbr{ S \cup S' }
            }^{
                \phi \sbr{ S \cup S'  } 
            }
                x_I
        }
    \\
    &=
        2^{
            \card{
                (S \cup S') \cap \phi^{-1}\sbr{ \set*{i} } 
            }
            - \card{
                \phi \sbr{ S \cup S' } \cap \set*{i} 
            }
        }
        \polyset[ \Psi( A \cup A' )
                \setminus \phi \sbr{ S \cup S' } ]{ \phi \sbr{ S \cup S' } }(\bx),
    \end{align*}
    where $\phi$ and $\Psi$ are as in the mentioned lemma.
    Fix such $S$ and $S'$.
    We have three possible cases:
    $i \in S$, 
    $i \notin S$ with $S' = \emptyset$,
    and
     $i \notin S$ with $S' \ne \emptyset$.
    We will show that 
    $
        \polyset{S \cup S'}( \by ) 
        = 2^{k-m} \polyset{ \phi[ S \cup S' ] }(\bx)
    $
    in each case.
    \begin{enumerate}[label=Case \Roman*:, leftmargin=2em, itemindent=2em]
        \item
        $i \in S$.
        
        In this case, 
        $\phi[ S \cup S' ] = S$ 
        and
        $
            \Psi( A \cup A' )
                \setminus \phi[ S \cup S' ]
            =
            A 
        $
        for all $A \subseteq [m] \setminus S$ and 
        $A' \subseteq [k] \setminus ([m] \cup S')$.
        Thus,
        $
            \card{
                \phi[ S \cup S' ] \cap \set*{i} 
            }
            =
            1
        $
        and
        $
            \card{
                (S \cup S') \cap \phi^{-1}[ \set*{i} ] 
            }
            =
            \card{S'}+1.
        $
        Hence,
        \begin{align*}
            \polyset{S \cup S'}( \by )
        &=
            \sum_{A'= \emptyset}^{[k] \setminus (S' \cup [m])}
                \sum_{A= \emptyset}^{[m] \setminus S}
                    \polyset[A \cup A']{S \cup S'}( \by )
        \\
        &=
            \sum_{A'= \emptyset}^{[k] \setminus (S' \cup [m])}
                \sum_{A= \emptyset}^{[m] \setminus S}
                    2^{ \abs{S'} }   
                    \polyset[ A ]{ S }(\bx)
            =
            2^{ 
                \abs{S'} 
                + \card{ [k] \setminus (S' \cup [m]) }
            }  
            \polyset{ S }(\bx)
            =
            2^{k-m} \polyset{ S }(\bx).
        \end{align*}
        Thus,
        $
            \polyset{S \cup S'}( \by ) 
            = 2^{k-m} \polyset{ S }(\bx)
            = 2^{k-m} \polyset{ \phi[ S \cup S' ] }(\bx).
        $
        \item 
        $i \notin S$ with $S' = \emptyset$. 
        
        In this case we have $\phi[ S \cup S' ] = S$ 
        and
        \begin{equation*}
            \hspace{-0.3cm}
            \forall A \subseteq [m] \setminus S
            \quad  
            \forall A' \subseteq [k] \setminus ([m] \cup S')
            \quad \, \, 
            \Psi( A \cup A' )
                \setminus \phi[ S \cup S' ]
            =
            \begin{cases}
                A &\text{if } \card{A'} \text{ is even},\\
                A \div \set*{i} & \text{otherwise}.
            \end{cases}
        \end{equation*}
        Also,
        $
            \card{
                \phi[ S \cup S' ] \cap \set*{i} 
            }
            =
            0
        $
        and
        $
            \card{
                (S \cup S') \cap \phi^{-1}[ \set*{i} ] 
            }
            =
            0,
        $
        hence
        \begin{equation*}
            2^{
                \card{
                    (S \cup S') \cap \phi^{-1}\sbr{ \set*{i} } 
                }
                - \card{
                    \phi \sbr{ S \cup S' } \cap \set*{i} 
                }
            }
            =
            2^{0-0}
            =
            1.
        \end{equation*}
        Therefore, by 
        \Cref{lem::between_A_and_B_cardinalities},
        \begin{align*}
            \polyset{S \cup S'}( \by )
        	=
            \sum_{A'= \emptyset}^{[k] \setminus (S' \cup [m])}
                \sum_{A= \emptyset}^{[m] \setminus S}
                    \polyset[A \cup A']{S \cup S'}( \by )
        \phantom{
        	+
			\sum_{  
			 \substack{
			     A'= \emptyset
			     \\ \card{A'} \text{ is odd}
			 }
			}^{[k] \setminus [m]}
			 \sum_{A= \emptyset}^{[m] \setminus S}
			     \polyset[A \cup A']{S \cup S'}( \by )
        }
        &
        \\
        	=
            \sum_{  
                \substack{
                    A'= \emptyset
                    \\ \card{A'} \text{ is even}
                }
            }^{[k] \setminus [m]}
                \sum_{A= \emptyset}^{[m] \setminus S}
                    \polyset[A \cup A']{S \cup S'}( \by )
            +
            \sum_{  
                \substack{
                    A'= \emptyset
                    \\ \card{A'} \text{ is odd}
                }
            }^{[k] \setminus [m]}
                \sum_{A= \emptyset}^{[m] \setminus S}
                    \polyset[A \cup A']{S \cup S'}( \by )
        &
        \\
			=
            \sum_{  
                \substack{
                    A'= \emptyset
                    \\ \card{A'} \text{ is even}
                }
            }^{[k] \setminus [m]}
                \sum_{A= \emptyset}^{[m] \setminus S}
                    \polyset[A]{S}( \bx )
            +
            \sum_{  
                \substack{
                    A'= \emptyset
                    \\ \card{A'} \text{ is odd}
                }
            }^{[k] \setminus [m]}
                \sum_{A= \emptyset}^{[m] \setminus S}
                    \polyset[A \div \set*{i}]{S}( \bx )
        &
        \\
        	=
            2^{k - m - 1}
            \del*{
                \sum_{A= \emptyset}^{[m] \setminus S}
                    \polyset[A]{S}( \bx )
                +
                \sum_{A= \emptyset}^{[m] \setminus S}
                    \polyset[A \div \set*{i}]{S}( \bx )
            }
        &=
            2^{k - m}
            \polyset{S}( \bx ),
        \end{align*}
        where in the last equality we used the fact that since
        $A \mapsto A \div \set*{i}$
        is a bijection from 
        $2^{ [m] \setminus S }$ 
        to itself, 
        we have
        $
            \sum_{A= \emptyset}^{[m] \setminus S}
                \polyset[A \div \set*{i}]{S}( \bx )
            =
            \sum_{A= \emptyset}^{[m] \setminus S}
                \polyset[A]{S}( \bx )
            =
            \polyset{S}( \bx ).
        $
        Therefore, we again have
        $
            \polyset{S \cup S'}( \by ) 
            = 2^{k-m} \polyset{ S }(\bx)
            = 2^{k-m} \polyset{ \phi[ S \cup S' ] }(\bx).
        $
        \item 
        $i \notin S$ with $S' \ne \emptyset$.
        
        In the last case we have $\phi[ S \cup S' ] = S \cup \set*{i}$ 
        and
        $
            \Psi( A \cup A' )
                \setminus \phi[ S \cup S ]
            =
            A \setminus \set*{i} 
        $
        for all $A \subseteq [m] \setminus S$ and 
        $A' \subseteq [k] \setminus ([m] \cup S')$.        
        Also,
        $
            \card{
                \phi[ S \cup S' ] \cap \set*{i} 
            }
            =
            1
        $
        and
        $
            \card{
                (S \cup S') \cap \phi^{-1}[ \set*{i} ] 
            }
            =
            \card{S'},
        $
        hence 
        \begin{equation*}
            2^{
                \card{
                    (S \cup S') \cap \phi^{-1}\sbr{ \set*{i} } 
                }
                - \card{
                    \phi \sbr{ S \cup S' } \cap \set*{i} 
                }
            }
            =
            2^{ \abs{S'} - 1}.
        \end{equation*}
        Thus,
        \begin{multline*}
            \polyset{S \cup S'}( \by )
           =
            \sum_{A'= \emptyset}^{[k] \setminus (S' \cup [m])}
                \sum_{A= \emptyset}^{[m] \setminus S}
                    \polyset[A \cup A']{S \cup S'}( \by )
        	=
            \sum_{A'= \emptyset}^{[k] \setminus (S' \cup [m])}
                \sum_{A= \emptyset}^{[m] \setminus S}
                    2^{ \abs{S'} - 1}
                    \polyset[A \setminus \set*{i} ]{S \cup \set*{i}}( \bx )
        \\
        	=
            2^{ 
                \card{ [k] \setminus (S' \cup [m]) }
            }
            \t
            2^{
                \abs{S'} - 1
            }
            \t 
            2
            \polyset{S \cup \set*{i} }(\bx)
        	=
            2^{k-m} \polyset{S \cup \set*{i}}(\bx),
        \end{multline*}
        where we used the fact that
        \begin{align*}
            \sum_{A = \emptyset }^{ [m] \setminus S }
                \polyset[A \setminus \set*{i}]{ S \cup \set*{i} }(\bx)
        &=
            \sum_{A = \emptyset }^{ [m] \setminus (S \cup \set*{i}) }
                \del*{
                    \polyset[A \setminus \set*{i}]{ S \cup \set*{i} }(\bx)
                    +
                    \polyset[(A \cup \set*{i}) \setminus \set*{i}]{ S \cup \set*{i} }(\bx)
                }
    	\\
        &=
            2 \sum_{A = \emptyset }^{ [m] \setminus (S \cup \set*{i}) }
                \polyset[A \setminus \set*{i}]{ S \cup \set*{i} }(\bx)
        	=
        	2\sum_{A = \emptyset }^{ [m] \setminus (S \cup \set*{i}) }
        		\polyset[A]{ S \cup \set*{i} }(\bx)
        	=
            2\polyset{S \cup \set*{i}}( \bx).
        \end{align*}
        Therefore,
        $
            \polyset{S \cup S'}( \by ) 
            = 2^{k-m} \polyset{ S \cup \set*{i} }(\bx)
            = 2^{k-m} \polyset{ \phi[ S \cup S' ] }(\bx),
        $
        as claimed.
    \end{enumerate}
    We have shown that, 
    for all $S \subseteq[m]$ and 
    $S' \subseteq [k] \setminus [m]$
    such that $S \cup S'$ is nonempty, 
    we have
    $
        \polyset{S \cup S'}(\by)
        =
        2^{k-m}\polyset{ \phi[S \cup S']}( \bx).
    $
    Since we have 
    $ \phi[S \cup S'] = S $
    or $ \phi[S \cup S'] = S \cup \set*{i}$,
    using \Cref{lem::normed_if_ScS'_then_PS'<=PS}
    we get
    $
        \polyset{S \cup S'}(\by)
        \le
        2^{k-m}\polyset{ S}( \bx).
    $    
    Also, if 
    $S = \emptyset$, 
    then we have 
    $S' \ne \emptyset$, 
    so
    $
        \polyset{S \cup S'}(\by)
        =
        2^{k-m}\polyset{ \set*{i} }( \bx).
    $   
    Next, fix $\cP' \in \partitions{[k]}$
    and define $\mathcal{A}, \mathcal{B}$ 
    as follows:
    \begin{equation*}
        \mathcal{A}
        \coloneq
        \setc*{
            S' \in \cP'
        }{
            S' \cap [m] \ne \emptyset
        }
        \quad \text{and} \quad 
        \mathcal{B} 
        \coloneq
        \cP' \setminus \mathcal{A}.
    \end{equation*}
    Then, define
    $
        \cP
        \coloneq 
        \setc*{
            S' \cap [m]
        }{
            S' \in \mathcal{A}
        }.
    $
    Since $\cP'$ is a partition of $[k]$ and $[m] \subseteq [k]$,
    it follows that $\cP$ is~a~partition of $[m]$.
    Moreover,
    \begin{align*}
        \polyfam{\cP'}(\by)
        =
        \prod_{S' \in \cP'}
            \polyset{S'}(\by)
    &=
        \prod_{S' \in \mathcal{A}}
            \polyset{S'}(\by)
        \cdot 
        \prod_{S' \in \mathcal{B}}
            \polyset{S'}(\by)
    \\
    &\le
        \prod_{S \in \cP}
            \del*{
                2^{k-m}
                \polyset{S}(\bx)
            }
        \cdot 
        \prod_{S' \in \mathcal{B}}
            \del*{
                2^{k-m}
                \polyset{\set*{i}}(\bx)
            }
    \\
    &=
        2^{ \abs{\cP'}(k-m)}
            \polyset{\cP}(\bx)
            \polyset{\set*{i}}(\bx)^{\card{\mathcal{B}}}
        \le 
        2^{k(k-m)}
        \polygen{ \gen{m} }(\bx)
        \del*{
            1 + \polyset{\set*{i}}(\bx)^{k-m}
        },
    \end{align*}
    where in the last inequality we have used 
    \Cref{ex::special_values_of_P_s_k::item::polyen_of_k}
    and the fact that since
    $\mc{B}$ is a disjoint family of~nonempty subsets of 
    $[k] \setminus [m]$, we have 
    $ \card*{ \mc{B} } \le k-m$.
    Therefore, for all $j \in [k-1]$,
    \begin{align*}
        \polynum{ j }(\by)
        =
        \sum_{ \cP' \in \partitions[j]{[k]}}
            \polyfam{\cP'}(\by)
    &\le 
        \sum_{ \cP' \in \partitions[j]{[k]}}
            2^{k(k-m)}
            \polygen{ \gen{m} }(\bx)
            \del*{
                1 + \polyset{\set*{i}}(\bx)^{k-m}
            }
    \\
    &=
        2^{k(k-m)}
        \card*{ \partitions[j]{[k]}}
        \polygen{ \gen{m} }(\bx)
            \del*{
                1 + \polyset{\set*{i}}(\bx)^{k-m}
            }.
    \end{align*}
    Also, since 
    $
        \polyset{\set*{j}}(\by) 
        \le 
        2^{k-m} \polyset{\set*{j}}(\bx)
    $ 
    for $j \in [m]$
    and 
    $
        \polyset{\set*{j}}(\by) 
        =
        2^{k-m} \polyset{\set*{i}}(\bx)
    $ 
    for 
    $j \in [k] \setminus [m]$, 
    \begin{align*}
        \polynum{ s }(\by)
    &=
        \polyset{\set*{k}}(\by)^{s-k+1}
        \prod_{ j = 1 }^{ k- 1}
            \polyset{\set*{j}}(\by),
    \intertext{
    	which, since $m \le k-1$,
    }
    &=
        \polyset{\set*{k}}(\by)^{s-k+1}
        \del*{
            \prod_{ j = 1 }^{ m}
                \polyset{\set*{j}}(\by)
        }
        \del*{
            \prod_{ j = m+1 }^{ k- 1}
                \polyset{\set*{j}}(\by)
        }
    \\
    &\le 
        \del*{ 2^{k-m} \polyset{\set*{i}}(\bx) }^{s-k+1}
        \del*{
            \prod_{ j = 1 }^{ m}
                2^{k-m} \polyset{\set*{j}}(\bx) 
        }
        \del*{
            \prod_{ j = m+1 }^{ k- 1}
                2^{k-m} \polyset{\set*{i}}(\bx) 
        }
    \\
    &=
        2^{s(k-m)}
        \polynum{ m }(\bx)
        \polyset{\set*{i}}(\bx)^{s-m}
    \\
    &\le 
        2^{k(k-m)}
        \polygen{ \gen{m} }(\bx)
        \del*{
            1 + \polyset{\set*{i}}(\bx)^{k-m}
        }
        =
        2^{k(k-m)}
        \card*{ \partitions[k]{[k]} }
        \polygen{ \gen{m} }(\bx)
        \del*{
            1 + \polyset{\set*{i}}(\bx)^{k-m}
        }.
    \end{align*}
    Finally, we get
    \begin{align*}
        \polygen{ \gen{s} }(\by)
    ={}&
        \polynum{ s }(\by)
        +
        \sum_{j=1}^{k-1}
            \polynum{ j }(\by)      
    \\    
    \le{}& 
        2^{k(k-m)}
        \card*{\partitions[k]{[k]}}
        \polygen{ \gen{m} }(\bx)
        \del*{
            1 + \polyset{\set*{i}}(\bx)^{k-m}
        }
    \\
    &+
        \sum_{j=1}^{k-1}
            2^{k(k-m)}
            \card*{ \partitions[j]{[k]}}
            \polygen{ \gen{m} }(\bx)
            \del*{
                1 + \polyset{\set*{i}}(\bx)^{k-m}
            }   
    \\
    ={}&
        2^{k(k-m)}
        \del*{
            \sum_{j=1}^{k}
                \card*{ \partitions[j]{[k]}}
        }
        \polygen{ \gen{m} }(\bx)
        \del*{
            1 + \polyset{\set*{i}}(\bx)^{k-m}
        }   
    \\
    ={}&
        2^{k (k-m)}
        \partitions{k} \polygen{ \gen{m} }( \bx )
        \del*{
            1 + \polyset{ \set*{ i} }( \bx )^{ k-m }
        },
    \end{align*} 
    where the last equality follows from the fact that
    $
        \sum_{j=1}^{k}
            \card*{ \partitions[j]{[k]}}
        =
        \card*{ \partitions{ [k] } }
        =
        \partitions{k}.
    $
\end{proof}
\begin{lemma}\label{lem::estimate_cube_with_points_and_vectors}
    Let $(V, \| \cdot \|)$ be a normed space. 
    Let $k \in \bN$
    and $S \subseteq [k]$ be nonempty.
    Fix
    $
        \bx = \set*{ x_I }_{ I \subseteq [k+1] }
        \subseteq 
        V
    $
    and define 
    $
        \set*{ v_J }_{ J \subseteq [k] }
    $
    by~the~formula
    $
        v_J \coloneqq x_{ J \cup \set*{ k+1 } } - x_J
    $
    for all $J \subseteq [k]$.
    Finally, define
    $
        \by = \set*{ y_L }_{ L \subseteq [k] }
    $
    by~the~formula
    $
        y_L = v_{L \cap S} + x_{L \cup S}
    $
    for all $L \subseteq [k]$.
    Then 
    \begin{gather}
    \label{eqlem::lem::estimate_cube_with_points_and_vectors::ell}
        \ell( \by )
        \le 
        2^k
        \ell( \bx )
	\shortintertext{and}
    \label{eqlem::lem::estimate_cube_with_points_and_vectors::poly}
        \forall j \in [k]
    \qquad 
        \polynum{ j }( \by ) 
        \le 
        2^{k^2} \partitions{k} 
        \polynum{ j }(\bx).
    \end{gather}
    %
    %
    %
\end{lemma}
\begin{proof}
    We will start the proof by showing that
    \begin{align}
        \label{proofeq::lem::estimate_cube_with_points_and_vectors::S'cS}
        \forall S' \subseteq [k] 
    &\qquad 
    	\del*{
	        S' \ne \emptyset
	        \text{ and }
	        S' \subseteq S
	    }
    \implies 
        \polyset{ S' }( \by )
        \le 
        2^k
        \polyset{S'\cup \set*{k+1}}(\bx)
    \intertext{
        and
    }
        \label{proofeq::lem::estimate_cube_with_points_and_vectors::S'outsideS}
        \forall S' \subseteq [k] 
    &\qquad 
    	\del*{
	        S' \ne \emptyset
	        \text{ and }
	        S' \subseteq [k] \setminus S
	    }
    \implies 
        \polyset{ S' }( \by )
        \le 
        2^k
       \polyset{S' }(\bx).
    \end{align}
    \begin{itemize}[leftmargin=1em, itemindent=1em]
        \item
        Suppose first that
        $ S' \subseteq S$ is nonempty and 
        fix $A \subseteq [k] \setminus S'$.
        Then
        \begin{equation*}
            \diff_{ L \upmapsto A}^{S'}
                x_{L \cup S}
            =
            \sum_{ L = \emptyset }^{ S' }
                (-1)^{ \card{S' \setminus L } }
                x_{L \cup A \cup S}
            =
            \sum_{ L = \emptyset }^{ S' }
                (-1)^{ \card{S' \setminus L } }
                x_{A \cup S}
            =
            \diff_{ L \upmapsto A }^{ S' }
            	x_{A \cup S }
            =
            0,
        \end{equation*}
        where the final equality follows from
        \Cref{lem::diff_constant_term}.
        Next, note that since $S' \subseteq S$,
        \begin{equation*}
            \diff_{ L \upmapsto A}^{S'}
                v_{L \cap S}
            =
            \sum_{ L = \emptyset }^{ S' }
                (-1)^{ \card{S' \setminus L } }
                v_{(L \cup A) \cap S}
            =
            \sum_{ L = \emptyset }^{ S' }
                (-1)^{ \card{S' \setminus L } }
                v_{L \cup \del*{ A \cap S} }
            =
            \diff_{ L \upmapsto A \cap S}^{S'}
                v_{L}.
        \end{equation*}
        Hence, 
        \begin{align*}
            \polyset[A]{S'}(\by)
            =
            \norm*{
                \diff_{ L \upmapsto A }^{S'}
                    y_L
            }
            =
            \norm*{
                \diff_{ L \upmapsto A }^{S'}
                    \del*{ v_{L \cap S} + x_{L \cup S} }
            }
        &
        \\=
            \norm*{
                \diff_{ L \upmapsto A \cap S }^{S'}
                    v_{L} 
            }
        &=
            \norm*{
                \diff_{ L \upmapsto A \cap S }^{S'}
                    \del*{
                        x_{L \cup \set*{k+1}}
                        - x_{L }
                    }
            }
        \\
    	&=
            \norm*{
                \diff_{ L \upmapsto A \cap S }^{S' \cup \set*{k+1}}
                    x_L
            }
            =
            \polyset[A \cap S]{S'\cup \set*{k+1}}(\bx),
        \end{align*}
        where we used the
        \say{moreover} part of
        \Cref{lem::diff_B_diff_C_to_diff_BuC} in the penultimate equality.
        Therefore,
        \begin{multline*}
            \polyset{S'}(\by)
            =
            \sum_{ A = \emptyset }^{ [k] \setminus S' }
                \polyset[A]{S'}(\by)
            =
            \sum_{ A = \emptyset }^{ [k] \setminus S' }
                \polyset[A \cap S]{S'\cup \set*{k+1}}(\bx)
        \\
            =
            \sum_{ A = \emptyset }^{ [k] \setminus S }
                \sum_{ B = \emptyset }^{S \setminus S' }
                    \polyset[(A\cup B) \cap S]{S'\cup \set*{k+1}}(\bx)
            =
            \sum_{ A = \emptyset }^{ [k] \setminus S }
                \sum_{ B = \emptyset }^{S \setminus S' }
                    \polyset[B]{S'\cup \set*{k+1}}(\bx)
            =
            2^{k - \card{S} }
            \sum_{ B = \emptyset }^{S \setminus S' }
                \polyset[B]{S'\cup \set*{k+1}}(\bx)   
        \\
            \le 
            2^k
            \sum_{ B = \emptyset }^{[k+1] \setminus (S' \cup \set*{k+1} ) }
                \polyset[B]{S'\cup \set*{k+1}}(\bx) 
            =
            2^k
            \polyset{S'\cup \set*{k+1}}(\bx),
        \end{multline*}
        where in the third equality we
        used the fact that, 
        since $S' \subseteq S \subseteq [k]$,
		the map
        $(A, B) \mapsto A \cup B$ is a bijection between
        $2^{[k] \setminus S} \times 2^{S \setminus S'}$
        and $2^{ [k] \setminus S' }$.
        Thus, we have
        $
            \polyset{S'}(\by)
            \le 
            2^k
            \polyset{S'\cup \set*{k+1}}(\bx),
        $
        as claimed.
        %
        \item 
        Next, suppose that $S' \subseteq [k] \setminus S$ is nonempty and fix
        $A \subseteq [k] \setminus S'$.
        This time, by \Cref{lem::diff_constant_term} we have
        \begin{equation*}
            \diff_{ L \upmapsto A}^{S'}
                v_{L \cap S}
            =
            \sum_{ L = \emptyset }^{ S' }
                (-1)^{ \card{S' \setminus L } }
                v_{(L \cup A) \cap S}
            =
            \sum_{ L = \emptyset }^{ S' }
                (-1)^{ \card{S' \setminus L } }
                v_{A \cap S}
            =
            \diff_{ L \upmapsto A }^{ S' }
            	v_{A \cap S }
            =
            0.
        \end{equation*}
        Next, note that
        \begin{equation*}
            \diff_{ L \upmapsto A}^{S'}
                x_{L \cup S}
            =
            \sum_{ L = \emptyset }^{ S' }
                (-1)^{ \card{S' \setminus L } }
                x_{L \cup A \cup S}
            =
            \diff_{ L \upmapsto A \cup S}^{S'}
                x_{L}.
        \end{equation*}
        Hence,
        \begin{equation*}
            \polyset[A]{S'}(\by)
            =
            \norm*{
                \diff_{ L \upmapsto A }^{S'}
                    y_L
            }
            =
            \norm*{
                \diff_{ L \upmapsto A }^{S'}
                    \del{ v_{L \cap S} + x_{L \cup S} }
            }
            =
            \norm*{
                \diff_{ L \upmapsto A \cup S }^{S'}
                    x_{L} 
            }
            =
            \polyset[A \cup S]{S'}(\bx).
        \end{equation*}
        Therefore, 
        \begin{multline*}
            \polyset{S'}(\by)
            =
            \sum_{A = \emptyset}^{ [k] \setminus S'}
                \polyset[A]{S'}(\by)
            =
            \sum_{A = \emptyset}^{ [k] \setminus S'}
                \polyset[A \cup S]{S'}(\bx)
        \\
            =
            \sum_{A = \emptyset}^{ [k] \setminus (S' \cup S) }
                \sum_{ B = \emptyset }^{S}
                    \polyset[A \cup B \cup S]{S'}(\bx)
            =
            \sum_{A = \emptyset}^{ [k] \setminus (S' \cup S) }
                \sum_{ B = \emptyset }^{S}
                    \polyset[A \cup S]{S'}(\bx)
            =
            2^{ \card{S} }
            \sum_{A = \emptyset}^{ [k] \setminus (S' \cup S) }
                \polyset[A \cup S]{S'}(\bx)
        \\
            =
            2^{ \card{S} }
            \sum_{A = S}^{ [k] \setminus S'}
                \polyset[A]{S'}(\bx)
            \le 
            2^k
            \sum_{A = \emptyset}^{ [k+1] \setminus S'}
                \polyset[A]{S'}(\bx)
            =
            2^k 
            \polyset{S'}(\bx),
        \end{multline*}
        where in the third equality we have used the fact that 
        since $S'$ and $S$ are disjoint subsets of $[k]$,
        the mapping
        $(A, B) \mapsto A \cup B$ is a bijection 
        from 
        $
            2^{[k] \setminus (S' \cup S)} 
            \times 
            2^{S}
        $
        to 
        $
            2^{[k] \setminus S'}.
        $
        Thus, we have
        $
            \polyset{S'}(\by)
            \le  
            2^k
            \polyset{S'}(\bx),
        $
        as claimed.
    \end{itemize}

    We will now prove \eqref{eqlem::lem::estimate_cube_with_points_and_vectors::ell}. 
    Let us notice that
    \begin{align*}
        \ell( \by )
    &=
        \sum_{ j = 1 }^{ k }
            \polyset{ \set*{j} }( \by )
    =
        \sum_{ \substack{ 
            j = 1 \\ j \in S } }^{ k }
            \polyset{ \set*{j} }( \by )
        +
        \sum_{ \substack{ 
            j = 1 \\ j \notin S } }^{ k }
            \polyset{ \set*{j} }( \by )
    \\
        &\le 
        \sum_{ \substack{ 
            j = 1 \\ j \in S } }^{ k }
            2^k
            \polyset{ \set*{j} \cup \set*{k+1} }( \bx )
        +
        \sum_{ \substack{ 
            j = 1 \\ j \notin S } }^{ k }
            2^k
            \polyset{ \set*{j} }( \bx )
    \\
        &\le 
        \sum_{ \substack{ 
            j = 1 \\ j \in S } }^{ k }
            2^k
            \polyset{ \set*{j} }( \bx )
        +
        \sum_{ \substack{ 
            j = 1 \\ j \notin S } }^{ k }
            2^k
            \polyset{ \set*{j} }( \bx )
        =
        2^k  \! \hspace{-0.25pt}
        \sum_{ j = 1 }^{ k }
            \polyset{ \set*{j} }( \bx )
        \le 
        2^k \! \hspace{-0.25pt}
        \sum_{ j = 1 }^{ k+1 } 
            \polyset{ \set*{j} }( \bx )
        =
        2^k \ell( \bx ),
    \end{align*}
    where in the first inequality we have used 
    \eqref{proofeq::lem::estimate_cube_with_points_and_vectors::S'cS}
    and 
    \eqref{proofeq::lem::estimate_cube_with_points_and_vectors::S'outsideS},
    while in the second, 
    \Cref{lem::normed_if_ScS'_then_PS'<=PS}.
    
    Next, we will prove \eqref{eqlem::lem::estimate_cube_with_points_and_vectors::poly}.
    Fix $j \in [k]$.
    We will now show that 
    $
        \polyfam{\cP}(\by)
        \le 
        2^{k^2} \polynum{ j }(\bx)
    $
    for all $\cP \in \partitions[j]{[k]}$.
    Fix $\cP \in \partitions[j]{[k]}$.
    We have two cases to consider:
	\begin{itemize}[leftmargin=1em, itemindent=1em, listparindent=\parindent, itemsep=0em, parsep=0em]
        \item
        Suppose there is $S' \in \cP$ such that 
        $S' \cap S \ne \emptyset$ 
        and 
        $S' \cap \del*{  [k] \setminus S } \ne \emptyset$.
        Clearly, $\set*{S, [k] \setminus S}$ is a partition of $[k]$.
        For $J \subseteq [k] \setminus S$, 
        define
        $
            w_J \coloneqq x_{J \cup S}.
        $
        Then, for all $L \subseteq [k]$,
        \begin{equation*}
            y_L
            = 
            v_{L \cap S} + x_{L \cup S}
            =
            v_{L \cap S} + x_{(L \cap ([k] \setminus S ))\cup S }
            =
            v_{L \cap S} + w_{ L \cap ([k] \setminus S ) }.
        \end{equation*}
        Hence, $\polyfam{\cP}(\by) = 0$ by \Cref{lem::poly_partitioned_gives_0}.
        Thus,
        $
            \polyfam{\cP}(\by)
            \le 
            2^{k^2} \polynum{ j }(\bx)
        $
        in this case.
        \item
        Suppose there is no $S' \in \cP$ such that 
        $S' \cap S \ne \emptyset$ 
        and 
        $S' \cap \del*{ [k] \setminus S } \ne \emptyset$.
        Hence, for all $S' \in \cP$ we either have
        $S' \subseteq S$ 
        or 
        $S' \subseteq [k] \setminus S$.
        Denote
        \begin{equation*}
            \mathcal{A} 
            \coloneq
            \setc*{
                S' \in \cP 
            }{
                S' \subseteq S
            }
            \quad \text{and} \quad 
            \mathcal{B}
            \coloneq
            \setc*{
                S' \in \cP 
            }{
                S' \subseteq [k] \setminus S
            }.
        \end{equation*}
        By the assumption of this case,
        $\cP = \mathcal{A} \cup \mathcal{B}$.
        Moreover, since $S \ne \emptyset$, we have
        $\mathcal{A} \ne \emptyset$.
        
        Fix $S'' \in \mathcal{A}$
        and define
        \begin{equation*}
            \cP'
            \coloneq 
            \set*{ S'' \cup \set*{k+1} }
            \cup \del{ \mathcal{A} \setminus \set*{S''} }
            \cup \mathcal{B}.
        \end{equation*}
        Since $\cP \in \partitions[j]{[k]}$
        and 
        $\cP = \mathcal{A} \cup \mathcal{B}$, 
        we have
        $\cP' \in \partitions[j]{[k+1]}$.
        
        From 
        \eqref{proofeq::lem::estimate_cube_with_points_and_vectors::S'cS}, 
        for all 
        $S' \in \mathcal{A}$ we have
        $ 
            \polyset{S'}(\by)
            \le 
            2^k \polyset{S' \cup \set*{k+1}}(\bx).
        $
        From 
        \eqref{proofeq::lem::estimate_cube_with_points_and_vectors::S'outsideS},
        it~follows that for all 
        $S' \in \mathcal{B}$ we have
        $
            \polyset{S'}(\by)
            \le
            2^k \polyset{S'}(\bx).
        $     
        By \Cref{lem::normed_if_ScS'_then_PS'<=PS}, 
        for all 
        $S' \in \mathcal{A}$ 
        we also have
        \begin{equation*}
            \polyset{S'}(\by)
            \le 
            2^k \polyset{S' \cup \set*{k+1}}(\bx)
            \le 
            2^k \polyset{S'}(\bx).
        \end{equation*}
        Therefore,
        \begin{multline*}
            \polyfam{\cP}(\by)
            =
            \prod_{S' \in \cP}
                \polyset{S'}(\by)
            =
            \polyset{S''}(\by)
            \del*{
                \prod_{
                    \substack{
                        S' \in \mathcal{A} \\ 
                        S' \ne S''
                    }
                }   
                    \polyset{S'}(\by)
            }
            \del*{
                \prod_{ S' \in \mathcal{B}}
                    \polyset{S'}(\by)
            }
        \\
            \le 
            2^{k} 
            \polyset{S'' \cup \set*{k+1}}(\bx)
            \del*{
                \prod_{
                    \substack{
                        S' \in \mathcal{A} \\ 
                        S' \ne S''
                    }
                }   
                    2^k\polyset{S'}(\bx)
            }
            \del*{
                \prod_{ S' \in \mathcal{B}}
                    2^k\polyset{S'}(\bx)
            }
        \\
            =
            2^{kj}
            \prod_{S' \in \cP'}
                \polyset{S'}(\bx)
            =
            2^{kj}
            \polyfam{\cP'}(\bx)
            \le 
            2^{k^2}
            \polynum{ j }(\bx),
        \end{multline*}
        where the last inequality follows from 
        the fact that $j \le k$ and 
        $\cP' \in \partitions[j]{[k+1]}$.
        Thus, we again have 
        $
            \polyfam{\cP}(\by)
            \le 
            2^{k^2}
            \polynum{ j }(\bx).
        $
    \end{itemize}
    Finally, we end the proof by noting that
    \begin{equation*}
        \polynum{ j }(\by)
        =
        \sum_{ \cP \in \partitions[j]{[k]} }
            \polyfam{\cP} (\by)
        \le 
        \sum_{ \cP \in \partitions[j]{[k]} }
            2^{k^2}
            \polynum{ j }(\bx)
    	=
        2^{k^2} \card*{ \partitions[j]{[k]} }
        \polynum{ j }(\bx)
    	\le 
        2^{k^2} \partitions{k}
        \polynum{ j }(\bx).
        \tag*{\qedhere}
    \end{equation*}
\end{proof}
\begin{lemma}
\label{lem::k+1_direction_<=1_gives_poly(s)<=poly(t)_for_t<=s}
    Let $k \in \bN_0$ and $s, t \in \intoc{k,k+1}$
    be such that $t \le s$.
    Let $(V, \| \cdot \|)$ be a normed space
    and 
    $
        \bx = \set*{ x_I }_{ I \subseteq [k + 1] }
        \subseteq 
        V.
    $
    If 
    $\polyset{ \set*{k+1} }\del*{ \bx } \le 1$,
    then
    $
    	\polygen{ \gen{s} }\del*{ \bx } 
    	\le 
    	\polygen{ \gen{t} }\del*{ \bx }.
    $
\end{lemma}
\begin{proof}
    Since $\polyset{ \set*{k+1} }(\bx) \le 1$,
    we have
    $
        \polyset{ \set*{k+1} }(\bx)^{ s-k}
        \le 
        \polyset{ \set*{k+1} }(\bx)^{ t-k}.
    $
    Hence,
    \begin{equation*}
        \polynum{ s }( \bx )
        =
        \polyset{ \set*{k+1} }( \bx )^{s-k}
        \prod_{ j = 1 }^{ k }
            \polyset{ \set*{j} }( \bx )
        \le 
        \polyset{ \set*{k+1} }( \bx )^{ t-k}
        \prod_{ j = 1 }^{ k }
            \polyset{ \set*{j} }( \bx )
        =
        \polynum{ t }(\bx).
    \end{equation*}
    Therefore, 
    \begin{equation*}
        \polygen{ \gen{s} }( \bx )
        =
        \polynum{ s }( \bx )
        +
        \sum_{ j = 1 }^{ k }
            \polynum{ j }( \bx )
        \le 
        \polynum{ t }( \bx )
        +
        \sum_{ j = 1 }^{ k }
            \polynum{ j }( \bx )
        =
        \polygen{ \gen{t} }( \bx ).
        \tag*{\qedhere}
    \end{equation*}
    %
    %
\end{proof}

\section{The Main Results}
\label{sec::characterizations}
\subsection{Results for higher-order local Sobolev spaces}
In this section, 
we will prove a rather general characterization of higher-order 
Sobolev
spaces. 
We will first show some partial results in terms of local Sobolev spaces and locally integrable functions.
Thanks to these results, 
there will effectively be no additional cost in expressing our characterization in terms of Sobolev spaces based on Banach function spaces 
rather than just in terms of the usual $\sobolev{k,p}$ spaces as in \Cref{thm::sobolev_characterization}.

We will assume that all norms on $\bR^n$ used in this section are the Euclidean ones, unless specified otherwise.
However, since on finite-dimensional normed spaces all norms are equivalent, 
a~different choice of norms would only change the constants present in some inequalities.

We begin our discussion by recalling a rather well-known inequality whose version appeared in \cite{Hajlasz}.
Let us note that the inequality is often cited as being satisfied almost everywhere and with $\mathsf{M}_{\norm{x-y}}$ instead~of~%
$\mathsf{M}_{2\norm{x-y}}$ as we do below 
(and was already done, for example, in~%
\cite{hajlas_bojarski_inequality_with_2}).
For us, it will be indispensable for the inequality to~be~satisfied for all Lebesgue points of $f$.
Because of this, as well as for the sake of completeness, we have decided to~include the proof of this admittedly weaker form of the mentioned inequality.
Let us also note that other related results appeared in the literature
\cite{bojarski,Bojarski_Hajlasz, boundary_Hajlasz,sobolev_met_poincare}; in particular, \cite[Theorem 3.2]{sobolev_met_poincare} can be considered a generalization of the below inequality, and has a rather similar proof.
\begin{lemma}
    \label{lem::bojarski_inequality}
    Let $n \in \bN$.
    Then for every 
    $f \in \locSobolev{1,1}(\bR^n)$,
    we have
    \begin{equation*}
        \forall x, y \in \lebesguePoints{f}
        \qquad 
        \abs{
            f(x) - f(y)
        }
        \le 
        C_M \,
        \norm{ x - y }
        \del{
            \restMaxFun{2 \norm{x-y} }{\nabla f}(x)
            + \restMaxFun{2 \norm{x-y} }{\nabla f}(y)
        },
    \end{equation*}
    where 
    $
        C_M
        \coloneq 
        2^{n+3} C_P,
    $
    where $C_P$ is the constant from the Poincar\'{e} inequality for a ball.
\end{lemma}
\begin{proof}
    Fix $f \in \locSobolev{1,1}(\bR^n)$.
    First, let us notice that the inequality is true 
    if $x, y \in \lebesguePoints{f}$
    are such that
    $x = y$. 
    Indeed, in this case we have
    $
    	\abs*{ f(x)  -f(y) } = 0
    $
    and the inequality follows.
    
    Next, we will first show that, 
    for all $x \in \bR^n$ and 
    all $R > 0$, 
    we have 
    \begin{equation} \label{proofeq::lem::bojarski_characterization::difference_of_two_averages}
        \abs{
            \ballAv{R}{f}(x) - \ballAv{2R}{f}(x)
        }
        \le 
        2^{n+1} C_P R \,
        \restMaxFun{2R}{\nabla f}(x).
    \end{equation}
    Indeed, we have
    \begin{align*}
        \abs{
            \ballAv{R}{f}(x) - \ballAv{2R}{f}(x)
        }
        &=
        \abs*{
            \sintegral{
                \ball{ x, R}
            }{f(y)}{y} - \ballAv{2R}{f}(x)
        }
        \\
        &\le 
        \sintegral{
            \ball{ x, R}
        }{
            \abs{ f(y) - \ballAv{2R}{f}(x) }
        }{y}
        \\
        &\le
        \frac{
            \abs{
                \ball{x, 2R}
            }
        }{
            \abs{
                \ball{x, R}
            }
        }
        \sintegral{
            \ball{ x, 2R}
        }{
            \abs{ f(y) - \ballAv{2R}{f}(x) }
        }{y}
        \\
        &=
        2^n
        \sintegral{
            \ball{ x, 2R}
        }{
            \abs{ f(y) - \ballAv{2R}{f}(x) }
        }{y}
        \\
        &\le 
        2^{n+1} C_P R 
        \sintegral{
            \ball{ x, 2R}
        }{
            \norm{ \nabla f (y) }
        }{y}
        \\
        &\le 
        2^{n+1} C_P R \,
        \restMaxFun{2R}{ \nabla f }(x),
    \end{align*}
    where we have used the Poincar\'{e} inequality in the penultimate inequality.
    Therefore, for all $x \in \lebesguePoints{f}$ and $R > 0$,
    we have
    \begin{align*}
        \abs{
            f(x) - \ballAv{R}{f}(x)
        }
        &=
        \abs*{
            \lim_{m \to \infty}
            	\del*{
	                \ballAv{2^{-m} R}{f}(x)
	                - \ballAv{R}{f}(x)
	            }
        }
    \\
        &=
        \abs*{
            \lim_{m \to \infty}
            \sum_{j = 1 }^{m}
                \del*{
                    \ballAv{2^{-j} R}{f}(x)
                    - \ballAv{2^{-j + 1} R}{f}(x)
                }
        }
    \\
        &\le
        \sum_{j=1}^{ \infty} 
        \abs*{
            \ballAv{2^{-j} R}{f}(x)
                - \ballAv{2^{-j + 1} R}{f}(x)
        }
    \\
        &\!\!\stackrel{\eqref{proofeq::lem::bojarski_characterization::difference_of_two_averages}}{\le}
        \sum_{j=1}^{ \infty} 
            2^{n+1} C_P 2^{-j} R \,
            \restMaxFun{ 2^{-j + 1} R }{ \nabla f }(x)
    \\
        &\le 
        \sum_{j=1}^{ \infty} 
        2^{n+1} C_P 2^{-j} R \,
        \restMaxFun{R}{ \nabla f }(x)
    \\
        &=
        2^{n+1} C_P R \,
        \restMaxFun{R}{ \nabla f }(x).
    \end{align*}
    Hence,
    \begin{equation} \label{proofeq::lem::bojarski_characterization::difference_of_value_and_average}
        \abs{
            f(x) - \ballAv{R}{f}(x)
        }
        \le 
        2^{n+1} C_P R \,
        \restMaxFun{R}{ \nabla f }(x).
    \end{equation}  

    Now, fix $x, y \in \lebesguePoints{f}$ such that
    $ x \ne y $ 
    and set
    $R \coloneqq \norm{x - y}$. 
    Then
    \begin{align*}
        \abs{
            \ballAv{2R}{f}(x) - \ballAv{R}{f}(y)
        }
        &=
        \abs*{
            \ballAv{2R}{f}(x) 
            - \sintegral{
                \ball{y, R}
            }{f(z)}{z}
        }
        \\
        &\le 
        \sintegral{
            \ball{y, R}
        }{
            \abs{
                \ballAv{2R}{f}(x) 
                - f(z)
            }
        }{z}
        \\
        &\le 
        \frac{
            \abs{ \ball{x, 2R}}
        }{
            \abs{ \ball{y, R}}
        }
        \sintegral{
            \ball{x, 2R}
        }{
            \abs{
                \ballAv{2R}{f}(x) 
                - f(z)
            }
        }{z}
        \\
        &=
        2^n
        \sintegral{
            \ball{x, 2R}
        }{
            \abs{
                \ballAv{2R}{f}(x) 
                - f(z)
            }
        }{z}
        \\
        &\le 
        2^{n+1}C_P  R 
        \sintegral{
            \ball{x, 2R}
        }{
            \norm{
                \nabla f(z)
            }
        }{z}
        \\
        &\le 
        2^{n+1} C_P R \, 
        \restMaxFun{2R}{ \nabla f}(x),
    \end{align*}
    where we used the Poincar\'{e} inequality in the penultimate inequality.
    
    Therefore,
    \begin{equation} \label{proofeq::lem::bojarski_characterization::offcenter_difference_of_averages}
        \abs{
            \ballAv{2R}{f}(x) - \ballAv{R}{f}(y)
        }
        \le 
        2^{n+1} C_P R \, 
        \restMaxFun{2R}{ \nabla f}(x).
    \end{equation}
    Finally, combining \eqref{proofeq::lem::bojarski_characterization::difference_of_value_and_average} with \eqref{proofeq::lem::bojarski_characterization::offcenter_difference_of_averages}
    gives us
    \begin{align*}
        \abs{
            f(x) - f(y)
        }
        &\le 
        \abs{
            f(x) - \ballAv{2R}{f}(x)
        }
        + \abs{
            \ballAv{2R}{f}(x) - \ballAv{R}{f}(y)
        }
        + \abs{
            \ballAv{R}{f}(y) - f(y)
        }
    \\
        & \!      \stackrel{\eqref{proofeq::lem::bojarski_characterization::difference_of_value_and_average}}{\le}
        \!
        2^{n+1} C_P (2R) \,
        \restMaxFun{2R}{ \nabla f }(x)
        + \abs{
            \ballAv{2R}{f}(x) - \ballAv{R}{f}(y)
        }
        + 2^{n+1} C_P R \,
        \restMaxFun{R}{ \nabla f }(y)
    \\ 
        & \!     \stackrel{\eqref{proofeq::lem::bojarski_characterization::offcenter_difference_of_averages}}{\le}
        \!
        2^{n+1} C_P (2R) \,
        \restMaxFun{2R}{ \nabla f }(x)
        + 2^{n+1} C_P R \, 
        \restMaxFun{2R}{ \nabla f}(x)
        + 2^{n+1} C_P R \,
        \restMaxFun{R}{ \nabla f }(y)
    \\
        &\le 
        2^{n+3} C_P R \del{
            \restMaxFun{2R}{ \nabla f }(x)
            + \restMaxFun{2R}{ \nabla f }(y)
        }
    \\
        &=
        2^{n+3} C_P \,
        \norm{x - y}
        \del{
            \restMaxFun{2\norm{x - y}}{ \nabla f }(x)
            + \restMaxFun{2\norm{x - y}}{ \nabla f }(y)
        },
    \end{align*}
    which is the desired inequality.
\end{proof}
\begin{corollary}\label{cor:bojarski_inequality_vector_version}
    Let $n, k \in \bN$.
    Then for every 
    $f \in \locSobolev{1,1}(\bR^n; \bR^k)$,
    we have
    \begin{equation*}
        \forall x, y \in \lebesguePoints{f}
        \qquad 
        \norm{
            f(x) - f(y)
        }
        \le 
        C_M k \,
        \norm{x - y}
        \del{
            \restMaxFun{2\norm{ x-y }}{\nabla f}(x)
            + \restMaxFun{2\norm{ x-y }}{\nabla f}(y)
        },
    \end{equation*}
    where 
    $
        C_M
    $
    is as in 
    \Cref{lem::bojarski_inequality}.
\end{corollary}
\begin{proof}
    Fix $f \in \locSobolev{1,1}(\bR^n; \bR^k)$. 
    Then there are functions 
    $
        f_1, \, \ldots, \, f_k \in \locSobolev{1,1}(\bR^n)
    $
    such that ${f = (f_i)_{ i = 1}^k }$.
    For~every $i \in [k]$,
    by \Cref{lem::bojarski_inequality}
    we have that, 
    for all $x, y \in \lebesguePoints{f_i}$,
    \begin{align*}
        \abs{
            f_i(x) - f_i(y)
        }
    &\le 
        C_M \,
        \norm{x-y}
        \del{
            \restMaxFun{2\norm{ x-y }}{\nabla f_i}(x)
            + \restMaxFun{2\norm{ x-y }}{\nabla f_i }(y)
        }
    \\
    &\le 
        C_M \,
        \norm{x-y}
        \del{
            \restMaxFun{2\norm{ x-y }}{\nabla f}(x)
            + \restMaxFun{2\norm{ x-y }}{\nabla f }(y)
        }.
    \end{align*}
    Now, note that if $x, y \in \lebesguePoints{f}$
    then 
    $x, y \in \lebesguePoints{f_i}$
    for all $i \in [k]$.
    Hence, for all 
    $x, y \in \lebesguePoints{f}$,
    \begin{align*}
        \norm{ f(x) - f(y) }^2
        &=
        \sum_{ i = 1 }^k 
            \abs{ f_i(x) - f_i(y)}^2
        \\ 
        &\le 
        \sum_{ i = 1 }^k 
            \del{ 
                C_M \,
                \norm{x-y}
                \del{
                    \restMaxFun{2\norm{ x-y }}{\nabla f}(x)
                    + \restMaxFun{2\norm{ x-y }}{\nabla f }(y)
                }
            }^2
        \\
        &=
        k \del{ 
                C_M \,
                \norm{x-y}
                \del{
                    \restMaxFun{2\norm{ x-y }}{\nabla f}(x)
                    + \restMaxFun{2\norm{ x-y }}{\nabla f }(y)
                }
            }^2.
    \end{align*}
    Taking square roots and noting that $\sqrt{k} \le k$ for all $k \in \bN$ implies the claim.
\end{proof}
%
Hajłasz \cite{Hajlasz} 
used an inequality similar to the one present in 
\Cref{lem::bojarski_inequality}
to characterize first-order Sobolev spaces via Hajłasz gradients. 
This notion was then extended by 
Yang \cite{dachun_fractional}
to the \emph{Hajłasz $s$-gradients}.
Below, we recall their definition and introduce another, related concept that we will use in 
\Cref{thm::multipointwise_bound_for_W^k_loc} --- the first step in~our characterization.
\begin{definition}[(restricted) Hajłasz $s$-gradients]
\label{def::Hajlasz_s-gradients}
    Let $s \in \intoc{0,1}$, $(X, \metricAlone)$ be a metric space, and $\measureAlone$ be a measure on $X$. 
    Let $f \colon X \to \bR$ be measurable.
    \begin{itemize}[leftmargin=1em]
    \item
		We will say that a measurable function $g \colon X \to \intcc{0,\infty}$ is a \emph{Hajłasz} $s$\emph{-gradient} of $f$ if 
	    \begin{equation*}
	        \measureAlone\forall 
	            x, y \in X
	        \qquad 
	        \abs{ f(x) - f(y) }
	        \le 
	        \norm{x-y}^{s}
	        \del{ g(x) + g(y) }.
	    \end{equation*}
	    We will denote the family of all Hajłasz $s$-gradients of $f$ by $D^s_{\measureAlone}(f)$.
	\item
	    We will say that a measurable function $g \colon \intco{0, \infty} \times X \to \intcc{0,\infty}$ is a 
	    \emph{restricted Hajłasz} $s$\emph{-gradient} of $f$ if 
	    \begin{equation*}
	        \forall t \in \intco{0, \infty}
	        \, \, \,  
	        \measureAlone\forall 
	            x, y \in X
	        \quad \, \,
	        \norm{x-y} \le t
	        \implies 
	        \abs{ f(x) - f(y) }
	        \le 
	        \norm{x-y}^{s}
	        \del{ g_t(x) + g_t(y) }.
	    \end{equation*}
	    (Note that the first argument of $g$ will be placed in the lower index and that the set of full measure on~which the inequality is satisfied may depend on it.)
	    We will denote the family of all restricted Hajłasz $s$-gradients of $f$ by 
	    $D^s_{\measureAlone, \mathrm{res}}(f)$.	 
    \end{itemize}
\end{definition}
For $n \in \bN$ and open $\Omega \subseteq \bR^n$,
let
$ \holder{k}\del*{ \Omega }$ denote the space of $k$-times continuously differentiable functions 
$
	f \colon \Omega \to \bR.
$
We will now try to motivate the first step of our characterization by proving a rather classical result for elements of
$\holder{k}\del*{ \bR^n }$. 
First, let us recall the following ``higher-order'' version of the fundamental theorem of calculus.

\begin{lemma}
\label{lem::higher_order_fundamental_theorem_of_calculus}
	Let $n, k \in \bN$, $f \in \holder{k}\del{\bR^n}$,
	and 
	$
		\bx = \set*{x_I}_{I \subseteq [k] }
		\subseteq 
		\bR^n.
	$
	Let $c_{\bx}$ be as in 
	\Cref{prop::interpolation_estimate}.
	Then
	\begin{equation*}
		\diff_{  I = \emptyset }^{ [k] }
			f(x_I)
		=
		\integral{[0,1]^k}{ \partial^{[k]}\del*{ f\circ c_{\bx}}(t)}{t}.
	\end{equation*}
\end{lemma}
\begin{proof}
	We will prove the lemma using induction over $k$. 
	The claim is correct when $k = 1$ since,
	by the fundamental theorem of calculus,
	we have
	\begin{equation*}
		\diff_{  I = \emptyset }^{ [1] }
			f(x_I)
		=
		f\del{ x_{\set*{1}}} - f\del*{ x_{\emptyset}}
		=
		\integral{[0,1]}{\partial_t \del*{f \circ c_{\bx} }(t) }{t},
	\end{equation*}
	where
	$
		c_{\bx}(t) = (1-t)x_{\emptyset} + tx_{\set*{1}}
	$
	for all $t \in [0,1]$.
	Now, fix $k \in \bN$ such that $k \ge 2$ and
	suppose that the claim is correct for $k-1$.
	For all $t_k \in [0,1]$,
	let 
	$
		\by(t_k) = \set*{ y_I(t_k)}_{I \subseteq [k-1] }
	$
	be defined by the formula 
	$   
		y_I(t_k)
		\coloneq 
		(1-t_k)x_{I  } + t_k x_{I \cup \set*{k}}
	$
	for all $I \subseteq [k-1]$.
	Let us notice that, 
	for all 
	$ t = \del*{ t_1, \, \ldots, \, t_k } \in [0,1]^k$,
	we have
	$
		c_{ \by(t_k) }(t') = c_{\bx}(t),
	$
	where we denote $t' = \del*{ t_1, \, \ldots, \, t_{k-1} }$.
	Indeed, 
	by 
	\Cref{prop::interpolation_estimate},
	for all such $t$ we have
	\begin{align*}
		c_{ \by(t_k) }(t')
	&=
		\sum_{ I = \emptyset}^{ [k-1] }
			\alpha_I(t')\beta_{ [k-1] \setminus I }(t') y_I(t_k)
	\\
	&=
		\sum_{ I = \emptyset}^{ [k-1] }
			\alpha_I(t')\beta_{ [k-1] \setminus I }(t')
			\del*{ (1-t_k)x_{I  } + t_k x_{I \cup \set*{k}} }
	\\
	&=
		\sum_{ I = \emptyset}^{ [k-1] }
			\del*{
				\alpha_I(t) \beta_{ [k] \setminus I }(t) x_I
				+ \alpha_{ I \cup \set*{k} } 
					\beta_{ [k] \setminus \del*{ I \cup \set*{ k} } }(t)
					x_{I \cup \set*{k}}
			}
	\\
	&=
		\sum_{ I = \emptyset}^{ [k] }
			\alpha_I(t)\beta_{ [k] \setminus I }(t) x_I
		=
		c_{\bx}(t). 
	\end{align*}
	Therefore, using the 
	induction hypothesis
	and 
	\Cref{lem::diff_B_diff_C_to_diff_BuC}, we have
	\begin{align*}
		\diff_{  I = \emptyset }^{ [k] }
			f(x_I)
	&=
		\diff_{  I = \emptyset }^{ [k-1] }
			f\del*{ x_{I \cup \set*{k}}}
		-
		\diff_{  I = \emptyset }^{ [k-1] }
			f\del*{ x_{I}}
	\\
	&=
		\integral{[0,1]^{k-1}}{
			\partial^{[k-1]}
			\del*{
				f \circ c_{ \by(1) }
			}(t')
		}{t'}
		- 
		\integral{[0,1]^{k-1}}{
			\partial^{[k-1]}
			\del*{
				f \circ c_{ \by(0) }
			}(t')
		}{t'}
	\\
	&=
		\integral{[0,1]^{k-1}}{
			\del*{
				\partial^{[k-1]}
				\del*{
					f \circ c_{ \by(1) }
				}(t')
				-
				\partial^{[k-1]}
				\del*{
					f \circ c_{ \by(0) }
				}(t')
			}
		}{t'}
	\\
	&=
		\integral{[0,1]^{k-1}}{
			\integral{[0,1]}{
				\partial_{x_k}
				\del*{
					\partial^{[k-1]}
					\del*{
						f \circ c_{ \by(t_k) }
					}(t')
				}
			}{t_k}
		}{t'}
		=
		\integral{[0,1]^{k}}{
			\partial^{[k]}
				\del*{
					f \circ c_{ \bx }
			}(t)
		}{t},
	\end{align*}
	as claimed.
\end{proof}
We will also need the following special case of the Fa\`{a} di Bruno formula for higher-order derivatives of the composition of two functions.
This result can also be viewed as a generalization of
\cite[Proposition 1]{Faa_di_Bruno_formula}
to~the~case when the inner function is vector-valued.
\begin{lemma}
\label{lem::Faa_di_Bruno}
	Let $n, k \in \bN$, 
	$ f \in \holder{k}\del*{ \bR^n } $,
	and 
	$ g \in \holder{k}\del*{ \intoo{0,1}^k ; \bR^n }$.
	Then for all $m \in [k]$, 
	\begin{equation}
	\label{eqlem::lem::Faa_di_Bruno::formula}
		\forall t \in \intoo{0,1}^k 
	\qquad 
		\partial^{ [m] }
			\del*{ f \circ g }(t)
		=
		\sum_{ j = 1 }^{ m }
			\sum_{ \cP \in \partitions[j]{ [m] } }
				\nabla^j f\del*{ g(t) }
				\bm{\t}
				\del*{ 
					\partial^s g( t )
				}_{ S \in \cP },			
	\end{equation}
	where we use $\bm{\t}$ to denote that $\nabla^j f(g(t))$, a $j$-linear operator, is evaluated on the tuple 
	$
		\del*{
			\partial^S g(t)
		}_{ S \in \cP };
	$%
	\footnote{Note that since the operator is symmetric, the order of the tuple's elements is not important.}
	more precisely, if
	$ g = \del*{ g_1, \, \ldots, \, g_n } $
	and 
	$
		\cP
		=
		\set*{
			S_1, \, \ldots, \, S_j	
		},
	$
	then
	\begin{equation*}
		\nabla^j f\del*{ g(t) }
		\bm{\t}
		\del*{ 
			\partial^s g( t )
		}_{ S \in \cP }
		=
		\sum_{
			\del*{ v_i }_{ i = 1 }^{ j }
			\in 
			[n]^{j}
		}
			\del*{
				\partial_{ x_{v_1} }
				\cdots 
				\partial_{ x_{v_j} }
			f 
			}
			\del*{ g(t) }
			\prod_{ b = 1 }^{ j }
				\partial^{S_b} g_{v_b}(t).
	\end{equation*}
\end{lemma}
%
%
\begin{proof}
	We will prove the lemma by induction over 
	$m\in [k]$.
	For $m=1$, the result is effectively a~rewriting of the usual chain rule in a slightly different notation. Indeed, for all 
	$t \in \intoo{0,1}^k$ we have
	\begin{align*}
		\partial^{[1]}\del*{ f \circ g}(t)
		=
		\partial_{x_1}\del*{ f \circ g}(t)
	&=
		\nabla f (g(t))
		\bm{ \t }
		\partial_{x_1} g(t)
	\\
	&=
		\nabla f (g(t))
		\bm{ \t }
		\partial^{[1]} g(t)
		=
		\sum_{ j = 1 }^{ 1 }
			\sum_{ \cP \in \partitions[j]{ [1] } }
				\nabla^j f\del*{ g(t) }
				\bm{\t}
				\del*{ 
					\partial^s g( t )
				}_{ S \in \cP },
	\end{align*}
	where we use the fact that $\set*{ [1] }$
	is the unique partition of $[1]$.

	Next, fix $m \in [k]$ such that $m < k$
	and assume that 
	\eqref{eqlem::lem::Faa_di_Bruno::formula}
	is true for $m$.
	Let us introduce a~partial order on 
	$ 2^{ [m] } $ by the following rule:
	\begin{equation*}
		\forall S, T \in 2^{ [m] }
	\qquad 
		S \le T 
		\iff 
		\min S \le \min T.
	\end{equation*}
	For all 
	$ j \in [m]$ and 
	$ \cP \in \partitions[j]{ [m] } $,
	there are 
	$ S_1, \, \ldots, \, S_j \in 2^{ [m] }$
	such that
	$ \cP = \setc*{ S_i }{ i \in [j] }$
	and, for all 
	$a, b \in [j]$,
	we~have
	$  a \le b $
	if and only if 
	$ S_a \le S_b $.
	(Note that the existence of such an ordering of elements of $\cP$ is guaranteed since $\cP$, 
	as a partition of $[m]$, has pairwise disjoint elements.)
	Then, for such 
	$j$, $\cP$, 
	and all
	$ \ell \in [j] $, 
	we~denote
	\begin{equation*}
		\cP_\ell
		\coloneq 
		\setc*{ S_i }{ 
			i \in [j] \setminus \set*{ \ell } 
		}
		\cup 
		\set*{ S_\ell \cup \set*{m+1} }.		
	\end{equation*}
	Furthermore, for all $j \in [m]$
	and $\ell \in [j]$, we define
	\begin{equation}
	\label{proofeq::lem::Faa_di_Bruno::partitions_j_k_definition}
		\partitions[j,\ell]{[m]}
		\coloneq 
		\setc*{
			\cP_{\ell}
		}{ \cP \in \partitions[j]{ [m] }}
	\quad \text{and} \quad 
		\partitions[j,j+1]{[m]}
		\coloneq 
		\setc*{
			\cP \cup \set*{ \set*{m+1} }
		}{ \cP \in \partitions[j]{ [m] }}.			
	\end{equation}
	Let us notice that
	\begin{equation}
	\label{proofeq::lem::Faa_di_Bruno::partitions_j_k_used_to_partition_j+1_first_part}
		\partitions[1]{ [m+1] }
		=
		\partitions[1,1]{ [m] },	
	\qquad     
		\partitions[m+1]{ [m+1] }
		=
		\partitions[m,m+1]{ [m] },	
	\end{equation}
	and
	\begin{equation}
	\label{proofeq::lem::Faa_di_Bruno::partitions_j_k_used_to_partition_j+1_second_part}
		\forall j \in [m] \setminus \set{1} 
	\qquad 
		\partitions[j]{[m+1]}
		=
		\partitions[j-1,j]{[m]}	
		\cup
		\bigcup_{ \ell = 1 }^{ j }
			\partitions[j,\ell]{[m]}.
	\end{equation}

		Next, let us notice that for all
		$ t \in \intoo{0,1}^k$,
		$j \in [m]$,
		and
		$
			\cP \in \partitions[j]{[m]},
		$
		we have
		\begin{align*}
			\partial_{x_{m+1}}
			\del*{
				\nabla^j f\del*{ g(\blank ) }
				\bm{\t}
				\del*{ 
					\partial^s g( \blank )
				}_{ S \in \cP }	
			}(t)
		\hspace{-6cm}&
	\\
	&=
			\partial_{x_{m+1}}
			\del*{
				\sum_{
					\del*{ v_i }_{ i = 1 }^{ j }
					\in 
					[n]^{j}
				}
					\del*{
						\partial_{ x_{v_1} }
						\cdots 
						\partial_{ x_{v_j} }
					f 
					}
					\del*{ g(\blank) }
					\prod_{b = 1}^{ j }
						\partial^{ S_b}
						g_{
							v_{b}
						}(\blank)
			}(t)			 
		\\
		&=
		\sum_{
			\del*{ v_i }_{ i = 1 }^{ j }
			\in 
			[n]^{j}
		}
			\partial_{x_{m+1}}
			\del*{
				\del*{
					\partial_{ x_{v_1} }
					\cdots 
					\partial_{ x_{v_j} }
				f 
				}
				\del*{ g(\blank) }
				\prod_{b = 1}^{ j }
					\partial^{ S_b}
					g_{
						v_{b}
					}(\blank)
			}(t),
		\intertext{
			which, 
			denoting
			$
				S_{b, b}
				\coloneq 
				S_{b} \cup \set{ m+1 }
			$
			and
			$
				S_{b, \ell}
				\coloneq
				S_{b}	
			$
			for all 
			$b \in [j]$
			and
			$\ell \in [j] \setminus \set{ b }$,
			by the Leibniz rule,
		}
		&=
		\sum_{
			\del*{ v_i }_{ i = 1 }^{ j }
			\in 
			[n]^{j}
		}
			\left(
				\partial_{x_{m+1}}
				\del*{
					\del*{
						\partial_{ x_{v_1} }
						\cdots 
						\partial_{ x_{v_j} }
					f 
					}
					\del*{ g(\blank) }
				}
				(t)
				\prod_{b = 1}^{ j }
					\partial^{ S_b}
					g_{
						v_{b}
					}(t)
			\right.
		\\
		& \hspace{3.5cm} \left.+
			\sum_{\ell = 1 }^{ j }
				\del*{
					\partial_{ x_{v_1} }
					\cdots 
					\partial_{ x_{v_j} }
				f 
				}
				\del*{ g(t) }
					\prod_{b = 1}^{ j }
					\partial^{ S_{b, \ell} }
					g_{
						v_{b}
					}(t)
			\right)		
		\\
		&=
		\sum_{
			\del*{ v_i }_{ i = 1 }^{ j }
			\in 
			[n]^{j}
		}
			\left(
				\del*{
					\sum_{w = 1}^{n}
					\del*{
						\partial_{ x_{w} }
						\del*{
							\partial_{ x_{v_1} }
							\cdots 
							\partial_{ x_{v_j} }
						f 
						}
						\del*{ g(t) }
					}
					\partial^{ \set{m+1} }
					g_{w}(t)
				}
				\prod_{b = 1}^{ j }
					\partial^{ S_b}
					g_{
						v_{b}
					}(t)
			\right.
		\\
		& \hspace{6.2cm} \left.+
			\sum_{\ell = 1 }^{ j }
				\del*{
					\partial_{ x_{v_1} }
					\cdots 
					\partial_{ x_{v_j} }
				f 
				}
				\del*{ g(t) }
					\prod_{b = 1}^{ j }
					\partial^{ S_{b, \ell} }
					g_{
						v_{b}
					}(t)
			\right)
		\\
		&=
		\sum_{
			\del*{ v_i }_{ i = 1 }^{ j }
			\in 
			[n]^{j}
		}
			\sum_{w = 1}^{n}
				\del*{
					\del*{
						\partial_{ x_{v_1} }
						\cdots 
						\partial_{ x_{v_j} }
						\partial_{ x_{w} }
					f 
					}
					\del*{ g(t) }
					\partial^{ \set{m+1} }
					g_{w}(t)
				}
				\prod_{b = 1}^{ j }
					\partial^{ S_b}
					g_{
						v_{b}
					}(t)
		\\
		& \hspace{1.3cm} 
			+
			\sum_{\ell = 1 }^{ j }
			\sum_{
				\del*{ v_i }_{ i = 1 }^{ j }
				\in 
				[n]^{j}
			}
				\del*{
					\partial_{ x_{v_1} }
					\cdots 
					\partial_{ x_{v_j} }
				f 
				}
				\del*{ g(t) }
					\prod_{b = 1}^{ j }
					\partial^{ S_{b, \ell} }
					g_{
						v_{b}
					}(t),
		\intertext{
			which, since 
			$
				\del*{
					\del*{
						v_i
					}_{i=1}^{j},
					w
				}
				\mapsto 
				\del*{
					v_1, \, \ldots, 
					\,
					v_j,
					w
				}
			$
			is a bijection from
			$
				[n]^j \times [n]				
			$
			to 
			$
				[n]^{j+1},
			$		
			and 
			denoting
			$
				S_{j+1, j+1}
				\coloneq 
				\set{ m+1 }
			$ 
			and
			$ 
				S_{b, j+1}
				\coloneq 
				S_{b}
			$
			for all 
			$ b \in [j] $,
		}
		&=
		\sum_{
			\del*{ w_i }_{ i = 1 }^{ j+1 }
			\in 
			[n]^{j+1}
		}
				\del*{
					\del*{
						\partial_{ x_{w_1} }
						\cdots 
						\partial_{ x_{w_{j+1}} }
					f 
					}
					\del*{ g(t) }
				}
				\prod_{b = 1}^{ j+1 }
					\partial^{ S_{b,j+1} }
					g_{
						w_{b}
					}(t)
		\\
		& \hspace{1.2cm} 
			+
			\sum_{\ell = 1 }^{ j }
				\sum_{
					\del*{ v_i }_{ i = 1 }^{ j }
					\in 
					[n]^{j}
				}
				\del*{
					\partial_{ x_{v_1} }
					\cdots 
					\partial_{ x_{v_j} }
				f 
				}
				\del*{ g(t) }
					\prod_{b = 1}^{ j }
					\partial^{ S_{b, \ell} }
					g_{
						v_{b}
					}(t)
		\\
		&=
		\nabla^j f\del*{ g(t) }
		\bm{\t}
		\del*{ 
			\partial^s g( t )
		}_{ S \in \cP \cup \set*{ \set{m+1} } }
			+
			\sum_{\ell = 1 }^{ j }
				\nabla^j f\del*{ g(t) }
				\bm{\t}
				\del*{ 
					\partial^s g( t )
				}_{ S \in \cP_{\ell} }.
		\end{align*}
		That is, for all
		$ t \in \intoo{0,1}^k$,
		$j \in [m]$,
		and
		$
			\cP \in \partitions[j]{[m]},
		$
		we have
		\begin{multline}
		\label{proofeq::lem::Faa_di_Bruno::consequence_of_Leibniz_rule}
			\partial_{x_{m+1}}
			\del*{
				\nabla^j f\del*{ g(\blank ) }
				\bm{\t}
				\del*{ 
					\partial^s g( \blank )
				}_{ S \in \cP }	
			}(t)
		\\
			=
			\nabla^j f\del*{ g(t) }
			\bm{\t}
			\del*{ 
				\partial^s g( t )
			}_{ S \in \cP \cup \set*{ \set{m+1} } }
				+
				\sum_{\ell = 1 }^{ j }
					\nabla^j f\del*{ g(t) }
					\bm{\t}
					\del*{ 
						\partial^s g( t )
					}_{ S \in \cP_{\ell} }.			
		\end{multline}

	Therefore, for all $t \in \intoo{0,1}^k$,
	\begin{align*}
		\partial^{ [m+1] }
			\del*{ f \circ g}(t)
	\hspace{-1.8cm}
	&\hspace{1.8cm}=
		\partial_{x_{m+1}}
			\partial^{ [m] }
				\del*{ f \circ g}(t)
	\\
	&=
		\partial_{x_{m+1}}
			\del*{
				\sum_{ j = 1 }^{ m }
					\sum_{ \cP \in \partitions[j]{ [m] } }
						\nabla^j f\del*{ g(\cdot ) }
						\bm{\t}
						\del*{ 
							\partial^s g( \cdot )
						}_{ S \in \cP }
			}(t)
	\\
	&=
		\sum_{ j = 1 }^{ m }
			\sum_{ \cP \in \partitions[j]{ [m] } }
				\partial_{x_{m+1}}
				\del*{
					\nabla^j f\del*{ g(\cdot ) }
					\bm{\t}
					\del*{ 
						\partial^s g( \cdot )
					}_{ S \in \cP }	
				}(t),
	\intertext{
		which, by 
		\eqref{proofeq::lem::Faa_di_Bruno::consequence_of_Leibniz_rule},
	}
	&=
		\sum_{ j = 1 }^{ m }
			\sum_{ \cP \in \partitions[j]{ [m] } }
				\del*{
					\nabla^{j+1} f\del*{ g( t ) }
					\bm{\t}
						\del*{ 
							\partial^s g( t )
						}_{ 
							S \in 
							\cP \cup \set*{ \set*{m+1} }
						}
					+	
					\sum_{ \ell = 1 }^{ j }
					\nabla^{j} f\del*{ g( t ) }
					\bm{\t}
						\del*{ 
							\partial^s g( t )
						}_{ S \in \cP_\ell }				
				},
	\intertext{
		which, by \eqref{proofeq::lem::Faa_di_Bruno::partitions_j_k_definition},
	}
	&=
		\sum_{ j = 1 }^{ m }
			\del*{
				\sum_{ \cR \in \partitions[j,j+1]{ [m] } }
					\nabla^{j+1} f\del*{ g( t ) }
					\bm{\t}
						\del*{ 
							\partial^s g( t )
						}_{ 
							S \in 
							\cR
						}
					+	
				\sum_{ \ell = 1 }^{ j }
					\sum_{ \cR \in \partitions[j,\ell]{ [m] } }
					\nabla^{j} f\del*{ g( t ) }
					\bm{\t}
						\del*{ 
							\partial^s g( t )
						}_{ S \in \cR }				
				}
	\\
	&=
		\sum_{ 
			\cR \in \partitions[1,1]{ [m] } 
		}
			\nabla^{1} f\del*{ g( t ) }
			\bm{\t}
				\del*{ 
					\partial^s g( t )
				}_{ S \in \cR }		
		+
		\sum_{ \cR \in \partitions[m,m+1]{ [m] } }
			\nabla^{m+1} f\del*{ g( t ) }
			\bm{\t}
				\del*{ 
					\partial^s g( t )
				}_{ 
					S \in 
					\cR
				}
	\\
	&\qquad 	
		+
		\sum_{ j = 2 }^{ m }
			\del*{
				\sum_{ \cR \in \partitions[j-1,j]{ [m] } }
					\nabla^{j} f\del*{ g( t ) }
					\bm{\t}
						\del*{ 
							\partial^s g( t )
						}_{ 
							S \in 
							\cR
						}
				+	
				\sum_{ \ell = 1 }^{ j }
					\sum_{ 
						\cR \in \partitions[j,\ell]{ [m] } 
					}
					\nabla^{j} f\del*{ g( t ) }
					\bm{\t}
						\del*{ 
							\partial^s g( t )
						}_{ S \in \cR }	
			}
	\intertext{
		which, by \eqref{proofeq::lem::Faa_di_Bruno::partitions_j_k_used_to_partition_j+1_first_part}
		and 
		\eqref{proofeq::lem::Faa_di_Bruno::partitions_j_k_used_to_partition_j+1_second_part}
		,
	}
	&=	
		\sum_{ j = 1 }^{ m+1 }
			\sum_{ \cP \in \partitions[j]{ [m+1] } }
				\nabla^j f\del*{ g(t) }
				\bm{\t}
				\del*{ 
					\partial^s g( t )
				}_{ S \in \cP },		
	\end{align*}
	proving the inductive step.
	Thus, the claim follows by induction.
\end{proof}
\begin{proposition}
\label{prop::higher_order_fundamental_theorem_of_calculus_inequality}
	Let $n, k \in \bN$
	and 
	$f \in \holder{k}\del{\bR^n}$.
	Then for all 
	$
		\bx = \set*{x_I}_{I \subseteq [k] }
		\subseteq 
		\bR^n,
	$
	we have
	\begin{equation}
	\label{eqprop::prop::higher_order_fundamental_theorem_of_calculus_inequality::inequality}
		\abs*{
			\diff_{ I= \emptyset}^{ [k] }
				f(x_I)
		}
		\le 
		\sum_{ j = 1 }^{ k }
			\polynum{j}(\bx) 
			\integral{[0,1]^k}{
				\norm{ \nabla^j f \del*{ c_{\bx}(t) }}_{\mathrm{op}}
			}{t},
	\end{equation}
	where $c_{\bx} \colon [0,1]^k \to \bR^n$
	is as in \Cref{prop::interpolation_estimate},
	and 
	$\normAlone[\mathrm{op}]$ is the operator norm of a multilinear bounded operator, that is,
	if 
	$
		\del*{ V, \normAlone[V] }
	$
	and 
	$
		\del*{ W, \normAlone[W] }
	$
	are normed spaces
	and
	$
		T \colon V^j \to W
	$ 
	is~a~$j$-linear operator, $j\in \bN$,
	then
	\begin{equation*}
		\norm{ T }_{ \mathrm{op} }
		\coloneq 
		\inf 
		\setc*{
			C \ge 0
		}{
			\norm*{ 
				T\del*{ 
					v_1, \, \ldots, \, v_j 
				}
			}_{W}
			\le 
			C
			\prod_{ i = 1 }^{ j }
				\norm{ v_j }_{V}
		\quad
			\text{for all }
			v_1, \, \ldots, v_j
			\in 
			V 
		}.
	\end{equation*}
\end{proposition}
\begin{proof}
	Fix $t \in \intoo{0,1}^k$.
	First, note that
	by 
	\Cref{prop::interpolation_estimate},
	for all 
	$j \in [k]$
	and
	$ \cP \in \partitions[j]{ [k] } $,
	we have
	\begin{align*}
		\abs*{
			\nabla^j f\del*{ c_{\bx}(t) }
			\bm{\t} 
			\del*{
				\partial^S c_{\bx}(t)
			}_{ S \in \cP }
		}
		\le 
		\norm*{ 
			\nabla^j f\del*{ c_{\bx}(t) } 
		}_{ \mathrm{op} }
		\prod_{ S \in \cP }
			\norm*{
				\partial^S c_{\bx}(t)
			}
	&
	\\
		\le
		\norm*{ 
			\nabla^j f\del*{ c_{\bx}(t) } 
		}_{ \mathrm{op} }
		\prod_{S \in \cP}
			\polyset{S}(\bx)
	&=
		\norm*{ 
			\nabla^j f\del*{ c_{\bx}(t) } 
		}_{ \mathrm{op} }
		\polyfam{ \cP }(\bx),
	\end{align*}
	Thus, using the version of the Fa\`{a} di Bruno formula that is present in \Cref{lem::Faa_di_Bruno}, we have
	\begin{align*}
		\abs*{
			\partial^{ [k] } \del*{ f \circ c_{\bx} }(t)
		}
		=
		\abs*{
			\sum_{ j = 1 }^{ k }
				\sum_{ \cP \in \partitions[j]{[k]} }
					\nabla^j f\del*{ c_{\bx}(t) }
					\bm{\t} 
					\del*{
						\partial^S c_{\bx}(t)
					}_{ S \in \cP }
		}
	&\\ 
		\le 
		\sum_{ j = 1 }^{ k }
			\sum_{ \cP \in \partitions[j]{[k]} }
				\abs*{
					\nabla^j f\del*{ c_{\bx}(t) }
					\bm{\t} 
					\del*{
						\partial^S c_{\bx}(t)
					}_{ S \in \cP }
				}
	\\
		\le 
		\sum_{ j = 1 }^{ k }
			\sum_{ \cP \in \partitions[j]{[k]} }
				\norm*{ 
					\nabla^j f\del*{ c_{\bx}(t) } 
				}_{\mathrm{op}}
				\polyfam{ \cP }(\bx)		
	&=
		\sum_{ j = 1 }^{ k }
			\polynum{ j }( \bx )
			\norm*{ 
				\nabla^j f\del*{ c_{\bx}(t) } 
			}_{\mathrm{op}}.
	\end{align*}
	Therefore, since the resulting inequality is valid for all 
	$ t \in \intoo{0,1}^k$,
	by 
	\Cref{lem::higher_order_fundamental_theorem_of_calculus}
	we have
	\begin{align*}
		\abs*{
			\diff_{ I= \emptyset}^{ [k] }
				f(x_I)
		}
	&=
		\abs*{
			\integral{[0,1]^k}{
				\partial^{ [k] } \del*{ f \circ c_{\bx} }(t)	
			}{t}
		}
	\\
	&\le 
		\integral{[0,1]^k}{
			\abs*{ \partial^{ [k] } \del*{ f \circ c_{\bx} }(t)}
		}{t}
	\\
	&\le
		\integral{[0,1]^k}{
			\sum_{ j = 1 }^{ k }
				\polynum{ j }( \bx )
				\norm*{ 
					\nabla^j f\del*{ c_{\bx}(t) } 
				}_{\mathrm{op}}
		}{t}
		=
		\sum_{ j = 1 }^{ k }
			\polynum{ j }( \bx )
			\integral{[0,1]^k}{
				\norm*{ 
					\nabla^j f\del*{ c_{\bx}(t) } 
				}_{\mathrm{op}}
			}{t},
	\end{align*}
	as claimed.
\end{proof}
\begin{remark}
\label{rem::motivation_for_the_characterization}
	Using the notation from \Cref{prop::higher_order_fundamental_theorem_of_calculus_inequality},
	let us consider the case when 
	$ k = 1$,
	and we have
	$ x_{\emptyset} = x$ and $x_{\set*{1}} = y$
	for some $x, y \in \bR^n$.
	Then 
	$
		\polynum{1}(\bx) = \norm{x-y}
	$
	and
	$
		c_{\bx}(t) = (1-t)x + ty
	$
	for all $t \in [0,1]$.
	Thus, inequality 
	\eqref{eqprop::prop::higher_order_fundamental_theorem_of_calculus_inequality::inequality}
	can be equivalently written as 
	\begin{equation*}
		\abs*{ f(x) - f(y) }
		\le 
		\norm{x-y}
		\integral{ [0,1] }{ 
			\norm{ 
				\nabla f \del*{ (1-t)x+ty } 
			}_{\mathrm{op}}
		}{ t }.
	\end{equation*}
	We can compare the above inequality with the one present within the definition of Hajłasz gradients:
	\begin{equation*}
		\abs*{ f(x) - f(y) }
		\le 
		\norm{x-y}
		\del*{
			g(x) + g(y)
		}.
	\end{equation*}
	The two inequalities have a somewhat similar form;
	the latter one looks as if it could be ``obtained'' from the former by changing the integral into a sum of values of some function $g$ at 
	$x$ and $y$.
	Note that $x = c_{\bx}(0)$ and 
	$ y = c_{\bx}(1)$, 
	so these points correspond to the values of~$c_{\bx}$ at the ends of the domain of integration.
	
	Suppose now that $k$ is a natural number greater than $1$.
	The above observation suggests that for a sufficiently regular $f$, we might hope to find functions $G_j$, $j \in [k]$, such that from the~inequality
	\begin{equation*}
		\abs*{
			\diff_{ I= \emptyset}^{ [k] }
				f(x_I)
		}
		\le 
		\sum_{ j = 1 }^{ k }
			\polynum{j}(\bx) 
			\integral{[0,1]^k}{
				\norm{ 
					\nabla^j f \del*{ c_{\bx}(t) }
				}_{\mathrm{op}}
			}{t},
	\end{equation*}
	we can ``obtain'' the following one:
	\begin{equation*}
		\abs*{
			\diff_{ I= \emptyset}^{ [k] }
				f(x_I)
		}
		\le 
		\sum_{ j = 1 }^{ k }
			\polynum{j}(\bx) 
			\sum_{ I = \emptyset }^{ [k] }
				G_j(x_I).
	\end{equation*}
\end{remark}
As we will see in  
\Cref{thm::multipointwise_bound_for_W^k_loc},
the above discussion, though obviously non-rigorous, has pointed us in the right direction.
\begin{theorem} 
\label{thm::multipointwise_bound_for_W^k_loc}
    Let $n, k \in \bN$. 
    Then there exists a constant $C_{n,k} \ge 0$ such that, for all 
    ${f \in \locSobolev{k,1}(\bR^n)}$ 
    and 
    $
        \bx = \set*{ x_I }_{ I \subseteq [k] }
            \subseteq \lebesguePoints{f},
    $
    we have
    \begin{align} 
        \label{thmeq::thm::multipointwise_bound_for_W^k_loc::eq::just_loc}
    \hspace{-0.25cm}
        \abs*{
            \diff_{ I = \emptyset }^{ [k] }
                f(x_I)
        }
        &\le 
        C_{n,k}
        \sum_{ j = 1 }^{ k }
            \polynum{  j  }( \bx )
            \!
            \smashoperator{
            \sum_{ I = \emptyset }^{ [k] }
            }
            \!
                \restMaxFun[k+j-1]{R(\bx) }{\nabla^j f}(x_I)
        ,
    \intertext{
        where 
        $
            R(\bx) 
            \coloneq 
            2^{k^2} \ell(\bx).
        $
        Moreover, 
        there exists a constant $\widehat{C}_{n,k+1} \ge 0$
        such that, 
        for every ${s \in \intoc{k, k+1}}$,
        $f \in \locSobolev{k,1}(\bR^n)$, 
        $g \in D_{\lambda, \mathrm{res}}^{s-k}( \nabla^k f )$, and all
        $
            \bx = \set*{ x_I }_{ I \subseteq [k+1] }
                \subseteq \lebesguePoints{f},
        $
    }
        \label{thmeq::thm::multipointwise_bound_for_W^k_loc::eq::with_s-k_gradient}
    \hspace{-0.25cm}
        \abs*{
            \diff_{ I = \emptyset }^{ [k+1] }
                f(x_I)
        }
        &\le 
        \widehat{C}_{n,k+1}
        \!
        \del*{
            \sum_{ j = 1 }^{ k }
                \polynum{ j }( \bx )
            \!
            \smashoperator{
                \sum_{ I = \emptyset }^{ [k+1] }
            }
            \!
                    \restMaxFun[k+j]{R(\bx)}{\nabla^j f}(x_I)
            + 
            \polynum{ s }(\bx)
            \!
            \smashoperator{
            \sum_{ I = \emptyset }^{ [k+1] }
	        }
	        \!
                    \restMaxFun[2k]{R(\bx)}{g_{\ell(\bx)} }(x_I)
        }
    \! 
    	,
    \end{align}
    where this time 
    $ R(\bx) = 2^{(k+1)^2} \ell(\bx)$,
    since now $\bx$ has $2^{k+1}$ elements.
\end{theorem}
\begin{proof}
    The proof will be inductive over the value of $k$.
    Assuming that \eqref{thmeq::thm::multipointwise_bound_for_W^k_loc::eq::just_loc} 
    is satisfied for a~given $k$,
    we will show that 
    \eqref{thmeq::thm::multipointwise_bound_for_W^k_loc::eq::with_s-k_gradient} 
    is true for all $s \in \intoc{k, k+1}$;
    then,
    \eqref{thmeq::thm::multipointwise_bound_for_W^k_loc::eq::just_loc} 
    for $k+1$ will quickly follow.

    First of all, let us note that 
    \eqref{thmeq::thm::multipointwise_bound_for_W^k_loc::eq::just_loc}
    for $k = 1$ follows directly from 
    \Cref{lem::bojarski_inequality}.
    Indeed, fix 
    $f \in \locSobolev{1,1}( \bR^n)$
    and 
    $
        \bx = 
        \set*{ x_{I} }_{ I \subseteq [1] }
        \subseteq 
        \lebesguePoints{f}.
    $
    Denote $x \coloneq x_{\emptyset}$
    and $y \coloneq x_{ \set*{1} }$. 
    Then
    \begin{equation*}
        \abs{ f(x) - f(y) }
        =
        \abs*{
            \diff_{ I = \emptyset }^{ [1] }
                f(x_I)
        },
        \quad 
        \polynum{ 1 }( \bx )
        =
        \ell( \bx)
        =
        \norm{x-y}.
    \end{equation*}
    Therefore, by
    \Cref{lem::bojarski_inequality}
    we have
    \begin{align*}
        \abs*{
            \diff_{ I = \emptyset }^{ [1] }
                f(x_I)
        }
        =
        \abs{ f(x) - f(y) }
        &\le 
        C_M \norm{x-y}
        \del*{
            \restMaxFun{2\norm{x-y}}{ \nabla f}(x)
            + \restMaxFun{2\norm{x-y}}{ \nabla f}(y)
        }
        \\
        &=
        C_M \polynum{ 1 }( \bx) 
        \sum_{ I = \emptyset }^{ [1] }
            \restMaxFun{ 2\ell(\bx) }{ \nabla f}(x_I),    
    \end{align*}
    and \eqref{thmeq::thm::multipointwise_bound_for_W^k_loc::eq::just_loc}
    for $k = 1$ is satisfied with 
    $C_{n,1} = C_M$.

	In the next part of the proof, 
	we will show some auxiliary results that we will refer to as ``sublemmas.''
	The~common assumptions of~the~sublemmas are listed in \Cref{asu::thm::multipointwise_bound_for_W^k_loc::inductive_assumptions}.
    \begin{assumptions}
    \label{asu::thm::multipointwise_bound_for_W^k_loc::inductive_assumptions}
        For some $k \in \bN$,
    	the inequality in
        \eqref{thmeq::thm::multipointwise_bound_for_W^k_loc::eq::just_loc}
        is true for every $f \in \locSobolev{k,1}( \bR^n)$
        and all 
        $
            \bx = \set*{ x_I }_{ I \subseteq [k] }
            \subseteq 
            \lebesguePoints{f}.
        $
        Fix $s \in \intoc{k, k+1}$,
	    $f \in \locSobolev{k,1}( \bR^n)$,
	    and
	    $g \in D^{s-k}_{ \lambda, \mathrm{res} }( \nabla^k f )$. 
    \end{assumptions}
    \begin{lemmaInThm}
        \label{sublem::thm::multipointwise_bound_for_W^k_loc::eq::sublemma_1}
        Suppose  \Cref{asu::thm::multipointwise_bound_for_W^k_loc::inductive_assumptions}
        are true.
        Let 
        $
            \bx = \set*{ x_I }_{ I \subseteq [k+1] }
            \subseteq 
            \lebesguePoints{f}
        $
        be such that
        $ \polyset{\set*{k+1}}(\bx)  > 0 $.
        Fix $w \in \clball{ 0,  \polyset{\set*{k+1}}(\bx) }$
        and denote
        \begin{equation*}
            A \coloneq
            \bigcap_{I \subseteq [k] }
                \del*{
                    \lebesguePoints{f} 
                    - x_I + x_{\emptyset}
                }
            \cap 
            \bigcap_{I \subseteq [k] }
                \del*{
                    \lebesguePoints{f} 
                    - x_I + x_{\emptyset} - w
                }.
        \end{equation*}
        %
        Fix $ a \in A $ and define  
        $
            \set*{ a_I }_{ I \subseteq [k+1] }
        $
        as follows:
        \begin{equation*}
            \forall I \subseteq [k]
        \qquad 
            a_I \coloneq a + x_I - x_{\emptyset}
            \quad \text{and} \quad 
            a_{I \cup \set*{k+1} } \coloneq a_I + w.
        \end{equation*}
        Then, 
        \begin{equation*}
            \abs*{
                \diff_{ I = \emptyset }^{ [k+1] }
                    f(a_I)
            }
            \le 
            C_{n,k}'
            \del*{
                \sum_{ j = 2 }^{k}
                    \polynum{ j }( \bx)
                    \sum_{ I = \emptyset }^{ [k+1] }
                        \restMaxFun[ k+j -1 ]{ R( \bx ) }{ 
                            \nabla^j f
                        }( a_I )
                +
                \polynum{ s }( \bx )
                \sum_{ I = \emptyset }^{ [k+1] }
                    \restMaxFun[2k-1]{ R( \bx ) }{ g_{\ell(\bx)} }( a_I )
            }
        \end{equation*}
        for some $C_{n,k}' \ge 0$ that depends only on $n$ and $k$.
    \end{lemmaInThm}
    \begin{proof}[Proof of \Cref{sublem::thm::multipointwise_bound_for_W^k_loc::eq::sublemma_1}]
        First, let us denote
        \begin{equation*}
            \bx' \coloneq \set*{ x_I }_{ I \subseteq [k] }
        \quad
        \text{and}
        \quad
            a - x_{\emptyset} + \bx'
            =
            \set*{ a - x_{\emptyset} + x_{I} }_{ I \subseteq [k] }
            =
            \set*{ a_I }_{ I \subseteq [k] }.	        
        \end{equation*}
        Since 
        $ y \mapsto a - x_{\emptyset} + y$ is an isometry of $\bR^n$,
        by \Cref{lem::invariant_transformations_of_poly},
        for all $j \in [k]$ we have 
        \begin{gather}
        \label{thmeq::thm::multipointwise_bound_for_W^k_loc::eq::sublemma_1::P_for_x'_to_x}
			\polynum{ j }\del*{
				a 
				- 
				x_{\emptyset} 
				+ 
				\bx' 
			}
        	=
            \polynum{ j }(\bx')        	
		\shortintertext{and}
        \label{thmeq::thm::multipointwise_bound_for_W^k_loc::eq::sublemma_1::l_for_x'_to_x}
		R\del*{
			a 
			-
			x_{\emptyset} 
			+
		\bx' 
		}
		=
		2^{k^2}\ell\del*{
			a 
			-
			x_{\emptyset} 
			+ 
			\bx' 
		}
		=
		2^{k^2}\ell( \bx' )
		\le 
		2^{(k+1)^2}\ell( \bx )
		=
		R(\bx),	        
        \end{gather}
        where the inequality follows from 
        \Cref{lem::adding_k+1_direction_to_poly_j_gives_poly_j+1}.
        
        Let us define function $f_w$ by the formula 
        $
            f_w(x) = f(x + w) - f(x).
        $
        Since $f \in \locSobolev{k,1}(\bR^n)$, we have
        $
            f_w \in \locSobolev{k,1}(\bR^n)
        $
        and, 
        for all $j \in [k]$,
        we have
        $
            \nabla^j f_w(x)
            =
            \nabla^j f(x+w)
            - \nabla^j f(x)
        $
        for almost every $x$.
        For all $j \in [k]$,
        let us denote the~set on which this equality is satisfied by $E_j$.
        
        Now, notice that by~%
        \Cref{lem::diff_B_diff_C_to_diff_BuC},
        \begin{equation*}
            \abs*{
                \diff_{ I = \emptyset }^{ [k+1] }
                    f(a_I)
            }
            =
            \abs*{
                \diff_{ I = \emptyset }^{ [k] }
                    \del*{
                        f\del*{a_{I\cup \set*{k+1} }}
                        - f(a_I)
                    }
            }
            =
            \abs*{
                \diff_{ I = \emptyset }^{ [k] }
                    \del*{
                        f(a_I + w)
                        - f(a_I)
                    }
            }
            =
            \abs*{
                \diff_{ I = \emptyset }^{ [k] }
                    f_w(a_I)
            }.
        \end{equation*}
        Since $a \in A$,
        we have
        $
            \set*{ a_I }_{ I \subseteq [k] }
            \subseteq 
            \lebesguePoints{ f_w }.
        $
        Therefore, by 
        \eqref{thmeq::thm::multipointwise_bound_for_W^k_loc::eq::just_loc}
        we have
        \begin{align}
            \abs*{
                \diff_{ I = \emptyset }^{ [k+1] }
                    f(a_I)
            }
            =
            \abs*{
                \diff_{ I = \emptyset }^{ [k] }
                    f_w(a_I)
            }
            &\le 
            C_{n,k}
            \sum_{ j = 1 }^{ k }
                \polynum{ j }( a - x_{\emptyset} + \bx' )
                \sum_{ I = \emptyset }^{ [k] }
                    \restMaxFun[k+j-1]{R\del*{a - x_{\emptyset} + \bx'}}{\nabla^j f_w}(a_I)
            \notag
        \\
            &\le
            C_{n,k}
            \sum_{ j = 1 }^{ k }
                \polynum{ j }( \bx' )
                \sum_{ I = \emptyset }^{ [k] }
                    \restMaxFun[k+j-1]{R(\bx)}{\nabla^j f_w}(a_I), 
            \label{thmeq::thm::multipointwise_bound_for_W^k_loc::eq::sublemma_1::after_inductive_step}
        \end{align}
        where the last inequality follows from \eqref{thmeq::thm::multipointwise_bound_for_W^k_loc::eq::sublemma_1::P_for_x'_to_x}
        and
        \eqref{thmeq::thm::multipointwise_bound_for_W^k_loc::eq::sublemma_1::l_for_x'_to_x}.

        Now, fix $j \in [k-1]$. 
        By \Cref{cor:bojarski_inequality_vector_version} there exists a set 
        $D'_j$ of full measure such that if~$x, y \in D'_j$, then
        \begin{equation*}
            \norm*{
                \nabla^j f(x) - \nabla^j f(y)
            }
            \le 
            C_M n^j
            \norm{x-y}
            \del*{
                \restMaxFun{ 2\norm{x-y} }{ \nabla^{j+1} f}(x)
                + \restMaxFun{ 2\norm{x-y} }{ \nabla^{j+1} f}(y)
            }.
        \end{equation*}
        Denote $D_j \coloneq E_j \cap D'_j \cap \del*{ D'_j - w}$.
        Then, for all $d \in D_j$, we have
        \begin{align*}
            \norm*{\nabla^j f_w(d)}
        &=
            \norm*{
                \nabla^j f(d+w) - \nabla^j f(d)
            }
        \\
        &\le 
            C_M n^j
            \norm{w }
            \del*{
                \restMaxFun{ 2 \norm{w} }{ \nabla^{j+1} f}(d+w)
                + \restMaxFun{ 2\norm{w} }{ \nabla^{j+1} f}(d)
            }.
        \end{align*}
        Since $D_j$ is of full measure as an intersection of sets of full measure, for all $I \subseteq [k]$ we have
        \begin{align*}
            \restMaxFun[k+j-1]{ R(\bx) }{
                \nabla^j f_w
            }(a_I)
        \hspace{-2cm}& \hspace{2cm} 
        \le
            \restMaxFun[k+j-1]{ R(\bx) }{
                C_M n^j
                \norm{ w }
                \del*{
                    \restMaxFun{ 2 \norm{w} }{ \nabla^{j+1} f}\del*{
                    	\blank {+} w 
                    }
                    + \restMaxFun{ 2\norm{w} }{ \nabla^{j+1} f}
                }
            }(a_I),
        \intertext{
        	which, since $2\norm{w} \le 2\polyset{\set*{k+1}}(\bx) \le 2\ell(\bx) \le R(\bx)$ and $n^j \le n^k$,
        }
            &\le 
            \restMaxFun[k+j-1]{ R(\bx) }{
                C_M n^k
                \polyset{\set*{k+1}}(\bx)
                \del*{
                    \restMaxFun{ R(\bx) }{ \nabla^{j+1} f}\del*{ \blank +w }
                    + \restMaxFun{ R(\bx) }{ \nabla^{j+1} f}
                }
            }(a_I),
        \intertext{which, by \Cref{lem::shifting_operators},
        }
            &\le
            C_M n^k
            \polyset{\set*{k+1}}(\bx) 
            \del*{
                \restMaxFun[k+j]{ R(\bx) }{ \nabla^{j+1} f}( a_I +w)
                + \restMaxFun[k+j]{ R(\bx) }{ \nabla^{j+1} f}( a_I) 
            }
        \\
            &=
            C_M n^k
            \polyset{\set*{k+1}}(\bx) 
            \del*{
                \restMaxFun[k+j]{ R(\bx) }{ \nabla^{j+1} f}\del*{ a_{I \cup \set*{k+1} }}
                + \restMaxFun[k+j]{ R(\bx) }{ \nabla^{j+1} f}( a_I) 
            }.
        \end{align*}
        Therefore,
        \begin{align*}
            \polynum{ j }(\bx')
            \sum_{ I = \emptyset }^{ [k] }
                \restMaxFun[k+j-1]{ R(\bx) }{
                    \nabla^j f_w
                }(a_I)
            \hspace{-4.5cm}&
        \\
            &\le 
            \polynum{ j }(\bx')
            \sum_{ I = \emptyset }^{ [k] }
                C_M n^k
                \polyset{\set*{k+1}}(\bx) 
                \del*{
                    \restMaxFun[k+j]{ R(\bx) }{ \nabla^{j+1} f}\del*{ a_{I \cup \set*{k+1} }}
                    + \restMaxFun[k+j]{ R(\bx) }{ \nabla^{j+1} f}( a_I) 
                }
        \\
            &=
            C_M n^k
            \polynum{ j }(\bx')
            \polyset{\set*{k+1}}(\bx) 
            \sum_{ I = \emptyset }^{ [k+1] }
                \restMaxFun[k+j]{ R(\bx) }{ \nabla^{j+1} f}( a_I) 
        \\
            &\le 
            C_M n^k
            \polynum{ j+1 }(\bx)
            \sum_{ I = \emptyset }^{ [k+1] }
                \restMaxFun[k+j]{ R(\bx) }{ \nabla^{j+1} f}( a_I),
        \end{align*}
        where the last inequality follows from 
        \Cref{lem::adding_k+1_direction_to_poly_j_gives_poly_j+1}.
        Thus, for all $j \in [k-1]$ we have
        \begin{equation}
        \label{thmeq::thm::multipointwise_bound_for_W^k_loc::eq::sublemma_1::case_j<k}
            \polynum{ j }\del*{ \bx' }
            \sum_{ I = \emptyset }^{ [k] }
                \restMaxFun[k+j-1]{ R(\bx) }{
                    \nabla^j f_w
                }(a_I)
            \le
            C_M n^k
            \polynum{ j+1 }(\bx)
            \sum_{ I = \emptyset }^{ [k+1] }
                \restMaxFun[k+j]{ R(\bx) }{ \nabla^{j+1} f}( a_I).
        \end{equation}

        Next, since $ g \in D^{s-k}_{ \lambda, \mathrm{res}}( \nabla^k f )$,
        there exists a set $D'_k$ of full measure such that if 
        $x, y \in D'_k$ are such that $\norm{x-y} \le \ell(\bx)$, then 
        \begin{equation*}
            \norm*{
                \nabla^k f(x)
                - \nabla^k f(y)
            }
            \le 
            \norm{ x- y }^{ s-k }
            \del*{
                g_{ \ell(\bx) }(x) + g_{ \ell(\bx) }(y)
            }.
        \end{equation*}
        Denote $D_k \coloneq E_k \cap D'_k \cap (D'_k - w)$.
        Then, since $\norm{w} \le \ell(\bx)$,
        for $d \in D_k$ we have
        \begin{equation*}
            \norm*{\nabla^k f_w(d)}
            =
            \norm*{
                \nabla^k f(d+w) - \nabla^k f(d)
            }
            \le 
            \norm{w }^{ s-k }
            \del*{
                g_{ \ell(\bx) }(d+w) 
                + g_{ \ell(\bx) }(d)
            }.
        \end{equation*}
        Since $D_k$ is of full measure as an intersection of sets of full measure, for all $I \subseteq [k]$ we have
        \begin{align*}
            \restMaxFun[2k-1]{ R(\bx) }{
                \nabla^k f_w
            }(a_I)
            &\le 
            \restMaxFun[2k-1]{ R(\bx) }{
                \norm{w }^{ s-k }
                \del*{
                    g_{ \ell(\bx) }\del*{ \blank + w } 
                    + g_{ \ell(\bx) }
                }
            }(a_I),
        \intertext{
            which, since 
            $\norm{w} \le \polyset{\set*{k+1}}(\bx)$,
        }
            &\le 
            \restMaxFun[2k-1]{ R(\bx) }{
                \polyset{\set*{k+1}}(\bx)^{ s-k }
                \del*{
                    g_{ \ell(\bx) }\del*{ \blank + w } 
                    + g_{ \ell(\bx) }
                }
            }(a_I),
        \intertext{which, by \Cref{lem::shifting_operators},
        }
            &\le
            \polyset{ \set*{k+1} }(\bx)^{ s-k }
            \del*{
                \restMaxFun[2k-1]{ R(\bx) }{ g_{ \ell(\bx) } }( a_I +w)
                + \restMaxFun[2k-1]{ R(\bx) }{ g_{ \ell(\bx) } }( a_I) 
            }
        \\
            &=
            \polyset{ \set*{k+1} }(\bx)^{ s-k }
            \del*{
                \restMaxFun[2k-1]{ R(\bx) }{ g_{ \ell(\bx) } }\del*{ a_{I\cup \set*{k+1} }}
                + \restMaxFun[2k-1]{ R(\bx) }{ g_{ \ell(\bx) } }( a_I) 
            }.
        \end{align*}
        Therefore,
        \begin{align*}
            \polynum{ k }(\bx')
            \sum_{ I = \emptyset }^{ [k] }
                \restMaxFun[2k-1]{ R(\bx) }{
                    \nabla^k f_w
                }(a_I)
        \hspace{-4cm}&
        \\
        &\le 
            \polynum{ k }(\bx')
            \sum_{ I = \emptyset }^{ [k] }
                \polyset{ \set*{k+1} }(\bx)^{ s-k }
                \del*{
                    \restMaxFun[2k-1]{ R(\bx) }{ g_{ \ell(\bx) } }\del*{ a_{I\cup \set*{k+1} }}
                    + \restMaxFun[2k-1]{ R(\bx) }{ g_{ \ell(\bx) } }( a_I) 
                }
        \\
        &=
            \polynum{ k }(\bx')
            \polyset{\set*{k+1}}(\bx)^{ s-k }
            \sum_{ I = \emptyset }^{ [k+1] }
                \restMaxFun[2k-1]{ R(\bx) }{ g_{ \ell(\bx) }}( a_I) 
        \\
        &\le 
            \polynum{ s }(\bx)
            \sum_{ I = \emptyset }^{ [k+1] }
                \restMaxFun[2k-1]{ R(\bx) }{ g_{ \ell(\bx) }}( a_I), 
        \end{align*}
        where the last inequality follows from 
        \Cref{lem::adding_k+1_direction_to_poly_j_gives_poly_j+1}.
        Thus, 
        \begin{equation}
        \label{thmeq::thm::multipointwise_bound_for_W^k_loc::eq::sublemma_1::case_j=k}
            \polynum{ k }(\bx')
            \sum_{ I = \emptyset }^{ [k] }
                \restMaxFun[2k-1]{ R(\bx) }{
                    \nabla^k f_w
                }(a_I)
            \le
            \polynum{ s }(\bx)
            \sum_{ I = \emptyset }^{ [k+1] }
                \restMaxFun[2k-1]{ R(\bx) }{ g_{ \ell(\bx) }}( a_I) 
        \end{equation}

        Finally, combining \eqref{thmeq::thm::multipointwise_bound_for_W^k_loc::eq::sublemma_1::case_j<k} and \eqref{thmeq::thm::multipointwise_bound_for_W^k_loc::eq::sublemma_1::case_j=k} with \eqref{thmeq::thm::multipointwise_bound_for_W^k_loc::eq::sublemma_1::after_inductive_step}, we get
        \begin{align*}
            \abs*{
                \diff_{ I = \emptyset }^{ [k+1] }
                    f(a_I)
            }
        \hspace{-1.2cm}&\hspace{1.2cm}\le
            C_{n,k}
            \sum_{ j = 1 }^{ k }
                \polynum{ j }( \bx' )
                \sum_{ I = \emptyset }^{ [k] }
                    \restMaxFun[k+j-1]{ R(\bx)}{\nabla^j f_w}(a_I), 
        \\
            &\le 
            C_{n,k}
            \del*{
                \sum_{j=1}^{k-1}
                    C_M n^k
                    \polynum{ j+1 }(\bx)
                    \sum_{ I = \emptyset }^{ [k+1] }
                        \restMaxFun[k+j]{ R(\bx) }{ \nabla^{j+1} f}( a_I)
                +
                \polynum{ s }(\bx)
                \sum_{ I = \emptyset }^{ [k+1] }
                    \restMaxFun[2k-1]{ R(\bx) }{ g_{ \ell(\bx) }}( a_I) 
            }
        \\
            &\le 
            C'_{n,k}
            \del*{
                \sum_{j=1}^{k-1}
                    \polynum{ j+1 }(\bx)
                    \sum_{ I = \emptyset }^{ [k+1] }
                        \restMaxFun[k+j]{ R(\bx) }{ \nabla^{j+1} f}( a_I)
                +
                \polynum{ s }(\bx)
                \sum_{ I = \emptyset }^{ [k+1] }
                    \restMaxFun[2k-1]{ R(\bx) }{ g_{ \ell(\bx) }}( a_I) 
            }
        \\
            &=
            C'_{n,k}
            \del*{
                \sum_{j=2}^{k}
                    \polynum{ j }(\bx)
                    \sum_{ I = \emptyset }^{ [k+1] }
                        \restMaxFun[k+j-1]{ R(\bx) }{ \nabla^{j} f}( a_I)
                +
                \polynum{ s }(\bx)
                \sum_{ I = \emptyset }^{ [k+1] }
                    \restMaxFun[2k-1]{ R(\bx) }{ g_{ \ell(\bx) }}( a_I) 
            },
        \end{align*}
        where we set 
        $
            C'_{n,k}
            \coloneq 
            C_{n,k}
            \max \del*{ C_M n^k, 1}.
        $
        This completes the proof
        of 
        \Cref{sublem::thm::multipointwise_bound_for_W^k_loc::eq::sublemma_1}.
    \end{proof}
    \begin{lemmaInThm}
    \label{sublem::thm::multipointwise_bound_for_W^k_loc::eq::sublemma_2}
        Under \Cref{asu::thm::multipointwise_bound_for_W^k_loc::inductive_assumptions},
        for all 
        $
            \bx = \set*{ x_I }_{ I \subseteq [k+1] }
            \subseteq 
            \lebesguePoints{f}
        $
        such that
        $ \polyset{ \set*{k+1} }(\bx)  > 0 $,
        \begin{align*}
        &
            \abs*{ 
                \diff_{ I = \emptyset }^{ [k] }
                    \del*{
                        \ballAv{
                            \polyset{ \set*{k+1} }(\bx)
                        }{
                            f
                        }(x_I)
                        - f(x_I)
                    }
            }
        \\
        &\hspace{2cm}
            \le 
            C_{n,k}''
            \del*{
                \sum_{ j = 2 }^{k}
                    \polynum{ j }( \bx)
                    \sum_{ I = \emptyset }^{ [k + 1] }
                        \restMaxFun[k+j]{R( \bx) }{ 
                            \nabla^j f
                        }( x_I )
                +
                \polynum{ s }( \bx )
                \sum_{ I = \emptyset }^{ [k + 1] }
                    \restMaxFun[2k]{R( \bx) }{ g_{\ell(\bx)} }( x_I )
            }
        \end{align*}
        for some $C_{n,k}'' \ge 0$ that depends only on $n$ and $k$.
    \end{lemmaInThm}
    \begin{proof}[Proof of \Cref{sublem::thm::multipointwise_bound_for_W^k_loc::eq::sublemma_2}]
        Fix 
        $
            \bx = \set*{ x_I }_{ I \subseteq [k+1] }
            \subseteq 
            \lebesguePoints{f}
        $
        such that
        $ 
        	\polyset{ \set*{k+1} }(\bx)  > 0 
        $,
        and denote
        \begin{equation*}
            B
            \coloneq
            \bigcap_{I \subseteq [k]}
                \del{
                    \lebesguePoints{f} - x_I + x_{\emptyset}
                }.
        \end{equation*}
        Let us note that if $b \in B$, 
        then $b + x_I - x_{\emptyset} \in \lebesguePoints{f}$
        for all $I \subseteq [k] $.
        Fix 
        $
            b \in
            B
            \cap 
            \ball{x_{\emptyset}, \polyset{ \set*{k+1} }(\bx) }
        $ 
        and denote 
        $
            w \coloneq b - x_{\emptyset}.
        $
        Then 
        $
            w \in \clball{0, \polyset{ \set*{k+1} }(\bx) }
        $
        and 
        \begin{equation*}
            x_{\emptyset}
            =
            b - w
            \in 
            B - w 
            =
            \del*{
                \bigcap_{I \subseteq [k]}
                \del{
                    \lebesguePoints{f} - x_I + x_{\emptyset}
                }
            } - w
            =
            \bigcap_{I \subseteq [k]}
                \del{
                    \lebesguePoints{f} - x_I + x_{\emptyset} - w
                }.
        \end{equation*}
        Moreover, for all $I \subseteq [k]$ we have
        $
            x_I \in \lebesguePoints{f},
        $
        hence
        \begin{equation*}
            x_{\emptyset}
            =
            x_I - x_I + x_{\emptyset}
            \in 
            \lebesguePoints{f} - x_I + x_{\emptyset}.
        \end{equation*}
        Therefore, 
        \begin{equation*}
            x_{\emptyset}
            \in 
            \bigcap_{I \subseteq [k] }
                \del*{
                    \lebesguePoints{f} 
                    - x_I + x_{\emptyset}
                }
            \cap 
            \bigcap_{I \subseteq [k] }
                \del*{
                    \lebesguePoints{f} 
                    - x_I + x_{\emptyset} - w
                }
            =
            A,
        \end{equation*}
        where $A$ is as in \Cref{sublem::thm::multipointwise_bound_for_W^k_loc::eq::sublemma_1}.

        Define $\set*{ b_I }_{ I \subseteq [k+1] }$
        as follows:
        \begin{equation*}
            \forall I \subseteq [k]
        \qquad 
            b_I \coloneq x_I
        \quad \text{and} \quad
            b_{I \cup \set*{k+1} }
            \coloneq x_I + w. 
        \end{equation*}
        Note that for all $I \subseteq [k]$
        we have
        \begin{equation*}
            b_I
            =
            x_{\emptyset} + x_I - x_{\emptyset}
            \quad \text{and} \quad
            b_{I \cup \set*{k+1} }
            = b_I + w.
        \end{equation*}
        Hence, since $x_{\emptyset} \in A,$
        by \Cref{sublem::thm::multipointwise_bound_for_W^k_loc::eq::sublemma_1}
        we have
        \begin{equation}
        \label{thmeq::thm::multipointwise_bound_for_W^k_loc::eq::sublemma_2::using_sublemma_1}
            \abs*{
                \diff_{ I = \emptyset }^{ [k+1] }
                    \mspace{-3mu}
                    f(b_I)
            }
            \! \le \! 
            C_{n,k}'
            \!
            \del*{
                \sum_{ j = 2 }^{k}
                	\mspace{-3mu}
                    \polynum{ j }( \bx)
                    \! \mspace{-3mu}
                    \sum_{ I = \emptyset }^{ [k+1] }
                    \! \mspace{-3mu}
                        \restMaxFun[k+j-1]{R( \bx) }{ 
                            \nabla^j f
                        }\del*{ b_{I } }
                +
                \polynum{ s }( \bx )
                \! \mspace{-3mu}
                \sum_{ I = \emptyset }^{ [k+1] }
                \! \mspace{-3mu}
                    \restMaxFun[2k-1]{ R( \bx) }{ g_{\ell(\bx)} }( b_I )
            \!
            }.
        \end{equation}
        Next, note that by 
        \Cref{lem::diff_B_diff_C_to_diff_BuC}
        we have
        \begin{equation*}
            \abs*{
                \diff_{ I = \emptyset }^{ [k+1] }
                    f(b_I)
            }
            =
            \abs*{ 
                \diff_{ I = \emptyset }^{ [k] }
                    f\del*{ b_{I \cup \set*{k+1} } }
                - \diff_{ I = \emptyset }^{ [k] }
                    f(b_I)
            }
            =
            \abs*{ 
                \diff_{ I = \emptyset }^{ [k] }
                    f\del*{ x_I + b - x_{\emptyset} }
                - \diff_{ I = \emptyset }^{ [k] }
                    f(x_I)
            }.
        \end{equation*}
        Also, we have
        \begin{align*}
            \sum_{ I = \emptyset }^{ [k+1] }
                \restMaxFun[2k-1]{ R( \bx) }{ g_{\ell(\bx)} }( b_I )
            &=
            \sum_{ I = \emptyset }^{ [k] }
                \restMaxFun[2k-1]{ R( \bx) }{ g_{\ell(\bx)} }\del*{
                    b_{ I \cup \set*{k+1} }
                }
            + \sum_{ I = \emptyset }^{ [k] }
                \restMaxFun[2k-1]{ R( \bx) }{ g_{\ell(\bx)} }( b_I )
        \\
            &=
            \sum_{ I = \emptyset }^{ [k] }
                \restMaxFun[2k-1]{ R( \bx) }{ g_{\ell(\bx)} }
                \del*{ x_I + b - x_{\emptyset} }
            + \sum_{ I = \emptyset }^{ [k] }
                \restMaxFun[2k-1]{ R( \bx) }{ g_{\ell(\bx)} }
                ( x_I ),
        \end{align*}
        and, for all $j \in [k]$,
        \begin{align*}
            \sum_{ I = \emptyset }^{ [k+1] }
                \restMaxFun[k+j-1]{ R( \bx) }{ 
                    \nabla^j f
                }( b_I )
                &=
            \sum_{ I = \emptyset }^{ [k] }
                \restMaxFun[k+j-1]{ R( \bx) }{ 
                    \nabla^j f
                }( b_{I \cup \set*{k+1} } )
            + \sum_{ I = \emptyset }^{ [k] }
                \restMaxFun[k+j-1]{ R( \bx) }{ 
                    \nabla^j f
                }( b_I )
        \\
                &=
            \sum_{ I = \emptyset }^{ [k] }
                \restMaxFun[k+j-1]{ R( \bx) }{ 
                    \nabla^j f
                }( x_I + b - x_{\emptyset} )
            + \sum_{ I = \emptyset }^{ [k] }
                \restMaxFun[k+j-1]{ R( \bx) }{ 
                    \nabla^j f
                }( x_I ).
        \end{align*}
        Thus, we can express \eqref{thmeq::thm::multipointwise_bound_for_W^k_loc::eq::sublemma_2::using_sublemma_1} equivalently as
        \begin{align}
            \abs*{ 
                \diff_{ I = \emptyset }^{ [k] }
                    f( x_I + b - x_{\emptyset} )
                - \diff_{ I = \emptyset }^{ [k] }
                    f(x_I)
            }
            \hspace{-5.4cm}&
            \notag
        \\
            &\le 
            C_{n,k}'
            \left(
                \sum_{ j = 2 }^{k}
                    \polynum{ j }( \bx)
                    \del*{
                        \sum_{ I = \emptyset }^{ [k] }
                            \restMaxFun[k+j-1]{ R( \bx) }{ 
                                \nabla^j f
                            }( x_I + b - x_{\emptyset} )
                        + \sum_{ I = \emptyset }^{ [k] }
                            \restMaxFun[k+j-1]{ R( \bx) }{ 
                                \nabla^j f
                            }( x_I )
                    }
            \right.
            \notag
        \\
            &\hspace{2cm}
            \left.
                +
                \polynum{ s }( \bx )
                \del*{
                    \sum_{ I = \emptyset }^{ [k] }
                        \restMaxFun[2k-1]{ R( \bx) }{ g_{\ell(\bx)} }
                        ( x_I + b - x_{\emptyset} )
                    + \sum_{ I = \emptyset }^{ [k] }
                        \restMaxFun[2k-1]{ R( \bx) }{ g_{\ell(\bx)} }
                        ( x_I )
                }
            \right).
            \label{thmeq::thm::multipointwise_bound_for_W^k_loc::eq::sublemma_2::equivalent_form_to_result_from_sublemma_1}
        \end{align}

        Since 
        $B$ is of full measure as an intersection of sets of full measure, 
        $
            B 
            \cap \ball{x_{\emptyset}, \polyset{\set*{k+1}}(\bx)}
        $
        is of full measure in 
        $
            \ball{x_{\emptyset}, \polyset{\set*{k+1}}(\bx)}
        $.
        Therefore, since \eqref{thmeq::thm::multipointwise_bound_for_W^k_loc::eq::sublemma_2::equivalent_form_to_result_from_sublemma_1} is true for all 
        $ 
            b \in B \cap \ball{x_{\emptyset}, \polyset{\set*{k+1}}(\bx)}
        $,
        it is true for almost all 
        $ 
            b \in \ball{x_{\emptyset}, \polyset{\set*{k+1}}(\bx)}
        $.

        In~consequence, 
        denoting 
        $
        	r
        	\coloneq 
        	\polyset{\set{k+1}}\del*{ \bx },
        $
        %
        \begin{align*}
            \ballAv{
                \polyset{\set*{k+1}}(\bx)
            }{
                \abs*{ 
                    \diff_{ I = \emptyset }^{ [k] }
                        f( x_I + \cdot - x_{\emptyset} )
                    - \diff_{ I = \emptyset }^{ [k] }
                        f(x_I)
                }
            }( x_{ \emptyset} )
            \hspace{-8.3cm}&
        \\
        	&=
            \ballAv{
                r
            }{
                \abs*{ 
                    \diff_{ I = \emptyset }^{ [k] }
                        f( x_I + \cdot - x_{\emptyset} )
                    - \diff_{ I = \emptyset }^{ [k] }
                        f(x_I)
                }
            }( x_{ \emptyset} )        	
        \\
            &\le 
            C_{n,k}'
            \mathsf{B}_{ 
                r
            }
            \mspace{-4mu}
            \left(
                \sum_{ j = 2 }^{k}
                    \polynum{ j }( \bx)
                    \del*{
                        \sum_{ I = \emptyset }^{ [k] }
                            \restMaxFun[k+j-1]{ R( \bx) }{ 
                                \nabla^j f
                            }\del*{ 
                            	x_I 
                            +
                            	\blank 
                            - 
                            	x_{\emptyset}
                            }
                        + \sum_{ I = \emptyset }^{ [k] }
                            \restMaxFun[k+j-1]{ R( \bx) }{ 
                                \nabla^j f
                            }( x_I )
                    }
            \right.
        \\
            &\hspace{2.4cm}
            \left.
                + \,
                \polynum{ s }( \bx )
                \del*{
                    \sum_{ I = \emptyset }^{ [k] }
                        \restMaxFun[2k-1]{ R( \bx) }{ g_{\ell(\bx)} }
                        \del*{ 
                        	x_I 
                        +
                        	\blank 
                        - 
                        	x_{\emptyset} 
                        }
                    + \sum_{ I = \emptyset }^{ [k] }
                        \restMaxFun[2k-1]{ R( \bx) }{ g_{\ell(\bx)} }
                        ( x_I )
                }
            \mspace{-4mu}
            \right)
            \mspace{-5mu}
            \del*{ x_{ \emptyset} },
        \intertext{
        	which, by the linearity of 
        	$
		        \mathsf{B}_{ 
                    r
                },
        	$
        }
            &=
            C_{n,k}' 
            \mspace{-4mu}
            \left(
                \sum_{ j = 2 }^{k}
                    \polynum{ j }( \bx)
                \mspace{-2mu}
                    \del*{
                        \sum_{ I = \emptyset }^{ [k] }
                            \ballAv{
                                r 
                            }{
                                \restMaxFun[k+j-1]{ R( \bx) }{ 
                                    \nabla^j f
                                }\del*{ 
                                	x_I
                                +
                                	\blank 
                                -
                                	x_{\emptyset} 
                                }
                            }(x_{\emptyset})
                    +
                   		\sum_{ I = \emptyset }^{ [k] }
                            \restMaxFun[k+j-1]{ R( \bx) }{ 
                                \nabla^j f
                            }( x_I )
                        \mspace{-4mu}
                    }
            \right.
        \\
            &\hspace{1.85cm}
            \left.
                + \,
                \polynum{ s }( \bx )
            \mspace{-2mu}
                \del*{
                    \sum_{ I = \emptyset }^{ [k] }
                        \ballAv{
                            r 
                        }{
                            \restMaxFun[2k-1]{ R( \bx) }{ g_{\ell(\bx)} }
                            \del*{ 
                            	x_I
                           	+
                           		\blank 
                           	- 
                           		x_{\emptyset}
                           	}
                        }(x_{\emptyset})
                +
                    \sum_{ I = \emptyset }^{ [k] }
                        \restMaxFun[2k-1]{ R( \bx) }{ g_{\ell(\bx)} }
                        ( x_I )
                    \mspace{-2mu}
                }
            \mspace{-3mu}
            \right)
            \mspace{-4mu},
        \intertext{which, by \Cref{lem::shifting_operators},}
            &\le
            C_{n,k}'
            \left(
                \sum_{ j = 2 }^{k}
                    \polynum{ j }( \bx)
                    \del*{
                        \sum_{ I = \emptyset }^{ [k] }
                            \ballAv{
                                r
                            }{
                                \restMaxFun[k+j-1]{ R( \bx) }{ 
                                    \nabla^j f
                                }
                            }(x_{I})
                        + \sum_{ I = \emptyset }^{ [k] }
                            \restMaxFun[k+j-1]{ R( \bx) }{ 
                                \nabla^j f
                            }( x_I )
                    }
            \right.
        \\
            &\hspace{1.87cm}
            \left.
                {}+
                \polynum{ s }( \bx )
                \del*{
                    \sum_{ I = \emptyset }^{ [k] }
                        \ballAv{
                            r
                        }{
                            \restMaxFun[2k-1]{ R( \bx) }{ g_{\ell(\bx)} }
                        }(x_{I})
                    + \sum_{ I = \emptyset }^{ [k] }
                        \restMaxFun[2k-1]{ R( \bx) }{ g_{\ell(\bx)} }
                        ( x_I )
                }
            \right),
        \intertext{
            which, since 
            $
            	r
            	=
                \polyset{\set*{k+1}}(\bx) 
                \le 2\ell(\bx) \le R(\bx)
            $,
        }
            &\le
            C_{n,k}'
            \left(
                \sum_{ j = 2 }^{k}
                    \polynum{ j }( \bx)
                    \del*{
                        \sum_{ I = \emptyset }^{ [k] }
                            \restMaxFun[k+j]{ R( \bx) }{ 
                                \nabla^j f
                            }(x_{I})
                        + \sum_{ I = \emptyset }^{ [k] }
                            \restMaxFun[k+j-1]{ R( \bx) }{ 
                                \nabla^j f
                            }( x_I )
                    }
            \right.
        \\
            &\hspace{1.8cm}
            \left.
                {}+
                \polynum{ s }( \bx )
                \del*{
                    \sum_{ I = \emptyset }^{ [k] }
                        \restMaxFun[2k]{ R( \bx) }{ g_{\ell(\bx)} }(x_{I})
                    + \sum_{ I = \emptyset }^{ [k] }
                        \restMaxFun[2k-1]{ R( \bx) }{ g_{\ell(\bx)} }
                        ( x_I )
                }
            \right),
        \intertext{
            and, by \Cref{lem::composing_maximal_functions},
        }
            &\le 
            2C_{n,k}'
            \del*{
                \sum_{ j = 2 }^{k}
                \polynum{ j }( \bx)
                \sum_{ I = \emptyset }^{ [k] }
                        \restMaxFun[k+j]{ R( \bx) }{ 
                            \nabla^j f
                        }(x_{I})
                +
                \polynum{ s }( \bx )
                \sum_{ I = \emptyset }^{ [k] }
                        \restMaxFun[2k]{ R( \bx) }{ g_{\ell(\bx)} }(x_{I})
            }
        \\
            &\le 
            2C_{n,k}'
            \del*{
                \sum_{ j = 2 }^{k}
                \polynum{ j }( \bx)
                \sum_{ I = \emptyset }^{ [k + 1] }
                        \restMaxFun[k+j]{ R( \bx) }{ 
                            \nabla^j f
                        }(x_{I})
                +
                \polynum{ s }( \bx )
                \sum_{ I = \emptyset }^{ [k + 1] }
                        \restMaxFun[2k]{ R( \bx) }{ g_{\ell(\bx)} }(x_{I})
            }.
        \end{align*}
        %
        Thus,
        \begin{align}
        &
            \ballAv{
                \polyset{\set*{k+1}}(\bx)
            }{
                \abs*{ 
                    \diff_{ I = \emptyset }^{ [k] }
                        f\del*{ x_I + \blank - x_{\emptyset} }
                    - \diff_{ I = \emptyset }^{ [k] }
                        f(x_I)
                }
            }( x_{ \emptyset} )
            \notag
        \\
        &\hspace{0.5cm}
            \le 
            2C_{n,k}'
            \del*{
                \sum_{ j = 2 }^{k}
                \polynum{ j }( \bx)
                \sum_{ I = \emptyset }^{ [k + 1] }
                        \restMaxFun[k+j]{ R( \bx) }{ 
                            \nabla^j f
                        }(x_{I})
                +
                \polynum{ s }( \bx )
                \sum_{ I = \emptyset }^{ [k + 1] }
                        \restMaxFun[2k]{ R( \bx) }{ g_{\ell(\bx)} }(x_{I})
            }.
            \label{thmeq::thm::multipointwise_bound_for_W^k_loc::eq::sublemma_2::estimate_with_average_over_abs}
        \end{align}
        It remains to note that by 
        \Cref{lem::shifting_operators}
        we have
        \begin{align*}
            \abs*{ 
                \diff_{ I = \emptyset }^{ [k] }
                    \del*{
                        \ballAv{
                            \polyset{\set*{k+1}}(\bx)
                        }{
                            f
                        }(x_I)
                        - f(x_I)
                    }
            }
            &=
            \abs*{ 
                \diff_{ I = \emptyset }^{ [k] }
                    \del*{
                        \ballAv{
                            \polyset{\set*{k+1}}(\bx)
                        }{
                            f\del*{ x_I + \blank - x_{\emptyset} }
                        }(x_{\emptyset})
                        - f(x_I)
                    }
            }
        \\
            &=
            \abs*{ 
                \ballAv{
                    \polyset{\set*{k+1}}(\bx)
                }{
                    \diff_{ I = \emptyset }^{ [k] }
                         f\del*{ x_I + \blank - x_{\emptyset} }
                        - f(x_I)
                }(x_{\emptyset})
            }
        \\
            &\le 
            \ballAv{
                \polyset{\set*{k+1}}(\bx)
            }{
                \abs*{ 
                    \diff_{ I = \emptyset }^{ [k] }
                        f\del*{ x_I + \blank - x_{\emptyset} }
                    - \diff_{ I = \emptyset }^{ [k] }
                        f(x_I)
                }
            }( x_{ \emptyset} ),
        \end{align*}
        which combined with \eqref{thmeq::thm::multipointwise_bound_for_W^k_loc::eq::sublemma_2::estimate_with_average_over_abs} gives
        \begin{align*}
        &\abs*{ 
                \diff_{ I = \emptyset }^{ [k] }
                    \del*{
                        \ballAv{
                            \polyset{\set*{k+1}}(\bx)
                        }{
                            f
                        }(x_I)
                        - f(x_I)
                    }
            }
        \\
        &\hspace{2cm}\le 
            2C_{n,k}'
            \del*{
                \sum_{ j = 2 }^{k}
                \polynum{ j }( \bx)
                \sum_{ I = \emptyset }^{ [k+ 1] }
                        \restMaxFun[k+j]{ R( \bx) }{ 
                            \nabla^j f
                        }(x_{I})
                +
                \polynum{ s }( \bx )
                \sum_{ I = \emptyset }^{ [k+ 1] }
                        \restMaxFun[2k]{ R( \bx) }{ g_{\ell(\bx)} }(x_{I})
            }.
        \end{align*}
        Therefore, we get the desired inequality with a constant
        $
            C_{n,k}'' \coloneq 2C_{n,k}'.
        $
    \end{proof}
    \begin{lemmaInThm}
    \label{sublem::thm::multipointwise_bound_for_W^k_loc::eq::sublemma_3}
        Under \Cref{asu::thm::multipointwise_bound_for_W^k_loc::inductive_assumptions},
        for all 
        $
            \bx = \set*{ x_I }_{ I \subseteq [k+1] }
            \subseteq 
            \lebesguePoints{f}
        $
        such that
        $ \polyset{ \set*{k+1} }(\bx)  > 0 $,
        \begin{align*}
        &\abs*{ 
                \diff_{ I = \emptyset }^{ [k] }
                    \del*{
                        \ballAv{
                            \polyset{\set*{k+1}}(\bx)
                        }{
                            f
                        }(x_I + x_{\set*{k+1}} - x_{\emptyset})
                        -
                        \ballAv{
                            \polyset{\set*{k+1}}(\bx)
                        }{
                            f
                        }(x_I)                      
                    }
            }
        \\
        &\hspace{2cm}\le 
            C_{n,k}'''
            \del*{
                \sum_{ j = 2 }^{k}
                    \polynum{ j }( \bx)
                    \sum_{ I = \emptyset }^{ [k + 1] }
                        \restMaxFun[k+j]{ R( \bx) }{ 
                            \nabla^j f
                        }( x_I )
                +
                \polynum{ s }( \bx )
                \sum_{ I = \emptyset }^{ [k + 1] }
                    \restMaxFun[2k]{ R( \bx) }{ g_{\ell(\bx)} }( x_I )
            }
        \end{align*}
        for some $C_{n,k}''' \ge 0$ that depends only on $n$ and $k$.
    \end{lemmaInThm}
    \begin{proof}[Proof of \Cref{sublem::thm::multipointwise_bound_for_W^k_loc::eq::sublemma_3}]
        Fix 
        $
            \bx = \set*{ x_I }_{ I \subseteq [k+1] }
            \subseteq 
            \lebesguePoints{f}
        $
        such that
        $ 
        	\polyset{ \set*{k+1} }(\bx)  > 0. 
        $
        Let us denote $w \coloneq x_{ \set*{k+1} } - x_{ \emptyset }$
        and
        \begin{equation*}
            A
            \coloneq
            \bigcap_{I \subseteq [k] }
                \del*{
                    \lebesguePoints{f} 
                    - x_I + x_{\emptyset}
                }
            \cap 
            \bigcap_{I \subseteq [k] }
                \del*{
                    \lebesguePoints{f} 
                    - x_I + x_{\emptyset} - w
                }.
        \end{equation*}
        Fix $a \in A$ and define $\set*{ a_I }_{ I \subseteq [k+1] }$
       as follows:
        \begin{equation*}
            \forall I \subseteq [k] 
           \qquad
            a_I \coloneq a + x_I - x_{ \emptyset }
        \quad \text{and} \quad 
            a_{I \cup \set*{k+1} } \coloneq a_I + w.
        \end{equation*}
        Let us note that since
        $
            \norm{w} \le \polyset{ \set*{k+1}}( \bx ),
        $
        $A$,
        and $\set*{a_I}_{ I \subseteq [k+1]}$ are as in 
        \Cref{sublem::thm::multipointwise_bound_for_W^k_loc::eq::sublemma_1}.
        Therefore, 
        %
        \begin{equation}
        \label{thmeq::thm::multipointwise_bound_for_W^k_loc::eq::sublemma_3::using_sublemma_1}
            \hspace{-0.25cm}
            \abs*{
                \diff_{ I = \emptyset }^{ [k+1] }
                    f(a_I)
            }
            \le 
            C_{n,k}'
            \del*{
            \!
                \sum_{ j = 2 }^{k}
                    \polynum{ j }( \bx)
                \!
                \smashoperator{
                    \sum_{ I = \emptyset }^{ [k+1] }
                }
                \!
                        \restMaxFun[k+j-1]{ R( \bx) }{ 
                            \nabla^j f
                        }( a_I )
                +
                \polynum{ s }( \bx )
            \!
            \smashoperator{
                \sum_{ I = \emptyset }^{ [k+1] }
            }
            \!
                    \restMaxFun[2k-1]{ R( \bx) }{ g_{\ell(\bx)} }( a_I )
            \!
            }
            .
        \end{equation}

        Next, note that by 
        \Cref{lem::diff_B_diff_C_to_diff_BuC}
        we have
        \begin{equation*}
            \abs*{
                \diff_{ I = \emptyset }^{ [k+1] }
                    f(a_I)
            }
            \!
            =
            \!
            \abs*{ 
                \diff_{ I = \emptyset }^{ [k] }
                    f\del*{ a_{I \cup \set*{k+1} } }
                - \diff_{ I = \emptyset }^{ [k] }
                    f(a_I)
            }
            \!
            =
            \!
            \abs*{ 
                \diff_{ I = \emptyset }^{ [k] }
                    f\del*{ a + x_I - x_{\emptyset} + w}
                - \diff_{ I = \emptyset }^{ [k] }
                    f\del*{ a + x_I - x_{\emptyset} }
            }.
        \end{equation*}
        Also, we have
        \begin{multline*}
            \sum_{ I = \emptyset }^{ [k+1] }
                \restMaxFun[2k-1]{ R( \bx) }{ 
                    g_{\ell(\bx)} 
                }(a_I)
        	=
            \sum_{ I = \emptyset }^{ [k] }
                \restMaxFun[2k-1]{ R( \bx) }{  
                    g_{\ell(\bx)} 
                }\del*{
                    a_{I \cup \set*{k+1} }
                }
            + \sum_{ I = \emptyset }^{ [k] }
                \restMaxFun[2k-1]{ R( \bx) }{ g_{\ell(\bx)} }\del*{
                	a_I 
                }
        \\
        	=
            \sum_{ I = \emptyset }^{ [k] }
                \restMaxFun[2k-1]{ R( \bx) }{ g_{\ell(\bx)} }
                \del*{ a + x_I - x_{\emptyset} + w}
            + \sum_{ I = \emptyset }^{ [k] }
                \restMaxFun[2k-1]{ R( \bx) }{ g_{\ell(\bx)} }
                \del*{ a + x_I - x_{\emptyset} },
        \end{multline*}
        and, for all $j \in [k]$,
        \begin{multline*}
            \sum_{ I = \emptyset }^{ [k+1] }
                \restMaxFun[k+j-1]{ R( \bx) }{ 
                    \nabla^j f
                }( a_I )
        	=
            \sum_{ I = \emptyset }^{ [k] }
                \restMaxFun[k+j-1]{ R( \bx) }{ 
                    \nabla^j f
                }\del*{
                	a_{I \cup \set*{k+1} } 
                }
            + \sum_{ I = \emptyset }^{ [k] }
                \restMaxFun[k+j-1]{ R( \bx) }{ 
                    \nabla^j f
                }\del*{
                	a_I 
                }
        \\
        	=
            \sum_{ I = \emptyset }^{ [k] }
                \restMaxFun[k+j-1]{ R( \bx) }{ 
                    \nabla^j f
                }\del*{ a + x_I - x_{\emptyset} + w}
            + \sum_{ I = \emptyset }^{ [k] }
                \restMaxFun[k+j-1]{ R( \bx) }{ 
                    \nabla^j f
                }\del*{ a + x_I - x_{\emptyset} }.
        \end{multline*}
        Thus, we can express \eqref{thmeq::thm::multipointwise_bound_for_W^k_loc::eq::sublemma_3::using_sublemma_1} equivalently as
        \begin{align}
            &
            \hspace{-0.25cm}
            \abs*{ 
                \diff_{ I = \emptyset }^{ [k] }
                    f\del*{ a + x_I - x_{\emptyset} + w}
                - \diff_{ I = \emptyset }^{ [k] }
                    f\del*{ a + x_I - x_{\emptyset} }
            }
            \notag
        \\
            &
            \hspace{-0.35cm}
            \le 
            C_{n,k}'\! \!
            \left(
            \!
                \smashoperator{
                \sum_{ j = 2 }^{k}
                }
                    \polynum{ j }( \bx)
                    \!
                    \del*{
                    	\smashoperator{
                        \sum_{ I = \emptyset }^{ [k] }
                        }
                    \!
                            \restMaxFun[k+j-1]{ R( \bx) }{ 
                                \nabla^j f
                            }
                        \!
                            \del*{ 
                            	a 
                           	{+} 
                           		x_I 
                       		{-} 
                       			x_{\emptyset} 
                   			{+} 
                   				w}
                {+} 
                \!
                        \smashoperator{
                        \sum_{ I = \emptyset }^{ [k] }
                        }
                    \!
                            \restMaxFun[k+j-1]{ R( \bx) }{ 
                                \nabla^j f
                            }
                        \!
                            \del*{ 
                            	a 
                           	{+} 
                           		x_I 
                       		{-} 
                       			x_{\emptyset} 
                       		}
                       		\!
                    }
            \right.
            \notag
        \\
            &
            \hspace{-0.35cm}
            \hspace{1.85cm}
            \left.
                \!  + \,
                \polynum{ s }( \bx )
            \!
                \del*{
                \!
                	\smashoperator{
                    \sum_{ I = \emptyset }^{ [k] }
                    }
                \!
                        \restMaxFun[2k-1]{ R( \bx) }{ g_{\ell(\bx)} }
                    \!
                        \del*{ 
                        	a 
                       	{+} 
                       		x_I 
                   		{-} 
                   			x_{\emptyset} 
               			{+} 
               				w
               			}
            {+} 
            \!
                    \smashoperator{
                    \sum_{ I = \emptyset }^{ [k] }
                    }
                \!
                        \restMaxFun[2k-1]{ R( \bx) }{ g_{\ell(\bx)} }
                    \!
                        \del*{ 
                        	a 
                       	{+} 
                       		x_I 
                   		{-} 
                   			x_{\emptyset} 
                   		}
                } \!
            \! 
            \right)
        \hspace{-3pt} .
            \label{thmeq::thm::multipointwise_bound_for_W^k_loc::eq::sublemma_3::equivalent_form_to_result_from_sublemma_1}
        \end{align}

        Since 
        $A$ is of full measure as an intersection of sets of full measure, 
        $
            A 
            \cap 
            \ball{
                x_{\emptyset}, 
                \polyset{ \set*{k+1} }(\bx)
            }
        $
        is of full measure in 
        $
            \ball{
                x_{\emptyset}, 
                \polyset{\set*{k+1}}(\bx)
            }
        $.
        Therefore, since \eqref{thmeq::thm::multipointwise_bound_for_W^k_loc::eq::sublemma_3::equivalent_form_to_result_from_sublemma_1} is true for all 
        $ a \in A$,
        it is true for almost all 
        $ 
            a 
            \in 
            \ball{
                x_{\emptyset}, 
                \polyset{\set*{k+1}}(\bx)
            }
        $.
        In consequence, 
        denoting
        $
        	r
        	\coloneq 
        	\polyset{ \set{k+1} }\del*{ \bx },
        $
        \begin{align*}
            \ballAv{
                \polyset{\set*{k+1}}(\bx)
            }{
                \abs*{ 
                    \diff_{ I = \emptyset }^{ [k] }
                        f\del*{ \blank + x_I - x_{\emptyset} + w}
                    - \diff_{ I = \emptyset }^{ [k] }
                        f\del*{ \blank + x_I - x_{\emptyset} }
                }
            }( x_{ \emptyset} )
            \hspace{-10.7cm}&
        \\    
        	&=
            \ballAv{
                r
            }{
                \abs*{ 
                    \diff_{ I = \emptyset }^{ [k] }
                        f\del*{ \blank + x_I - x_{\emptyset} + w}
                    - \diff_{ I = \emptyset }^{ [k] }
                        f\del*{ \blank + x_I - x_{\emptyset} }
                }
            }( x_{ \emptyset} )       	 
        \\
            &\le
            C_{n,k}'
            \mathsf{B}_{
                r
            }
            \mspace{-4mu}
            \left(
            \mspace{-3mu}
                \sum_{ j = 2 }^{k}
                \mspace{-3mu}
                    \polynum{ j }( \bx)
                    \mspace{-4mu}
                    \del*{
                    	\mspace{-3mu}
                        \sum_{ I = \emptyset }^{ [k] }
                        \mspace{-3mu}
                            \restMaxFun[k+j-1]{ R( \bx) }{ 
                                \nabla^j f
                            }
                            \mspace{-2mu}
                            \del*{ 
                            \blank {+} x_I {-} x_{\emptyset} {+} w}
                    +
	                    \mspace{-3mu}
                        \sum_{ I = \emptyset }^{ [k] }
                        \mspace{-3mu}
                            \restMaxFun[k+j-1]{ R( \bx) }{ 
                                \nabla^j f
                            }
                            \mspace{-2mu}
                            \del*{ \blank {+} x_I {-} x_{\emptyset} }
                    \mspace{-5mu}
                    }
            \right.
        \\
            &\hspace{2.2cm}
            \left.
                + \,
                \polynum{ s }( \bx )
                \mspace{-4mu}
                \del*{
                	\mspace{-5mu}
                    \sum_{ I = \emptyset }^{ [k] }
                    \mspace{-3mu}
                        \restMaxFun[2k-1]{ R( \bx) }{ g_{\ell(\bx)} }
                        \mspace{-2mu}
                        \del*{ \blank {+} x_I {-} x_{\emptyset} + w}
                +
                	\mspace{-3mu}
                    \sum_{ I = \emptyset }^{ [k] }
                    \mspace{-3mu}
                        \restMaxFun[2k-1]{ R( \bx) }{ g_{\ell(\bx)} }
                        \mspace{-2mu}
                        \del*{ \blank {+} x_I {-} x_{\emptyset} }
                    \mspace{-7mu}
                }
                \mspace{-9mu}
            \right) \mspace{-7mu} ( x_{ \emptyset} ),
        \intertext{
            which, by the linearity of 
            $
                \mathsf{B}_{ r }
            $ 
            and \Cref{lem::shifting_operators},
        }
            &\le
            C_{n,k}'
            \left(
                \sum_{ j = 2 }^{k}
                    \polynum{ j }( \bx)
                    \del*{
                        \sum_{ I = \emptyset }^{ [k] }
                            \ballAv{ 
                                r
                            }{
                                \restMaxFun[k+j-1]{ R( \bx) }{ 
                                    \nabla^j f
                                }
                            }(x_I + w)
                    + 
                        \sum_{ I = \emptyset }^{ [k] }
                            \ballAv{ 
                                r
                            }{
                                \restMaxFun[k+j-1]{ R( \bx) }{ 
                                    \nabla^j f
                                }
                            }(x_I)
                    }
            \right.
        \\
            &\hspace{2cm}
            \left.
                + \,
                \polynum{ s }( \bx )
                \del*{
                    \sum_{ I = \emptyset }^{ [k] }
                        \ballAv{ 
                            r
                        }{
                            \restMaxFun[2k-1]{ R( \bx) }{ g_{\ell(\bx)} }
                        }(x_I + w)
                +
                    \sum_{ I = \emptyset }^{ [k] }
                        \ballAv{ 
                            r
                        }{
                            \restMaxFun[2k-1]{ R( \bx) }{ g_{\ell(\bx)} }
                        }(x_I)
                }
            \right),
        \intertext{
            which, by \Cref{lem::centered_max_fun_bounds_shifted_average} and the fact that we have
           	$
           		r = \polyset{\set*{k+1}}(\bx),
           	$
            $
            	\norm{w}
            	\le
                \polyset{\set*{k+1}}(\bx) 
                \le 
                \ell(\bx),
            $ 
            and 
            $2\ell(\bx) \le R(\bx)$,
        }
            &\le
            C_{n,k}'
            \left(
                \sum_{ j = 2 }^{k}
                    \polynum{ j }( \bx)
                    \del*{
                        \sum_{ I = \emptyset }^{ [k] }
                            2^n
                            \restMaxFun[k+j]{ R( \bx) }{ 
                                \nabla^j f
                            }(x_{I})
                        + \sum_{ I = \emptyset }^{ [k] }
                            \restMaxFun[k+j]{ R( \bx) }{ 
                                \nabla^j f
                            }( x_I )
                    }
            \right.
        \\
            &\hspace{2cm}
            \left.
                + \,
                \polynum{ s }( \bx )
                \del*{
                    \sum_{ I = \emptyset }^{ [k] }
                        2^n
                        \restMaxFun[2k]{ R( \bx) }{ g_{\ell(\bx)} }(x_{I})
                    + \sum_{ I = \emptyset }^{ [k] }
                        \restMaxFun[2k]{ R( \bx) }{ g_{\ell(\bx)} }
                        ( x_I )
                }
            \right),
        \intertext{
            which, by setting 
            $C_{n,k}''' \coloneq C_{n,k}'(2^n+1)$, 
        }
            &= 
            C_{n,k}'''
            \del*{
                \sum_{ j = 2 }^{k}
                \polynum{ j }( \bx)
                \sum_{ I = \emptyset }^{ [k] }
                        \restMaxFun[k+j]{ R( \bx) }{ 
                            \nabla^j f
                        }(x_{I})
                +
                \polynum{ s }( \bx )
                \sum_{ I = \emptyset }^{ [k] }
                        \restMaxFun[2k]{ R( \bx) }{ g_{\ell(\bx)} }(x_{I})
            }
        \\
            &\le 
            C_{n,k}'''
            \del*{
                \sum_{ j = 2 }^{k}
                \polynum{ j }( \bx)
                \sum_{ I = \emptyset }^{ [k + 1] }
                        \restMaxFun[k+j]{ R( \bx) }{ 
                            \nabla^j f
                        }(x_{I})
                +
                \polynum{ s }( \bx )
                \sum_{ I = \emptyset }^{ [k + 1] }
                        \restMaxFun[2k]{ R( \bx) }{ g_{\ell(\bx)} }(x_{I})
            }.
        \end{align*}
        Thus,
        \begin{align}
        &
            \ballAv{
                \polyset{\set*{k+1}}(\bx)
            }{
                \abs*{ 
                    \diff_{ I = \emptyset }^{ [k] }
                        f\del*{ \cdot + x_I - x_{\emptyset} + w}
                    - \diff_{ I = \emptyset }^{ [k] }
                        f\del*{ \cdot + x_I - x_{\emptyset} }
                }
            }( x_{ \emptyset} )
        \notag
        \\
        &\hspace{1cm}
            \le 
            C_{n,k}'''
            \del*{
                \sum_{ j = 2 }^{k}
                \polynum{ j }( \bx)
                \sum_{ I = \emptyset }^{ [k+1] }
                        \restMaxFun[k+j]{ R( \bx) }{ 
                            \nabla^j f
                        }(x_{I})
                +
                \polynum{ s }( \bx )
                \sum_{ I = \emptyset }^{ [k+1] }
                        \restMaxFun[2k]{ R( \bx) }{ g_{\ell(\bx)} }(x_{I})
            }.
            \label{thmeq::thm::multipointwise_bound_for_W^k_loc::eq::sublemma_3::estimate_with_average_over_abs}
        \end{align}

        It remains to note that by 
        \Cref{lem::shifting_operators}
        we have
        \begin{align*}
            &
            \abs*{ 
                \diff_{ I = \emptyset }^{ [k] }
                    \del*{
                        \ballAv{
                            \polyset{\set*{k+1}}(\bx)
                        }{
                            f
                        }\del*{ x_I + x_{\set*{k+1}} - x_{\emptyset} }
                    -
                        \ballAv{
                            \polyset{\set*{k+1}}(\bx)
                        }{
                            f
                        }(x_I)
                    }
            }
            \\
            &=
            \abs*{ 
                \diff_{ I = \emptyset }^{ [k] }
                    \del*{
                        \ballAv{
                            \polyset{ \set*{k+1} }(\bx)
                        }{
                            f\del*{ \blank + x_I - x_{\emptyset} + w }
                        }(x_{\emptyset})
                        - 
                        \ballAv{
                            \polyset{ \set*{k+1} }(\bx)
                        }{
                            f\del*{ \blank + x_I - x_{\emptyset} }
                        }(x_{\emptyset})
                    }
            }
        \\
            &=
            \abs*{ 
                \ballAv{
                    \polyset{ \set*{k+1} }(\bx)
                }{
                    \diff_{ I = \emptyset }^{ [k] }
                        \del*{
                            f\del*{ \blank + x_I - x_{\emptyset} + w }
                            - f\del*{ \blank + x_I - x_{\emptyset} }
                        }
                }(x_{\emptyset})
            }
        \\
            &\le 
            \ballAv{
                \polyset{ \set*{k+1} }(\bx)
            }{
                \abs*{ 
                    \diff_{ I = \emptyset }^{ [k] }
                        \del*{
                            f\del*{ \blank + x_I - x_{\emptyset} + w }
                            - f\del*{ \blank + x_I - x_{\emptyset} }
                        }
                }
            }( x_{ \emptyset} ),
        \end{align*}
        which combined with \eqref{thmeq::thm::multipointwise_bound_for_W^k_loc::eq::sublemma_3::estimate_with_average_over_abs} gives
        \begin{align*}
        &
            \abs*{ 
                \diff_{ I = \emptyset }^{ [k] }
                    \del*{
                        \ballAv{
                            \polyset{ \set*{k+1} }(\bx)
                        }{
                            f
                        }\del*{ 
                        	x_I + x_{\set*{k+1}} - x_{\emptyset}
                        }
                    -
                        \ballAv{
                            \polyset{ \set*{k+1} }(\bx)
                        }{
                            f
                        }(x_I)
                    }
            }
        \\
        &\hspace{2cm}
            \le 
            C_{n,k}'''
            \del*{
                \sum_{ j = 2 }^{k}
                \polynum{ j }( \bx)
                \sum_{ I = \emptyset }^{ [k+1] }
                        \restMaxFun[k+j]{ R( \bx) }{ 
                            \nabla^j f
                        }\del*{ x_{I} }
                +
                \polynum{ s }( \bx )
                \sum_{ I = \emptyset }^{ [k+1] }
                        \restMaxFun[2k]{ R( \bx) }{ g_{\ell(\bx)} }\del*{ x_{I} }
            },
        \end{align*}
        as claimed.
    \end{proof}
    \begin{lemmaInThm}
    \label{sublem::thm::multipointwise_bound_for_W^k_loc::eq::sublemma_4}
        Suppose \Cref{asu::thm::multipointwise_bound_for_W^k_loc::inductive_assumptions}
        are true.
        Let 
        $
            \bx = \set*{ x_I }_{ I \subseteq [k+1] }
            \subseteq 
            \lebesguePoints{f}
        $
        be such that
        $ \polyset{\set*{k+1}}(\bx)  > 0 $.
        For $I \subseteq [k]$ denote
        $ v_I \coloneq x_{ I \cup \set*{k+1}} - x_{I}$.
        Let $S \subseteq [k]$ be nonempty and for 
        $I,J \subseteq [k]$ define 
        $ 
            x_{I,J} \coloneq v_I + x_J.
        $
        Then,
        for some $\wave{C}_{n,k} \ge 0$ that depends only on $n$ and $k$,
        \begin{equation*}
            \abs*{ 
                \diff_{  I = \emptyset }^{ S }
                    \diff_{ J = S }^{ [k] }
                        \ballAv{ 
                            \polyset{\set*{k+1} } ( \bx) 
                        }{
                            f
                        }( x_{I,J} )
            }
            \le 
            \wave{C}_{n,k}
            \sum_{ j = 1 }^{k}
                \polynum{ j }( \bx)
                \sum_{ I = \emptyset }^{ [k+1] }
                    \restMaxFun[k+j]{R( \bx) }{ 
                        \nabla^j f
                    }(x_{I}).
        \end{equation*}
        %
        %
    \end{lemmaInThm}
    \begin{proof}[Proof of \Cref{sublem::thm::multipointwise_bound_for_W^k_loc::eq::sublemma_4}.]
        First of all, let us define
        $
            \by = \set*{ y_I }_{ I \subseteq [k] }
        $
        as follows:
        \begin{equation*}
            \forall I \subseteq [k]
        \qquad 
            y_I \coloneq v_{I \cap S} + x_{I \cup S}.
        \end{equation*}
        Note that by 
        \Cref{lem::estimate_cube_with_points_and_vectors}
        for all $j \in [k]$ we have
        \begin{equation}
        \label{thmeq::thm::multipointwise_bound_for_W^k_loc::eq::sublemma_4::using_estimate_cube_with_points_and_vectors}
            \polynum{ j }(\by)
            \le 
            2^{k^2} \partitions{k} \polynum{ j }( \bx )
            \quad \text{and} \quad 
            \ell( \by )
            \le 
            2^k \ell( \bx ).
        \end{equation}

        Next, let us define
        \begin{equation*}
            D \coloneq 
            \bigcap_{ I \subseteq [k] }
                \del*{
                    \lebesguePoints{f} - y_I + y_{\emptyset}
                }.
        \end{equation*}
        Fix $d \in D$ and define 
        $\bd = \set*{d_I}_{ I \subseteq [k] }$ 
        as follows:
        \begin{equation*}
            \forall I \subseteq [k] 
        \qquad 
            d_I \coloneq d + y_I - y_{\emptyset}.
        \end{equation*}
        Then for all $I \subseteq [k]$ we have 
        $
            d_I
            =
            d + y_I - y_{\emptyset}
            \in 
            \lebesguePoints{f}.
        $
        Hence,
        by
        \Cref{asu::thm::multipointwise_bound_for_W^k_loc::inductive_assumptions},
        \begin{equation*}
            \abs*{
                \diff_{ I = \emptyset }^{ [k] }
                    f(d_I)
            }
            \le 
            C_{n,k}
            \sum_{ j = 1 }^{ k }
                \polynum{ j }( \bd )
                \sum_{ I = \emptyset }^{ [k] }
                    \restMaxFun[k+j-1]{2^{k^2}\ell(\bd)}{
                        \nabla^j f
                    }(d_I).
        \end{equation*}

        Since 
        $ y \mapsto d + y - y_{\emptyset} $ is an isometry of $\bR^n$,
        by \Cref{lem::invariant_transformations_of_poly},
        for all $j \in [k]$ we have 
        \begin{equation*}
            \polynum{ j }( \bd )
            =
            \polynum{ j }( d + \by - y_{\emptyset} )
            =
            \polynum{ j }(\by)
            \le 
            2^{k^2} \partitions{k} \polynum{ j }( \bx )
        \,\,\,\, \text{and} \,\,\,\, 
            \ell( \bd) 
            =
            \ell( d + \by - y_{\emptyset} )
            =
            \ell( \by )
            \le 
            2^k \ell( \bx ),
        \end{equation*}
        where the inequalities follow from 
        \eqref{thmeq::thm::multipointwise_bound_for_W^k_loc::eq::sublemma_4::using_estimate_cube_with_points_and_vectors}.
        Hence,
        \begin{equation*}
            2^{k^2} \ell( \bd )
            \le 
            2^{k^2} 2^k \ell( \bx )
            \le 
            2^{(k+1)^2} \ell( \bx )
            = 
            R(\bx).
        \end{equation*}
        In consequence, 
        \begin{align*}
            \abs*{
                \diff_{ I = \emptyset }^{ [k] }
                    f(d_I)
            }
            &\le 
            C_{n,k}
            \sum_{ j = 1 }^{ k }
                \polynum{ j }( \bd )
                \sum_{ I = \emptyset }^{ [k] }
                    \restMaxFun[k+j-1]{2^{k^2}\ell(\bd)}{
                        \nabla^j f
                    }(d_I)
        \\
            &\le 
            C_{n,k}
            \sum_{ j = 1 }^{ k }
                2^{k^2} \partitions{k} \polynum{ j }( \bx )
                \sum_{ I = \emptyset }^{ [k] }
                    \restMaxFun[k+j-1]{R(\bx)}{
                        \nabla^j f
                    }(d_I).
        \end{align*}
        Therefore, since
        $d_I = d + y_I - y_{\emptyset}$
        for all $I \subseteq [k]$, 
        we have
        \begin{equation}
        \label{thmeq::thm::multipointwise_bound_for_W^k_loc::eq::sublemma_4::estimate_for_the_elements_of_D}
            \abs*{
                \diff_{ I = \emptyset }^{ [k] }
                    f( d + y_I - y_{\emptyset} )
            }
            \le 
            \wave{C}'_{n,k}
            \sum_{ j = 1 }^{ k }
                \polynum{ j }( \bx )
                \sum_{ I = \emptyset }^{ [k] }
                    \restMaxFun[k+j-1]{R(\bx)}{
                        \nabla^j f
                    }( d + y_I - y_{\emptyset} ),
        \end{equation}
        where we set
        $
            \wave{C}'_{n,k}
            \coloneq
            2^{k^2} \partitions{k}
            C_{n,k} .
        $
        
        Next, let us note that $D$ is of full measure as an intersection of sets of full measure.
        Therefore, 
        $
            D 
            \cap 
            \ball{
                y_{\emptyset}, 
                \polyset{ \set*{k+1}}( \bx ) 
            }
        $
        is of full measure in 
        $
            \ball{
                y_{\emptyset}, 
                \polyset{ \set*{k+1}}( \bx ) 
            }
        $.
        Therefore, 
        \eqref{thmeq::thm::multipointwise_bound_for_W^k_loc::eq::sublemma_4::estimate_for_the_elements_of_D}
        is true for almost all 
        $
            d 
            \in 
            \ball{
                y_{\emptyset}, 
                \polyset{ \set*{k+1}}( \bx ) 
            }
        $.
        Hence,
        \begin{align*}
            \ballAv{ 
                \polyset{ \set*{k+1}}( \bx ) 
            }{
                \abs*{
                    \diff_{ I = \emptyset }^{ [k] }
                        f\del*{ \blank + y_I - y_{\emptyset} }
                }
            }(y_{\emptyset})
        \hspace{-5cm}&
        \\
        &\le 
            \wave{C}'_{n,k}
            \ballAv{ 
                \polyset{ \set*{k+1}}( \bx ) 
            }{
                \sum_{ j = 1 }^{ k }
                    \polynum{ j }( \bx )
                    \sum_{ I = \emptyset }^{ [k] }
                        \restMaxFun[k+j-1]{R(\bx)}{
                            \nabla^j f
                        }\del*{ \blank + y_I - y_{\emptyset} }
            }(y_{\emptyset}),
        \intertext{
            which, by the linearity of 
            $\mathsf{B}_{ \polyset{ \set*{k+1}}(\bx) } $,
        }
        &= 
            \wave{C}'_{n,k}
            \sum_{ j = 1 }^{ k }
                \polynum{ j }( \bx )
                \sum_{ I = \emptyset }^{ [k] }
                    \ballAv{ 
                        \polyset{ \set*{k+1} }( \bx ) 
                    }{
                        \restMaxFun[k+j-1]{R(\bx)}{
                            \nabla^j f
                        }\del*{ \blank + y_I - y_{\emptyset} }
                    }(y_{\emptyset}),
        \intertext{
            which, by 
            \Cref{lem::shifting_operators},
        }
        &= 
            \wave{C}'_{n,k}
            \sum_{ j = 1 }^{ k }
                \polynum{ j }( \bx )
                \sum_{ I = \emptyset }^{ [k] }
                    \ballAv{ 
                        \polyset{ \set*{k+1}}( \bx ) 
                    }{
                        \restMaxFun[k+j-1]{R(\bx)}{
                            \nabla^j f
                        }
                    }(y_{I}),
        \intertext{
            hence, recalling that 
            $y_I = v_{I \cap S} + x_{I \cup S}$
            for all $I \subseteq [k]$,
        }
        &\le 
            \wave{C}'_{n,k}
            \sum_{ j = 1 }^{ k }
                \polynum{ j }( \bx )
                \sum_{ I = \emptyset }^{ [k] }
                    \ballAv{ 
                        \polyset{ \set*{k+1} }( \bx ) 
                    }{
                        \restMaxFun[k+j-1]{R(\bx)}{
                            \nabla^j f
                        }
                    }\del*{ v_{I \cap S} + x_{I \cup S} },
        \intertext{
            which, by \Cref{lem::centered_max_fun_bounds_shifted_average} 
            as well as
            the facts that 
            $ 
                \norm*{ v_{I \cap S} } 
                = 
                \polyset[I\cap S]{\set*{k+1}}(\bx)
                \le 
                \polyset{ \set*{k+1} }(\bx) 
                \le 
                \ell(\bx)
            $ 
            and
            $   
                2\ell(\bx) \le R(\bx)
            $,
        }
        &\le 
            \wave{C}'_{n,k}
            \sum_{ j = 1 }^{ k }
                \polynum{ j }( \bx )
                \sum_{ I = \emptyset }^{ [k] }
                    2^n
                    \restMaxFun[k+j]{R(\bx)}{
                        \nabla^j f
                    }( x_{I \cup S} ),
        \intertext{
            which, by 
            \Cref{lem::sum_over_x_IuS_gives_2^|S|_before_the_sum},
        }
        &= 
            \wave{C}'_{n,k}
            \sum_{ j = 1 }^{ k }
                \polynum{ j }( \bx )
                2^{ \card{S} }
                \sum_{ I = S }^{ [k] }
                    2^n
                    \restMaxFun[k+j]{R(\bx)}{
                        \nabla^j f
                    }( x_{I} ),
        \intertext{
            which, since $S \subseteq [k]$,
        }
        &\le 
            \wave{C}'_{n,k}
            \sum_{ j = 1 }^{ k }
                \polynum{ j }( \bx )
                2^k
                \sum_{ I = \emptyset }^{ [k] }
                    2^n
                    \restMaxFun[k+j]{R(\bx)}{
                        \nabla^j f
                    }( x_{I} ),
        \\
        &\le 
            \wave{C}'_{n,k}
            \sum_{ j = 1 }^{ k }
                \polynum{ j }( \bx )
                2^k
                \sum_{ I = \emptyset }^{ [k+1] }
                    2^n
                    \restMaxFun[k+j]{R(\bx)}{
                        \nabla^j f
                    }( x_{I} ).
        \end{align*}
        Therefore, setting 
        $
            \wave{C}_{n,k}
            \coloneq
            2^{n+k}
            \wave{C}'_{n,k},
        $
        we have
        \begin{equation}
        \label{thmeq::thm::multipointwise_bound_for_W^k_loc::eq::sublemma_4::estimate_with_average_outside_abs}
            \ballAv{ \polyset{ \set*{k+1}}( \bx ) } {
                \abs*{
                    \diff_{ I = \emptyset }^{ [k] }
                        f\del*{ \blank + y_I - y_{\emptyset} }
                }
            }(y_{\emptyset})
            \le 
            \wave{C}_{n,k}
            \sum_{ j = 1 }^{ k }
                \polynum{ j }( \bx )
                \sum_{ I = \emptyset }^{ [k+1] }
                    \restMaxFun[k+j]{R(\bx)}{
                        \nabla^j f
                    }( x_{I } ).
        \end{equation}

        Next, notice that,
        using the linearity of 
        $\mathsf{B}_{ \polyset{ \set*{k+1}}(\bx) } $,
        we have
        \begin{align*}
            \ballAv{ 
                \polyset{ \set*{k+1} }( \bx ) 
            }{
                \abs*{
                    \diff_{ I = \emptyset }^{ [k] }
                        f\del*{ \blank + y_I - y_{\emptyset} }
                }
            }(y_{\emptyset})
            &\ge 
            \abs*{
                \ballAv{ 
                    \polyset{ \set*{k+1} }( \bx ) 
                }{
                    \diff_{ I = \emptyset }^{ [k] }
                        f\del*{ \blank + y_I - y_{\emptyset} }
                    }(y_{\emptyset})
            }
        \\  
            &=
            \abs*{
                \diff_{ I = \emptyset }^{ [k] }
                    \ballAv{ 
                        \polyset{ \set*{k+1} }( \bx ) 
                    }{
                        f\del*{ \blank + y_I - y_{\emptyset} }
                    }(y_{\emptyset})
            },
        \intertext{
            which, by \Cref{lem::shifting_operators},
        }
            &=
            \abs*{
                \diff_{ I = \emptyset }^{ [k] }
                    \ballAv{ 
                        \polyset{ \set*{k+1} }( \bx ) 
                    }{
                        f
                    }(y_I)
            },
        \intertext{
            which, by \Cref{lem::difference_to_S_and_difference_from_S_into_single_one},
        }
            &=
            \abs*{
                \diff_{  I = \emptyset }^{ S }
                    \diff_{ J = S }^{ [k] }
                        \ballAv{ 
                            \polyset{ \set*{k+1} }( \bx) 
                        }{
                            f
                        }( x_{I,J} )
            }.
        \end{align*}
        Therefore, 
        \begin{equation*}
            \abs*{
                \diff_{  I = \emptyset }^{ S }
                    \diff_{ J = S }^{ [k] }
                        \ballAv{ 
                            \polyset{\set*{k+1} } ( \bx) 
                        }{
                            f
                        }( x_{I,J} )
            }
            \le 
            \ballAv{ 
                \polyset{ \set*{k+1} }( \bx ) 
            }{
                \abs*{
                    \diff_{ I = \emptyset }^{ [k] }
                        f( \cdot + y_I - y_{\emptyset} )
                }
            }(y_{\emptyset}),
        \end{equation*}
        which combined with \eqref{thmeq::thm::multipointwise_bound_for_W^k_loc::eq::sublemma_4::estimate_with_average_outside_abs}, gives
        \begin{equation*}
            \abs*{
                \diff_{  I = \emptyset }^{ S }
                    \diff_{ J = S }^{ [k] }
                        \ballAv{ 
                            \polyset{\set*{k+1} } ( \bx) 
                        }{
                            f
                        }( x_{I,J} )
            }
            \le 
            \wave{C}_{n,k}
            \sum_{ j = 1 }^{ k }
                \polynum{ j }( \bx )
                \sum_{ I = \emptyset }^{ [k+1] }
                    \restMaxFun[k+j]{R(\bx)}{
                        \nabla^j f
                    }( x_{I } ),
        \end{equation*}
        as claimed.
    \end{proof}

    Let us move to the proof of \eqref{thmeq::thm::multipointwise_bound_for_W^k_loc::eq::with_s-k_gradient} for all $s \in \intoc{k, k+1}$. 
    Suppose
    \Cref{asu::thm::multipointwise_bound_for_W^k_loc::inductive_assumptions}
    are true.
    Fix~%
    $
        \bx = \set*{ x_I }_{ I \subseteq [k+1] }
        \subseteq
        \lebesguePoints{f}.
    $
    
    Let us first suppose that 
    $ \polyset{ \set*{k+1} }( \bx ) = 0$. 
    Then for all 
    $I \subseteq [k]$ 
    we have 
    $x_I = x_{ I \cup \set*{k+1} }$, 
    hence 
    $f(x_I) = f\del*{x_{ I \cup \set*{k+1} } }$.
    Since 
    $\bx \subseteq \lebesguePoints{f}$, 
    we have 
    $\abs{ f(x_I) } < \infty $ 
    for all 
    $I \subseteq [k+1]$.
    Thus,
    $ f\del*{x_{ I \cup \set*{k+1} } } - f(x_I) = 0 $ for all 
    $I \subseteq [k]$.
    Therefore, by \Cref{lem::diff_B_diff_C_to_diff_BuC},
    \begin{equation*}
        \abs*{
            \diff_{ I = \emptyset }^{ [k+1] }
                f(x_I)
        }
        =
        \abs*{
            \diff_{ I = \emptyset }^{ [k] }
                \del*{
                    f\del*{x_{ I \cup \set*{k+1} } }
                    - f(x_I)
                }
        }
        =
        \abs*{
            \diff_{ I = \emptyset }^{ [k] }
                0
        }
        =
        0.
    \end{equation*}
    Since the right-hand side of \eqref{thmeq::thm::multipointwise_bound_for_W^k_loc::eq::with_s-k_gradient}
    is nonnegative, we see that \eqref{thmeq::thm::multipointwise_bound_for_W^k_loc::eq::with_s-k_gradient}
    is satisfied in this case.

    Next, let us assume that 
    $ \polyset{ \set*{k+1} }( \bx ) > 0$. 
    Denote 
    $r \coloneq \polyset{ \set*{k+1} }( \bx )$
    and 
    $v_I \coloneq x_{ I \cup \set*{k+1} } - x_I$
    for~all~$I \subseteq [k]$.
    Then, for all $I, J \subseteq [k] $, define
    $ x_{I,J} \coloneq v_I + x_J$. 
    Note that we have
    $ x_{I,I} = x_{ I \cup \set*{k+1} }$
    for all $I \subseteq [k]$.

    We will now split the expression at the left-hand side of \eqref{thmeq::thm::multipointwise_bound_for_W^k_loc::eq::with_s-k_gradient}
    into several terms.
    \begin{align*}
        \abs*{
            \diff_{ I = \emptyset }^{ [k+1] }
                f(x_I)
        }
        &\le
        \underbrace{
            \abs*{
                \diff_{ I = \emptyset }^{ [k] }
                    \del*{
                        f(x_I)
                        - \ballAv{r}{f}(x_I)
                    }
            }
        }_{ \eqcolon D_1 }
        +
        \underbrace{
            \abs*{
                \diff_{ I = \emptyset }^{ [k] }
                    \del*{
                        f\del*{
                        	x_{ I \cup \set*{k+1} }
                        }
                        - 
                        \ballAv{r}{f}\del*{
                        	x_{ I \cup \set*{k+1} }
                        }
                    }
            }
        }_{ \eqcolon D_2 }
    \\
    &\quad \,\,\, + 
        \abs*{
            \diff_{ I = \emptyset}^{ [k+1] }
                \ballAv{r}{f}(x_I)
        }
        =
        D_1 + D_2 
        + 
        \abs*{
            \diff_{ I = \emptyset}^{ [k] }
                \ballAv{r}{f}(x_I)
            -
            \diff_{ I = \emptyset}^{ [k] }
                \ballAv{r}{f}(x_{I,I})
        },
    \intertext{
        which, using \Cref{lem::diag_diff_into_row_and_col_diffs},
    }
        &=
        D_1 + D_2 
        + 
        \abs*{
            \diff_{ I = \emptyset}^{ [k] }
                \ballAv{r}{f}(x_I)
            -
            \sum_{
                S \subseteq [k]
            }
                \diff_{ I = \emptyset }^{ S }
                    \diff_{ J = S }^{ [k] }
                        \ballAv{r}{f}( x_{I, J} )
        }
        \\
        &\le
        D_1 + D_2 
        + 
        \abs*{
            \diff_{ I = \emptyset}^{ [k] }
                \ballAv{r}{f}(x_I)
            -
            \diff_{ I = \emptyset }^{ \emptyset }
                \diff_{ J = \emptyset }^{ [k]}
                    \ballAv{r}{f}( x_{I, J} )
        }
        +
        \abs*{
            \sum_{ 
                \substack{ 
                    S \subseteq [k] \\ S \ne \emptyset 
                } 
            }
                \diff_{ I = \emptyset }^{ S }
                    \diff_{ J = S }^{ [k] }
                        \ballAv{r}{f}( x_{I, J} )
        }
        \\
        &\le
        D_1 + D_2 
        + 
        \underbrace{
            \abs*{
                \diff_{ I = \emptyset}^{ [k] }
                    \ballAv{r}{f}(x_I)
                -
                \diff_{ J = \emptyset }^{ [k]}
                    \ballAv{r}{f}( x_{\emptyset, J} )
            }
        }_{ \eqqcolon D_3 }
        +
        \sum_{ 
            \substack{ 
                S \subseteq [k] \\ S \ne \emptyset 
            } 
        }
        \underbrace{
            \abs*{
                \diff_{ I = \emptyset }^{ S }
                    \diff_{ J = S }^{ [k]  }
                        \ballAv{r}{f}( x_{I, J} )
            }
        }_{ \eqqcolon D_{4,S} }
        \\
        &=
        D_1 + D_2 + D_3
        +
        \sum_{ 
            \substack{ 
                S \subseteq [k] \\ S \ne \emptyset 
            } 
        }
            D_{4,S}.
    \end{align*}
    Thus,
    \begin{equation}
    \label{thmeq::thm::multipointwise_bound_for_W^k_loc::eq::estimates_for_diff_with_Di's_and_D_4S's}
        \abs*{
            \diff_{ I = \emptyset }^{ [k+1] }
                f(x_I)
        }
        \le 
        D_1 + D_2 + D_3
        +
        \sum_{ 
            \substack{ 
                S \subseteq [k] \\ S \ne \emptyset 
            } 
        }
            D_{4,S}.
    \end{equation}

    Let us now estimate 
    $D_i$, where $i \in [3]$, and $D_{4,S}$, where 
    $ S\subseteq [k]$
    is nonempty,
    keeping in~mind that $r = \polyset{ \set*{k+1} }( \bx) > 0$.
    \begin{enumerate}
        \item[$\textbf{D}_{\bm{1}}$:]
        Here we have
        \begin{equation*}
            D_1
            =
            \abs*{
                \diff_{ I = \emptyset }^{ [k] }
                    \del*{
                        f(x_I)
                        - \ballAv{r}{f}(x_I)
                    }
            }
            =
            \abs*{
                \diff_{ I = \emptyset }^{ [k] }
                    \del*{
                        \ballAv{ 
                            \polyset{ \set*{k+1} }( \bx) 
                        }{f}(x_I)
                        - f(x_I)
                    }
            },
        \end{equation*}
        so, by \Cref{sublem::thm::multipointwise_bound_for_W^k_loc::eq::sublemma_2},
        \begin{equation}
        \label{thmeq::thm::multipointwise_bound_for_W^k_loc::eq::estimating_D1}
            D_1
            \le 
            C_{n,k}''
            \del*{
                \sum_{ j = 2 }^{k}
                    \polynum{ j }( \bx)
                \smashoperator{
                    \sum_{ I = \emptyset }^{ [k + 1] }
                }
                        \restMaxFun[k+j]{R( \bx) }{ 
                            \nabla^j f
                        }( x_I )
                +
                \polynum{ s }( \bx )
            \smashoperator{
                \sum_{ I = \emptyset }^{ [k + 1] }
            }
                    \restMaxFun[2k]{R( \bx) }{ g_{\ell(\bx)} }( x_I )
            }.
        \end{equation}
        \item[$\textbf{D}_{\bm{2}}$:]
        Define 
        $
            \wave{\bx}
            =
            \set*{ \wave{x}_I }_{ I \subseteq [k+1] }
        $
        by the formula 
        $
            \wave{x}_I
            \coloneq
            x_{I \div \set*{k+1}}
        $
        for $I \subseteq [k+1]$.
        Then by \Cref{lem::invariant_transformations_of_poly} we have
        $
            \polyset{ \set*{k+1} }( \wave{\bx} )
            =
            \polyset{ \set*{k+1} }( \bx) > 0.
        $
        Therefore, since
        \begin{align*}
            D_2
        &=
            \abs*{
                \diff_{ I = \emptyset }^{ [k] }
                    \del*{
                        f\del*{x_{ I \cup \set*{k+1} } }
                        - \ballAv{r}{f}\del*{
                        	x_{ I \cup \set*{k+1} } 
                        }
                    }
            }
        \\
        &=
            \abs*{
                \diff_{ I = \emptyset }^{ [k] }
                    \del*{
                        \ballAv{ 
                            \polyset{ \set*{k+1} }( \bx) 
                        }{f}\del*{
                        	x_{ I \div \set*{k+1} } 
                        }
                        - f\del*{ x_{ I \div \set*{k+1} }}
                    }
            }
            =
            \abs*{
                \diff_{ I = \emptyset }^{ [k] }
                    \del*{
                        \ballAv{ 
                            \polyset{ \set*{k+1} }\del*{ \wave{ \bx }} 
                        }{f}\del*{ \wave{x}_I  }
                        - f\del*{ \wave{x}_I  }
                    }
            },
        \end{align*}
        by 
        \Cref{sublem::thm::multipointwise_bound_for_W^k_loc::eq::sublemma_2}
        we have
        \begin{align*}
            D_2
            \hspace{-18pt}&\hspace{18pt}\le 
            C_{n,k}''
            \del*{
                \sum_{ j = 2 }^{k}
                    \polynum{ j }( \wave{\bx})
                    \sum_{ I = \emptyset }^{ [k + 1] }
                        \restMaxFun[k+j]{R( \wave{\bx} ) }{ 
                            \nabla^j f
                        }( \wave{x}_I )
                +
                \polynum{ s }( \wave{\bx} )
                \sum_{ I = \emptyset }^{ [k + 1] }
                    \restMaxFun[2k]{R( \wave{\bx} ) }{ 
                    	g_{\ell\del*{ \wave{ \bx }} } 
                    }( \wave{x}_I )
            },
        \intertext{
            which, by \Cref{lem::invariant_transformations_of_poly} 
            and the fact that
            $
                R( \wave{\bx} )
                =
                2^{(k+1)^2} \ell( \wave{\bx} )
                =
                2^{(k+1)^2} \ell( \bx )
                =
                R( \bx)
            $,
        }
            &\le 
            C_{n,k}''
            \del*{
                \sum_{ j = 2 }^{k}
                    \polynum{ j }( \bx )
                    \sum_{ I = \emptyset }^{ [k + 1] }
                        \restMaxFun[k+j]{R( \bx) }{ 
                            \nabla^j f
                        }\del*{ \wave{x}_I }
                +
                \polynum{ s }( \bx  )
                \sum_{ I = \emptyset }^{ [k + 1] }
                    \restMaxFun[2k]{R( \bx) }{ g_{\ell(\bx)} }\del*{ \wave{x}_I }
            },
        \intertext{
            which, since 
            $\wave{x}_I = x_{I \div \set*{k+1} }$ 
            for all 
            $I \subseteq [k+1]$,
        }
            &=
            C_{n,k}''
            \del*{
                \sum_{ j = 2 }^{k}
                    \polynum{ j }( \bx )
                    \smashoperator{
                    \sum_{ I = \emptyset }^{ [k + 1] }
                    }
                        \restMaxFun[k+j]{R( \bx) }{ 
                            \nabla^j f
                        }\del*{ x_{I \div \set*{k+1} } }
            {+}
                \polynum{ s }( \bx  )
                \smashoperator{
                \sum_{ I = \emptyset }^{ [k + 1] }
                }
                    \restMaxFun[2k]{R( \bx) }{ g_{\ell(\bx)} }\del*{ x_{I \div \set*{k+1} } }
            },
        \intertext{
            which, since 
            $I \mapsto I \div \set*{k+1}$ 
            is a bijection from $2^{[k+1]}$ to itself,
        }
            &= 
            C_{n,k}''
            \del*{
                \sum_{ j = 2 }^{k}
                    \polynum{ j }( \bx )
                    \sum_{ I = \emptyset }^{ [k + 1] }
                        \restMaxFun[k+j]{R( \bx) }{ 
                            \nabla^j f
                        }( x_{I } )
                +
                \polynum{ s }( \bx  )
                \sum_{ I = \emptyset }^{ [k + 1] }
                    \restMaxFun[2k]{R( \bx) }{ g_{\ell(\bx)} }( x_{I } )
            }.
        \end{align*}
        Thus, we have
        \begin{equation}
        \label{thmeq::thm::multipointwise_bound_for_W^k_loc::eq::estimating_D2}
            D_2
            \le 
            C_{n,k}''
            \del*{
                \sum_{ j = 2 }^{k}
                    \polynum{ j }( \bx)
                \smashoperator{
                    \sum_{ I = \emptyset }^{ [k + 1] }
                }
                        \restMaxFun[k+j]{R( \bx) }{ 
                            \nabla^j f
                        }( x_I )
                +
                \polynum{ s }( \bx )
            \smashoperator{
                \sum_{ I = \emptyset }^{ [k + 1] }
            }
                    \restMaxFun[2k]{R( \bx) }{ g_{\ell(\bx)} }( x_I )
            }.
        \end{equation}
        \item[$\textbf{D}_{\bm{3}}$:]
        Here, since 
        $
            x_{\emptyset, I} = v_{\emptyset} + x_I = 
            x_I + x_{ \set*{k+1} } - x_{\emptyset}
        $
        for all $I \subseteq [k]$, we have
        \begin{align*}
            D_3
        &=
            \abs*{
                \diff_{ I = \emptyset}^{ [k] }
                    \ballAv{r}{f}(x_I)
                -
                \diff_{ I = \emptyset }^{ [k]}
                    \ballAv{r}{f}( x_{\emptyset, I} )
            }
        \\
        &=
            \abs*{ 
                \diff_{ I = \emptyset }^{ [k] }
                    \del*{
                        \ballAv{
                            \polyset{\set*{k+1}}(\bx)
                        }{
                            f
                        }(x_I + x_{\set*{k+1}} - x_{\emptyset})
                        -
                        \ballAv{
                            \polyset{\set*{k+1}}(\bx)
                        }{
                            f
                        }(x_I)                      
                    }
            },
        \end{align*}
        hence, by \Cref{sublem::thm::multipointwise_bound_for_W^k_loc::eq::sublemma_3},
        \begin{equation}
        \label{thmeq::thm::multipointwise_bound_for_W^k_loc::eq::estimating_D3}
            D_3
            \le 
            C_{n,k}'''
            \del*{
                \sum_{ j = 2 }^{k}
                    \polynum{ j }( \bx)
                \smashoperator{
                    \sum_{ I = \emptyset }^{ [k + 1] }
                }
                        \restMaxFun[k+j]{ R( \bx) }{ 
                            \nabla^j f
                        }( x_I )
                +
                \polynum{ s }( \bx )
            \smashoperator{
                \sum_{ I = \emptyset }^{ [k + 1] }
            }
                    \restMaxFun[2k]{ R( \bx) }{ g_{\ell(\bx)} }( x_I )
            }.
        \end{equation}
        \item[$\textbf{D}_{\bm{4,}\textbf{S}}$:]
        Here, since $S \subseteq [k]$ is nonempty, 
        by \Cref{sublem::thm::multipointwise_bound_for_W^k_loc::eq::sublemma_4} we have
        \begin{equation}
        \label{thmeq::thm::multipointwise_bound_for_W^k_loc::eq::estimating_D4}
            D_{4,S}
            =
            \abs*{ 
                \diff_{  I = \emptyset }^{ S }
                    \diff_{ J = S }^{ [k] }
                        \ballAv{ 
                            \polyset{\set*{k+1} }( \bx) 
                        }{
                            f
                        }( x_{I,J} )
            }
            \le 
            \wave{C}_{n,k}
            \sum_{ j = 1 }^{k}
                \polynum{ j }( \bx)
                \sum_{ I = \emptyset }^{ [k+1] }
                    \restMaxFun[k+j]{R( \bx) }{ 
                        \nabla^j f
                    }(x_{I}).
        \end{equation}
    \end{enumerate}

    Therefore, using the estimates from \eqref{thmeq::thm::multipointwise_bound_for_W^k_loc::eq::estimating_D1},
    \eqref{thmeq::thm::multipointwise_bound_for_W^k_loc::eq::estimating_D2},
    \eqref{thmeq::thm::multipointwise_bound_for_W^k_loc::eq::estimating_D3},
    and 
    \eqref{thmeq::thm::multipointwise_bound_for_W^k_loc::eq::estimating_D4}
    in  \eqref{thmeq::thm::multipointwise_bound_for_W^k_loc::eq::estimates_for_diff_with_Di's_and_D_4S's}, we have
    \begin{align*}
        \abs*{
            \diff_{ I = \emptyset }^{ [k+1] }
                f(x_I)
        }
    \hspace{-2cm}&\hspace{2cm}\le
        D_1 + D_2 + D_3
        +
        \sum_{ 
            \substack{ 
                S \subseteq [k] \\ S \ne \emptyset 
            } 
        }
            D_{4,S}
    \\
    &\le 
        \del*{2C_{n,k}'' + C_{n,k}'''}
        \del*{
                \sum_{ j = 2 }^{k}
                    \polynum{ j }( \bx)
                    \sum_{ I = \emptyset }^{ [k + 1] }
                        \restMaxFun[k+j]{ R( \bx) }{ 
                            \nabla^j f
                        }( x_I )
                +
                \polynum{ s }( \bx )
                \sum_{ I = \emptyset }^{ [k + 1] }
                    \restMaxFun[2k]{ R( \bx) }{ g_{\ell(\bx)} }( x_I )
            }
    \\
    &\hspace{1.15cm} +
        \sum_{ 
            \substack{ 
                S \subseteq [k] \\ S \ne \emptyset 
            } 
        }
        \wave{C}_{n,k}
            \sum_{ j = 1 }^{k}
                \polynum{ j }( \bx)
                \sum_{ I = \emptyset }^{ [k+1] }
                    \restMaxFun[k+j]{R( \bx) }{ 
                        \nabla^j f
                    }(x_{I}),
    \intertext{
        which, since there are $2^k-1$ subsets $S$ of $[k]$ such that $S \ne \emptyset$,
    }
    &\le 
    \del*{
    	2C_{n,k}'' 
   	{+} 
   		C_{n,k}''' 
   	{+} 
   		2^k \wave{C}_{n,k}
    }
        \!
        \del*{
        \!
                \smashoperator{
                \sum_{ j = 1 }^{k}
                }
                    \polynum{ j }( \bx)
                    \smashoperator{
                    \sum_{ I = \emptyset }^{ [k + 1] }
                    }
                    \!
                        \restMaxFun[k+j]{ R( \bx) }{ 
                            \nabla^j f
                        }( x_I )
            +
                \polynum{ s }( \bx )
                \smashoperator{
                \sum_{ I = \emptyset }^{ [k + 1] }
                }
                \!
                    \restMaxFun[2k]{ R( \bx) }{ g_{\ell(\bx)} }( x_I )
        \!    
        }.
    \end{align*}
    Therefore, setting 
    $
        \widehat{C}_{n,k+1}
        \coloneq 
        2C_{n,k}'' + C_{n,k}''' + 2^k \wave{C}_{n,k},
    $
    we get
    \begin{equation*}
        \abs*{
            \diff_{ I = \emptyset }^{ [k+1] }
                f(x_I)
        }
        \le 
        \widehat{C}_{n,k+1}
        \del*{
            \sum_{ j = 1 }^{ k }
                \polynum{ j }( \bx )
                \sum_{ I = \emptyset }^{ [k+1] }
                    \restMaxFun[k+j]{R(\bx)}{\nabla^j f}(x_I)
            + 
            \polynum{ s }(\bx)
            \sum_{ I = \emptyset }^{ [k+1] }
                    \restMaxFun[2k]{R(\bx)}{g_{\ell(\bx)} }(x_I)
        },
    \end{equation*}
    which proves \eqref{thmeq::thm::multipointwise_bound_for_W^k_loc::eq::with_s-k_gradient} for $s \in \intoc{k,k+1}$.

    Finally, let us prove \eqref{thmeq::thm::multipointwise_bound_for_W^k_loc::eq::just_loc} for $k+1$
    in order to finish the inductive step.
    To that end, suppose that 
    $f \in \locSobolev{k+1, 1}( \bR^n)$.
    Note that by 
    \Cref{cor:bojarski_inequality_vector_version},
    the function
    $
        (R,x) 
        \mapsto 
        C_M n^{k}
        \restMaxFun{2R}{ \nabla^{k+1} f }(x)
    $
    is an element of 
    $D^{1}_{\lambda, \mathrm{res}}\del*{ \nabla^k f }$,
    so 
    by \eqref{thmeq::thm::multipointwise_bound_for_W^k_loc::eq::with_s-k_gradient} for $s = k+1$,
    we have 
    \begin{align*}
        \abs*{
            \diff_{ I = \emptyset }^{ [k+1] }
                f(x_I)
        }
    \hspace{-2cm}&
    \\
        &\le 
        \widehat{C}_{n,k+1}
        \del*{ \!
        	\smashoperator{
            \sum_{ j = 1 }^{ k }
            }
                \polynum{ j }( \bx ) \!
                \smashoperator{
                \sum_{ I = \emptyset }^{ [k+1] }
                }
                \!
                    \restMaxFun[k+j]{R(\bx)}{\nabla^j f}(x_I)
        {+} 
            \polynum{ k+1 }(\bx) \!
            \smashoperator{
            \sum_{ I = \emptyset }^{ [k+1] }
            }
            \!
                    \restMaxFun[2k]{R(\bx)}{
                        C_M n^{k} \restMaxFun{2 \ell(\bx)}{ \nabla^{k+1} f }
                    }(x_I)
        \! \! }
    \\  
        &\le 
        \max\del*{ 1, C_M n^{k}}
        \widehat{C}_{n,k+1}
        \del*{ \!
        	\smashoperator{
            \sum_{ j = 1 }^{ k }
            }
                \polynum{ j }( \bx ) \!
                \smashoperator{
                \sum_{ I = \emptyset }^{ [k+1] }
                }
                \!
                    \restMaxFun[k+j]{R(\bx)}{\nabla^j f}(x_I)
        {+} 
            \polynum{ k+1 }(\bx) \!
            \smashoperator{
            \sum_{ I = \emptyset }^{ [k+1] }
            }
            \!
                    \restMaxFun[2k+1]{R(\bx)}{ \nabla^{k+1} f }(x_I)
        \! \! }
    \\
        &=
        C_{n, k+1}
        \sum_{ j = 1 }^{ k + 1 }
            \polynum{ j }( \bx )
            \sum_{ I = \emptyset }^{ [k+1] }
                \restMaxFun[k+j]{R(\bx)}{\nabla^j f}(x_I),
    \end{align*}
    where we set 
    $
        C_{n, k+1}
        \coloneq
        \max\del*{ 1, C_M n^{k} }
        \widehat{C}_{n,k+1}.
    $
    Thus, we proved \eqref{thmeq::thm::multipointwise_bound_for_W^k_loc::eq::just_loc} for $k+1$ and the~theorem follows by~induction.
\end{proof}
\begin{definition}[Families $\bD_{\measureAlone}^s(f)$ and $\mF_{\measureAlone}^s(f, G)$]
\label{def::s-hypergradient_normed}
    Let 
    $\del{V, \normAlone[V]}$,
    $\del{W, \normAlone[W]}$
    be normed spaces and $X \subseteq V$. Let $s \ge 0$ and $k \in \bN_0$ be such that $s \in \intoc{k-1, k}$.
    Let $\measureAlone$ be a measure on $X$. 
    For a given measurable 
    function $f \colon X \to W$,
	we will denote by 
    $ \bD_{\measureAlone}^s(f)$ the family of measurable functions
    $G \colon X \to \intcc{0,\infty}$ 
    such that
    \begin{equation} \label{eqdef::def::s-hypergradient_normed::m-upper}
        \measureAlone
        \forall 
            \bx = \set*{x_{I}}_{ I \subseteq [k] }
            \subseteq X
        \qquad 
        \norm*{
            \diff_{I = \emptyset }^{ [k] }
                f(x_I)
        }_{W}
        \le 
        \polygen{ \gen{s} }(\bx)
        \sum_{ I = \emptyset }^{ [k] }
            G(x_I).
    \end{equation}
    For a given $G \in \bD_{ \measureAlone }^{ s }(f)$, 
    we will denote by 
    $\mF_{\measureAlone}^{s}(f,G)$
    the family of all subsets 
    of $X$ of
    full measure on~which the~inequality in \eqref{eqdef::def::s-hypergradient_normed::m-upper} is satisfied.
    Also, 
    when 
    $ \measureAlone $
    is the counting measure $\#$,
    we~might write 
    $ \bD^s(f) $
    instead of 
    $ \bD^s_{\#}(f) $.
\end{definition}
%
%
\begin{remark}
    ${}$
    \begin{enumerate}[listparindent=\parindent]
    \item
        Since the only set with a cardinality of $0$ is the empty set,
    	if $\measureAlone = \#$
    	and 
    	$G \in \bD^s(f) = \bD^s_{\measureAlone}(f)$,
        then the inequality in
        \eqref{eqdef::def::s-hypergradient_normed::m-upper}
        has to be satisfied everywhere
        and
        $ \mF_{\#}^s(f,G) = \set*{ X } $.
    \item 
        When $s = 0$, then 
        $\polygen{ \gen{s} }(\bx) = 1$
        and the condition present in 
        \eqref{eqdef::def::s-hypergradient_normed::m-upper}
        is equivalent to stating that 
        we have
        $\norm{f}_W \le G$
        $\measureAlone$-almost everywhere.
    \item 
        When $s \in \intoc{0,1}$, the notions of 
        Hajłasz $s$-gradients from
        \Cref{def::Hajlasz_s-gradients}
        and the elements of
        $ \bD_{\measureAlone}^s(f) $ 
        in~the~sense of 
        \Cref{def::s-hypergradient_normed} 
        coincide exactly.
        Indeed, this follows from the fact that 
        for all such $s$, we~have
        $
        	\diff_{I=\emptyset}^{[k]}
        		f(x_I)
        	=
        	f\del*{ x_{\set*{1}}}
        	- f\del*{ x_{\emptyset}}
        $
        and
		$
			\polygen{ \gen{s} }(\bx)
			=
			\norm{ x_{\set*{1}} - x_{\emptyset} }^s.
		$
    \item
    	Let us also note that
    	we can define families
    	$ \bD_{\measureAlone}^s(f)$
    	and
    	$
    		\mF_{\measureAlone}^{s}\del*{f, G}
    	$
    	for measurable functions
    	$ f \colon X \to \bRExt $
    	by replacing the normed space
    	$
    		\del*{
    			W, \normAlone[W]
    		}
    	$
    	with the extended real line 
    	$ 
    		\del*{
    			\bRExt,
    			\absAlone	
    		}
    	$
    	in 
    	\Cref{def::s-hypergradient_normed}.
    	However, if we do so, then the increase in generality is effectively negligible.
    	Indeed, suppose we have such an $f$
        and 
        $ s > 0$ such that
        $
        	\bD_{\measureAlone}^{s}(f)
        	\ne 
        	\emptyset.
        $
        Then $f$ is finite almost everywhere, 
        hence it has an everywhere finite representative.
        As such, with only minimal loss of generality, 
        we can restrict ourselves to only functions
        $ f \colon X \to \bR $.
        
        Indeed, 
        let us fix
        $ G \in \bD_{\measureAlone}^{s}(f)$
        and
        $ F \in \mF_{\measureAlone}^{s}(f, G)$.
        Suppose we have
        $\abs{f(x)} = \infty$ 
        for some 
        $x \in F$.
        Let $k \in \bN$ be such that 
        $ s \in \intoc{ k-1, k }$
        and define 
        $\bx = \set*{x_I}_{ I \subseteq [k] }$
        by the formula 
        $ x_I \coloneq x $
        for all 
        $ I \subseteq [k] $.
        Then, by~our~convention, 
        $
            \abs*{
                \diff_{ I = \emptyset }^{ [k] }
                    f(x_I)
            }
            =
            \infty
        $
        since 
        $
        	\abs*{
        		f\del{ x_\emptyset }
        	}
        	=
        	\abs*{ f(x) }
        	=
        	\infty.
        $
        Also, by
        \Cref{cor::value_of_Poly_when_tuple_is_a_hyperparallelogram},
        we have
        $
            \polygen{ \gen{s} }(\bx) = 0,
        $
        so by the convention $0 \t \infty = 0$ we have
        $
            \polygen{ \gen{s} }(\bx)
            \sum_{ I = \emptyset }^{ [k] }
                G(x_I)
            =
            0.
        $
        In~consequence, $x \notin F$, contradicting our assumption.
        Thus,
        $ f \rvert_{F} $
        is finite everywhere,
        hence
        $ f $ is finite almost everywhere.
    \end{enumerate}
\end{remark}
\begin{corollary}
\label{cor::hypergradients_from_multipointwise_bound_for_W^k_loc}
    Let $n, k \in \bN$ and $s \in \intoc{k, k+1}$. Let $C_{n,k}$ and $\widehat{C}_{n,k+1}$ be the constants from \Cref{thm::multipointwise_bound_for_W^k_loc}. 
    Let $f \in \locSobolev{k, 1}( \bR^n)$ and $g \in \bD^{s-k}_{\lambda}( \nabla^k f)$. Then
    \begin{align*}
        G_k
        &\coloneq
        C_{n,k}
        \sum_{j=1}^k 
            \maxFun[k+j-1]{ \nabla^j  f}
        \in \bD^{k}_{\lambda}(f)
    \intertext{
    	and
    }
        \widehat{G}_s
        &\coloneq
        \widehat{C}_{n,k+1}
        \del*{
            \sum_{j=1}^k 
                \maxFun[k+j]{ \nabla^j  f}
            +
                \maxFun[2k]{g}
        }
        \in \bD^{s}_{\lambda}(f).
    \end{align*}
    Moreover, 
    $
        \lebesguePoints{f} 
        \in 
        \mF^k_{\lambda}(f, G_k)
        \cap 
        \mF^s_{\lambda}\del*{
        	f, \widehat{G}_s
        }.
    $
\end{corollary}
\begin{proof}
    Both statements follow immediately from \Cref{thm::multipointwise_bound_for_W^k_loc}
    and the fact that 
    $\restMaxFun{R}{h} \le \maxFun{h}$ everywhere for all $R \ge 0$ and every
    measurable function $h$.
\end{proof}
\begin{definition}[H\"{o}lder space $\holder{m,\alpha}\del*{ \Omega } $]
\label{def::higher_order_holder_spaces}
	Let $n \in \bN$, $m \in \bN_0$,
	and $\alpha \in \intoc{0,1}$.
	Let 
	$ \Omega \subseteq \bR^n$ be open.
	We define the family 
	\begin{equation*}
		\holder{m,\alpha}\del*{ \Omega } 
		\coloneq 
		\setc*{
			f \in \holder{m}\del*{ \Omega} 
		}{
			\text{
			$ 
				\nabla^{j} f 
				\in 
				\continuousAndBounded\del*{ \Omega }
			$ 
			for all $j \in [m] $
			and }
			\seminorm*{
				\nabla^m f
			}_{ 
				\holder{0,\alpha}\del*{ \Omega }
			}
			<
			\infty
		},
	\end{equation*}
	where
	\begin{equation*}
		\seminorm*{ \nabla^m f }_{ 
			\holder{0,\alpha}\del*{ \Omega } 
		}		
		\coloneq 
		\sup_{
			\substack{
				x, y \in \Omega 
			\\
				x \ne y
			}
		}
		\frac{ 
			\norm{ 
				\nabla^m f (x) - \nabla^m f(y)
			}
		}{
			\norm{x-y}^{\alpha}
		}.
	\end{equation*}
	We endow this family with the seminorm
	$
		\seminormAlone[{
			\holder{m,\alpha}\del*{ \Omega }
		}]
	$
	that is defined for 
	$
		f \in \holder{m}\del*{ \Omega }
	$
	as follows:
	\begin{equation*}
		\seminorm*{ f }_{ 
			\holder{m,\alpha}\del*{ \Omega } 
		}
		\coloneq 
		\sum_{ j = 1 }^{ m }
			\norm*{ \nabla^j f }_{
				\infty
			}
		+
		\seminorm*{ \nabla^m f }_{ 
			\holder{0,\alpha}\del*{ \Omega } 
		}.
	\end{equation*}
\end{definition}
\begin{corollary}
\label{cor::multipointwise_bound_for_higher_order_holder}
    Let $n \in \bN$, $m \in \bN_0$, $\alpha \in \intoc{0,1}$, and 
    $
        f \in \holder{m, \alpha}\del*{ \bR^n }.
    $
    Then
    we have that
    $
    	{
        \widehat{C}_{n,m+1} 
        \seminorm{ f }_{ 
        	\holder{m, \alpha}\del*{
	        	\bR^n 
	        } 
	    }
        \in \bD^{ m + \alpha }( f ),
        }
    $
    where 
    $ \widehat{C}_{n,m+1}  $
    is the constant from 
    \Cref{thm::multipointwise_bound_for_W^k_loc}
    for $m > 0$ and 
    $ \widehat{C}_{n,m+1} = 1$ if~$m = 0$. 
\end{corollary}
\begin{proof}
    Let us first consider the case when $m = 0.$
    Then for all $x, y \in \bR^n$ we have the inequality
    $
        \abs*{ f(x) - f(y) }
        \le 
        \seminorm*{ f }_{ 
        	\holder{0,\alpha}\del*{
	        	\bR^n 
	        } 
        }
        \norm*{ x - y }^{ \alpha },
    $
    so 
    $
        \seminorm{ f }_{ 
        	\holder{0,\alpha}\del*{ 
        		\bR^n 
        	} 
        }
        \in 
        \bD^{ \alpha }(f),
    $
    as claimed.
    Next, assume that $m > 0$. 
    Let us note that if~%
    $
        f 
        \in 
        \holder{m, \alpha}\del*{ \bR^n },
    $
    then
    $
        f \in \locSobolev{m, \infty}\del*{ \bR^n },
    $
    so also 
    $
        f \in \locSobolev{m, 1}\del*{ \bR^n }.
    $
    Moreover, for all $x, y \in \bR^n$, we have
    \begin{equation*}
        \norm*{
            \nabla^{ m } f(x)
            - \nabla^{ m }f (y)
        }
        \le 
        \norm{x-y}^{ \alpha }
        \seminorm*{ \nabla^{m} f }_{
            \holder{0,\alpha}\del*{ \bR^n }
        }.
    \end{equation*}
    Hence,
    $
        \seminorm{ \nabla^{k} f }_{
            \holder{0,\alpha}\del*{ \bR^n } 
        }
        \in 
        \bD_{\lambda}^{ \alpha }( 
            \nabla^{m} f
        ).
    $
    Therefore, by 
    \Cref{cor::hypergradients_from_multipointwise_bound_for_W^k_loc},
    \begin{equation*}
        G
        \coloneq 
        \widehat{C}_{n,m+1}
        \del*{
            \sum_{ j = 1 }^{ m }
                \maxFun[m+j]{ \nabla^j f}
            +
            \seminorm{ \nabla^{m} f }_{
                \holder{0,\alpha}\del*{ \bR^n }
            }
        }
        \in 
        \bD_{\lambda}^{m+\alpha}(f)
    \end{equation*}
    and 
    $
        \lebesguePoints{f}
        \in 
        \mF_{\lambda}^{m+\alpha}(f, G).
    $
    Note that here we have used the fact that
    $
        \maxFun{C} = C 
    $
    for all constant functions $C \ge 0.$
    Also, since $f$ is continuous, we have
    $
        \lebesguePoints{f} = \bR^n.
    $
    Hence, in fact, 
    $
        G \in \bD^{m+\alpha}(f).
    $

    Next, note that for all $j \in [m]$ we have
    $
        \norm*{
            \nabla^j f(x)
        }
        \le 
        \norm*{
            \nabla^j f
        }_{ 
            \infty
        }
    $
    for all $x \in \bR^n$.
    In consequence, for all $j \in [m]$
    and $x \in \bR^n$,
    \begin{equation*}
        \maxFun[m+j]{
            \nabla^j f
        }(x)
        \le 
        \maxFun[m+j]{
            \norm*{
                \nabla^j f
            }_{ 
                \infty
            }
        }(x)
        =
        \norm*{
            \nabla^j f
        }_{ 
            \infty
        }.
    \end{equation*}
    In consequence, we have
    \begin{equation*}
        G
        \le 
        \widehat{C}_{n,m+1}
        \del*{
            \sum_{ j = 1 }^{ m }
                \norm*{
                    \nabla^j f
                }_{ \infty }
            +
            \seminorm{ \nabla^{m} f }_{
                \holder{0,\alpha}\del*{ \bR^n }
            }
        }
        =
        \widehat{C}_{n,m+1}
        \seminorm{ f }_{
            \holder{k,\alpha}\del*{ \bR^n }
        }
    \end{equation*}
    everywhere in $\bR^n$.
    Thus, we have
    $
        \widehat{C}_{n,m+1}
        \seminorm{ f }_{
            \holder{k,\alpha}\del*{ \bR^n }
        }
        \in 
        \bD^{m+\alpha}(f),
    $
    as claimed.
\end{proof}
So far, we have found elements of 
$ \bD_{\lambda}^s(f)$ or $ \bD^s(f)$
under the assumption that sufficiently many weak or~classical derivatives of $f$ exist.
Our next goal will be to answer a somewhat opposite question.
Let $\Omega \subseteq \bR^n$, $s > 0$
and suppose that 
$f, G \in \locIntegrable{1}(\Omega)$
are such that
$ G \in \bD_{\lambda}^s(f)$.
What can we then say about the weak derivatives of $f$?
As we will see in~%
\Cref{thm::local_hypergradients_imply_local_sobolev},
for such $f$ and $G$, we have
$\partial^{\alpha} f \in 
\locIntegrable{1}(\Omega)$
for all $\alpha \in \bN_0^n$ such that
$ 0 < \abs{ \alpha } \le s $.
Moreover, 
$ 
	\abs{ \partial^{\alpha} f }
	\le 
	2^{k^2} \partitions{k} G
$
almost everywhere, where $k \in \bN$ is such that
$ s \in \intoc{k-1,k}$.
However, before we will state and prove the mentioned theorem, we will need to obtain some auxiliary results.
\begin{proposition}
\label{prop::embedding_into_lower_orders_normed} 
    Let 
    $ \del*{ V, \normAlone[V] }$
    and
    $ \del*{ W, \normAlone[W] }$
    be normed spaces, 
    $X \subseteq V$ 
    be nonempty, and 
    $\measureAlone$ 
    be~a~measure on $X$. 
    Let $s \ge 0$  and $k \in \bN_0$ be~such that 
    $ s \in \intoc{k-1, k}$.
    Let
    $
        f \colon X \to W
    $
    and 
    $
        G \colon X \to \intcc{0, \infty}
    $
    be measurable and such that 
    $
        G \in \bD^{s}_{\measureAlone}(f) .
    $
    Then, 
    \begin{equation*}
        \forall t \in \intcc{0,s}
    \qquad 
        2^{k^2} \partitions{k} ( \norm{f}_W + G)
        \in \bD^{t}_{\measureAlone}(f)
        \quad \text{and} \quad 
        \mF_{\measureAlone}^s(f, G)
        \subseteq
        \mF_{\measureAlone}^t\del*{ 
            f, 
            2^{k^2} \partitions{k} \del*{ 
                \norm{ f }_W + G 
            } 
        }.
    \end{equation*}
\end{proposition}
\begin{proof}
    Using induction on $m \in \bN_0$ satisfying $m \le k$, 
    we will show that for all such $m$
    the statement is correct for every 
    $
    	t \ge 0
    $
    such that
    $ 
        t \in \intoc{ m -1, \min(m, s) }.
    $ 
    Let us notice that this is true for $t = 0$.
    Hence, the base case of induction, $m = 0$, is proved.
    
    Now, let us suppose that the claim is true for some $m \in \bN_0$,
    and we have $m+1 \le k$. Fix 
    $
        t \in \intoc{ m, \min(m + 1,s)},
    $
    $F \in \mF_{\measureAlone}^s(f,G) $, 
    and
    $
        \bx = \set*{ x_I }_{ I \subseteq [m+1] }
        \subseteq
        F.
    $
    We have three cases to~consider:
    \begin{enumerate}[label=Case \Roman*:,leftmargin=3.5em]
        \item
        $
            \polyset{ \set*{m+1} }( \bx ) \le 1
        $
        and 
        $
            m+1 = k,
        $
        \item 
        $
            \polyset{ \set*{m+1} }( \bx ) \le 1
        $
        and 
        $
            m+1 < k,
        $
        \item 
        $
            \polyset{ \set*{m+1} }( \bx ) > 1.
        $
    \end{enumerate}
    We will consider these cases in the above order.
    \begin{enumerate}[label=Case \Roman*:,leftmargin=0em, itemindent=3.5em]
        \item 
        Here, since 
        $F \in \mF_{\measureAlone}^s(f,G)$,
        $
            \bx = \set*{ x_I }_{ I \subseteq [m+1] }
            \subseteq
            F,
        $
        and $m+1 = k$,
        we have
        \begin{align*}
			\norm*{
                \diff_{ I = \emptyset }^{ [m+1] }
                    f(x_I)
            }_W
            =
            \norm*{
                \diff_{ I = \emptyset }^{ [k] }
                    f(x_I)
            }_W
        &\le 
            \polygen{ \gen{s} }( \bx )
            \sum_{ I = \emptyset }^{ [k] }
                G( x_I )
        \\
        &\le 
            \polygen{ \gen{t} }( \bx )
            \sum_{ I = \emptyset }^{ [k] }
                G( x_I )
            =
            \polygen{ \gen{t} }( \bx )
            \sum_{ I = \emptyset }^{ [m+1] }
                G( x_I ),
        \end{align*}
        where the second inequality follows by \Cref{lem::k+1_direction_<=1_gives_poly(s)<=poly(t)_for_t<=s}.
        Therefore, 
        \begin{align}
            \norm*{
                \diff_{ I = \emptyset }^{ [m+1] }
                    f(x_I)
            }_W
        &\le 
            \polygen{ \gen{t} }( \bx )
            \sum_{ I = \emptyset }^{ [m+1] }
                \del*{ 
                    \norm*{f( x_{I}) }_W
                    + G(x_{I} )
                }
        \notag
        \\
        &\le 
            \polygen{ \gen{t} }( \bx )
            \sum_{ I = \emptyset }^{ [m+1] }
                2^{k^2} \partitions{k}
                \del*{ 
                    \norm*{f( x_{I}) }_W
                    + G(x_{I} )
                }.
            \label{eqprop::prop::embedding_into_lower_orders_normed::end_case_m+1=k}
        \end{align}
        \item
        Define $\by = \set*{ y_J }_{ J \subseteq [k] }$ by 
        \begin{equation*}
            \forall I \subseteq [m+1]
            \quad \forall L \subseteq [k] \setminus [m+1] 
            \qquad 
            y_{ I \cup L }
            \coloneq 
            \begin{cases}
                x_I 
                &\text{if } \card{L} \text{ is even},
            \\
                x_{I \div \set*{m+1} }   &\text{otherwise}.
            \end{cases}
        \end{equation*}
        By \Cref{lem::folding_hypercubes_in_normed_spaces} we have
        \begin{align*}
            \polygen{ \gen{s} }(\by) 
        &\le 
            2^{k (k-m-1)}
            \partitions{k} \polygen{ \gen{m+1} }( \bx )
            \del*{
                1 + \polyset{ \set*{ m+1 } }( \bx )^{ k-m -1 }
            }
        \\
        &\le 
            2^{k(k-m-1)+1}
            \partitions{k} 
            \polygen{ \gen{ m+1 } }( \bx )
        	\le 
        	2^{k^2} \partitions{k} 
                \polygen{ \gen{ m+1 } }(\bx )
            \le 
            2^{k^2}
            \partitions{k} 
            \polygen{ \gen{ t } }( \bx ),
        \end{align*}
        where the last inequality follows from
        \Cref{lem::k+1_direction_<=1_gives_poly(s)<=poly(t)_for_t<=s}.
        Also, by \Cref{lem::folding_hypercubes_in_groups}
        we have
        \begin{equation*}
            \norm*{
                \diff_{ J = \emptyset }^{ [k] }
                    f(y_J)
            }_W
            =
            \norm*{
                (-1)^{k-m-1}
                2^{k-m-1}
                \diff_{ I = \emptyset }^{ [m+1] }
                    f(x_I)
            }_W
            =
            2^{k-m-1}
            \norm*{
                \diff_{ I = \emptyset }^{ [m+1] }
                    f(x_I)
            }_W.
        \end{equation*}
        Finally, by \Cref{lem::folding_hypercubes_with_a_sum},
        \begin{equation*}
            \sum_{ J = \emptyset }^{ [k] }
                G(y_J)
            =
            2^{k-m-1}
            \sum_{ I = \emptyset }^{ [m+1] }
                G(x_I).
        \end{equation*}

        Because $\bx \subseteq F$, we have $\by \subseteq F$, hence
        \begin{align*}
            2^{k-m-1}
            \norm*{
                \diff_{ I = \emptyset }^{ [m+1] }
                    f(x_I)
            }_W
        &=
            \norm*{
                \diff_{ J = \emptyset }^{ [k] }
                    f(y_J)
            }_W
        \\
        &\le 
            \polygen{ \gen{s} }(\by)
            \sum_{ J = \emptyset }^{ [k] }
                G( y_J )
            \le 
            2^{k^2}
            \partitions{k} 
            \polygen{ \gen{t} }( \bx )
            \sum_{ J = \emptyset }^{ [k] }
                G( y_J )
        \\
        &
        \phantom{
        {}\le
        	\polygen{ \gen{s} }(\by)
	        \sum_{ J = \emptyset }^{ [k] }
	            G( y_J )
	    }
            =
            2^{k^2}
            \partitions{k} 
            \polygen{ \gen{t} }( \bx )
            \del*{
                2^{k-m-1}
                \sum_{ I = \emptyset }^{ [m+1] }
                    G(x_I)
            }.
        \end{align*}
        Therefore,
        \begin{align}
            \norm*{
                \diff_{ I = \emptyset }^{ [m+1] }
                    f(x_I)
            }_W
        &\le 
            2^{k^2}
            \partitions{k} 
            \polygen{ \gen{t} }( \bx )
            \sum_{ I = \emptyset }^{ [m+1] }
                G(x_I)     
        \notag
        \\
        &\le
            \polygen{ \gen{t} }( \bx )
            \sum_{ I = \emptyset }^{ [m+1] }
                2^{k^2}
                \partitions{k}
                \del*{
                    \norm*{f(x_I)}_W
                    + G(x_I)
                }.
            \label{eqprop::prop::embedding_into_lower_orders_normed::end_case<=1}
        \end{align}
        \item 
        Denote
        $
            \bx'
            \coloneq 
            \set*{ x_I }_{ I \subseteq [m] }
        $
        and 
        $
            \bx''
            \coloneq 
            \set*{ x_{I \cup \set*{m+1}} }_{ I \subseteq [m] }.
        $
        By \Cref{lem::adding_k+1_direction_to_poly_j_gives_poly_j+1},
        for all $j \in [m]$, we have
        \begin{align*}
            \polyset{ \set*{m+1} }( \bx )
            \del*{
                \polynum{ j }( \bx' )
                + \polynum{ j }( \bx'' )
            }
            &\le 
            \polynum{ j+1 }( \bx )
        \shortintertext{and}
            \polyset{ \set*{m+1} }( \bx)^{ t-m } 
            \del*{
                \polynum{ m }( \bx' )
                + \polynum{ m }( \bx'' )
            }
            &\le 
            \polynum{ t }( \bx ).
        \end{align*}
        Since 
        $\polyset{ \set*{m+1} }(\bx) > 1$,
        then also
        $ \polyset{ \set*{m+1} }( \bx )^{t-m} > 1$, 
        so
        \begin{align*}
            \polygen{ \gen{m} }( \bx' )
            + \polygen{ \gen{m} }( \bx'' )
        \hspace{-2cm}&\hspace{2cm}=
            \sum_{ j = 1 }^{ m }
                \del*{
                    \polynum{ j }( \bx' )
                    + \polynum{ j }( \bx'' )
                }
        \\
            &\le 
            \polyset{ \set*{m+1} }( \bx )^{t-m}
            \del*{
                    \polynum{ m }( \bx' )
                    + \polynum{ m }( \bx'' )
                }
            +
            \sum_{ j = 1 }^{ m-1 }
                \polyset{ \set*{m+1} }(\bx)
                \del*{
                    \polynum{ j }( \bx' )
                    + \polynum{ j }( \bx'' )
                }
        \\
            &\le 
            \polynum{ t }( \bx )
            + 
            \sum_{ j = 1 }^{ m-1 }
                \polynum{ j+1 }( \bx )
        \\
            &= 
            \polynum{ t }( \bx )
            + 
            \sum_{ j = 2 }^{ m }
                \polynum{ j }( \bx )
            \le 
            \polygen{ \gen{t} }( \bx ).
        \end{align*}
        Hence,
        \begin{equation}
        \label{eqprop::prop::embedding_into_lower_orders_normed::bounding_sum_of_poly_m_by_poly_t}
            \polygen{ \gen{m} }( \bx' )
            + \polygen{ \gen{m} }( \bx'' )
            \le 
            \polygen{ \gen{t} }( \bx ).
        \end{equation}

        Since $\bx \subseteq F$, by the inductive assumption we have
        \begin{align*}
            \norm*{
                \diff_{ I = \emptyset }^{ [m] }
                    f(x_I)
            }_W
            &\le 
            \polygen{ \gen{m} }( \bx' )
            \sum_{ I = \emptyset }^{ [m] }
                2^{k^2} \partitions{k}
                \del*{ \norm*{f(x_I)}_W + G(x_I) }
        \shortintertext{and}
            \norm*{
                \diff_{ I = \emptyset }^{ [m] }
                    f\del*{ x_{I \cup \set*{m+1} } }
            }_W
            &\le 
            \polygen{ \gen{m} }( \bx'' )
            \sum_{ I = \emptyset }^{ [m] }
                2^{k^2} \partitions{k}
                \del*{ 
                    \norm*{f \del*{ x_{I \cup \set*{m+1} } } }_W
                    + G \del*{ x_{I \cup \set*{m+1} } } 
                }.
        \end{align*}
        Therefore, by the 
        \say{in particular}
        part of 
        \Cref{lem::diff_B_diff_C_to_diff_BuC},
        \begin{align*}
            \norm*{
                \diff_{ I = \emptyset }^{ [m+1] }
                    f(x_I)
            }_W
            &=
            \norm*{
                \diff_{ I = \emptyset }^{ [m] }
                	\del*{
	                    f\del*{ x_{ I \cup \set*{m+1} } }
	                -
	                    f(x_I)
	                }
            }_W
        \\  
            &\le 
            \norm*{
                \diff_{ I = \emptyset }^{ [m] }
                    f(x_I)
            }_W
            +
            \norm*{
                \diff_{ I = \emptyset }^{ [m] }
                    f\del*{ x_{ I \cup \set*{m+1} } }
            }_W
        \\  
            &\le 
            \polygen{ \gen{m} }( \bx' )
            \sum_{ I = \emptyset }^{ [m] }
                2^{k^2} \partitions{k}
                \del*{ 
                	\norm*{f(x_I)}_W 
                	+ G(x_I) 
                }
        \\
            &\quad +
            \polygen{ \gen{m} }( \bx'' )
            \sum_{ I = \emptyset }^{ [m] }
                2^{k^2} \partitions{k}
                \del*{ 
                    \norm*{f\del*{ x_{I \cup \set*{m+1} } } }_W
                    + G\del*{ x_{I \cup \set*{m+1} } }
                }
        \\
            &\le
            \del{
                \polygen{ \gen{m} }( \bx' )
                +  \polygen{ \gen{m} }( \bx'' )
            }
            \sum_{ I = \emptyset }^{ [m+1] }
                2^{k^2} \partitions{k}
                \del*{ 
                    \norm{f(x_I)}_W + G(x_I) 
                },
        \intertext{
            which, by \eqref{eqprop::prop::embedding_into_lower_orders_normed::bounding_sum_of_poly_m_by_poly_t},
        }
            &\le             
            \polygen{ \gen{t} }(\bx)
            \sum_{ I = \emptyset }^{ [m+1] }
                2^{k^2} \partitions{k}
                \del*{ 
                    \norm{f(x_I)}_W + G(x_I) 
                }.
        \end{align*}
        Therefore, we have 
        \begin{equation}
        \label{eqprop::prop::embedding_into_lower_orders_normed::end_case>1}
            \norm*{
                \diff_{ I = \emptyset }^{ [m+1] }
                    f(x_I)
            }_W
            \le 
            \polygen{ \gen{t} }(\bx)
            \sum_{ I = \emptyset }^{ [m+1] }
                2^{k^2} \partitions{k} 
                \del*{ \norm{f(x_I)}_W + G(x_I) }
        \end{equation}
        in this case.
    \end{enumerate}
    Combining 
    \eqref{eqprop::prop::embedding_into_lower_orders_normed::end_case_m+1=k}, \eqref{eqprop::prop::embedding_into_lower_orders_normed::end_case<=1},
    and \eqref{eqprop::prop::embedding_into_lower_orders_normed::end_case>1}, 
    we have 
    \begin{equation*}
        \norm*{
            \diff_{ I = \emptyset }^{ [m+1] }
                f(x_I)
        }_W
        \le 
        \polygen{ \gen{t} }(\bx)
        \sum_{ I = \emptyset }^{ [m+1] }
            2^{k^2} \partitions{k}
            \del*{ \norm{f(x_I)}_W + G(x_I) }
    \end{equation*}
    in all cases.
    
    Since 
    $
        \bx = \set*{ x_I }_{ I \subseteq [m+1] }
        \subseteq
        F
    $
    is arbitrary, we see that 
    \begin{equation*}
        2^{k^2} \partitions{k}
        \del*{ \norm{f}_W + G }
        \in
        \bD_{\measureAlone}^{t}( f )    
    \quad \text{and} \quad 
        F \in \mF_{\measureAlone}^{ t }\del*{
            f, 
            2^{k^2} \partitions{k}
            \del*{ \norm{f}_W + G } 
        },
    \end{equation*}
    which ends the proof of the inductive step, since
    $ F \in \mF_{\measureAlone}^{ t }(f, G)$
    is arbitrary.
    Therefore, by~induction, the proposition is true for all $t \in \intcc{0,s}$. 
\end{proof}
\begin{lemma} 
\label{lem::from_m_upper_to_upper(normed)}
    Let 
    $ \del*{ V, \normAlone[V]} $
    and
    $ \del*{ W, \normAlone[W]} $
    be normed spaces, 
    $X \subseteq V$ be nonempty, 
    and $\measureAlone$ be a measure on $X$. 
    Let $s \ge 0$. 
    Then for every measurable functions
    $
        f \colon X \to W
    $
    and 
    $
        G \in \bD^{s}_{\measureAlone}(f),
    $
    there are functions
    $ f', G'$ such that 
    $ f = f'$ and $G = G'$
    $\measureAlone$-almost everywhere 
    and 
    $G' \in \bD^{s}(f')$.
\end{lemma}
\begin{proof}
    Fix $F \in \mF_{\measureAlone}^{s}(f, G)$. 
    Let us define $f'$ and $G'$ by the formulas
    \begin{equation*}
        f' \coloneq f \indicator{F}
        \quad\text{and}\quad
        G' \coloneq G \indicator{F} + \infty \indicator{X \setminus F}.
    \end{equation*}
    Then $f = f'$ and $G = G'$ 
    $\measureAlone$-almost everywhere.
    Moreover, we have $G' \in \bD^{s}(f')$.
    Indeed, let 
    $ k \in \bN_0 $
    be such that
    $ s \in \intoc{k-1,k}$
    and let us fix 
    $
        \bx = \set*{ x_I }_{ I \subseteq [k] }
        \subseteq 
        X.
    $
    If, on the one hand,
    $
        \bx \subseteq F,
    $
    then
    \begin{equation*}
        \norm*{
            \diff_{ I = \emptyset }^{ [k] }
                f'(x_I)
        }_W
        =
        \norm*{
            \diff_{ I = \emptyset }^{ [k] }
                f(x_I)
        }_W
        \le 
        \polygen{ \gen{s} }( \bx )
        \sum_{ I = \emptyset }^{ [k] }
            G(x_I)
        =
        \polygen{ \gen{s} }( \bx )
        \sum_{ I = \emptyset }^{ [k] }
            G'(x_I).
    \end{equation*}
    On the other hand, 
    if there exists $J \subseteq [k]$
    such that $x_J \notin F$,
    then we have two possibilities:
    \begin{itemize}[leftmargin=0.5cm]
        \item
        If 
        $
            \diff_{ I = \emptyset }^{ [k] }
                f'(x_I)
            =
            0,
        $
        then
        \begin{equation*}
            \norm*{
                \diff_{ I = \emptyset }^{ [k] }
                    f'(x_I)
            }_W
            =
            0
            \le 
            \polygen{ \gen{s} }( \bx )
            \sum_{ I = \emptyset }^{ [k] }
                G'(x_I).
        \end{equation*}
        \item 
        If
        $
            \diff_{ I = \emptyset }^{ [k] }
                f'(x_I)
            \ne 
            0,
        $
        then $ k > 0 $ and $s > 0$.
        Thus, by 
        \Cref{cor::If_P^set(j)(x)=0_for_any_j_then_diff_f(x)=0},
        we have
        $
            \polyset{ \set*{j} }( \bx ) > 0
        $
        for all $ j \in [k] $. 
        
        In~consequence,
        \begin{equation*}
            \polygen{ \gen{s} }( \bx )
            \ge 
            \polynum{ s }( \bx )
            =
            \polyset{ \set*{k} }(\bx)^{ s-k+1 }
            \prod_{ j = 1 }^{ k-1 }
                \polyset{ \set*{j} }( \bx )
            >
            0.
        \end{equation*}
        Since $x_J \notin F$, we have
        $
            G'(x_J) = \infty.
        $
        Noting that $f'$~is finite everywhere, we have
        \begin{equation*}
            \norm*{
                \diff_{ I = \emptyset }^{ [k] }
                    f'(x_I)
            }_W
            \le 
            \infty 
            =
            \polygen{ \gen{s} }( \bx ) G'(x_J)
            \le 
            \polygen{ \gen{s} }( \bx )
            \sum_{ I = \emptyset }^{ [k] }
                G'(x_I).
        \end{equation*}
    \end{itemize}
    Thus, for all 
    $ 
        \bx = \set*{ x_I }_{ I \subseteq [k] }
        \subseteq
        X
    $
    we have
    $
        \norm*{
            \diff_{ I = \emptyset }^{ [k] }
                f(x_I)
        }_W
        \le 
        \polygen{ \gen{s} }( \bx )
        \sum_{ I = \emptyset }^{ [k] }
            G'(x_I).
    $
    Therefore,
    $ G' \in \bD^s( f' )$,
    ending the~proof.
\end{proof}
\begin{lemma}
\label{lem::diff_for_smooth_converges_to_derivative}
    Let $n, m \in \bN$, 
    $\Omega \subseteq \bR^n$ be open,
    $
        \alpha 
        = 
        \del{ \alpha_1, \, \ldots, \, \alpha_n }
        \in 
        \bN_0^n
    $ 
    be such that 
    $ \abs{ \alpha } = m$.
    Let 
    $   
        \set*{v_{\ell}}_{ \ell =1 }^{ m }
    $
    be a tuple of unit vectors in $\bR^n$ such that, 
    for all $i \in [n]$, 
    exactly 
    $ \alpha_i $
    of them equal $\be_i$, 
    where $\be = \del{ \be_i }_{ i=1 }^n$ 
    is~the~canonical basis of $\bR^n$.
    Define
    $
        \set*{ v_L }_{ L \subseteq [m] }
    $
    by the formula
    $
        v_L 
        =
        \sum_{ \ell \in L }
            v_{ \ell }.
    $
    (Note that we have $v_{\emptyset} = 0$.)
    Then for all
    $ \psi \in \smoothCompactSupport( \Omega )$
    and all $y \in \Omega$, 
    we have
    \begin{equation*}
        h^{-m}
        \diff_{ L = \emptyset }^{ [m] }
            \psi( y + hv_L )
        \xrightarrow{ h \to 0 }
        \partial^{ \alpha }
            \psi(y).
    \end{equation*}
\end{lemma}
\begin{proof}
    Fix $y \in \Omega$
    and let 
    $
        D 
        \coloneq 
        \min\del*{
        	\dist\del{ y, \partial \Omega },
        	1
        }.
    $
    By the 
    \say{in particular}
    part of
    \Cref{lem::diff_B_diff_C_to_diff_BuC},
    for all $h \in \intoo{ -D/m, D/m}$,
    we have 
    \begin{align*}
        \diff_{L = \emptyset }^{ [m] }
            \psi( y + hv_L )
    \hspace{-2cm}&\hspace{2cm}=
        \diff_{L = \emptyset}^{ [m-1] }
            \del*{
                \psi( y + hv_L + hv_m )
                - \psi( y + hv_L )              
            },
    \intertext{
    	which, by the fundamental theorem of calculus,
    }
        &=
        \diff_{L = \emptyset}^{ [m-1] }
        \integral{ [0,1] }{
            \partial_{v_m}
            \psi( y + hv_L + ht_m v_m )
            h
        }{t_m},
    \intertext{
    	which, by repeatedly applying the previous two steps for $j = 2, \, \ldots, \,m$,
    }
        &=
        \diff_{L = \emptyset}^{ [m-j] }
        \integral{ [0,1]^j }{
            \del*{
                \partial_{v_{m-j+1} }
                \,
                \cdots 
                \,
                \partial_{v_m}
            }
            \psi\del*{
                y + hv_L
                +
                \sum_{i = m-j+1}^m
                    ht_i v_i
            } 
            h^j
        }{\del*{ t_{m-j+1}, \, \ldots, \, t_m } },
    \intertext{
    	which, since 
    	$
    		\partial_{v_1} \, \cdots \, \partial_{v_m} 
    		\psi
    		= 
    		\partial^{\alpha}
    		\psi,
    	$
    	for 
    	$j=m$,
    }
        &=
        \integral{ [0,1]^m }{
            \partial^{\alpha} 
            \psi \del*{
                y 
                + 
                h \sum_{ i = 1 }^{ m }
                    t_i v_i
            }
            h^m
        }{t}.
    \end{align*}
    As
    $\partial^{\alpha} \psi$ is continuous,
    for all 
    $
        t = \del{ t_i }_{i=1}^m
        \in 
        \intcc{ 0,1 }^m
    $
    we have
    $
    	{
        \partial^{\alpha} 
        \psi \del*{
            y 
            + 
            h \sum_{ i = 1 }^{ m }
                t_i v_i
        }
        \to 
        \partial^{\alpha} 
        \psi(y)
        }
    $
    as~%
    $
        h \to 0.
    $
    Thus, as 
    $
        \partial^{\alpha} 
        \psi
    $
    is bounded, 
    by the Lebesgue dominated convergence theorem,
    we have
    \begin{equation*}
        h^{-m}
        \diff_{ L = \emptyset }^{ [m] }
            \psi( y + h v_L )
        =
        \integral{ [0,1]^m }{
            \partial^{\alpha} 
            \psi \del*{
                y 
                + h \sum_{ i = 1 }^{ m }
                    t_i v_i
            }
        }{t}
        \xrightarrow{ h \to 0 }
        \partial^{ \alpha }
            \psi(y).
        \tag*{\qedhere}
    \end{equation*}
\end{proof}
\begin{lemma}
\label{lem::diff(f)g=fdiff(g)_under_integral}
    Let $n, m \in \bN$, 
    $ \Omega \subseteq \bR^n $ 
    be open,
    $ f \in \locIntegrable{1}( \Omega )$,
    and 
    $ g \in \locIntegrable{\infty}( \Omega )$.
    Fix 
    $ 
        \set*{ v_{ \ell } }_{ \ell = 1 }^{ m } 
        \subseteq 
        \bR^n
    $
    and~define
    $
        \set*{ v_L }_{ L \subseteq [m] }
    $
    by the formula
    $
        v_L
        =
        \sum_{ \ell \in L }
            v_{ \ell }
    $
    for all $L \subseteq [m] $.    
    Fix $x \in \Omega$
    and $r > 0$ such that 
    $\clball{x,r} \subseteq \Omega$,
    and~suppose that 
    $
        \supp(g)
        \subseteq 
        \ball{ x, r - D },
    $
    where 
    $
        D
        \ge 
        \sum_{ \ell = 1 }^{ m }
            \norm{ v_{ \ell } }.
    $
    Then
    \begin{equation*}
        \integral{\ball{x, r}}{
            \del*{
                \diff_{ L = \emptyset }^{ [m] }
                    f\del*{ y + v_L }
            }
            g(y)
        }{y}
        =
        \integral{\ball{x, r}}{
            f(y)
            \del*{
                \diff_{ L = \emptyset }^{ [m] }
                    g\del*{ y - v_L }
            }
        }{y}.
    \end{equation*}
\end{lemma}
\begin{proof}
    Let us notice that for all 
    $L \subseteq [m]$ we have
    $
        \norm{ v_L } \le D,
    $
    so 
    $
        \ball{x, r-D} + v_L
        \subseteq 
        \ball{x,r}.
    $
    In~consequence, for every $L \subseteq [m] $,
    \begin{align*}
        \integral{ \ball{x,r} }{
            f(y + v_L)
            g(y)
        }{y}
        &=
        \integral{ \ball{x,r-D} }{
            f(y + v_L)
            g(y)
        }{y}
    \\
        &=
        \integral{ \ball{x,r-D} + v_L }{
            f(y)
            g(y- v_L)
        }{y}
        =
        \integral{ \ball{x,r} }{
            f(y)
            g(y- v_L)
        }{y}.
    \end{align*}
    Therefore,
    \begin{multline*}
        \integral{\ball{x, r}}{
            \del*{
                \diff_{ L = \emptyset }^{ [m] }
                    f\del*{ y + v_L }
            }
            g(y)
        }{y}
        =
        \diff_{ L = \emptyset }^{ [m] }
            \integral{\ball{x, r}}{
                f\del*{ y + v_L }
                g(y)
            }{y}
    \\
        =
        \diff_{ L = \emptyset }^{ [m] }
            \integral{\ball{x, r}}{
                f(y)
                g\del*{ y - v_L }
            }{y}
        =
        \integral{\ball{x, r}}{
            f(y)
            \del*{
                \diff_{ L = \emptyset }^{ [m] }
                    g\del*{ y - v_L }
            }
        }{y},
    \end{multline*}
    as claimed.
\end{proof}
\begin{lemma}
\label{lem::bounded_finite_differences_implies_pointwise_bound_on_derivative}
    Let $n, m \in \bN$, $\Omega \subseteq \bR^n$ be open, and 
    $
        \alpha = \del{ \alpha_1, \, \ldots, \, \alpha_n }
        \in \bN_0^n
    $ 
    satisfy 
    $ \abs{ \alpha } = m$.
    Let~%
    $   
        \set*{v_{\ell}}_{ \ell =1 }^{ m }
    $
    be a tuple of unit vectors in $\bR^n$ such that, 
    for all $i \in [n]$, 
    exactly 
    $ \alpha_i $
    of them equal $\be_i$, 
    where 
    $\del{ \be_i }_{ i=1 }^n$ 
    is~the~canonical basis of 
    $\bR^n$.
    Define 
    $
        \set*{ v_L }_{ L \subseteq [m] }
    $
    by the formula
    $
        v_L 
        =
        \sum_{ \ell \in L }
            v_{ \ell }.
    $
    Fix $f \in \locSobolev{m, 1}(\Omega)$
    and~let~%
    $G \in \locIntegrable{1}(\Omega)$
    be nonnegative.
    Fix 
    $
        x \in
        \lebesguePoints{ \partial^{\alpha} f }
        \cap \lebesguePoints{ G }
    $ 
    and let 
    $R > 0$ be such that 
    $
        \clball{x, 2R } \subseteq \Omega.
    $
    Fix~$r \in \intoc{0,R}$ 
    and 
    $\eps \in \intoo*{0,\frac{r}{m+1}}$.
    If 
    \begin{equation*}
        \forall h \in \intoc{ 0, \eps }
        \quad 
        \lambda \forall y \in \ball{x,r}
        \qquad 
        \abs*{
            \diff_{L = \emptyset}^{ [m] }
               f\del*{y + hv_L }
        }
        \le 
        h^m 
        \sum_{ L = \emptyset }^{ [m] }
            G\del*{y + hv_L },
    \end{equation*}
    then
    \begin{equation*}
        \abs{
            \partial^{ \alpha } f(x)
        }
        \le 
        2^m G(x).
    \end{equation*}
\end{lemma}
\begin{proof}
    First of all, let us notice that
    \begin{equation}
	\label{proofeq::lem::bounded_finite_differences_implies_pointwise_bound_on_derivative::estimate_for_hvL}
        \forall 
        h \in \intoc{0, \eps }
    \qquad 
        \norm{
            hv_L
        }
        =
        \norm*{
            h \sum_{ \ell \in L } v_{ \ell }
        }
        \le 
        h \sum_{ \ell \in L } \norm{ v_{\ell} }
        =
        h \card{ L }
        \le 
        hm.
    \end{equation}
    Next, fix  
    $
        \psi \in \smoothCompactSupport\del*{
            \ball{x, r - m\eps }
        }.
    $
    By 
    \Cref{lem::diff_for_smooth_converges_to_derivative},
    for all 
    $ y \in \ball{x, r}$
    we have
    \begin{equation*}
        h^{-m}
        \diff_{L = \emptyset}^{ [m] }
            \psi\del*{ y - hv_L }
        =
        (-1)^{m}
        (-h)^{-m}
        \diff_{L = \emptyset }^{ [m] }
            \psi\del*{ y + (-h)v_L }
        \xrightarrow{ h \to 0^+ }
        (-1)^{ m }
        \partial^{ \alpha }
            \psi(y).
    \end{equation*}
    Also, note that $\psi$ can be treated as an element of
    $\smoothCompactSupport(\bR^n)$ by
    putting
    $
        \psi \rvert_{ 
            \bR^n \setminus 
            \ball{x, r - m\eps }
        }
        \equiv
        0.
    $
    In~this sense, we have
    $
        \psi 
        \in 
        \holder{ m-1, 1}\del*{ \bR^n }.
    $
    Thus, by   \Cref{cor::multipointwise_bound_for_higher_order_holder},
    for almost all 
    $
        y \in \ball{x, r}
    $
    and all 
    $h \in \intoc{0, \eps}$,
    we have
    \begin{align*}
        \abs*{
            f(y)
            \del*{
                \diff_{L = \emptyset}^{ [m] }
                    \psi\del*{ y - hv_L }
            }
        }
        \le
        \abs{ f(y) }
        \del*{
            \polygen{ \gen{m} }\del*{ y - h\bv }
            2^m
            \widehat{C}_{n,m}
            \seminorm{ \psi }_{ 
                \holder{m-1, 1}\del*{ \bR^n } 
            }
        }
    &
	\\
        =
        h^m
        2^{m^2}
        \widehat{C}_{n,m}
        \abs{ f(y) }
        \seminorm{ \psi }_{ 
            \holder{m-1, 1}\del*{ \bR^n } 
        },
    &
    \end{align*}
    where we used the fact by
    \Cref{lem::invariant_transformations_of_poly}
    and 
    \Cref{cor::value_of_Poly_when_tuple_is_a_hyperparallelogram},
    we have
    \begin{equation*}
        \polygen{ \gen{m} }(y - h\bv )
        =
        \polygen{ \gen{m} }( h\bv )
        =
        2^{m(m-1)}
        \prod_{ \ell = 1 }^{ m }
            \norm{ h v_{\ell } }
        =
        h^m
        2^{m(m-1)}
    \end{equation*}
    since
    $
        a \mapsto y - a
    $
    is an isometry of $\bR^n$.
	Therefore, by the
    Lebesgue dominated convergence theorem, 
    \begin{multline*}
        h^{-m}
        \integral{ \ball{x, r} }{
            \del*{
                \diff_{L = \emptyset}^{ [m] }
                    f \del*{ y + h v_L }
            }
            \psi(y)
        }{ y }
    \\
        =
        h^{-m}
        \integral{ \ball{x, r} }{
            f(y)
            \del*{
                \diff_{L = \emptyset}^{ [m] }
                    \psi\del*{ y - h v_L }
            }
        }{ y }
        \xrightarrow{h \to 0^+}
        (-1)^{m}
        \integral{ \ball{x, r} }{
            f(y)
            \partial^{\alpha} \psi(y)
        }{ y }
    \\
    	=
        \integral{ \ball{x, r} }{
            \partial^{\alpha} f(y)
            \psi(y)
        }{ y },
    \end{multline*}
    where, by 
    \Cref{lem::diff(f)g=fdiff(g)_under_integral}
    (with $D \coloneq \eps m$),
    the first equality is true for all $h \in \intoc{0,\eps}$
    since we have
    $
    	\sum_{ \ell = 1 }^{ m }
    		\norm{ h v_{\ell} }
    	=
    	hm
    	\le 
    	\eps m.
    $    
    Thus, 
    %
    \begin{multline}
    \label{proofeq::lem::bounded_finite_differences_implies_pointwise_bound_on_derivative::convergence_on_test_functions}
        \forall 
            \psi \! \in \! \smoothCompactSupport\del*{
                \ball{x, r - m\eps }
            }
	\\ 
        h^{-m} 
        \integral{ \ball{x, r} }{
            \del*{
                \diff_{L = \emptyset}^{ [m] }
                    f\del*{ y + h v_L }
            }
            \psi(y)
        }{ y }
        \xrightarrow{h \to 0^+}
        \integral{ \ball{x, r} }{ 
            \partial^{\alpha} f(y)
            \psi(y)
        }{ y }.
    \end{multline}

    Therefore, since thanks to
    \eqref{proofeq::lem::bounded_finite_differences_implies_pointwise_bound_on_derivative::estimate_for_hvL} 
    we have
    $
        \ball{x + h v_L, r}
        \subseteq
        \ball{x, r + hm}
    $
    for all
    $ h \in \intoc{0, \eps},$ 
    \begin{align}
        \integral{ \ball{x, r} }{
            \sum_{ L = \emptyset }^{ [m] }
                G\del*{ y + h v_L }
        }{ y }
    	=
        \sum_{ L = \emptyset }^{ [m] }
            \integral{ \ball{x, r} }{
                G\del*{ y + h v_L }
            }{y}
    &
    \notag
    \\
    	=
        \sum_{ L = \emptyset }^{ [m] }
            \integral{ \ball{x + hv_L, r  } }{
                G
            }{y}
    &\le 
        \sum_{ L = \emptyset }^{ [m] }
            \integral{ \ball{x, r + hm  } }{
                G
            }{y}
    \notag
    \\
    &=
        2^m
        \integral{ \ball{x, r + hm  } }{
            G
        }{y}.
        \label{proofeq::lem::bounded_finite_differences_implies_pointwise_bound_on_derivative::bound_on_integral_of_sum_G}
    \end{align}
    %

    Now, fix $\delta \in \intoo{0,\eps}$ and let
    $
        \psi_{ \delta }
        \coloneq 
        \del*{
            \sign\del{ \partial^{\alpha} f } 
            \indicator{ \ball{x, r-(m+1)\eps} }
        }_{ \delta },
    $
    where 
    $
        ( \cdot )_{\delta}
    $
    denotes the mollification by the standard mollifier 
    \cite[Remark C.18.(ii)]{leoni_first_course}.
    Then
    $ \psi_{ \delta } $ 
    is smooth and we have
    $
        \supp( \psi_{ \delta } )
        \subseteq 
        \ball{x, r- m \eps }.
    $
    Hence,
    $
        \psi_{ \delta } \in \smoothCompactSupport\del*{
            \ball{x, r - m\eps }
        }.
    $
    Moreover, since 
    \begin{equation*}
        \norm*{
            \sign\del*{ \partial^{\alpha} f } 
            \indicator{ \ball{x, r-(m+1)\eps} }
        }_{ \integrable{ \infty }\del*{
                \bR^n
            }
        }
        \le 
        1,
    \end{equation*}
    by Young convolution inequality 
    \cite[Theorem C.15]{leoni_first_course}
    we have
    $
        \norm*{
            \psi_{ \delta }
        }_{ \integrable{ \infty }\del*{
                \ball{x, r }
            }
        }
        \le 
        \norm*{
            \psi_{ \delta }
        }_{ \integrable{ \infty }( \bR^n ) }
        \le 
        1.
    $
    Therefore, for all $h \in \intoc{0, \eps}$, we have
    \begin{multline*}
        h^{-m} \!
        \integral{ \ball{x, r} }{
            \del*{
                \diff_{L = \emptyset}^{ [m] }
                    f\del*{y + h v_L }
            }
            \psi_{ \delta }(y)
        }{ y }
    \\
    \le 
        h^{-m}
        \integral{ \ball{x, r} }{
            \abs*{
                \diff_{L = \emptyset}^{ [m] }
                    f\del*{y + h v_L }
            }
        }{ y }
    \le 
        \integral{ \ball{x, r} }{
            \sum_{ L = \emptyset }^{ [m] }
                G\del*{y + h v_L }
        }{y}
    	\le 
        2^m \!
        \integral{ \ball{x, r+hm} }{G}{y},
    \end{multline*}
    where the last inequality follows from 
    \eqref{proofeq::lem::bounded_finite_differences_implies_pointwise_bound_on_derivative::bound_on_integral_of_sum_G}.
    Since
    $
        G \in \locIntegrable{1}( \Omega )
    $
    and 
    $
        \ball{x, 2r} \subseteq \Omega ,
    $
    we have
    \begin{equation*}
        \integral{ \ball{x, r+hm} }{G}{y}
        \xrightarrow{ h \to 0^+}
        \integral{ \ball{x, r} }{G}{y}.
    \end{equation*}
    In consequence, because of 
    \eqref{proofeq::lem::bounded_finite_differences_implies_pointwise_bound_on_derivative::convergence_on_test_functions},
    \begin{equation}
    \label{proofeq::lem::bounded_finite_differences_implies_pointwise_bound_on_derivative::estimate_with_mollification}
        \integral{ \ball{x, r} }{
            \partial^{\alpha} f(y)
            \psi_{ \delta }(y)
        }{ y }
        \le 
        2^m \integral{ \ball{x, r} }{ G }{ y }.
    \end{equation}

    Next, 
    as 
    $
        \norm*{
            \psi_{ \delta }
        }_{ \integrable{ \infty }\del*{
                \ball{x, r }
            }
        }
        \le 
        1,
    $
    for almost all
    $
        y \in \ball{x, r}
    $
    we have
    $
        \abs*{
            \partial^{\alpha} f(y)
            \psi_{ \delta }(y)
        }
        \le 
        \abs*{
            \partial^{\alpha} f(y)
        }.  
    $
    Also, note that by the mollifier theorem
    \cite[Theorem C.19]{leoni_first_course},
    for almost every
    $
        y \in \ball{x, r}
    $
    we have
    \begin{equation*}
        \psi_{ \delta }(y)
        \xrightarrow{ \delta \to 0^+ }
        \sign\del{ \partial^{\alpha} f(y) } 
        \indicator{ \ball{x, r-(m+1)\eps} }(y).
    \end{equation*}
    Therefore, by 
    the Lebesgue dominated convergence theorem,
    \begin{align*}
        \integral{ \ball{x, r} }{
            \partial^{\alpha} f(y)
            \psi_{ \delta }(y)
        }{ y }
        \xrightarrow{ \delta \to 0^+ }
    &
        \integral{ \ball{x, r} }{
            \abs*{ \partial^{\alpha} f(y) }
            \indicator{ \ball{x, r-(m+1)\eps} }(y)
        }{ y }
    \\
    &=
        \integral{ \ball{x, r-(m+1)\eps} }{
            \abs*{ \partial^{\alpha} f(y) } 
        }{ y }.
    \end{align*}
    Thus, since the estimate in 
    \eqref{proofeq::lem::bounded_finite_differences_implies_pointwise_bound_on_derivative::estimate_with_mollification}
    is independent of 
    $\delta \in \intoo{0, \eps}$,
    we have
    \begin{equation}
    \label{proofeq::lem::bounded_finite_differences_implies_pointwise_bound_on_derivative::estimate_with_r_eps}
        \integral{ \ball{x, r-(m+1)\eps} }{
            \abs*{ \partial^{\alpha} f(y) } 
        }{ y }
        \le 
        2^m \integral{ \ball{x, r} }{ G }{ y }.
    \end{equation}
    %

    Recall that within the statement of the lemma the only assumptions on $ \eps$ and $r$ are
    \begin{itemize}[leftmargin=1.5em]
        \item
        $ r \in \intoc{0,R}$
        and 
        $ \eps \in \intoo*{0, \frac{r}{m+1} } $;
        \item 
        For all $h \in \intoc{0,\eps}$
        and almost all $y \in \ball{x, r}$,
        $
            \abs*{
                \diff_{L = \emptyset}^{ [m] }
                   f\del*{y + hv_L }
            }
            \le 
            h^m 
            \sum_{ L = \emptyset }^{ [m] }
                G\del*{y + hv_L }.
        $
    \end{itemize}
    Let us notice that both of these assumptions remain satisfied if we replace 
    $r$ with an arbitrary $r' \in \intoc{0,r}$
    and
    $\eps$ with an arbitrary 
    $
    	\eps' 
    	\in 
    	\intoo*{
    		0,
    		\min\del*{ \eps, \frac{r'}{m+1} }
    	}.
    $
    Therefore, as
    \eqref{proofeq::lem::bounded_finite_differences_implies_pointwise_bound_on_derivative::estimate_with_r_eps}
    remains valid if we replace
    $\eps$ and $r$ with such $\eps'$ and $r'$, respectively, we get that
    \begin{equation}
    \label{proofeq::lem::bounded_finite_differences_implies_pointwise_bound_on_derivative::estimate_after_changing_to_prims}
        \forall 
        	r' \! \in \! \intoc{0,r}
    \, \, \, \, 
    	\forall 
    		\eps' \! \in \!
	    	\intoo*{
	    		0,
	    		\min\del*{ \eps, \tfrac{r'}{m+1} }
	    	}
	\quad \, \,
        \integral{ \ball{x, r'-(m+1)\eps'} }{
            \abs*{ \partial^{\alpha} f(y) } 
        }{ y }
        \le 
        2^m 
        \! \!
        \integral{ \ball{x, r'} }{
        \! \! G }{ y }.
    \end{equation}

    Let us fix
    $ 
    	r' \in \intoc{0,r}
    $
    and note that,
    by the monotone convergence theorem,
    we have
    \begin{align*}
        &
        \integral{ \ball{x, r'-(m+1)\eps'} }{
            \abs*{ \partial^{\alpha} f(y) } 
        }{ y }  
    \\
    &\quad
        =
        \integral{ \ball{x, r' } }{
            \abs*{ \partial^{\alpha} f(y) } 
            \indicator{
                \ball{x, r'-(m+1)\eps'}
            }(y)
        }{ y }
        \xrightarrow{ \eps' \to 0^+ }
        \integral{ \ball{x, r' } }{
            \abs*{ \partial^{\alpha} f(y) } 
        }{ y }.
    \end{align*}
    Therefore, since the right-hand side of 
    the inequality in
    \eqref{proofeq::lem::bounded_finite_differences_implies_pointwise_bound_on_derivative::estimate_after_changing_to_prims}
    does not depend on $\eps'$, 
    we obtain the estimate
    \begin{equation*}
        \abs*{
            \integral{ \ball{x, r'} }{
                \partial^{ \alpha } f
            }{y}
        }
        \le 
        \integral{ \ball{x, r'} }{
            \abs*{ \partial^{\alpha} f(y) } 
        }{ y }
        \le 
        2^m \integral{ \ball{x, r'} }{ G }{ y }.
    \end{equation*}
    Finally, since $r' \in \intoc{0,r}$ is arbitrary, and we have
    $
        x \in 
        \lebesguePoints{ \partial^{ \alpha} f }
        \cap \lebesguePoints{ G },
    $ 
    \begin{equation*}
        \abs{ \partial^{ \alpha} f(x) }
        \xleftarrow{ r' \to 0^+ }
        \abs*{
            \sintegral{ \ball{x, r'} }{
                \partial^{ \alpha } f
            }{y}
        }
        \le 
        2^m \sintegral{ \ball{x, r'} }{ G }{ y }
        \xrightarrow{ r' \to 0^+ }
        2^m G(x),
    \end{equation*}
    giving us the desired estimate
    $
        \abs{ \partial^{ \alpha} f(x) }
        \le 
        2^m G(x).
    $
\end{proof}
\begin{lemma}
\label{lem::loc_hajlasz_implies_loc_sobolev}
    Let $n \in \bN$ and $\Omega \subseteq \bR^n$ be open. 
    If 
    $f \in \locHajlasz{1,1}\del*{ \Omega }$,
    i.e.,
    $
    	f \in \locIntegrable{1}\del*{ \Omega }
    $  
    and there exists 
    $
        g \in \locIntegrable{1}(\Omega)
    $
    such that
    $
        g \in \bD_{\lambda}^1(f),
    $
    then     
    $f \in \locSobolev{1,1}(\Omega)$.
\end{lemma}
\begin{proof}
    Fix $x_0 \in \Omega$ and define
    \begin{equation*}
        \Omega_m
        \coloneq 
        \setc*{
            x \in \Omega
        }{
            \dist\del*{
                x, \partial \Omega
            } 
            > \frac{1}{m}
        }
        \cap 
        \ball{x_0, m }.
    \end{equation*}
    It is easy to see that each $\Omega_m$ is open, that $\cl{ \Omega_m } \subseteq \Omega$, and that $\Omega_m \subseteq \Omega_{m'}$ if $m \le m'$.
    Fix $\ell \in \bN$.
    Since $f, g \in \locIntegrable{1}(\Omega)$, we have
    $
        f \rvert_{\Omega_\ell}, 
        g \rvert_{\Omega_\ell}
        \in 
        \integrable{1}(\Omega_\ell).
    $
    In consequence, 
    $
        f \rvert_{ \Omega_\ell }
        \in 
        \hajlasz{1,1}( \Omega_\ell).
    $
    Hence, by 
    \cite[5.13 Remark on the case $p=1$]{analysis_on_metric_spaces},
    we have 
    $
        f \rvert_{ \Omega_{\ell} }
        \in 
        \sobolev{1,1}( \Omega_{\ell}).
    $
    Therefore, for all $ i \in [n]$,
    $
        \partial_{x_i} \del*{
            f \rvert_{ \Omega_\ell }
        }
    $
    exists in a weak sense and
    $
        \partial_{x_i}\del*{
            f \rvert_{ \Omega_\ell }
        }
        \in 
        \integrable{1}( \Omega_\ell).
    $
    Moreover, for all such $i$ and all $m, m' \in \bN$ such that
    $m \le m'$, we have 
    $
        \partial_{x_i}\del*{
            f \rvert_{ \Omega_m }
        }
        =
        \partial_{x_i}\del*{
            f \rvert_{ \Omega_m' }
        } 
    $
    almost everywhere in 
    $\Omega_m$.

    Fix $i \in [n]$.
    Let us define
    $
       h \colon \Omega \to \bR
    $
    by the formula
    \begin{equation*}
        h
        \coloneq
        \sum_{m=1}^{\infty}
        \partial_{x_i} \del*{
            f \rvert_{ \Omega_{m+1} }
        }
        \indicator{
            \Omega_{m+1}
            \setminus 
            \Omega_{m}
        }
        +   
        \partial_{x_i} \del*{
            f \rvert_{ \Omega_1 }
        }
        \indicator{
            \Omega_{1}
        }. 
    \end{equation*}
    For all $m \in \bN$, we have
    $
        h \rvert_{ \Omega_m }
        =
        \partial_{x_i} \del*{
            f \rvert_{\Omega_m}
        }
    $
    almost everywhere in $\Omega_m$.
    We will show that 
    $
        h = \partial_{x_i} f
    $   
    in a weak sense.
    Fix $\phi \in \smoothCompactSupport(\Omega)$.
    There exists $m \in \bN$ such that
    $\supp(\phi) \subseteq \Omega_m$.
    Therefore, 
    \begin{align*}
        \integral{\Omega}{
            f 
            \partial_{x_i} \phi 
        }{x}
        =
        \integral{\Omega_m}{
            f 
            \partial_{x_i} \phi 
        }{x}
    &=
        \integral{\Omega_m}{
            \del*{
                f \rvert_{\Omega_m}
            }
            \partial_{x_i} \phi 
        }{x}
    \\
    &= 
    	-
        \integral{\Omega_m}{
            \partial_{x_i} \del*{
                f \rvert_{\Omega_m}
            }
            \phi 
        }{x}
        =
        -
        \integral{\Omega_m}{
            h
            \phi 
        }{x}
    	=
        -
        \integral{\Omega}{
            h
            \phi 
        }{x}.
    \end{align*}
    Thus, 
    $
        h = \partial_{x_i} f
    $   
    in a weak sense. 
    Moreover, 
    since 
    $
        h \rvert_{ \Omega_m }
        \in 
        \integrable{1}( \Omega_m ),
    $
	we have that
    $
        \partial_{x_i} f 
        \in 
        \locIntegrable{1}( \Omega ).
    $
    Since $i \in [n]$ is arbitrary, 
    we conclude that
    $
        f \in \locSobolev{1,1}(\Omega).
    $
\end{proof}
\begin{theorem}
\label{thm::local_hypergradients_imply_local_sobolev}
    Let $n, k \in \bN$, $s \in \intoc{ k-1, k}$, and $\Omega \subseteq \bR^n$ be open. 
    Let 
    $f \in \locIntegrable{1}( \Omega )$
    be such that there is 
    $ 
        G 
        \in 
        \bD^{s}_{ \lambda }(f)
        \cap 
        \locIntegrable{1}( \Omega ).
    $
    Then, for every $\alpha \in \bN_0^n$
    satisfying $0 < \abs{ \alpha } \le s$,
    \begin{enumerate}[label=(\alph*)]
        \item
        \label{thmitem::thm::local_hypergradients_imply_local_sobolev::item::existence}
        $\partial^{ \alpha } f$ exists in a weak sense;
        in particular,
        $
        	\partial^{ \alpha } f
        	\in 
        	\locIntegrable{1}\del*{ \Omega }
        $.
        \item 
        \label{thmitem::thm::local_hypergradients_imply_local_sobolev::item::bounded_by_G}
        $
            \abs{ \partial^{ \alpha } f } 
            \le 
            2^{k^2}
            \partitions{ k }
            G
        $
        almost everywhere.
        \item 
        \label{thmitem::thm::local_hypergradients_imply_local_sobolev::item::G_is_again_a_hypergradient}
        $
            2^{ 2k \abs{ \alpha } } 
            G
            \in 
            \bD^{s- \abs{ \alpha}}_{ \lambda}(
                \partial^{ \alpha } f
            ).
        $ 
    \end{enumerate}
\end{theorem}
\begin{proof}
    First of all, let us note that there is nothing to prove in the case when $s \in \intoo{0,1}$. 
    For~this reason, let~us assume that $s \ge 1$.

    We will prove the theorem by induction. 
    To that end, let us first note that properties 
    \ref{thmitem::thm::local_hypergradients_imply_local_sobolev::item::existence} and \ref{thmitem::thm::local_hypergradients_imply_local_sobolev::item::G_is_again_a_hypergradient}
    are also satisfied for $\alpha = 0$. 
    (%
    	The existence of $\partial^{0} f$ merely means that 
	    $\partial^{0}f = f$ is~an~element of $\locIntegrable{1}(\Omega)$.%
    )
    Now, suppose that 
    properties 
    \ref{thmitem::thm::local_hypergradients_imply_local_sobolev::item::existence} and \ref{thmitem::thm::local_hypergradients_imply_local_sobolev::item::G_is_again_a_hypergradient}
    are satisfied for all 
    $
        \alpha \in \bN_0^n
    $
    such that
    $
        0 \le \abs{ \alpha } \le m - 1
    $
    for some 
    $m \in \bN$ 
    with
    $m \le s$.
    We will prove that this implies that 
    properties 
    \ref{thmitem::thm::local_hypergradients_imply_local_sobolev::item::existence},
    \ref{thmitem::thm::local_hypergradients_imply_local_sobolev::item::bounded_by_G},
    and \ref{thmitem::thm::local_hypergradients_imply_local_sobolev::item::G_is_again_a_hypergradient}
    are satisfied for all 
    $
        \alpha \in \bN_0^n
    $
    such that
    $
        \abs{ \alpha } = m.
    $
    Let us note that within the proof,
    as a matter of convenience, whenever we will use a representative of 
    $ \partial^{ \alpha } f $,
    we will implicitly assume that the representative is finite everywhere.
    \begin{enumerate}[listparindent=\parindent    , leftmargin=0em, itemindent=1.75em, labelsep=0.5em]
        \item[\ref{thmitem::thm::local_hypergradients_imply_local_sobolev::item::existence}:]
        Let us fix 
        $
            \beta \in \bN_0^n
        $
        such that
        $
            0 \le \abs{ \beta } \le m - 1.
        $
        By the inductive hypothesis,
        we have
        $
        	2^{2k\card{\beta}}G
        	\in 
        	\bD_{\lambda}^{s-\card{\beta}}\del*{
        		\partial^{\beta} f 
        	}.
        $
        Since
        $ s- \abs{ \beta } \ge 1$,
        by \Cref{prop::embedding_into_lower_orders_normed}
        we have that 
		\begin{equation*}
            2^{k^2} \partitions{k}
            \del{
            	\abs{\partial^{\beta} f }
            	+
                2^{2k \card{\beta} } G
            }
            \in 
            \bD^1_{\lambda}\del*{
                \partial^{\beta} f
            }.			
		\end{equation*}
        Also, since
        $
            2^{k^2} \partitions{k}
            \del{
            	\abs{\partial^{\beta} f }
            	+
                2^{2k \card{\beta} }G
            }
            \in 
            \locIntegrable{1}(\Omega)
        $
        and 
        $
            \partial^{\beta} f
            \in 
            \locIntegrable{1}(\Omega),
        $
        we have that
        $
            \partial^{\beta} f 
            \in 
            \locHajlasz{1,1}(\Omega).
        $
        Therefore, by~%
        \Cref{lem::loc_hajlasz_implies_loc_sobolev}
        we conclude that
        $
            \partial^{\beta} f 
            \in 
            \locSobolev{1,1}(\Omega).
        $
        Since 
        $ \beta \in \bN_{0}^{n}$ 
        with
        $
            0 \le \abs{ \beta } \le m-1
        $
        is arbitrary, 
        we~conclude that
        $
            f 
            \in 
            \locSobolev{m,1}(\Omega)
        $
        and 
        $
            \partial^{\alpha} f
        $
        exists in a weak sense for all 
        $
            \alpha \in \bN_0^n
        $
        such that
        $
            \abs{ \alpha } \le m.
        $
        \item[\ref{thmitem::thm::local_hypergradients_imply_local_sobolev::item::bounded_by_G}:]
        Fix 
        $
            \alpha \in \bN_0^n
        $
        such that 
        $
            \abs{ \alpha } = m.
        $
        There are $\alpha_1, \, \ldots, \, \alpha_n \in \bN_0^n$ such that
        $
            \alpha = \del{ \alpha_1, \, \ldots, \, \alpha_n}.
        $
        Let 
        $
            \set*{ v_\ell }_{ \ell = 1 }^{ m }   
        $
        be the unit vectors in $\bR^n$ such that, 
        for all $i \in [n]$, exactly $\alpha_i$ of~%
        them are equal to $\be_i$, where 
        $
            \set*{ \be_i }_{ i=1 }^{ n }
        $
        is the canonical basis of $\bR^n$.
        Define
        $
            \bv = \set*{ v_L }_{ L \subseteq [m] }
        $
        by the formula
        \begin{equation*}
            \forall L \subseteq [m] 
            \qquad
            v_L
            \coloneq 
            \sum_{ \ell \in L }
                v_{ \ell }.
        \end{equation*}
        Fix $h > 0$.
        Note that
        \begin{equation}
        \label{proofeq::thm::local_hypergradients_imply_local_sobolev::part::b::bound_on_hv_L}
            \forall 
                L \subseteq [m] 
            \qquad
            \norm*{ hv_L }
            =
            \norm*{
                h
                \sum_{ \ell \in L }
                    v_{ \ell }
            }
            \le 
            h
            \sum_{ \ell \in L }
                \norm*{ v_{\ell } }
            =
            h
            \sum_{ \ell \in L }
                1
            =
            h \card{ L } 
            \le 
            hm.
        \end{equation}
        %
        By \Cref{cor::value_of_Poly_when_tuple_is_a_hyperparallelogram} we have
        \begin{equation}
        \label{proofeq::thm::local_hypergradients_imply_local_sobolev::part::b::value_of_P^m(hv)}
            \polygen{ \gen{m} }( h\bv )
            =
            2^{m(m-1)}
            \prod_{ \ell = 1 }^{ m }
                \norm{ h v_{ \ell } }
            =
            2^{m(m-1)} h^m.
        \end{equation}
        Moreover,
        \begin{align}
            \polyset{ \set*{m} }( h\bv )
        &=
            \sum_{ A = \emptyset }^{ [m-1] }
                \polyset[A]{\set*{m}}( h\bv )
        \notag
        \\
        &=
            \sum_{ A = \emptyset }^{ [m-1] }
                \norm*{
                    h
                    \sum_{\ell \in A \cup \set*{m} }
                        v_{\ell}
                    -
                    h
                    \sum_{\ell \in A }
                        v_{\ell}
                }
            =
            \sum_{ A = \emptyset }^{ [m-1] }
                h \norm{ v_m }
            =
            2^{m-1}h.
            \label{proofeq::thm::local_hypergradients_imply_local_sobolev::part::b::value_of_P^set(m)(hv)}
        \end{align}

        Next, define
        $
            \by
            =
            \set*{
                y_J
            }_{ J \subseteq [k] }
        $
        by the formula
        \begin{equation*}
            \forall L \subseteq [m]
            \quad 
            \forall I \subseteq [k] \setminus [m]
            \qquad 
            y_{L \cup I}
            \coloneq 
            \begin{cases}
                v_L &\text{if } \card{I} \text{ is even}, \\ 
                v_{L \div \set*{m} }
                &\text{otherwise}.
            \end{cases}
        \end{equation*}
        Note that 
        by 
        \eqref{proofeq::thm::local_hypergradients_imply_local_sobolev::part::b::bound_on_hv_L}
        we have
        \begin{equation}
        \label{proofeq::thm::local_hypergradients_imply_local_sobolev::part::b::bound_on_hv_J}
            \forall 
                J \subseteq [k]
            \qquad
            \norm{
                h y_J
            }
            \le 
            hm.
        \end{equation}
        Then, by 
        \Cref{lem::folding_hypercubes_in_normed_spaces}, 
        for $m < k$ we have
        \begin{equation*}
            \polygen{ \gen{s} }( h \by )
            \le 
            2^{k(k-m)}
            \partitions{k}
            \polygen{ \gen{m} }( h \bv )
            \del*{
                1 + \polyset{\set*{m}}( h \bv )^{k-m}
            }.
        \end{equation*}
        Therefore, if $m < k$ and $h \in \intoc{0,1}$, 
        thanks to 
        \eqref{proofeq::thm::local_hypergradients_imply_local_sobolev::part::b::value_of_P^m(hv)}
        and 
        \eqref{proofeq::thm::local_hypergradients_imply_local_sobolev::part::b::value_of_P^set(m)(hv)},
        we have
        \begin{align*}
            \polygen{ \gen{s} }(h \by)
            &\le
            2^{k(k-m)}
            \partitions{k}
            \polygen{ \gen{m} }(h \bv)
            \del*{
                1 + \polyset{ \set*{m} }(h \bv)^{k-m}
            }
        \\
            &=
            2^{k(k-m)}
            \partitions{k}
            \del*{
                2^{m(m-1)} h^m
            }
            \del*{
                1 
                + \del*{
                    2^{m-1} h    
                }^{k-m}
            }
        \\
            &\le 
            2^{k(k-m)}
            \partitions{k}
            \del*{
                2^{m(m-1)} h^m
            }
            \del*{
                2^{(m-1)(k-m)} 
                + 2^{(m-1)(k-m)} 
            }
        \\
            &=
            2^{k(k-m) + m(m-1) + (m-1)(k-m) + 1}
            \partitions{k}
            h^m
        \\
            &= 
            2^{k(k-1) + 1}
            \partitions{k}
            h^m
        \\
            &\le
            2^{k(k-1)+k-m}
            \partitions{k}
            h^m
        \\
            &=
            2^{k^2 - m}
            \partitions{k}
            h^m.
        \end{align*} 
        Also, if $m = k$, then $m = s$ and 
        $\by = \bv$, so by 
        \eqref{proofeq::thm::local_hypergradients_imply_local_sobolev::part::b::value_of_P^m(hv)}
        we have
        \begin{equation*}
            \polygen{ \gen{s} }(h \by)
            =
            \polygen{ \gen{m} }( h \bv)
            =
            2^{m(m-1)} h^m
            =
            2^{k^2 - m}
            h^m
            \le 
            2^{k^2 - m}
            \partitions{k}
            h^m,
        \end{equation*} 
        where the last inequality follows from the fact that  
        $
            1 \le \partitions{k}
        $
        since
        $1 \le m \le k$.
        Therefore, we see that as~long as 
        $m \le k$ and $h \in \intoc{0,1}$, we have
        \begin{equation}
        \label{proofeq::thm::local_hypergradients_imply_local_sobolev::part::b::estimate_for_P^s(hy)}
            \polygen{ \gen{s} }( h\by )
            \le 
            2^{k^2 - m}
            \partitions{k}
            h^m.
        \end{equation}

        Let $f'$ and $G'$ be such that
        $f = f'$ and $G = G'$ almost everywhere
        and 
        $G ' \in \bD^s(f')$. 
        (The existence of such functions is assured by 
        \Cref{lem::from_m_upper_to_upper(normed)}.)
        Note that
        $ f ' \in \locSobolev{m, 1}(\Omega)$
        since we have already proved in~part
        \ref{thmitem::thm::local_hypergradients_imply_local_sobolev::item::existence}
        that 
        $ f \in \locSobolev{m, 1}( \Omega)$.
        Also, we have $G' \in \locIntegrable{1}(\Omega)$.
        Let 
        \begin{equation*}
            F \coloneq 
            \lebesguePoints{ \partial^{ \alpha} f'}
            \cap \lebesguePoints{ G' }.	        
        \end{equation*}
        %
%
        Fix $x \in F$ and~let 
        $R > 0$ be such that 
        $
            \clball{x, 2R} \subseteq \Omega.
        $
        Fix $r \in \intoc{0,R}$,
        $ 
        	\eps \in 
        	\intoo*{ 
        		0,
        		\min\intoo*{ 1, \frac{r}{m+1} } 
        	},
        $
        and 
        $
            h \in 
            \intoc{
                0, \eps
            }.
        $
        Then,
        thanks to 
        \eqref{proofeq::thm::local_hypergradients_imply_local_sobolev::part::b::bound_on_hv_J},
        for every 
        $
            y \in \ball{x, r}
        $
        and $J \subseteq [k]$,
        we have
        \begin{equation*}
            \norm*{
                y + hy_J
                - x
            }
            \le 
            \norm*{ y - x }
            + \norm*{ hy_J }
            \le 
            r + hm 
            \le 
            2R.
        \end{equation*}
        In consequence, for all such $y$ and $J$
        we have
        $
            y + hy_J 
            \in 
            \clball{x, 2R}
            \subseteq
            \Omega.
        $
        
        Thus, using \Cref{lem::folding_hypercubes_in_groups},
        for all 
        $ y \in \ball{x, r}$ we have
        \begin{align*}
            2^{k-m}
            \abs*{
                \diff_{L = \emptyset}^{[m] }
                f' \del*{ y + h v_L}
            }
        	=
            \abs*{
                \diff_{ J = \emptyset }^{ [k] }
                    f' \del*{ y + h y_J}
            }
        &\le 
            \polygen{ \gen{s} }( y + h \by )
            \sum_{ J = \emptyset }^{ [k] } 
                G'  \del*{  y + h y_J }
        \\
        &=
            2^{k-m}
            \polygen{ \gen{s} }( y + h \by )
            \sum_{ L = \emptyset }^{ [m] } 
                G'  \del*{  y + h v_L },
        \end{align*}
        where the final equality follows from
        \Cref{lem::folding_hypercubes_with_a_sum}.
        Dividing the~resulting inequality by~$2^{k-m}$,
        and using the fact that 
        by \Cref{lem::invariant_transformations_of_poly}
        we have 
        $
            \polygen{ \gen{s} }( y + h \by )
            =
            \polygen{ \gen{s} }( h \by )
        $       
        since
        $a \mapsto y + a$ is an isometry of $\bR^n$,
        we get
        \begin{align*}
            \abs*{
                \diff_{L = \emptyset}^{[m] }
                f' \del*{ y + h v_L}
            }
            \le 
            \polygen{ \gen{s} }( y + h \by )
            \sum_{ L = \emptyset }^{ [m] } 
                G'  \del*{  y + h v_L } &
        \\
            =
            \polygen{ \gen{s} }( h \by )
            \sum_{ L = \emptyset }^{ [m] } 
                G'  \del*{  y + h v_L } &
            \le 
            2^{k^2 - m}
            \partitions{ k }
            h^m 
            \sum_{ L = \emptyset }^{ [m] } 
                G'  \del*{  y + h v_L },
        \end{align*}
        where the last inequality follows from 
        \eqref{proofeq::thm::local_hypergradients_imply_local_sobolev::part::b::estimate_for_P^s(hy)}.
        Recall that $x \in F$
        and we have
        \begin{equation*}
            F=
            \lebesguePoints{ \partial^{\alpha} f' }
            \cap \lebesguePoints{ G' }
            =
            \lebesguePoints{ \partial^{\alpha} f' }
            \cap \lebesguePoints{ 
                2^{k^2 - m} 
                \partitions{ k } G' 
            }.
        \end{equation*}
        Therefore, by applying 
        \Cref{lem::bounded_finite_differences_implies_pointwise_bound_on_derivative}
        to functions 
        $f'$ and 
        $ 
            2^{k^2 - m}
            \partitions{ k }
            G',
        $
        we obtain the estimate
        $
            \abs*{
                \partial^{ \alpha } f'(x)
            }
            \le 
            2^{k^2} 
            \partitions{ k }
            G'(x).
        $
        Since $x \in F$ is arbitrary, we have
        $
            \abs*{
                \partial^{ \alpha } f'
            }
            \le 
            2^{k^2} 
            \partitions{ k }
            G'
        $
        everywhere in $F$.
        
        Recall that 
        $f \in \locSobolev{m,1}(\Omega)$, 
        $G \in \locIntegrable{1}(\Omega) $,
        and we have
        $f = f'$ and $G = G'$
        almost everywhere.
        In consequence, $F$ is of full measure in 
        $\Omega$ and we have
        $
            \abs*{
                \partial^{ \alpha } f'
            }
            \le 
            2^{k^2} 
            \partitions{ k }
            G'
        $
        almost everywhere,
        hence also
        $
            \abs*{
                \partial^{ \alpha } f
            }
            \le 
            2^{k^2} 
            \partitions{ k }
            G
        $
        almost everywhere.
        Finally, note that 
        $ \alpha $ is an arbitrary element 
        of
        $ \bN_0^n $ 
        such that
        $ \abs{ \alpha } = m $,
        so~%
        \ref{thmitem::thm::local_hypergradients_imply_local_sobolev::item::bounded_by_G}~%
        is~true for all such $\alpha$.

        \item[\ref{thmitem::thm::local_hypergradients_imply_local_sobolev::item::G_is_again_a_hypergradient}:]
        Let 
        $ \alpha $, 
        $ \bv $, 
        $ f' $, $ G' $, 
        and $F$ be as in the above proof of 
        \ref{thmitem::thm::local_hypergradients_imply_local_sobolev::item::bounded_by_G}.
        Fix 
        \begin{equation*}
            \bx = \set*{x_I}_{ I \subseteq [k-m] }
            \subseteq 
            F
            =
            \lebesguePoints{ \partial^{\alpha} f'}
            \cap \lebesguePoints{ G' }.	        
        \end{equation*}
        Next, let
        \begin{equation*}
            D
            \coloneq 
            \min_{ I \subseteq [k-m] }
                \dist\del{ x_I, \partial \Omega }
            \quad 
            \text{and}
            \quad 
            R
            \coloneq
            \min\del*{
            	1, 
            	D/4
            }.
        \end{equation*}
        Note that
        \begin{equation*}
            \forall I \subseteq [k-m]
            \qquad 
            \clball{ x_I, 2R } \subseteq \Omega. 
        \end{equation*}
        Define functions
        $
            \wave{f}
            \colon 
            \ball{0, 2R}
            \to 
            \bR
        $
        and
        $
            \wave{G}
            \colon 
            \ball{0, 2R}
            \to 
            \intcc{0,\infty}        
        $
        by the formulas
        \begin{equation}
        \label{proofeq::thm::local_hypergradients_imply_local_sobolev::part::c::definition_of_wave_f_and_G}
            \forall w \in \ball{0, 2R}
        \qquad 
            \wave{f}(w)
            \coloneq 
            \diff_{ I = \emptyset }^{ [k-m] }
                f'\del*{ w + x_I }
        \, \, \, \text{and}\, \, \,  
            \wave{G}(w)
            \coloneq 
            \sum_{ I = \emptyset }^{ [k-m] }
                G'\del*{ w + x_I }.
        \end{equation}

        Since $f' = f$ almost everywhere and 
        $
            f \in \locSobolev{m, 1}( \Omega),
        $
        we have
        $
            \wave{f} \in \locSobolev{m, 1}\del*{
                \ball{0, 2R}
            }.
        $
        Since
        $ \abs{ \alpha } \le m $,
        we can define a~representative of~%
        $
        	\partial^{ \alpha } \wave{f}
        $
        by the formula
        \begin{equation*}
            \forall w \in \ball{0, 2R}
        \qquad 
            \partial^{ \alpha } \wave{f}(w)
            \coloneq 
            \diff_{ I = \emptyset }^{ [k-m] }
                \partial^{ \alpha } f'\del*{ w + x_I }.      
        \end{equation*}
        %
        Let us note that 
        if for a given
        $
        	w \in \ball{0, 2R}
        $
        we have
        $
        	w + x_I 
        	\in 
        	\lebesguePoints{
				\partial^{ \alpha} f'
			}
        $
        for all
        $ I \subseteq [k -m ] $,
        then
        $
        	w 
        	\in
			\lebesguePoints{
				\partial^{ \alpha} \wave{f}
			}.
        $
        In consequence, since
        $
            \bx = \set*{ x_I }_{ I \subseteq [k - m] }
            \subseteq 
            \lebesguePoints{ \partial^{\alpha} f' },
        $
        we have
        $
            0 \in \lebesguePoints{\partial^{\alpha}   \wave{f} }
        $
        and 
        \begin{equation}
        \label{proofeq::thm::local_hypergradients_imply_local_sobolev::part::c::partial_alpha_of_wave_f_at_0}
            \partial^{\alpha} \wave{f}(0)
            =
            \diff_{ I = \emptyset }^{ [ k-m ] }
                \partial^{ \alpha} f'(x_I).
        \end{equation}

        Similarly, as 
        $G' = G$ almost everywhere
        and $G \in \locIntegrable{1}(\Omega)$,
        we have
        $ 
            \wave{G} \in \locIntegrable{1}\del{
                \ball{0,2R}
            }.
        $
        Moreover, if for a given 
        $ w \in \ball{0,2R}$
        we have
        $
            w+x_I \in \lebesguePoints{G'}
        $
        for~all~$I \subseteq [k-m]$,
        then 
        $
            w \in \lebesguePoints{ \wave{G} }.
        $
        Thus, since 
        $
            \bx = \set*{ x_I }_{ I \subseteq [k-m] }
            \subseteq 
            \lebesguePoints{ G' },
        $
        we have
        $
            0 \in \lebesguePoints{ \wave{G} }.
        $
        Hence,
        \begin{equation*}
            \wave{f} \in \locSobolev{m,1}\del{
                \ball{0,2R}
            },
        \quad
            \wave{G} \in \locIntegrable{1}\del{
                \ball{0,2R}
            },
        \quad \text{and} \quad
            0 \in
            \lebesguePoints{ 
            	\partial^{\alpha} \wave{f} 
            }
            \cap 
            \lebesguePoints{ \wave{G} }.
        \end{equation*}

        Fix 
        $ r \in \intoc*{ 0, R }$,
        $ \eps \in \intoo*{0, \frac{r}{m+1}}$, 
        and
        $h \in \intoc*{0, \eps}$.
        Let us define
        $
            \by^h = \set*{ y_J^h }_{ 
            	J \subseteq [k] 
            }
        $
        by the formula
        \begin{equation*}
            \forall L \subseteq [m]
        \quad 
        	\forall I \subseteq [k-m]
        \qquad 
            y_{L \cup (m + I) }^h
            \coloneq 
            h v_L
            +
            x_I.
        \end{equation*}
        Then,
        by 
        \Cref{cor::tuple_part_tuple_part_vectors}
        and the fact that
        $\norm{ v_{\ell} } = 1$ 
        for all $\ell \in [m]$, 
        \begin{equation*}
            \polygen{ \gen{ s } }\del*{ \by^h }
            \le
            2^{ (2k-1)m } 
            \polygen{ \gen{ s-m } }( \bx) 
            \prod_{ \ell = 1 }^{ m } 
                \norm{ h v_{\ell} }
            =
            2^{ (2k-1)m } h^m
            \polygen{ \gen{ s-m } }( \bx). 
        \end{equation*}
        Therefore, since 
        $
            a \mapsto y + a
        $
        is an isometry of $\bR^n$,
        by \Cref{lem::invariant_transformations_of_poly} we have
        \begin{equation}
        \label{proofeq::thm::local_hypergradients_imply_local_sobolev::part::c::estimate_for_P^s(w+y^h)}
            \forall w \in \bR^n
        \qquad 
            \polygen{ \gen{s} }\del*{
                w + \by^h
            }
            \le
            2^{ (2k-1)m }
            h^m
            \polygen{ \gen{s-m} }(\bx).
        \end{equation}
        Note that, 
        for all 
        $ w \in \ball{0, r} $,
        $ L \subseteq [m] $,
        and
        $I \subseteq [k-m]$,
        we have
        \begin{align*}
            \norm*{
                w + y_{L \cup (m+I)}^h
                - x_I
            }
            =
            \norm*{
                w + hv_L + x_I
                - x_I
            }
        &
        \\
            \le 
            \norm{ w } 
            + \norm{ h v_L } 
        &\le 
            r + hm
        	\le
        	r + \eps m
        	<
        	2r
        	\le 
        	2R,
        \end{align*}
        so, for all such $w$, $L$, and $I$,
        \begin{equation*}
            w + y_{L \cup (m+I)}^h
            \in
            \clball{ x_I, 2R}
            \subseteq 
            \Omega
        \quad \text{and} \quad 
        	w + hv_L \in \ball{0,2R}.
        \end{equation*}

        Fix $w \in \ball{0,r}$.
        Since 
        $ G' \in \bD^s( f' )$, we have
        \begin{equation}
        \label{proofeq::thm::local_hypergradients_imply_local_sobolev::part::c::estimate_between_f'_and_G'}
            \abs*{
                \diff_{ J = \emptyset }^{ [k] }
                    f'\del*{ w + y_J^h }
            }
            \le 
            \polygen{ \gen{s} }\del*{ w + \by^h }
            \sum_{ J = \emptyset }^{ [k] }
                G'\del*{ w + y_J^h }.
        \end{equation}
        Let us notice that,
        by 
        \Cref{lem::diff_[k]_to_diff_[m]_diff_[k-m]},
        we have
        \begin{align*}
            \abs*{
                \diff_{ J = \emptyset }^{ [k] }
                    f'\del*{ w + y_J^h }
            }
        	=
            \abs*{
                \diff_{ L = \emptyset }^{ [m] }
                    \diff_{ I = \emptyset }^{ [k-m] }
                        f'\del*{ w + y_{L \cup (m+I)}^h }
            }
        &
        \\
        	=
            \abs*{
                \diff_{ L = \emptyset }^{ [m] }
                    \diff_{ I = \emptyset }^{ [k-m] }
                        f'\del*{ w + hv_L + x_I}
            }
        &=
            \abs*{
                \diff_{ L = \emptyset }^{ [m] }
                    \wave{f}\del*{ w + hv_L}
            },
        \end{align*}
        and, since 
        $ (L, I) \mapsto L \cup (m+I)$
        is a bijection from 
        $ 2^{ [m] } \times 2^{ [k-m] }$
        to 
        $ 2^{[k] }$,
        \begin{align*}
            \sum_{ J = \emptyset }^{ [k] }
                G'\del*{ w + y_J^h }
        	=
            \sum_{ L = \emptyset }^{ [m] }
                \sum_{ I = \emptyset }^{ [ k-m ] }
                    G'\del*{ w + y_{L \cup (m+I) }^h }
        &
        \\    
        	=
            \sum_{ L = \emptyset }^{ [m] }
                \sum_{ I = \emptyset }^{ [ k-m ] }
                    G'\del*{ w + hv_L + x_I}
        &=
            \sum_{ L = \emptyset }^{ [m] }
                \wave{G}\del*{ w + hv_L }.
        \end{align*}

        Therefore, 
		applying the estimate present in        
        \eqref{proofeq::thm::local_hypergradients_imply_local_sobolev::part::c::estimate_for_P^s(w+y^h)}
        to
        \eqref{proofeq::thm::local_hypergradients_imply_local_sobolev::part::c::estimate_between_f'_and_G'}
        and noting that
        $ h \in \intoc{0,\eps}  $
        as well as
        $ w \in \ball{0,r} $
        are arbitrary,
        we obtain that 
        \begin{equation*}
            \forall w \! \in \! \ball{0,r}
        \quad
            \forall h \! \in \! \intoc{0,\eps} 
        \qquad 
            \abs*{
                \diff_{ L = \emptyset }^{ [m] }
                    \wave{f}\del*{ w + hv_L}
            }
            \le 
            2^{(2k-1)m}
            h^m
            \polygen{ \gen{ s-m } }(\bx) 
            \sum_{ L = \emptyset }^{ [m] }
                \wave{G}\del*{ w + hv_L }.
        \end{equation*}

        Let us recall that
        $
            0 \in 
            \lebesguePoints{ \partial^{\alpha} \wave{f}}
            \cap \lebesguePoints{ \wave{G} }.
        $
        Hence, 
        since 
        $$
            \lebesguePoints{ \wave{G} }
            =
            \lebesguePoints{ 
                2^{(2k-1)m} 
                \polygen{ \gen{s-m} }(\bx) 
                \wave{G} 
            },
        $$
        by~%
        \Cref{lem::bounded_finite_differences_implies_pointwise_bound_on_derivative}
        we have 
        $
            \abs*{
                \partial^{\alpha} \wave{f}(0)
            }
            \le 
            2^{2km}
            \polygen{ \gen{s-m} }(\bx)
            \wave{G}(0).
        $
        Thus, thanks to 
        \eqref{proofeq::thm::local_hypergradients_imply_local_sobolev::part::c::partial_alpha_of_wave_f_at_0}
        and 
        \eqref{proofeq::thm::local_hypergradients_imply_local_sobolev::part::c::definition_of_wave_f_and_G},
        \begin{equation*}
            \abs*{
                \diff_{ I = \emptyset }^{ [k-m] }
                    \partial^{\alpha} f'(x_I)
            }
            =
            \abs*{ \partial^{\alpha} \wave{f}(0) }
            \le 
            2^{2km}
            \polygen{ \gen{s-m} }(\bx)
            \wave{G}(0)
            =
            2^{2km}
            \polygen{ \gen{s-m} }(\bx)
            \sum_{ I = \emptyset }^{ [k-m] }
                G'(x_I).
        \end{equation*}
        As
        $
            \bx = \set*{ x_I }_{ I \subseteq [k-m] }
            \subseteq 
            F
        $
        is arbitrary,
        the above inequality implies that
        $
            2^{2km}G'
            \in 
            \bD_{ \lambda }^{ s-m }\del{ 
                \partial^{ \alpha }( f' )
            }.
        $
        Hence, since $f = f'$
        and $G = G'$ almost everywhere, 
        we have
        $
            2^{2km}G
            \in 
            \bD_{ \lambda }^{ s-m }\del{ 
                \partial^{ \alpha }( f )
            }.
        $
        Finally, note that 
        $ \alpha $ is an arbitrary element 
        of
        $ \bN_0^n $ 
        such that
        $ \abs{ \alpha } = m $,
        so 
        \ref{thmitem::thm::local_hypergradients_imply_local_sobolev::item::G_is_again_a_hypergradient}
        is true for all such $\alpha$.
    \end{enumerate}
    We have proved that, 
    for all 
    $m \in \bN$ such that $m \le s$,
    if 
    \ref{thmitem::thm::local_hypergradients_imply_local_sobolev::item::existence}
    and \ref{thmitem::thm::local_hypergradients_imply_local_sobolev::item::G_is_again_a_hypergradient}
    are satisfied for all 
    $
        \alpha \in \bN_0^n
    $
    such that
    $
        \abs{ \alpha } \le m-1,
    $
    then 
    \ref{thmitem::thm::local_hypergradients_imply_local_sobolev::item::existence},
    \ref{thmitem::thm::local_hypergradients_imply_local_sobolev::item::bounded_by_G},
    and \ref{thmitem::thm::local_hypergradients_imply_local_sobolev::item::G_is_again_a_hypergradient}
    are satisfied for all 
    $
        \alpha \in \bN_0^n
    $
    such that
    $
        \abs{ \alpha } = m.
    $
    Since 
    \ref{thmitem::thm::local_hypergradients_imply_local_sobolev::item::existence}
    and \ref{thmitem::thm::local_hypergradients_imply_local_sobolev::item::G_is_again_a_hypergradient}
    are true when $\alpha = 0$,
    by induction we can conclude that 
    \ref{thmitem::thm::local_hypergradients_imply_local_sobolev::item::existence},
    \ref{thmitem::thm::local_hypergradients_imply_local_sobolev::item::bounded_by_G},
    and \ref{thmitem::thm::local_hypergradients_imply_local_sobolev::item::G_is_again_a_hypergradient}
    are true for all 
    $
        \alpha \in \bN_0^n
    $ 
    such that 
    $ 0 < \abs{ \alpha } \le s$,
    which ends the proof.
\end{proof}
%
%
\subsection{Characterizations of higher-order Sobolev spaces}
In this section, we will prove a characterization of higher-order Sobolev spaces based on Banach function spaces.
However, we will need some definitions in order to~state the result in full generality.
\begin{definition}[Banach function space]
	Let 
	$ \del*{ X, \Sigma, \measureAlone } $ 
	be a measure space.
	We will say that 
	$ 
		\del*{
			\mc{F}\del*{ X, \Sigma, \measureAlone },
			\normAlone[{\mc{F}\del*{ X, \Sigma, \measureAlone }}]
		}
	$
	is a \emph{Banach function space}
	(over $ \del*{ X, \Sigma, \measureAlone } $)
	if 
	\begin{itemize}[leftmargin=1.5em,itemsep=0.25em, parsep=0.25em,topsep=0.25em]
	\item 
		$ \mc{F}\del*{ X, \Sigma, \measureAlone }$
		is a (vector) subspace of
		$
			\integrable{0}\del*{
				X, \Sigma, \measureAlone
			},
		$
		the space of measurable functions 
		$ f \colon X \to \bR $ 
		(modulo the relation of equality
		$ \measureAlone $-almost everywhere).
	\item
		$ 
			\del*{
				\mc{F}\del*{ X, \Sigma, \measureAlone },
				\normAlone[{\mc{F}\del*{ X, \Sigma, \measureAlone }}]
			}
		$
		is a Banach space.
	\item 
		$
			\forall E \in \Sigma 
		\qquad 
			\measure{E} < \infty 
			\implies 
			\forall f \in \mc{F}(X)
		\quad 
			\integral{E}{ \abs{f} }{\measureAlone}
			<
			\infty.
		$
	\item 
		$ \mc{F}\del*{ X, \Sigma, \measureAlone }$
		satisfies the \emph{ideal property},
		that is, 
		for all
		$
			f, h \in \integrable{0}\del*{ X, \Sigma, \measureAlone },
		$
		we have
		\begin{equation*}
			\del*{
				\abs{ f } \le \abs{ h } \text{ $\measureAlone$-a.e.} \,\,\, \text{and} \,\,\,  h \in \mc{F}(X)
			}
		\, \, \implies \, \, 
			\del*{
				f \in \mc{F}(X) \,\,\, \text{and} \,\,\, 
				\norm{ f }_{ \mc{F}(X) }
				\le 
				\norm{ h }_{ \mc{F}(X) }
			}.
		\end{equation*}
	\end{itemize}
	%
	Much like we already did above,
	we will often abbreviate the notation by writing 
	$\mc{F}(X)$ 
	or 
	$\mc{F}\del*{X, \measureAlone}$
	instead of 
	$\mc{F}\del*{ X, \Sigma, \measureAlone }$
	if~the~measure $\measureAlone$ or the $\sigma$-field $\Sigma$ should be inferrable from the context 
	(or irrelevant to~the~discussion).
	Also, 
	if $\del*{ Y, \Sigma, \measureAlone} $
	is a measure space and 
	$ 
		\del*{ 
			X, \Sigma \rvert_X, \measureAlone \rvert_X 
		}
	$
	is its subspace,\footnote{
	i.e., 
	$ X $ is measurable and
	we have
	$
		\Sigma \rvert_X 
		\coloneq 
		\setc*{
			A \cap X
		}{
			A  \in \Sigma 
		}
	$ 
	and
	$
		\measureAlone \rvert_X\del*{ 
			A \cap X
		}
		=
		\measure{ A \cap X }
	$
	for all $A \in \Sigma $.
	}
	then by writing 
	$\mc{F}(X)$ and $\mc{F}(Y)$ we~imply that 
	for all $f \in \mc{F}(Y)$ we have 
	$ f \rvert_{X} \in \mc{F}(X)$
	with 
	$
		\norm{ f \rvert_{X} }_{ \mc{F}(X) }
		\le 
		\norm{ f }_{ \mc{F}(Y) }.
	$
	Moreover, 
	if 
	$\del*{ V, \normAlone[V] }$
	and 
	$\del*{ W, \normAlone[W] }$
	are normed spaces such that
	$ V \subseteq \mc{F}(X) $
	and
	$ W \subseteq \mc{F}(Y) $,
	then we will say that 
	$ \ms{E} \colon V \to W $
	is~a~\emph{bounded extension operator} if~%
	there exists $C > 0$ such that
	\begin{equation*}
		\forall f \in V
	\qquad 
		\norm{ \ms{E}(f) }_{W}
		\le 
		C \norm{ f }_{ V }
	\quad \text{and} \quad 
		\ms{E}(f) \rvert_{X} = f.
	\end{equation*}
	Note that we do not require that $\ms{E}$ is linear.
	Furthermore, if 
	$ \Omega \subseteq \bR^n$ is open, 
	then,
	unless specified otherwise,
	by~writing
	$
		\mc{F}\del*{ \Omega },
	$
	we will mean a Banach function space defined over the Lebesgue measure space
	$ \del*{ \Omega, \lambda} $.
	Finally, as a matter of convenience,
	we extend the norm
	$ \normAlone[\mc{F}\del*{X}]$
	to all elements of 
	$ \integrable{0}\del*{X}$
	by putting
	$
		\norm{ f }_{ \mc{F}\del*{X} }
		\coloneq 
		\infty 
	$
	for all
	$ f \in \integrable{0}\del*{X} \setminus \mc{F}\del*{X} $.
\end{definition}
\begin{definition}[$\multipoint{s}\mc{F}(X)$ space]
\label{def::multipointwise_Banach_function_space}
	Let $n \in \bN$, $s \ge 0$, 
	$ X \subseteq \bR^n$,
	$\measureAlone$ be a measure on $X$,
	and 
	$ 
		\mc{F}\del*{X, \measureAlone}
		=
		\mc{F}(X)
	$ 
	be~a~Banach function space.
	We introduce the family
	\begin{align*}
		\multipoint{s}\mc{F}(X)
	&\coloneq 
		\setc*{
			f \in \mc{F}(X)
		}{
			\text{%
				there exists 
				$ G \in \bD_{\measureAlone}^{s}(f) $
				such that
				$ G \in \mc{F}(X) $%
			}
		}.
	\end{align*}
	We endow it with the following seminorm and norm:
	\begin{equation*}
		\forall f \in \multipoint{s}(X)
	\qquad 
		\seminorm{
			f
		}_{ \multipoint{s}\mc{F}(X) }
		\coloneq 
		\inf_{ G \in \bD^{s}_{\measureAlone}(f) }
			\norm{ G }_{ \mc{F}(X) }
	\quad \text{and} \quad
		\norm{ f }_{ \multipoint{s}\mc{F}(X) }
		\coloneq 
		\norm{ f }_{ \mc{F}(X) }
		+
		\seminorm{ f }_{ \multipoint{s}\mc{F}(X) }.
	\end{equation*}
	Note that if 
	$ s = 0 $, 
	then
	$
		\multipoint{s}\mc{F}(X)
		=
		\mc{F}(X),
	$
	$
		\seminormAlone[
			{\multipoint{s}\mc{F}(X)}
		]
		=
		\normAlone[{\mc{F}(X)}],
	$
	and 
	$
		\normAlone[{\multipoint{s}\mc{F}(X)}]
		=
		2\normAlone[{\mc{F}(X)}].
	$
\end{definition}
\begin{definition}[$\sobolev{k}\mc{F}(\Omega)$ space]
\label{def::sobolev_Banach_function_space}
	Let $n\in \bN$, $k \in \bN_0$, $\Omega \subseteq \bR^n$
	be open, 
	and 
    $
    	\mc{F}\del*{ \Omega, \lambda}
    	=
    	\mc{F}(\Omega)
    $ 
	be a Banach function space.
	We introduce the family
	\begin{align*}
		\sobolev{k}\mc{F}(\Omega)
		&\coloneq 
		\setc*{ 
			f \in \locSobolev{k,1}(\Omega)
		}
		{ 
			\text{%
				$ 
					\partial^{\alpha} f
					\in 
					\mc{F}(\Omega)
				$ 
				for all 
				$\alpha \in \bN_0^n$
				such that
				$ \abs{ \alpha } \le k $%
			}
		}.
	\end{align*}
	We endow it with the following seminorm and norm:
	\begin{align*}
		\forall f \in \sobolev{k}\mc{F}(\Omega)
	\qquad 
		\seminorm{f}_{                
		    \sobolev{k}\mc{F}(\Omega) 
		}
		\coloneq 
		\sum_{ j = 1 }^{ k }
			\norm*{ \nabla^j f }_{ \mc{F}\del*{ \Omega} }
	\,\,\,\, \text{and} \,\,\,\,
		\norm{f}_{ 
			\sobolev{k}\mc{F}(\Omega) 
		}
		\coloneq 
		\norm{ f }_{ \mc{F}(\Omega) }
		+
		\seminorm{ f }_{
			\sobolev{k}\mc{F}(\Omega) 
		}.
	\end{align*}
	Note that if $k = 0$, then
	$
		\sobolev{k}\mc{F}(\Omega)
		=
		\mc{F}(\Omega),
	$
	$
		\seminormAlone[{
			\sobolev{k}\mc{F}(\Omega)
		}]
		=
		0,
	$
	and 
	$
		\normAlone[{
			\sobolev{k}\mc{F}(\Omega)
		}]
		=
		\normAlone[{ \mc{F}(\Omega) }].
	$
\end{definition}
\begin{definition}[$\hajlasz{s}\mc{F}(\Omega)$ space]
\label{def::hajlasz_Banach_function_space}
	Let $n \in \bN$, $s > 0$, 
	$ \Omega \subseteq \bR^n$ be open,
	and 
    $
    	\mc{F}\del*{ \Omega, \lambda}
    	=
    	\mc{F}(\Omega)
    $ 
    be a Banach function space.
	Let $k \in \bN$ be such that $ s \in \intoc{k-1, k}$.
	We introduce the family
	\begin{equation*}
		\hajlasz{s}\mc{F}(\Omega)
		\coloneq 
		\setc*{
			f \in \sobolev{k-1}\mc{F}( \Omega )
		}{
			\text{%
				there exists 
				$ 
					g \in \bD_{\lambda}^{s-k+1}\del*{
						\nabla^{k-1} f
					} 
				$
				such that
				$ g \in \mc{F}(X) $%
			}
		}.
	\end{equation*}
	We endow it with the following seminorm and norm:
	\begin{align*}
		\forall f \in \hajlasz{s}\mc{F}(\Omega)
	\qquad 
		\seminorm{
			f
		}_{ \hajlasz{s}\mc{F}(\Omega) }
	&\coloneq 
		\sum_{ j = 1 }^{ k-1 }
			\norm*{ \nabla^j f }_{ \mc{F}\del*{ \Omega} }
		+
		\inf_{ G \in \bD^{s-k+1}_{\lambda}\del*{ \nabla^{k-1} f } }
			\norm{ G }_{ \mc{F}(\Omega) }
	\intertext{and}
		\forall f \in \hajlasz{s}\mc{F}(\Omega)
	\qquad 
		\norm{ f }_{ \hajlasz{s}\mc{F}(\Omega) }
	&\coloneq 
		\norm{ f }_{ \mc{F}(\Omega) }
		+
		\seminorm{ f }_{ \hajlasz{s}\mc{F}(\Omega) }.
	\end{align*}
\end{definition}
Before we state our characterization, we will make a small remark listing the results that directly follow from the above-stated definitions.
\begin{remark}
\label{rem::up_to_first-order_characterization}
	Let $n \in \bN$.
    Let $ \Omega \subseteq \bR^n $
    be open and 
    $
    	\mc{F}\del*{ \Omega, \lambda}
    	=
    	\mc{F}(\Omega)
    $ 
    be a Banach function space.
    Then for~all $s \in \intoc{0,1}$ we have
    $
    	\multipoint{s}\mc{F}(\Omega)
		=
		\hajlasz{s}\mc{F}(\Omega)
    $
    with equality of both seminorms and norms.
	Also,  
	\begin{equation*}
		\multipoint{0}\mc{F}(\Omega)
		=
		\sobolev{0}\mc{F}(\Omega)
		=
		\mc{F}(\Omega)
	\quad \text{with} \quad
		\normAlone[{\multipoint{0}\mc{F}(\Omega)}]
		=
		2\normAlone[{\sobolev{0}\mc{F}(\Omega)}]
		=
		2\normAlone[{\mc{F}(\Omega)}].
	\end{equation*}
\end{remark}
\begin{theorem}
\label{thm::higher-order_characterization}
    Let $n \in \bN$ and $s \ge 1$.
    Let $ \Omega \subseteq \bR^n $
    be open and 
    $
    	\mc{F}\del*{ \Omega, \lambda}
    	=
    	\mc{F}(\Omega)
    $ 
    be~a~Banach function space. 
    Then
	\begin{enumerate}[label=(\alph*)]
	\item 
		$
			\multipoint{s}\mc{F}(\Omega)
			\hookrightarrow
			\hajlasz{s}\mc{F}(\Omega).
		$
	\item
		If the Hardy-Littlewood maximal function 
		$ \mathsf{M} $ is bounded on $\mc{F}\del*{ \bR^n }$ and
		there exists a bounded extension operator
		$
			\mathsf{E}
			\colon 
			\hajlasz{s}\mc{F}(\Omega)
			\to
			\hajlasz{s}\mc{F}\del*{ \bR^n },
		$
		then
		$
			\hajlasz{s}\mc{F}(\Omega)
			\hookrightarrow
			\multipoint{s}\mc{F}(\Omega).
		$
	\end{enumerate}
    Furthermore, if $s =k$ for some $k\in \bN$, 
    then the above statements are true with 
    $ \hajlasz{s} $ being replaced by 
    $ \sobolev{k} $.
\end{theorem}
\begin{remark}
\label{rem::thm::higher-order_characterization_case_s=1}
	Let us note that by
	\Cref{rem::up_to_first-order_characterization},
	when $s = 1$,
	we have
	$ 
		\multipoint{s}\mc{F}(\Omega) 
		= 
		\hajlasz{s}\mc{F}(\Omega)
	$
	with equality of norms
	for any $\Omega$ and $\mc{F}(\Omega)$.
	Nevertheless, we have included the case $s=1$
	in the above theorem since it~would produce correct statements when 
	$
		\hajlasz{1}\mc{F}(\Omega)
	$
	is replaced by~%
	$
		\sobolev{1}\mc{F}(\Omega).
	$
\end{remark}
\begin{proof}[Proof of \Cref{thm::higher-order_characterization}]
	In light of 
	\Cref{rem::thm::higher-order_characterization_case_s=1},
	we only need to prove the theorem for $s > 1$.
	\begin{enumerate}[label=(\alph*), listparindent=\parindent, leftmargin=0.5em, itemindent=1.5em,itemsep=0.25em, parsep=0.25em, topsep=0.25em]
	\item 
		Fix 
		$ f \in \multipoint{s}\mc{F}(\Omega) $
		and 
		$ G \in \bD_{\lambda}^s(f) $
		with
		$ G \in \mc{F}(\Omega) $.
		Let $k \in \bN$ be such that
		$ s \in \intoc{k-1, k}$.
		By 
		\Cref{thm::local_hypergradients_imply_local_sobolev}
		we know that, 
		for all 
		$ \alpha \in \bN_0^n$ such that 
		$ 0 < \abs{ \alpha } \le s $,
		we have that
		\begin{itemize}
		\item
			$ \partial^{ \alpha} f $ exists in a weak sense;
			in particular,
			$
				\partial^{ \alpha} f
				\in 
				\locIntegrable{1}(\Omega).
			$
		\item 
			$    
				\abs{ \partial^{\alpha} f }
				\le 
				2^{k^2} \partitions{k} G
			$
			almost everywhere.
		\item 
			$
				2^{k^2} \partitions{k} G
				\in 
				\bD_{\lambda}^{s - \abs{\alpha} }\del*{
					\partial^{ \alpha} f
				}.
			$
		\end{itemize}
		Therefore, 
		$ f \in \locSobolev{k-1,1}(\Omega)$.
		Also,
		for all $j \in [k-1]$, there is 
		$C_j > 0$ that depends solely on~$k$
		such that 
		$
			\norm*{
				\nabla^j f
			}_{ \mc{F}(\Omega) }
			\le 
			C_j
			\norm*{ G }_{ \mc{F}(\Omega) }.
		$
		Moreover, there is $\wave{C} > 0$ that only depends on~$k$
		and such that 
		$
			\wave{C}G
			\in 
			\bD_{\lambda}^{s-k+1}\del*{ \nabla^{k-1} f }.				
		$
		
		Thus,
		since 
		$ f \in \mc{F}(\Omega)$
		as
		$
			f \in \multipoint{s}\mc{F}(\Omega)
		$,
		it~follows that
		$ f \in \hajlasz{s}\mc{F}(\Omega) $
		and 
		\begin{align*}
			\seminorm{ f }_{
				\hajlasz{s}\mc{F}(\Omega) 
			}
			\le 
			\sum_{j=1}^{ k-1}
				\norm{ \nabla^j f }_{   
					\mc{F}(\Omega) 
				}
			+
			\inf_{ 
				g \in 
				\bD^{s-k+1}_{\lambda}\del*{
				    \nabla^{k-1} f 
				} 
			}
				\norm{ g }_{ \mc{F}(\Omega) }
		&\\
			\le
			\sum_{j=1}^{ k-1}
				C_j
				\norm{ G }_{ \mc{F}(\Omega) }
			+
			\norm*{ \wave{C}G }_{ \mc{F}(\Omega) }
		&=
			C
			\norm{ G }_{ \mc{F}(\Omega) },
		\end{align*}
		where 
		$
			C \coloneqq 
			\wave{C} 
			+ \sum_{j=1}^{ k-1}
				C_j. 
		$
		Passing in the resulting inequality
		$
			\seminorm{ f }_{
				\hajlasz{s}\mc{F}(\Omega) 
			}
			\le 
			C
			\norm{ G }_{ \mc{F}(\Omega) }
		$
		to~the~infimum over 
		all 
		$ G \in \bD^s_{\lambda}(f)$,		
		we get
		$
			\seminorm{ f }_{ 
				\hajlasz{s}\mc{F}(\Omega) 
			}
			\le 
			C
			\seminorm{ f }_{     
				\multipoint{s}\mc{F}(\Omega) 
			}.
		$
		Also, we~have that
		\begin{equation*}
			\norm{ f }_{
				\hajlasz{s}\mc{F}(\Omega) 
			}
			=
			\norm{ f }_{ \mc{F}(\Omega) }
			+
			\seminorm{ f }_{
				\hajlasz{s}\mc{F}(\Omega) 
			}
			\le 
			\norm{ f }_{ \mc{F}(\Omega) }
			+
			C
			\seminorm{ f }_{
				\multipoint{s}\mc{F}(\Omega) 
			}
			\le 
			\del*{1+C}
			\norm{ f}_{
				\multipoint{s}\mc{F}(\Omega) 
			}.
		\end{equation*}
		In consequence, 
		$
			\multipoint{s}\mc{F}(\Omega)
			\hookrightarrow
			\hajlasz{s}\mc{F}(\Omega).
		$
	\item
		Let us suppose that there is 
		a bounded extension operator
		$
			\mathsf{E}
			\colon 
			\hajlasz{s}\mc{F}(\Omega)
			\to 
			\hajlasz{s}\mc{F}\del*{ \bR^n }.
		$
		Then 
		there exists $C_1 > 0$ such that
		for all $f \in \hajlasz{s}\mc{F}(\Omega)$ 
		we have
		$
			\mathsf{E}(f) \rvert_{\Omega} = f
		$
		and the estimate
		$
			\norm{ \mathsf{E}(f) \rvert_{\Omega} }_{   
				\hajlasz{s}\mc{F}\del*{ \bR^n } 
			}
			\le 
			C_1 \norm{f}_{ \hajlasz{s}\mc{F}(\Omega) }.
		$
		We~also~assume that $\mathsf{M}$
		is a bounded operator on
		$ \mc{F}\del*{ \bR^n }$,
		so there exists $C_2  > 0$ such that for all $f \in \mc{F}\del*{ \bR^n }$
		we~have
		$
			\norm{ \maxFun{f} }_{   
				\mc{F}\del*{ \bR^n } 
			}
			\le 
			C_2 \norm{f}_{ \mc{F}\del*{ \bR^n } }.
		$
		
		Fix 
		$
			f \in \hajlasz{s}\mc{F}(\Omega)
		$,
		$ \eps > 0 $,
		and let 
		$ k \in \bN $
		be such that
		$ s \in \intoc{k-1,k}$.		
		Then, there exists 
		$
			g \in \bD_{ \lambda }^{s-k+1}\del*{ \nabla^{k-1} \mathsf{E}(f) }
		$
		such that
		$
			\norm{ g }_{ \mc{F}\del*{ \bR^n } }
			\le 
			\norm{ \mathsf{E}(f) }_{ 
				\hajlasz{s}\mc{F}\del*{ \bR^n }
			}
			+
			\eps.
		$
		By 
		\Cref{cor::hypergradients_from_multipointwise_bound_for_W^k_loc}
		we know that 
		$ G \colon \bR^n \to \intcc{0, \infty }$
		defined by~the~formula 
		\begin{equation*}
			\forall x \in \bR^n
		\qquad 
			G(x)
			\coloneq 
			\widehat{C}_{n,k}
			\del*{
				\sum_{j = 1}^{ k-1 }
					\maxFun[k+j-1]{
						\nabla^j \mathsf{E}(f)
					}(x)
				+
				\maxFun[2k-2]{
					g
				}(x)
			}
		\end{equation*}
		is an element of 
		$ \bD^s_{\lambda}\del*{ \ms{E}(f) }$.
		In consequence,
		$
			G \rvert_{\Omega} 
			\in 
			\bD_{\lambda}^s(f).
		$
		Also, from the definition of~$G$ and~%
		the~boundedness of 
		$\mathsf{M}$ on $\mc{F}\del*{ \bR^n }$ 
		we have 
		$ G \in \mc{F}\del*{ \bR^n }$,
		hence
		$ G \rvert_{\Omega} \in \mc{F}(\Omega)$.
		In~consequence,
		since 
		$ f \in \mc{F}(\Omega)$
		as~%
		$f \in \hajlasz{s}\mc{F}(\Omega)$,
		we get that
		$
			f \in \multipoint{s}\mc{F}(\Omega).
		$
		
		Moreover, since
		\begin{align*}
			\norm{ G\rvert_{\Omega} }_{ \mc{F}(\Omega) }
			\le
			\norm{ G }_{ \mc{F}\del*{ \bR^n } }
		&\le 
			\widehat{C}_{n,k}
			\del*{
				\sum_{j=1}^{k-1}
					\norm*{
						\maxFun[k+j-1]{ \nabla^j \mathsf{E}(f)}
					}_{ \mc{F}\del*{ \bR^n } }
				+
				\norm{
					\maxFun[2k-2]{ g }
				}_{ \mc{F}\del*{ \bR^n } }
			}
		\\
		&\le 
			\widehat{C}_{n,k}
			\del*{
				\sum_{j=1}^{k-1}
					C_2^{k+j-1}
					\norm*{
						\nabla^j \mathsf{E}(f)
					}_{ \mc{F}\del*{ \bR^n } }
				+
				C_2^{2k-2}
				\norm{
					g
				}_{ \mc{F}\del*{ \bR^n } }
			},
		\intertext{
			which, denoting 
			$
				C
				\coloneq 
				\widehat{C}_{n,k}
				\sum_{j=1}^{k-1}
					C_2^{k+j-1},
			$
		}
		&\le 
			C\del*{
				\sum_{j=1}^{k-1}
					\norm*{
						\nabla^j \ms{E}(f)
					}_{ \mc{F}\del*{ \bR^n } }
				+
				\norm{
					g
				}_{ \mc{F}\del*{ \bR^n } }
			}
		\\
		&\phantom{mmmmm}\le 
			C\del*{
				2 \norm{ \ms{E}(f) }_{ \hajlasz{s}\mc{F}\del*{ \bR^n } }
				+ \eps
			}
		\le 
			C\del*{
				2 C_1 
				\norm{ f }_{ \hajlasz{s}\mc{F}\del*{ \Omega } }
				+ \eps
			}.	
		\end{align*}
		Therefore,
		\begin{align*}
			\seminorm{ f }_{ \multipoint{s}\mc{F}(\Omega) }
			=
			\inf_{ G' \in \bD_{\lambda}^s(f) }
				\norm{ G' }_{ \mc{F}(\Omega) }
		&\le 
			\norm{ G \rvert_{\Omega} }_{ \mc{F}(\Omega) }
		\\
		&\le 
			C\del*{
				2 C_1 
				\norm{ f }_{ 
					\hajlasz{s}\mc{F}\del*{ \Omega }
				}
				+ \eps
			}
			\xrightarrow{\eps \to 0^+}
			2CC_1
			\norm{ f }_{
				\hajlasz{s}\mc{F}\del*{ \Omega } 
			}.
		\end{align*}
		In consequence,
		denoting $C' \coloneq 2C C_1$,
		\begin{equation*}
			\norm{ f }_{ \multipoint{s}\mc{F}(\Omega) }
			=
			\norm{ f }_{ \mc{F}(\Omega) }
			+
			\seminorm{ f }_{ 
				\multipoint{s}\mc{F}(\Omega) 
			}
			\le 
			\norm{ f }_{ \mc{F}(\Omega) }
			+
			C'
			\norm{ f }_{  
				\hajlasz{s}\mc{F}(\Omega) 
			}
			\le 
			\del*{1 + C'}
			\norm{ f }_{    
				\hajlasz{s}\mc{F}(\Omega) 
			}.
		\end{equation*}
		Thus,
		$
			\hajlasz{s}\mc{F}(\Omega)
			\hookrightarrow
			\multipoint{s}\mc{F}(\Omega).
		$
	\end{enumerate}
	Now, suppose that $s = k$, where $k \in \bN$. We will show that both claims of
	\Cref{thm::higher-order_characterization}
	remain true if we replace every 
	occurrence of 
	$ \hajlasz{s} $
	with 
	$ \sobolev{k} $.
	\begin{enumerate}[label=(\alph*), listparindent=\parindent, leftmargin=0.5em, itemindent=1.5em]
	\item 
		Fix 
		$ f \in \multipoint{k}\mc{F}(\Omega) $
		and 
		$ G \in \bD_{\lambda}^k(f) $
		such that 
		$ G \in \mc{F}(\Omega) $.
		By 
		\Cref{thm::local_hypergradients_imply_local_sobolev}
		we know that,
		for all 
		$ \alpha \in \bN_0^n$ such that 
		$ 0 < \abs{ \alpha } \le k $,
		we have 
		\begin{itemize}
		\item
			$ \partial^{ \alpha} f $ exists in a weak sense
			and 
			$
				\partial^{ \alpha} f
				\in 
				\locIntegrable{1}(\Omega);
			$
		\item 
			$    
				\abs{ \partial^{\alpha} f }
				\le 
				2^{k^2} \partitions{k} G
			$
			almost everywhere.
		\end{itemize}
		Therefore, 
		$ f \in \locSobolev{k,1}(\Omega)$.
		Also,
		for all $j \in [k]$, 
		there exists
		$C_j > 0$ that depends solely on $j$ and $k$
		and
		such that 
		$
			\norm*{
				\nabla^j f
			}_{ \mc{F}(\Omega) }
			\le 
			C_j
			\norm*{ G }_{ \mc{F}(\Omega) }.
		$
		Thus, 
		$ f \in \sobolev{k} \mc{F}(\Omega)$.
		Moreover,
		\begin{align*}
			\seminorm{ f }_{
				\sobolev{k}\mc{F}(\Omega) 
			}
			\le 
			\sum_{j=1}^{ k}
				\norm{ \nabla^j f }_{
					\mc{F}(\Omega) 
				}
			\le
			\sum_{j=1}^{ k}
				C_j
				\norm{ G }_{ \mc{F}(\Omega) }
		&=
			C
			\norm{ G }_{ \mc{F}(\Omega) }
			,
		\end{align*}
		where we denote
		$
			C 
			\coloneq
			\sum_{j=1}^{ k}
				C_j. 
		$
		Passing in the resulting inequality
		$
			\seminorm{ f }_{ 
				\sobolev{k}\mc{F}(\Omega) 
			}
			\le 
			C
			\norm{ G }_{ \mc{F}(\Omega) }
		$
		to~the~infimum over 
		all 
		$ G \in \bD^{k}_{\lambda}(f)$,		
		we get
		$
			\seminorm{ f }_{     
				\sobolev{k}\mc{F}(\Omega) 
			}
			\le 
			C
			\seminorm{ f }_{ 
				\multipoint{k}\mc{F}(\Omega) 
			}.
		$
		As~%
		$
			f \in \multipoint{k}\mc{F}(\Omega) 
		$
		is arbitrary, we get the desired comparison of seminorms.
		Furthermore,
		\begin{equation*}
			\norm{ f }_{
				\sobolev{k}\mc{F}(\Omega) 
			}
			=
			\norm{ f }_{ \mc{F}(\Omega) }
			+
			\seminorm{ f }_{
				\sobolev{k}\mc{F}(\Omega) 
			}
			\le 
			\norm{ f }_{ \mc{F}(\Omega) }
			+
			C
			\seminorm{ f }_{ 
				\multipoint{k}\mc{F}(\Omega) 
			}
			\le 
			\del*{1+C}
			\norm{ f}_{
				\multipoint{k}\mc{F}(\Omega) 
			}.
		\end{equation*}
		Thus, 
		$
			\multipoint{k}\mc{F}(\Omega)
			\hookrightarrow
			\sobolev{k}\mc{F}(\Omega).
		$
	\item
		Let us suppose that there is a bounded extension operator
		$
			\mathsf{E}
			\colon 
			\sobolev{k}\mc{F}(\Omega)
			\to 
			\sobolev{k}\mc{F}\del*{ \bR^n }.
		$
		Then there exists $C_1 > 0$ such that 
		for all $f \in \sobolev{k}\mc{F}(\Omega)$ 
		we have
		$
			\mathsf{E}(f) \rvert_{\Omega} = f
		$
		with the estimate
		$
			\norm{ \mathsf{E}(f) \rvert_{\Omega} }_{   
				\sobolev{k}\mc{F}\del*{ \bR^n } 
			}
			\le 
			C_1 \norm{f}_{ \sobolev{k}\mc{F}(\Omega) }.
		$
		We~also assume that $\mathsf{M}$
		is a bounded operator on~%
		$ \mc{F}\del*{ \bR^n }$,
		so there exists $C_2 > 0$ such that
		for all $f \in \mc{F}\del*{ \bR^n }$
		we have
		$
			\norm{ \maxFun{f} }_{   
				\mc{F}\del*{ \bR^n } 
			}
			\le 
			C_2 \norm{f}_{ \mc{F}\del*{ \bR^n } }.
		$
		
		Fix 
		$
			f \in \sobolev{k}\mc{F}(\Omega).
		$
		Then
		$
			\mathsf{E}(f) 
			\in
			\sobolev{k}\mc{F}\del*{ \bR^n }
		$
		and by
		\Cref{cor::hypergradients_from_multipointwise_bound_for_W^k_loc}
		we know that 
		\begin{equation*}
			G
			\coloneq 
			C_{n,k}
			\sum_{j = 1 }^{ k }
				\maxFun[k+j-1]{
					\nabla^j \mathsf{E}(f)
				}
			\in 
			\bD_{\lambda}^k \del*{ \mathsf{E}(f) }.
		\end{equation*}
		Therefore, 
		$ 
			G \rvert_{ \Omega}
			\in 
			\bD_{\lambda}^k(f).
		$
		Also, from the definition of $G$ and the boundedness of 
		$\mathsf{M}$ on~$\mc{F}\del*{ \bR^n }$ 
		we have 
		$ G \in \mc{F}\del*{ \bR^n }$,
		so 
		$ G \rvert_{\Omega} \in \mc{F}(\Omega)$.
		In consequence,
		since 
		$ f \in \mc{F}(\Omega)$
		as~%
		$f \in \sobolev{k}\mc{F}(\Omega)$,
		we get that
		$
			f \in \multipoint{k}\mc{F}(\Omega).
		$
		
		Moreover, since
		\begin{align*}
			\norm{ G\rvert_{\Omega} }_{ \mc{F}(\Omega) }
			\le
			\norm{ G }_{ \mc{F}\del*{ \bR^n } }
		&\le 
			C_{n,k}
			\sum_{j=1}^{k}
				\norm*{
					\maxFun[k+j-1]{ \nabla^j \mathsf{E}(f)}
				}_{ \mc{F}\del*{ \bR^n } }
		\\
		&\le 
			C_{n,k}
			\sum_{j=1}^{k}
				C_2^{k+j-1}
				\norm*{
					\nabla^j \mathsf{E}(f)
				}_{ \mc{F}\del*{ \bR^n } },
		\intertext{
			which, denoting 
			$
				C
				\coloneq 
				C_{n,k}
				\sum_{j=1}^{k}
					C_2^{k+j-1},
			$
		}
		&\le 
			C
			\sum_{j=1}^{k}
				\norm*{
					\nabla^j \mathsf{E}(f)
				}_{ \mc{F}\del*{ \bR^n } }
			=
			C
			\seminorm*{
				\mathsf{E}(f)
			}_{ \sobolev{k}\mc{F}( \bR^n ) }
			\le
			C'
			\seminorm*{
				f
			}_{ \sobolev{k}\mc{F}( \Omega ) },
		\end{align*}
		where we denote 
		$
			C' \coloneq CC_1.
		$
		In consequence,
		\begin{equation*}
			\seminorm{ f }_{ 
				\multipoint{s}\mc{F}(\Omega) 
			}
			=
			\inf_{ G' \in \bD_{\lambda}^s(f) }
				\norm{ G' }_{ \mc{F}(\Omega) }
			\le 
			\norm{ G \rvert_{\Omega} }_{ \mc{F}(\Omega) }
			\le 
			C'
			\seminorm*{
				f
			}_{ 
				\sobolev{k}\mc{F}( \Omega )
			}.
		\end{equation*}
		Therefore,
		\begin{equation*}
			\norm{ f }_{ 
				\multipoint{k}\mc{F}(\Omega) 
			}
			=
			\norm{ f }_{ \mc{F}(\Omega) }
			+
			\seminorm{ f }_{
				\multipoint{k}\mc{F}(\Omega) 
			}
			\le 
			\norm{ f }_{ \mc{F}(\Omega) }
			+
			C'
			\seminorm{ f }_{
				\sobolev{k}\mc{F}(\Omega) 
			}
			\le 
			\del*{1 + C'}
			\norm{ f }_{ 
				\sobolev{k}\mc{F}(\Omega) 
			}.
		\end{equation*}
		Thus,
		$
			\sobolev{k}\mc{F}(\Omega)
			\hookrightarrow
			\multipoint{k}\mc{F}(\Omega).
		$
		\qedhere
	\end{enumerate}
\end{proof}
\subsubsection{Examples}
We will now give some examples of Banach function spaces that could play the role of~%
$\mc{F}(\Omega)$ in  \Cref{thm::higher-order_characterization}.
\begin{example}[$\integrable{p}$ space]
	Let $p \in \intcc{1,\infty}$. 
	Then $\integrable{p}(\Omega)$ is a Banach function space
	over $\del*{ \Omega, \lambda}$.
	Moreover, if $ p \in \intoc{1,\infty}$, then
	$ \mathsf{M} $ is a bounded operator on
	$ \integrable{p}(\bR^n )$ \cite[Chapter 1, Theorem 1]{Stein}.
	Therefore, 
	\Cref{thm::sobolev_characterization} 
	is a direct consequence of 
	\Cref{thm::higher-order_characterization}
	when
	$ s \in \bN$ 
	and 
	$\mc{F}(\Omega) = \integrable{p}(\Omega)$. 
\end{example}
\begin{example}[$\integrable{p}(\omega)$ space, $\omega \in A_p$]
	Let 
	$ \omega \in \locIntegrable{1}( \bR^n )$
	be positive almost everywhere. 
	One can naturally associate with $\omega$ the~measure 
	defined by the formula
	$ A \mapsto \integral{A}{ \omega }{x} $
	for all Lebesgue measurable sets $A$;
	said measure is~often also denoted by $\omega$.
	Since
	$ \omega $ is positive almost everywhere, 
	for all Lebesgue measurable sets $A$,
	we have
	$
		\omega(A)
		=
		0
	$
	if and only if
	$
		\lambda\del{A} = 0.
	$
	Hence,
	$
		\integrable{0}\del*{ \Omega, \omega}
		=
		\integrable{0}\del*{ \Omega, \lambda}.
	$
	We~will say that
	\begin{itemize}
	\item
		$ \omega \in A_1$
		if there exists $C > 0$ such that
		$ \maxFun{\omega} \le C \omega$ 
		almost everywhere;
	\item 
		$ \omega \in A_p$,
		where $p \in \intoo{1,\infty}$,
		if there exists $C > 0$ such that
		for all balls $B \subseteq \bR^n$,
		we have
		\begin{equation*}
			\del*{
				\sintegral{B}{\omega}{x}
			}\del*{
				\sintegral{B}{\omega^{-q/p}}{x}
			}^{p/q}
			< C,
		\end{equation*}
		where $q$ is the H\"{o}lder conjugate of $p$;
	\item 
		$ \omega \in A_\infty$
		if 
		$ \omega \in A_p$ for some $p \in \intco{1, \infty}$.\footnote{This definition is motivated by the fact that
		by 
		\cite[Theorem 297]{Sawano_Morrey},
		for all $p \in \intoo{1,\infty}$ we have 
		$A_p = \bigcup_{q \in \intco{1,p}} A_q $.
		} 
	\end{itemize}
	The elements of the above-defined sets $A_p$, $p \in \intcc{1,\infty}$, are usually referred to as 
	\emph{Muckenhoupt $A_p$ weights} or \emph{$A_p$ weights}.
	Fix $p \in \intcc{1,\infty}$,
	$ \omega \in A_p$,
	and let 
	$ \Omega \subseteq \bR^n$ 
	be open.
	Then 
	$\integrable{p}\del*{ \Omega, \omega\rvert_{\Omega} }$
	is a Banach function space over $\del*{ \Omega, \lambda}$.
	Moreover,  
	$ \mathsf{M} $ is a bounded operator on
	$ \integrable{p}(\bR^n , \omega)$
	for all $ p \in \intoc{1,\infty} $ 
	\cite[Theorem 290]{Sawano_Morrey}.\footnote{%
		The cited theorem justifies the claim for $p \in \intoo{1, \infty}$.
		For $p = \infty$, the claim follows since
		$
			\integrable{\infty}\del*{
				\bR^n, \omega
			}
			=
			\integrable{\infty}\del*{
				\bR^n
			}.	
		$%
	}
\end{example}
\begin{example}[$\integrable{p,q}$ space]
	Let $n \in \bN$ and $\Omega \subseteq \bR^n$.
	Let $f \colon \Omega \to \bR$ be measurable.
	Let
	$
		\lambda_f
	$
	be the \emph{distribution function} of $f$, 
	that is, the function 
	$
		\lambda_f 
		\colon 
		\intco{0,\infty} \to \intcc{0,\infty}
	$
	defined by the formula
	\begin{equation*}
		\forall t \in \intco{0,\infty }
	\qquad 
		\lambda_f(t)
		\coloneq 
		\abs*{
			\setc*{
				x \in \Omega 
			}{
				\abs*{ f(x) } > t
			}
		}.
	\end{equation*}
	Also, we define
	the \emph{nonincreasing rearrangement} of $f$ 
	as the function 
	$ 
		f^* 
		\colon 
		\intco{0,\infty} \to \intcc{0,\infty} 
	$
	given by the formula
	\begin{equation*}
		\forall s \in \intco{0,\infty}
	\qquad 
		f^*(s)
		\coloneq 
		\inf
		\setc*{
			t \ge 0
		}{
			\lambda_f(t) \le s
		}.
	\end{equation*}

	Fix $p \in \intco{1,\infty}$ and $q \in \intcc{1,\infty}$. 
	Then the \emph{Lorentz space}
	$\integrable{p,q}(\Omega)$
	is defined as the space
	\begin{equation*}
		\integrable{p,q}(\Omega)
		\coloneq 
		\setc*{
			f \colon \Omega \to \bR
			\text{ measurable}
		}{
			\seminorm{ f }_{\integrable{p,q}(\Omega)} <\infty
		},
	\end{equation*}
	where
	$
		\seminormAlone[\integrable{p,q}(\Omega)]
	$
	is a quasinorm defined by the formula
	\begin{equation*}
		\seminorm{ f }_{\integrable{p,q}(\Omega)}
		\coloneq 
		\begin{cases}
			\del*{
				\int_0^{\infty}
					\del*{ t^{1/p} f^*(t) }^q
				\frac{ \mathrm{d}t}{t} 
			}^{1/q}
			& \text{if } q \in \intoc{1, \infty},
		\\
			\sup_{ t > 0 }
				t^{1/p}
				f^*(t)
			& \text{if } q = \infty.
		\end{cases}
	\end{equation*}
	It is known that $\integrable{p,q}(\Omega)$
	can be endowed with a norm
	$ \normAlone[\integrable{p,q}(\Omega)] $
	that is Lipschitz comparable to 
	$ \seminormAlone[\integrable{p,q}(\Omega)] $
	\cite[Corollary 39]{Sawano_Morrey}\footnote{The cited corollary justifies the claim for $p \in \intoo{1, \infty}$. For $p=1$, the claim follows since 
	$\seminormAlone[\integrable{p,q}]$ is a norm when $p=1$. In fact, this is the case for all $p$, $q$ such that $1 \le p \le q < \infty$.}.
	Then
	$\integrable{p,q}(\Omega)$ with this norm
	is a Banach function space 
	over $\del*{ \Omega, \lambda}$.
	Moreover, if $p \in \intoo{1,\infty}$,
	then $\mathsf{M}$
	is a bounded operator on~%
	$ \integrable{p,q}( \bR^n )$
	\cite{boundedness_of_maximal_function_on_Lorentz}.\footnote{%
		The cited reference contains a result for weighted Lorentz spaces.
		The Lorentz spaces 
		$ \integrable{p,q}\del*{ \bR^n } $
		as defined above, 
		correspond to the weight
		$ W(x) = x^{(p/q)-1} $.
		One can then check that for $p>1$, this weight satisfies the condition present in
		\cite[Theorem (1.10)]{boundedness_of_maximal_function_on_Lorentz} from which the boundedness of
		$ \mathsf{M}$ follows.
	}
\end{example}
\subsection{H\"{o}lder spaces}
As the very last result of this section, 
we will show that our approach allows for~a~characterization of not just higher-order Sobolev or Sobolev-like spaces, but also the following H\"{o}lder spaces.
\begin{definition}[H\"{o}lder space $\holderBounded{m,\alpha}
\del*{
	\cl{\Omega }
}$]
	Let $n \in \bN$, $m \in \bN_0$, 
	$ \alpha \in \intoc{ 0, 1 }$,
	and 
	$ \Omega \subseteq \bR^n$
	be open.
	We define the family
	\begin{equation*}
		\holderBounded{m,\alpha}
		\del*{
			\cl{\Omega }
		}
		\coloneq 
		\setc*{
		\!
			f \in 
			\holder{m,\alpha}\del*{ \Omega }
		\!
		}{
			\forall j \in [m] \cup \set{0},
			\text{%
				$\nabla^j f $
				can be extended to an element of
				$ 
					\continuousAndBounded\del*{
						\cl{ \Omega }
					}
				$%
			}\!
		},
	\end{equation*}
	where $\nabla^0 f \coloneq f $,
	and endow it with the following norm:
	\begin{equation*}
		\forall 
		f \in 
		\holderBounded{m,\alpha}\del*{
			\cl{\Omega }
		}
	\qquad 
		\norm*{
			f
		}_{
			\holderBounded{m,\alpha}\del*{
				\cl{\Omega }
			}
		}
		\coloneq 
		\norm{ f }_{ \infty }
		+
		\seminorm{
			f
		}_{
			\holder{m,\alpha}\del*{ \Omega }
		}.
	\end{equation*}
\end{definition}
Let us begin with the following lemma.
\begin{lemma}
\label{label::almost_fundamental_theorem_of_calculus_after_averaging}
	Let $n \in \bN$
	and
	$ \Omega \subseteq \bR^n $ be open.
	Fix 
	$ x \in \Omega $
	and let 
	$ R > 0 $ be such that
	$ \clball{x,2R} \subseteq \Omega $.
	Then for all 
	$ f \in \sobolev{1,1}\del*{ \ball{x,2R} } $,
	$ r \in \intoc{0,R} $,
	and
	$
		v \in \ball{0,r},
	$
	we have
	\begin{equation*}
		\sintegral{\ball{x,r}}{
			\del*{ f(z+v) - f(z) }
		}{z}
		=
		\sintegral{\ball{x,r}}{
			\integral[1]{0}{
				\nabla f(z+tv) \t v
			}{t}
		}{z}.
	\end{equation*}
\end{lemma}
\begin{proof}
	First, let us note that the lemma's statement is true if 
	$ f \in \holder{1}\del*{ \ball{x,2R} } $.
	Indeed, for~such $f$ this follows since for all $z \in \ball{x,r} $
	we have
	$
		f(z+v) - f(z)
		=
		\integral[1]{0}{
			\nabla f(z+tv) \t v
		}{t}		
	$
	by~the~fundamental theorem of calculus.
	
	Next, suppose that 
	$f \in \sobolev{1,1}\del*{ \ball{x,2R} }$.
	Then by 
	\cite[5.3.2. Theorem 2]{evans_partial}
	there is a sequence 
	$\del*{ f_m }_{m \in \bN }$ 
	of elements of 
	$ 
		\sobolev{1,1}\del*{ \ball{x,2R} } 
		\cap 
		\holder{1}\del*{ \ball{x,2R} }
	$
	such that
	$
		f_m \to f
	$
	in 
	$
		\sobolev{1,1}\del*{ \ball{x,2R} } 
	$
	as $m \to \infty$.
	Therefore,
	$
		\norm{ f_m - f }_{
			\integrable{1}\del*{ \ball{x,2R} }
		}
		\to 0
	$
	and 
	$
		\norm{ \nabla f_m - \nabla f }_{
			\integrable{1}\del*{ \ball{x,2R} }
		}
		\to 0
	$	
	as $m \to  \infty $.
	
	We have
	\begin{align*}
		\abs*{
			\sintegral{\ball{x,r}}{
				\del*{ 
					f_m(z+v) - f_m(z)
				}
			}{z}
			-
			\sintegral{\ball{x,r}}{
				\del*{ 
					f(z+v) - f(z)
				}
			}{z}		
		}
	\hspace{-8cm}&
	\\
	&\le 
		\sintegral{\ball{x,r}}{
			\abs*{ f_m(z+v) - f(z+v) } 
		}{z}
		+
		\sintegral{\ball{x,r}}{
			\abs*{ f_m- f } 
		}{z}
	\\		
	&=
		\frac{1}{\abs*{ \ball{x,r}  } }
		\del*{
			\integral{\ball{x+v,r}}{
				\abs*{ f_m - f} 
			}{z}
			+
			\integral{\ball{x,r}}{
				\abs*{ f_m - f} 
			}{z}			
		},
	\intertext{
		which, since both
		$
			\ball{x, r}
		$
		and
		$
			\ball{x + v, r} 
		$
		are subsets of 
		$
			\ball{x,2R},
		$
	}
	&\le 
		\frac{2}{\abs*{ \ball{x,r}  } }
		\integral{\ball{x,2R}}{
			\abs*{ f_m - f} 
		}{ z }
		=
		\frac{
			2 \norm{ f_m- f }_{ 
				\integrable{1}\del*{ \ball{x,2R } }
			}
		}{
			\abs*{ \ball{x,r}  } 
		}
		\xrightarrow{ m \to \infty}
		0.
	\end{align*}
	Therefore, 
	\begin{equation}
	\label{proofeq::label::almost_fundamental_theorem_of_calculus_after_averaging::convergence_scalar_part}
		\sintegral{\ball{x,r}}{
			f_m(z+v) - f_m(z)
		}{z}
		\xrightarrow{ m \to \infty}		
		\sintegral{\ball{x,r}}{
			f(z+v) - f(z)
		}{z}.
	\end{equation}

	Similarly,
	\begin{align*}
		\abs*{
			\sintegral{\ball{x,r}}{
				\integral[1]{0}{
					\nabla f_m(z+tv) \t v
				}{t}
			}{z}
			-
			\sintegral{\ball{x,r}}{
				\integral[1]{0}{
					\nabla f(z+tv) \t v
				}{t}
			}{z}			
		}
	\hspace{-11cm}&
	\\
	&\le 
		\sintegral{\ball{x,r}}{
			\integral[1]{0}{
				\abs*{
					\nabla f_m(z+tv) \t v
					-
					\nabla f(z+tv) \t v
				}
			}{t}
		}{z}	
	\\
	&\le 
		\sintegral{\ball{x,r}}{
			\integral[1]{0}{
				\norm{ \nabla f_m(z+tv) - \nabla f(z+tv)}
				\norm{v}
			}{t}
		}{z}	
	\\
	&=
		\frac{\norm{v}}{\abs{\ball{x,r}}}
		\integral[1]{0}{
			\integral{\ball{x+tv,r}}{
				\norm{ \nabla f_m (z) - \nabla f(z)}
			}{z}			
		}{t},
	\intertext{
		which, since 
		$
			\ball{x+tv,r} \subseteq \ball{x,2R} 
		$
		for all $t \in \intcc{0,1}$,
	}
	&\le 
		\frac{\norm{v}}{\abs{\ball{x,r}}}
		\integral[1]{0}{
			\norm{ \nabla f_m - \nabla f }_{
				\integrable{1}\del*{ \ball{x,2R} }
			}			
		}{t}
		=
		\frac{\norm{v}}{\abs{\ball{x,r}}}
		\norm{ \nabla f_m - \nabla f }_{
			\integrable{1}\del*{ \ball{x,2R} }
		}	
		\xrightarrow{m \to  \infty}
		0.
	\end{align*}
	Therefore,
	\begin{equation}
	\label{proofeq::label::almost_fundamental_theorem_of_calculus_after_averaging::convergence_gradient_part}
		\sintegral{\ball{x,r}}{
			\integral[1]{0}{
				\nabla f_m(z+tv) \t v
			}{t}
		}{z}
		\xrightarrow{m \to  \infty}
		\sintegral{\ball{x,r}}{
			\integral[1]{0}{
				\nabla f(z+tv) \t v
			}{t}
		}{z}.			
	\end{equation}

	Finally, thanks to
	\eqref{proofeq::label::almost_fundamental_theorem_of_calculus_after_averaging::convergence_scalar_part},
	\eqref{proofeq::label::almost_fundamental_theorem_of_calculus_after_averaging::convergence_gradient_part}, 
	and the fact that
	$
		f_m \in \holder{1}\del*{ \ball{x,2R } }
	$
	for all $m \in \bN$,
	\begin{multline*}
		\sintegral{\ball{x,r}}{
			\del*{ f(z+v) - f(z) }
		}{z}
		\xleftarrow{ m \to \infty}	
		\sintegral{\ball{x,r}}{
			\del*{ f_m(z+v) - f_m(z) }
		}{z}
	\\
		=
		\sintegral{\ball{x,r}}{
			\integral[1]{0}{
				\nabla f_m(z+tv) \t v
			}{t}
		}{z}
		\xrightarrow{m \to  \infty}
		\sintegral{\ball{x,r}}{
			\integral[1]{0}{
				\nabla f(z+tv) \t v
			}{t}
		}{z},	
	\end{multline*}
	and the claim follows.
\end{proof}
\begin{lemma}
\label{lem::everywhere_hypergradient_implies_holder_low_order}
	Let $n \in \bN$,
	$ \Omega \subseteq \bR^n$ be open,
	and $s \in \intoc{1,2}$.
	Let 
	$ 
		f \in 
		\integrable{\infty}\del*{
			\Omega, 2^{\Omega}, \#
		}
	$
	be such that there is~a~constant $G \ge 0$
	such that 
	$ G \in \bD^s(f)$.
	Then 
	$
		f 
		\in 
		\holderBounded{1,s-1}\del*{ 
			\cl{\Omega} 
		}.
	$
\end{lemma}
\begin{proof}
	First, let us notice that by 
	\Cref{prop::embedding_into_lower_orders_normed}
	we have that 
	$
		2^{2^2}
		\partitions{2}
		\del*{
			\abs{ f } + G
		}
		\in 
		\bD_{\#}^{1}(f)
		=
		\bD^{1}(f),
	$
	so also
	$
		2^{2^2}
		\partitions{2}
		\del*{
			\norm{ f }_{
				\infty
			} 
			+ G
		}
		\in 
		\bD^{1}(f).
	$
	Thus, $f$ is Lipschitz.
	Since $f$ is also bounded as~an~element of
	$ \integrable{\infty}\del*{\Omega, \#}$,
	it follows that 
	$ f \in \continuousAndBounded(\Omega)$
	and it can be extended to an element of~%
	$\continuousAndBounded\del*{ \cl{\Omega} }$.
	Moreover, by the Rademacher theorem 
	\cite[5.8, Theorem 5]{evans_partial},
	$f$ is differentiable almost everywhere
	and the strong gradient is almost everywhere equal to~the weak gradient of $f$.
	
	Note that we can treat both $f$ and $G$ as elements of 
	$ \locIntegrable{1}(\Omega)$.
	Thus, by 
	\Cref{thm::local_hypergradients_imply_local_sobolev},
	we~know that, for all 
	$ \alpha \in \bN_0^n$ such that
	$ \abs{\alpha} = 1$,
	we have
	\begin{itemize}
	\item
		$
			\abs*{
				\partial^{\alpha} f 
			}
			\le 
			2^{2^2}\partitions{2}G
		$
		almost everywhere;
	\item 
		$
			16G = 2^{2(2)(1)}G
			\in 
			\bD_{\lambda}^{s-1}\del*{
				\partial^{\alpha}f
			}.
		$
	\end{itemize}
	It follows that, for all such $\alpha$, 
	we have
	$
		\partial^{\alpha} f
		\in 
		\integrable{\infty}( \Omega ).
	$
	
	Moreover, there exists a set $F \subseteq \Omega$
	of full measure such that
	\begin{itemize}
	\item
		$f$ is strongly differentiable everywhere in $F$;
	\item 
		$
			\norm*{
				\nabla f
			}
			\le 
			C G
		$
		everywhere in $F$
		for some $C > 0$ that depends only on $n$;
	\item  
		for all $x, y \in F$ we have
		$
			\norm*{
				\nabla f(x)
				- 
				\nabla f(y)
			}
			\le 
			\norm{x-y}^{s-1}
			2\wave{C}G
		$
		for some $\wave{C} > 0$ that depends only on $n$.
	\end{itemize}
	As $F$ is of full measure in $\Omega$, 
	$ \del*{  \nabla f } \rvert_F $
	can be extended to an element of 
	$ \holderBounded{0,s-1}\del*{ \cl{ \Omega} }$,
	which we will denote by~$h$.
	Note that since the strong gradient of $f$ is almost everywhere equal to the weak gradient of $f$, 
	$h$~is~a~representative of the weak gradient of $f$.
	Also, we have
	$
		\wave{C}G \in \bD^{s-1}\del*{ h }.
	$  
	
	Next, we will show that $h$ is the strong differential of $f$ everywhere in $\Omega$.
	Fix $ x \in \Omega $.
	Let $R > 0$ be such that 
	$ \clball{x,2R} \subseteq \Omega $
	and fix $y \in \ball{x,R}$.
	Denote $v \coloneq y-x$ and 
	$ r \coloneq \norm{v} = \norm{y-x}$.
	Let~us~notice that
	\begin{align*}
		\abs*{
			f(y) - f(x) - h(x) \cdot (y-x) 
		}
	\hspace{-4.75cm}&
	\\
	&=
		\abs*{
			\sintegral{\ball{x,r}}{
				\del*{
					f(y) - f(x) 
					- \del*{ f\del*{z + v} - f(z)}
					+ \del*{ 
						f\del*{z + v} - f(z)
						- h(x) \cdot (y-x) 
					}
				}
			}{z}
		},
	\intertext{
		which, by 
		\Cref{label::almost_fundamental_theorem_of_calculus_after_averaging},
		since
		$ h $ is a representative of the weak gradient of $f$,
	}
	&=
		\abs*{
			\sintegral{\ball{x,r}}{
				\del*{
				f(y) - f(x) 
				- \del*{ f\del*{z + v} - f(z)}
				+ 
				\integral[1]{0}{
					\del*{
						h\del*{ z + tv } \cdot v
						- h(x) \cdot (y-x) 
					}
				}{t}
				}
			}{z}
		}	
	\\
	&\le 
	\underbrace{
		\sintegral{\ball{x,r}}{
			\abs*{
				f(y) - f(x) 
				- f\del*{z {+} v} + f(z)
			}
		}{z}
	}_{ \textstyle{ \eqcolon I_1 } }
		+
	\underbrace{
		\sintegral{\ball{x,r}}{
			\integral[1]{0}{
				\abs*{
					\del*{ 
						h\del*{ z {+} tv } 
						- h(x)
					}\cdot (y {-} x)
				}
			}{t}
		}{z}
	}_{ \textstyle{ \eqcolon I_2 } }.
	\end{align*}
	Thus, 
	\begin{equation}
	\label{proofeq::lem::everywhere_hypergradient_implies_holder_low_order::estimate_with_I1_and_I2}
		\abs*{
			f(y) - f(x) - h(x) \cdot (y-x) 
		}
		\le 
		I_1 + I_2.
	\end{equation}

	Let us estimate $I_1$ first.
	Fix $z \in \ball{x,r}$
	and let $\bz = \set*{z_I}_{I \subseteq [2]}$
	denote the tuple given by 
	\begin{equation*}
		z_{\emptyset} = x,
	\quad 
		z_{\set*{1}} = y,
	\quad 
		z_{\set*{2}} = z,
	\quad \text{and} \quad    
		z_{\set*{1,2}} = z+v.
	\end{equation*}
	Let us notice that,
	denoting 
	$ v_1 \coloneq v $
	and 
	$ v_2 \coloneq z-x$,
	we have
	$
		z_I
		=
		x + \sum_{ i \in I} v_i
	$
	for all $I \subseteq [2]$.
	Thus, by~%
	\Cref{cor::value_of_Poly_when_tuple_is_a_hyperparallelogram}
	we have
	\begin{equation*}
		\polygen{ \gen{s} }( \bz )
		=
		2^{s}
		\norm{ v_2 }^{s-1}
		\norm{ v_1 }
		=
		2^s 	
		\norm{ z - x }^{s-1}
		r
		\le 
		2^{s}
		r^s
		=
		2^{s} \norm{ y - x }^{s}.
	\end{equation*}
	Therefore, since 
	$ G \in \bD^{s}(f)$, 
	it follows that, 
	for all 
	$ z \in \ball{x, r}$,
	we have
	\begin{equation*}
		\abs*{
			f(y) - f(x) 
			- f\del*{z + v} + f(z)
		}
		=
		\abs*{
			\diff_{ I = \emptyset }^{ [2] }
				f\del*{z_I}
		}
		\le 
		\polygen{ \gen{s} }( \bz )
		\sum_{ I = \emptyset }^{ [2] }
			G
		\le 
		2^{2+s} \norm{ y - x }^{s} G.
	\end{equation*}
	In consequence,
	\begin{align}
		I_1
		=
		\sintegral{\ball{x,r}}{
			\abs*{
				f(y) - f(x) 
				- f\del*{z + v} + f(z)
			}
		}{z}
	&
		\notag
	\\
		\le 
		\sintegral{\ball{x,r}}{
			2^{2+s} \norm{ y - x }^{s} G
		}{z}	
	&=
		2^{2+s} \norm{ y - x }^{s} G.	
		\label{proofeq::lem::everywhere_hypergradient_implies_holder_low_order::estimate_for_I1}
	\end{align}
	
	Moving on to $I_2$,
	let us notice that, 
	for all 
	$ z \in \ball{x,r} $
	and every $t \in \intcc{0,1}$,
	we have
	\begin{equation*}
		\norm{ z + tv - x }
		\le 
		\norm{z-x} + t \norm{v}
		\le 
		r + \norm{v}
		= 
		2\norm{y-x}.		
	\end{equation*}
	Thus, since 
	$
		\wave{C}G \in \bD^{s-1}\del*{ h },
	$ 
	\begin{align*}
		\abs*{
			\del*{ 
				h\del*{ z + tv } 
				- h(x)
			}\cdot (y-x)
		}		
	&\le 
		\norm*{
			h\del*{ z + tv } 
			- h(x)
		}
		\norm{ 
			y- x 
		}
	\\
	&\le 
		\norm{ z + tv - x }^{s-1} 
			2 \wave{C} G
		\norm{ y - x }
		\le 
		2^{s} \wave{C} G\norm{ y-x }^{s}.
		%
	\end{align*}
	In consequence,
	\begin{align}
	\label{proofeq::lem::everywhere_hypergradient_implies_holder_low_order::estimate_for_I2}
		I_2
		=
		\sintegral{\ball{x,r}}{
			\integral[1]{0}{
				\abs*{
					\del*{ 
						h\del*{ z + tv } 
						- h(x)
					}\cdot (y-x)
				}
			}{t}
		}{z}
	\notag
	&
	\\
		\le 
		\sintegral{\ball{x,r}}{
			\integral[1]{0}{
				2^{s} \wave{C} G
				\norm{ y-x }^{s}
			}{t}
		}{z}
	&=
		2^{s} \wave{C} G
		\norm{ y-x }^{s}.
	\end{align}
	Thus, using 
	\eqref{proofeq::lem::everywhere_hypergradient_implies_holder_low_order::estimate_for_I1}
	and
	\eqref{proofeq::lem::everywhere_hypergradient_implies_holder_low_order::estimate_for_I2}
	to estimate the right-hand side of \eqref{proofeq::lem::everywhere_hypergradient_implies_holder_low_order::estimate_with_I1_and_I2},
	we get
	\begin{equation*}
		\abs*{
			f(y) - f(x) - h(x) \cdot (y-x) 
		}
		\le
		I_1 + I_2
		\le 
		\del*{
			2^{2+s} + 2^{s}\wave{C}
		} 
		G \norm{ y- x }^s.
	\end{equation*}
	Since $s > 1$, using the above estimate,
	we can conclude that 
	$f$ is differentiable at $x$ with $h(x)$ being the strong differential.
	
	We have shown that, for all 
	$ x \in \Omega $, 
	function $f$ is differentiable at $x$ 
	with $h(x)$ being the~strong differential.	
	Thus, $h = \nabla f$ everywhere in $\Omega$.
	Finally, recall that we have previously shown that 
	$
		f \in \continuousAndBounded\del*{ \cl{\Omega}}
	$
	and~we~have
	$
		h \in \holderBounded{0,s-1}\del*{ \cl{\Omega}}.
	$
	Since 
	$ h = \nabla f $ everywhere, 
	it follows that
	$
		f \in \holderBounded{1,s-1}\del*{ \cl{\Omega}},
	$
	as~claimed.
\end{proof}
Before stating the characterization of higher-order H\"{o}lder spaces, let us note that, much like in the case of the Hajłasz-like spaces (%
\Cref{rem::up_to_first-order_characterization}),
we have a much stronger result for the orders $s \in \intoc{0,1}$.
\begin{remark}
\label{rem::up_to_first-order_characterization_holder}
	Let $n \in \bN$, $\alpha \in \intoc{0,1}$,
	and $ \Omega \subseteq \bR^n $ be open.
	Then 
	$
		\multipoint{\alpha}
		\integrable{\infty}\del*{
			\Omega, \#
		}
		\cong
		\holderBounded{0,\alpha}\del*{ \cl{ \Omega } }
	$
	with 
	$
		2
		\seminormAlone[{
			\multipoint{\alpha}
			\integrable{\infty}\del*{
				\Omega, \#
			}
		}]
		=
		\seminormAlone[
			{\holder{0,\alpha}\del*{ \Omega }}
		].
	$
\end{remark} 
\begin{proof}
	Fix 
	$
		f \in \holderBounded{0,\alpha}\del*{ 
			\cl{ \Omega } 
		}.
	$
	Then 
	$
		f \in \continuousAndBounded\del*{
			\cl{\Omega}
		},
	$
	so 
	$
		f \in \integrable{\infty}\del*{ \Omega, \# }.
	$
	Moreover, for all $x, y \in \Omega$ we have
	$
		\abs*{ f(x) - f(y) }
		\le 
		\norm{ x - y }^{\alpha}
		\seminorm{ f }_{
			\holder{0,\alpha}\del*{ \Omega }
		},
	$
	so 
	$
		\frac{1}{2}
		\seminorm{ f }_{
			\holder{0,\alpha}\del*{ \Omega }
		}
		\in 
		\bD^{\alpha}(f).	
	$
	In consequence,
	$
		f 
		\in 
		\multipoint{\alpha}
		\integrable{\infty}\del*{
			\Omega, \#
		}
	$ 
	and 
	$
		2 \seminorm{ f }_{
			\multipoint{\alpha}
			\integrable{\infty}\del*{
				\Omega, \#
			}
		}
		\le 
		\seminorm{ f }_{
			\holder{0,\alpha}\del*{ \Omega }
		}.		
	$
	Thus, 
	$
		\holderBounded{0,\alpha}\del*{ \cl{ \Omega } }
		\hookrightarrow
		\multipoint{\alpha}
		\integrable{\infty}\del*{
			\Omega, \#
		}.
	$

	Next, fix 
	$
		f \in 
		\multipoint{\alpha}
		\integrable{\infty}\del*{
			\Omega, \#
		}.
	$
	Then we have
	$
		\abs*{ f(x) - f(y )}
		\le 
		\norm{x-y}^{\alpha}
		\del*{ 
			2 \seminorm{ f }_{
				\multipoint{\alpha}
				\integrable{\infty}\del*{
					\Omega, \#
				}
			}
		}
	$
	for~all $x, y \in \Omega$.
	Thus, $f$ is continuous on $\Omega$ and can be uniquely extended to a continuous function on 
	$ \cl{\Omega} $ 
	with
	$
		\norm{ f }_{
			\infty
		}
		=
		\norm{ f }_{
			\integrable{\infty}\del*{ \Omega, \# }
		}.
	$
	Moreover, we have that
	$
		\seminorm{ f }_{
			\holder{0,\alpha}\del*{ \Omega }
		}
		\le 
		2 \seminorm{ f }_{
			\multipoint{\alpha}
			\integrable{\infty}\del*{
				\Omega, \#
			}
		}.
	$
	Thus, 
	$
		\multipoint{\alpha}
		\integrable{\infty}\del*{
			\Omega, \#
		}
		\hookrightarrow
		\holderBounded{0,\alpha}\del*{ \cl{ \Omega } }.
	$
	
	Finally, let us note that we have 
	$
		\seminormAlone[{
			\holder{0,\alpha}\del*{ \Omega }
		}]
		= 
		2 \seminormAlone[{
			\multipoint{\alpha}
			\integrable{\infty}\del*{
				\Omega, \#
			}
		}]
	$	
	by combining the previously obtained estimates.
\end{proof}
\begin{theorem}
\label{thm::higher-order_characterization_of_holder_spaces}
    Let $n, k \in \bN,$ $\alpha \in \intoc{0,1}$,
    and
	$ \Omega \subseteq \bR^n $
    be open. Then
	\begin{enumerate}[label=(\alph*)]
	\item 
	\label{thm::higher-order_characterization_of_holder_spaces::part_from_multipointwise}
		$
			\multipoint{k+\alpha}
			\integrable{\infty}\del*{
				\Omega, \#
			}
			\hookrightarrow
			\holderBounded{k,\alpha}\del*{\cl{\Omega}}.
		$
	\item 
	\label{thm::higher-order_characterization_of_holder_spaces::part_to_multipointwise}
		If there exists a bounded extension operator
		$
			\mathsf{E}
			\colon 
			\holderBounded{k,\alpha}\del*{ \cl{\Omega}}
			\to
			\holderBounded{k,\alpha}\del*{ \cl{ \bR^n }},
		$
		then we have
		$
			\holderBounded{k,\alpha}\del*{\cl{\Omega}}
			\cong
			\multipoint{k+\alpha}
			\integrable{\infty}\del*{
				\Omega, \#
			}.
		$
	\end{enumerate}
\end{theorem}
\begin{proof}
	${}$
	\begin{enumerate}[label=(\alph*),listparindent=\parindent, leftmargin=0.5em, itemindent=1.5em]
	\item 
		The proof will be inductive over the value of $k \in \bN$.
		First, let us notice that 
		although the~theorem's statement is for $k \in \bN$,
		by 
		\Cref{rem::up_to_first-order_characterization_holder}
		the two claims are also true for~%
		$k = 0$.
		Next, let us assume that $k \in \bN$
		and 
		\ref{thm::higher-order_characterization_of_holder_spaces::part_from_multipointwise}
		is true for $k-1$.
	
		First, let us suppose that 
		$
			f \in 
			\multipoint{k+\alpha}
			\integrable{\infty}\del*{
				\Omega, \#
			}.
		$
		Then
		$
			\seminorm{
				f
			}_{ 
				\multipoint{k+\alpha}
				\integrable{\infty}\del*{
					\Omega, \#
				}
			}
			\in \bD^{k+\alpha}(f).
		$
		By~%
		\Cref{prop::embedding_into_lower_orders_normed}
		we have that 
		$
			2^{(k+1)^2}
			\partitions{k+1}
			\del*{
				\abs{ f } 
				+
				\seminorm{
					f
				}_{ 
					\multipoint{k+\alpha}
					\integrable{\infty}\del*{
						\Omega, \#
					}
				}
			}
			\in 
			\bD_{\#}^{k}(f)
			=
			\bD^{k}(f).
		$
		Hence, we also have
		$
			2^{(k+1)^2}
			\partitions{k+1}
			\norm{
				f
			}_{ 
				\multipoint{k+\alpha}
				\integrable{\infty}\del*{
					\Omega, \#
				}
			}
			\in 
			\bD^{k}(f).
		$
		Thus, 
		$
			f \in \multipoint{k}
			\integrable{\infty}\del*{
				\Omega, \#
			}.
		$ 
		From the induction hypothesis,
		it~follows that
		$
			f 
			\in 
			\holderBounded{k-1,1}\del*{ \cl{\Omega} }.
		$
		
		Fix $\beta \in \bN_0^n$ such that
		$ \abs{ \beta } = k-1$.
		Then
		$
			\partial^{\beta} f
			\in \continuousAndBounded\del*{ \Omega },
		$
		so also 
		$
			\partial^{\beta} 
			f
			\in
			\integrable{\infty}\del*{\Omega, \# }.
		$
		Furthermore, by~%
		\Cref{thm::local_hypergradients_imply_local_sobolev}
		we know that 
		$
			G
			\coloneq 
			2^{2(k+1)(k-1)}
			\seminorm{
				f
			}_{ 
				\multipoint{k+\alpha}
				\integrable{\infty}
				\del*{
					\Omega, \#
				}
			}
			\in 
			\bD_{\lambda}^{1+\alpha}\del*{
				\partial^{\beta} f
			}.
		$
		Fix~%
		$
			F \in 
			\mf{F}_{\lambda}^{1+\alpha}\del*{
				\partial^{\beta}f,
				G
			}.
		$
		Then
		\begin{equation*}
			\forall 
			\bx = \set*{ x_I }_{ I \subseteq [2] }
			\subseteq 
			F
		\qquad 
			\abs*{
				\diff_{I = \emptyset}^{ [2] }
					\partial^{\beta} f\del*{ x_I }
			}
			\le 
			\polygen{ \gen{1+\alpha} }(\bx)
			\sum_{ I = \emptyset }^{ [2] }
				G\del*{ x_I}.
		\end{equation*}
		Since 
		$
			\partial^{\beta} f,
			G
			\in 
			\continuousAndBounded\del*{ \Omega }
		$
		and 
		$
			\polygen{ \gen{1+\alpha} }
		$
		is continuous as a function from 
		$ \Omega^4 $, it follows that the~above inequality is also satisfied for all tuples
		$
			\bx = \set*{ x_I }_{ I \subseteq [2] }
			\subseteq 
			\Omega.
		$
		Thus, 
		$
			G \in \bD^{1+\alpha}\del*{
				\partial^{\beta} f
			}.
		$
		Therefore,
		$
			\partial^{\beta} f 
			\in 
			\multipoint{1+\alpha}\integrable{\infty}
			\del*{
				\Omega, \#
			}. 
		$
		By 
		\Cref{lem::everywhere_hypergradient_implies_holder_low_order}
		we have
		$
			\partial^{\beta} f
			\in 
			\holderBounded{1,\alpha}\del*{
				\cl{\Omega} 
			}.
		$
		
		We showed that
		$
			f \in \holderBounded{k-1,1}\del*{
				\cl{\Omega} 
			}
		$
		and that,
		for all 
		$\beta \in \bN_0^n$ 
		with
		$\abs{ \beta } = k - 1$,
		we~have
		$
			\partial^{\beta} f
			\in 
			\holderBounded{1,\alpha}\del*{
				\cl{\Omega} 
			}.
		$
		Hence, it~follows that
		$
			f \in 
			\holderBounded{k,\alpha}\del*{ 
				\cl{\Omega} 
			}.
		$
		It remains to prove that the~resulting embedding is bounded.
		
		Since 
		$
			\seminorm{ f }_{
				\multipoint{k+\alpha}\integrable{\infty}\del*{ \Omega, \# }
			}
			\in 
			\bD^{k+\alpha}(f),
		$
		by 
		\Cref{thm::local_hypergradients_imply_local_sobolev},
		for all 
		$ 
			\beta \in \bN_0^n
		$
		such that 
		$ 0 < \abs{ \beta } \le k $,
		we have
		$
			\abs*{ \partial^{ \beta } f }
			\le 
			2^{(k+1)^2}
			\partitions{k+1}
			\seminorm{ f }_{
				\multipoint{k+\alpha}\integrable{\infty}\del*{ \Omega, \# }
			}
		$
		almost everywhere in $\Omega$.
		Since 
		$
			\partial^{ \beta } f \in \continuousAndBounded\del*{ \cl{\Omega} },
		$
		this inequality is, in fact, satisfied everywhere in 
		$ \Omega $.		
		Thus, we have shown that
		\begin{equation}
		\label{proofeq::thm::higher-order_characterization_of_holder_spaces::estimate_for_supremum_norm_for_partial_beta}
			\forall 
				\beta \in \bN_0^n
		\qquad 
			0 < \abs{\beta} \le k
		\quad \implies \quad 
			\norm*{
				\partial^{ \beta } f	
			}_{
				\infty
			}
			\le 
			2^{(k+1)^2}
			\partitions{k+1}
			\seminorm{ f }_{
				\multipoint{k+\alpha}\integrable{\infty}\del*{ \Omega, \# }
			}.
		\end{equation} 
		In consequence,
		there exists
		$C_1 $
		(that only depends on $n$ and $k$)
		such that
		\begin{equation}
		\label{proofeq::thm::higher-order_characterization_of_holder_spaces::estimate_forfirst_part_of_holder_seminorm}
			\sum_{j = 1}^{k}
				\norm*{
					\nabla^j f
				}_{\infty}
			\le 
			C_1
			\seminorm{ f }_{
				\multipoint{k+\alpha}\integrable{\infty}\del*{ \Omega, \# }
			}.
		\end{equation}

		Moreover, 
		by 
		\Cref{thm::local_hypergradients_imply_local_sobolev}
		we also know that,
		for all 
		$ 
			\beta \in \bN_0^n
		$
		such that 
		$ 0 < \abs{ \beta } \le k $,
		we have
		$
			G_{\beta}
			\coloneq 
			2^{2(k+1)\abs{\beta}}
			\seminorm{ f }_{
				\multipoint{k+\alpha}\integrable{\infty}\del*{ \Omega, \# }
			}
			\in 
			\bD_{\lambda}^{k+\alpha-\abs{\beta}}\del*{
				\partial^{ \beta } f
			}.
		$
		Fix such $\beta$ and
		$
			F_{\beta} 
			\in 
			\mf{F}_{\lambda}^{k+\alpha-\abs{\beta}}
			\del*{
				\partial^{\beta} f, 
				G_{\beta}
			}.
		$
		Then
		\begin{equation*}
			\forall 
			\bx = \set*{ x_I }_{ 
				I \subseteq [k+1-\abs{\beta}]
			}
			\subseteq 
			F_{\beta}
		\qquad 
			\abs*{
				\diff_{ I = \emptyset}^{ 
					[k+1-\abs{\beta}]
				}
					\partial^{ \beta }f\del*{x_I}
			}
			\le 
			\polygen{
				\gen{k+\alpha - \abs{\beta} }
			}( \bx )
			\sum_{ I = \emptyset }^{ 
					[k+1-\abs{\beta}]			
			}
				G_{\beta}(x_I).
		\end{equation*}
		Since 
		$
			\partial^{ \beta } f,
			G_{\beta}
			\in 
			\continuousAndBounded(\Omega)
		$
		and 
		$
			\polygen{
				\gen{k+\alpha - \abs{\beta} }
			}
		$
		is continuous as a function from
		$
			\Omega^{
				2^{k+1 - \abs{\beta}}
			},
		$
		it~follows that the~above inequality is satisfied for all tuples 
		$
			\bx = \set*{ x_I }_{ 
				I \subseteq \sbr*{ k+1-\abs{\beta} }
			}
			\subseteq 
			\Omega.
		$
		Thus, we have
		$
			G_{\beta}
			\in 
			\bD^{k+\alpha-\abs{\beta}}\del*{
				\partial^{\beta} f
			}.
		$
		In consequence,
		since
		$
			\partial^{ \beta } f,
			G_{\beta}
			\in 
			\integrable{\infty}\del*{
				\Omega, \#
			}
		$	
		as
		$
			\partial^{ \beta } f,
			G_{\beta}
			\in 
			\continuousAndBounded(\Omega),
		$
		we~have
		$
			\partial^{\beta} f
			\in 
			\multipoint{k+\alpha-\abs{\beta}}
			\integrable{\infty}\del*{
				\Omega, \#
			}.
		$
		Moreover, 
		\begin{align}
			\seminorm*{
				\partial^{\beta} f
			}_{
				\multipoint{k+\alpha-\abs{\beta}}
				\integrable{\infty}\del*{
					\Omega, \#
				}
			}
			\le 
			\norm{ G_{\beta} }_{
				\integrable{\infty}\del*{
					\Omega, \#
				}
			}
			=
			2^{2(k+1)\abs{\beta}}
			\seminorm{
				f
			}_{
				\multipoint{k+\alpha}
				\integrable{\infty}\del*{
					\Omega, \#
				}
			}
			\notag
		&
		\\
			\le 
			2^{2k(k+1)}
			\seminorm{
				f
			}_{
				\multipoint{k+\alpha}
				\integrable{\infty}\del*{
					\Omega, \#
				}
			}			
		&.
		\label{proofeq::thm::higher-order_characterization_of_holder_spaces::estimate_for_multipoint_seminorm_of_partial_beta}
		\end{align}
		Note that the above inequality is true for all 
		$ \beta \in \bN_0^n $
		with 
		$
			0 < \abs{ \beta} \le k.
		$		
		Hence, this fact, combined with 
		\eqref{proofeq::thm::higher-order_characterization_of_holder_spaces::estimate_for_supremum_norm_for_partial_beta},
		implies that for all such $\beta$, 
		we have
		\begin{align*}
			\norm*{
				\partial^{\beta} f
			}_{
				\multipoint{k+\alpha-\abs{\beta}}
				\integrable{\infty}\del*{
					\Omega, \#
				}			
			}
		&=
			\norm*{ 
				\partial^{\beta} f
			}_{
				\integrable{\infty}
				\del*{
					\Omega, \#
				}
			}
			+
			\seminorm*{
				\partial^{\beta} f
			}_{
				\multipoint{k+\alpha-\abs{\beta}}
				\integrable{\infty}
				\del*{
					\Omega, \#
				}			
			}
		\\
		&\le 
			\del*{
				2^{(k+1)^2}
				\partitions{k+1}
				+
				2^{2k(k+1)}
			}
			\seminorm{ f }_{
				\multipoint{k+\alpha}
				\integrable{\infty}
				\del*{
					\Omega, \#
				}
			}.
		\end{align*}
		%
		
		Now, fix
		$
			\beta \in \bN_0
		$
		with 
		$ \abs{ \beta } = k $.
		From previous reasoning, it follows that 
		$
			\partial^{\beta} f
			\in 
			\multipoint{\alpha}
			\integrable{\infty}
			\del*{
				\Omega, \#
			},
		$
		so by 
		\Cref{rem::up_to_first-order_characterization_holder}
		we have that
		$
			\partial^{\beta} f
			\in 
			\holderBounded{0,\alpha}\del*{\cl{\Omega}}
		$		
		with 
		\begin{equation*}
			\seminorm*{
				\partial^{\beta} f
			}_{
				\holder{0,\alpha}\del*{ \Omega }			
			}
			=
			2 
			\seminorm*{
				\partial^{\beta} f		
			}_{
				\multipoint{\alpha}
				\integrable{\infty}
				\del*{
					\Omega, \#
				}
			}
			\le 
			2\t 
			2^{2k(k+1)}
			\seminorm{
				f			
			}_{
				\multipoint{k+\alpha}
				\integrable{\infty}
				\del*{
					\Omega, \#
				}
			},
		\end{equation*}
		where the estimate follows from
		\eqref{proofeq::thm::higher-order_characterization_of_holder_spaces::estimate_for_multipoint_seminorm_of_partial_beta}.
		Thus, there exists 
		$ C_2 > 0 $, 
		depending only on $n$ and $k$,
		such that
		$
			\seminorm*{
				\nabla^k f
			}_{
				\holder{0,\alpha}\del*{ \Omega }
			}
			\le 
			C_2
			\seminorm{
				f
			}_{
				\multipoint{k+\alpha}
				\integrable{\infty}
				\del*{
					\Omega, \#
				} 
			}.
		$
		This inequality,
		combined with 
		\eqref{proofeq::thm::higher-order_characterization_of_holder_spaces::estimate_forfirst_part_of_holder_seminorm},
		implies that
		\begin{equation*}
			%
			\seminorm{
				f
			}_{ 
				\holder{k,\alpha}\del*{
					\Omega  
				}
			}
			=
			\sum_{j = 1}^{k}
				\norm*{ \nabla^j f }_{
					\infty
				}
			+
			\seminorm*{
				\nabla^k f
			}_{
				\holder{0,\alpha}\del*{ \Omega }
			} 
			\le 
			C'
			\seminorm{ f }_{
				\multipoint{k+\alpha}
				\integrable{\infty}\del*{\Omega, \#}
			},
		\end{equation*}
		where we denote
		$
			C'
			\coloneq 
			C_1 + C_2.
		$
		In consequence,
		since
		$
			\norm{ f }_{ 
				\infty
			}
			=
			\norm{ f }_{
				\integrable{\infty}\del*{\Omega, \#}
			},
		$
		\begin{align*}
			\norm{
				f
			}_{
				\holderBounded{k,\alpha}\del*{ 
					\cl{\Omega}
				}
			}
		&=
			\norm{ f }_{ 
				\infty
			}
			+
			\seminorm{
				f
			}_{
				\holder{k,\alpha}\del*{ \Omega }
			}
		\\	
		&\le 
			\norm{ f }_{
				\integrable{\infty}\del*{\Omega, \#}
			}
			+
			C'
			\seminorm{ f }_{
				\multipoint{k+\alpha}\integrable{\infty}\del*{\Omega, \#}
			}	
			\le 
			\del*{1+C'}
			\norm{ f }_{
				\multipoint{k+\alpha}\integrable{\infty}\del*{\Omega, \#}		
			}.
		\end{align*}
		Thus, since
		$
			f \in 
			\multipoint{k+\alpha}
			\integrable{\infty}\del*{
				\Omega, \#
			}
		$ 
		is arbitrary, we have
		$
			\multipoint{k+\alpha}
			\integrable{\infty}\del*{
				\Omega, \#
			}	
			\hookrightarrow
			\holderBounded{k,\alpha}\del*{\cl{\Omega}}.
		$
	\item
		Suppose now that there exists a continuous extension operator
		$
			\mathsf{E}
			\colon 
			\holderBounded{k,\alpha}\del*{
				\cl{\Omega}
			}
			\to
			\holderBounded{k,\alpha}\del*{ 
				\cl{ \bR^n }
			}.
		$
		Then there exists $C > 0 $ such that,
		for all 
		$ 
			f \in 
			\holderBounded{k,\alpha}\del*{
				\cl{\Omega}
			}
		$,
		we have
		$
			\norm{
				\ms{E}(f)
			}_{
				\holderBounded{k,\alpha}\del*{ 
					\cl{\bR^n}
				}
			}
			\le 
			C
			\norm{
				f
			}_{
				\holderBounded{k,\alpha}\del*{ 
					\cl{\Omega}
				}
			}.			
		$
		Let us fix
		$
			f \in 
			\holderBounded{k,\alpha}\del*{
				\cl{\Omega}
			}.			
		$
		Then we have
		$
			\mathsf{E}(f)
			\in 
			\holderBounded{k,\alpha}\del*{ 
				\cl{\bR^n}
			}
		$
		and by 
		\Cref{cor::multipointwise_bound_for_higher_order_holder}
		we obtain that
		$
			\widehat{C}_{n,k+1}
			\seminorm*{
				\mathsf{E}(f)
			}_{
				\holder{k,\alpha}\del*{
					\bR^n
				}			
			}
			\in 
			\bD^{k+\alpha}\del*{ \mathsf{E}(f) },
		$
		where
		$
			\widehat{C}_{n,k+1}
		$
		is the constant from 
		\Cref{thm::multipointwise_bound_for_W^k_loc}.
		In~consequence,
		$
			\widehat{C}_{n,k+1}
			\seminorm{
				\mathsf{E}(f)
			}_{
				\holder{k,\alpha}\del*{
					\bR^n
				}			
			}
			\in 
			\bD^{k+\alpha}\del*{ f }.
		$
		Since 
		$
			f \in \integrable{\infty}\del*{
				\Omega, \#
			}
		$
		as~%
		$
			f \in \continuousAndBounded\del*{
				\cl{\Omega}
			},
		$
		it~follows that
		$
			f \in 
			\multipoint{k+\alpha}
			\integrable{\infty}
			\del*{
				\Omega, \#
			}
		$
		with 
		\begin{equation*}
			\seminorm{
				f
			}_{
				\multipoint{k+\alpha}
				\integrable{\infty}
				\del*{
					\Omega, \#
				}			
			}
			\le 
			\widehat{C}_{n,k+1}
			\seminorm{
				\mathsf{E}(f)
			}_{
				\holder{k,\alpha}\del*{ \bR^n }			
			}
			\le 
			\widehat{C}_{n,k+1}
			\norm{
				\mathsf{E}(f)
			}_{
				\holderBounded{k,\alpha}\del*{
					\cl{\bR^n}
				}			
			}
			\le
			C' 
			\norm{
				f		
			}_{
				\holderBounded{k,\alpha}\del*{
					\cl{\Omega}
				}			
			},
		\end{equation*}
		where we denote 
		$
			C'
			\coloneq 
			\widehat{C}_{n,k+1}
			C.
		$
		Thus,
		\begin{align*}
			\norm{
				f
			}_{
				\multipoint{k+\alpha}
				\integrable{\infty}
				\del*{
					\Omega, \#
				}			
			}
			=
			\norm{ f }_{
				\integrable{\infty}\del*{ \Omega, \# }
			}
			+
			\seminorm{
				f
			}_{
				\multipoint{k+\alpha}
				\integrable{\infty}
				\del*{
					\Omega, \#
				}			
			}
		&
		\\
			\le 
			\norm{ f }_{
				\infty
			}				
			+
			C'
			\norm{ f }_{
				\holderBounded{k,\alpha}\del*{ 
					\cl{\Omega}
				}
			}			
		&\le 
			\del*{ 1+C' }		
			\norm{ f }_{
				\holderBounded{k,\alpha}\del*{ 
					\cl{\Omega}
				}
			}.
		\end{align*}
		Since 
		$
			f \in \holderBounded{k,\alpha}\del*{
				\cl{\Omega}
			} 
		$
		is arbitrary, 
		we have
		$
			\holderBounded{k,\alpha}\del*{ 
				\cl{\Omega}
			}
			\hookrightarrow
			\multipoint{k+\alpha}
			\integrable{\infty}\del*{
				\Omega, \#
			},
		$
		as claimed. \qedhere
	\end{enumerate}

\end{proof}

\section{Higher-Order Sobolev Spaces on Metric Spaces}
\label{sec::metric_spaces}

In the last section of the paper, we will explain how we can use the results from the~previous section to define higher-order Sobolev and H\"{o}lder spaces on metric spaces.

Let us note that although the definition of 
$\multipoint{s} \mc{F}(X)$ spaces present in 
\Cref{def::multipointwise_Banach_function_space}
does not use gradients, 
it does involve norms to explain the meaning of 
$ \polygen{ \gen{s} }(\bx) $ 
present within the definition of 
$ \bD^s_{\measureAlone}(f)$.
In order to~be~able to define an analogous family on~an~arbitrary metric space, 
we need to explain what $\polygen{\gen{s}}(\bx)$ would mean in that context.
Assuming that the structure of functions 
$ \polyset[A]{S} $, 
$ \polyset{S} $,
$ \polyfam{ \cP }$,
$ \polynum{ j }$,
and 
$ \polygen{ \gen{s } }$
would be~the~same as in 
\Cref{def::poly_normed_space},
we only need to define $ \polyset[A]{S} $.
The following lemma will provide a suggestion for how to do so.
\begin{lemma}
\label{lem::P_A^S_in_inner_product_spaces_in_terms_of_distances}
    Let $\del*{ V, \inner{ \blank, \blank[:]} }$ be an inner product space 
    and let $\normAlone$ be the norm induced by~%
    the~inner product. 
    Let $k \in \bN$, $S \subseteq [k]$ be nonempty, 
    and $A \subseteq [k] \setminus S$. 
    Fix 
    $
        \bx = \set*{x_{I}}_{ I \subseteq [k] }
        \subseteq 
        V.
    $
    Then, denoting 
    $
    	\diff_{ I, J \upmapsto A}^S
    	\coloneq 
    	\diff_{ I \upmapsto A }^S
    		\diff_{J \upmapsto A }^S,
    $
    \begin{equation*}
        \poly[A]{S}(\bx)^2
        =
        \norm*{
            \diff_{ I \upmapsto A }^{ S }
                x_I
        }^2
        =
        -
        \frac{1}{2}
        \diff_{ I, J \upmapsto A }^S
            \norm{x_I - x_J}^2.
    \end{equation*}
\end{lemma}
\begin{proof}
    The first equality is just the definition of 
    $ \poly[A]{S}(\bx) $
    we have used so far.
    For the second equality,
    let us recall that for any $x, y \in V$, we have
    \begin{equation*}
    	\norm*{
    		x-y
    	}^2
    	=
    	\inner{
    		x-y, x-y
    	}
    	=
    	\inner{x,x} - \inner{x,y} - \inner{y,x} + \inner{y,y}
    	=
    	\norm{x}^2- \inner{x,y} - \inner{y,x} + \norm{y}^2.
    \end{equation*}
    Thus,
    \begin{equation}
    \label{proofeq::lem::P_A^S_in_inner_product_in_terms_of_distances::consequence_of_Law_of_Cosines}
    	\forall x, y \in V
    \qquad
    	\inner{x,y} + \inner{y,x}
    	=
    	\norm{x}^2 + \norm{y}^2 - \norm{x-y}^2.
    \end{equation}
   	Also, let us note that
   	\begin{equation}
   	\label{proofeq::lem::P_A^S_in_inner_product_in_terms_of_distances::changing_the_order_in_double_diff}
   		\diff_{ I, J \upmapsto A }^{ S }
	   		\inner*{ 
	   			x_I - x_A
	   			,    
	   			x_J - x_A
	   		}
 		=
   		\diff_{ J,I \upmapsto A }^{ S }
	   		\inner*{ 
	   			x_I - x_A
	   			,    
	   			x_J - x_A
	   		}
		=
   		\diff_{ I, J \upmapsto A }^{ S }
			\inner*{ 
				x_J - x_A
				,    
				x_I - x_A
			},
   	\end{equation}
   	where the first equality is just a change in the order of summation, and the second, a renaming of dummy variables $(I,J) \mapsto (J,I)$.
   	
    Now, let us note that by \Cref{lem::diff_constant_term},
    we have
    \begin{align*}
        \polyset[A]{S}( \bx )^2
    &= 
        \norm*{
            \diff_{ I \upmapsto A }^{ S }
                x_I
        }^2
    \\
    &=
        \inner*{ 
            \diff_{ I \upmapsto A }^{ S }
                x_I,
            \diff_{ J \upmapsto A }^{ S }
                x_J
        }
    	=
        \inner*{ 
            \diff_{ I \upmapsto A }^{ S }
                x_I
            - 
            \diff_{ I \upmapsto A }^{ S }
                x_A
            ,
            \diff_{ J \upmapsto A }^{ S }
                x_J
            -
            \diff_{ J \upmapsto A }^{ S }
                x_A
        }
	\\
    &
	\phantom{
		=\;
		\inner*{ 
			\diff_{ I \upmapsto A }^{ S }
			   x_I,
			\diff_{ J \upmapsto A }^{ S }
			   x_J
		}
	}   
    =
        \inner*{ 
            \diff_{ I \upmapsto A }^{ S }
                \del*{ x_I - x_A }
            ,
            \diff_{ J \upmapsto A }^{ S }
                \del*{ x_J - x_A }
        }
    	=
        \diff_{ I, J \upmapsto A }^{ S }
        \inner*{ 
            x_I - x_A
            ,    
            x_J - x_A
        }.
    \end{align*}
    Therefore, thanks to  \eqref{proofeq::lem::P_A^S_in_inner_product_in_terms_of_distances::changing_the_order_in_double_diff},
   	\begin{align*}
   		\polyset[A]{S}( \bx )^2
      	&=      	
      	\frac{1}{2}
      	\diff_{ I, J \upmapsto A }^{ S }
      		\del*{
		      	 \inner*{ 
		      	 	x_I - x_A
		      	 	,    
		      	 	x_J - x_A
		      	 } 
		      	 +
		      	 \inner*{ 
		      	 	x_J - x_A
		      	 	,    
		      	 	x_I - x_A
		      	 } 
		    }
        ,
    \intertext{
        which, by \eqref{proofeq::lem::P_A^S_in_inner_product_in_terms_of_distances::consequence_of_Law_of_Cosines},
    }
        &=
        \frac{1}{2}
        \diff_{ I, J \upmapsto A }^{ S }
            \del*{
                \norm*{x_{ I } - x_{ A } }^2
                + \norm*{x_{J} - x_{ A }}^2
                - \norm*{x_{ I } - x_{ J }}^2
            }
    \\
        &=
        \frac{1}{2} 
        \del*{
            \diff_{ I \upmapsto A }^{ S }
                \diff_{ J \upmapsto A }^{ S }
                    \norm*{x_{ I } - x_{ A } }^2
            +
            \diff_{ I \upmapsto A }^{ S }
                \diff_{ J \upmapsto A }^{ S }
                    \norm*{x_{ J } - x_{ A } }^2
            -
            \diff_{ I, J \upmapsto A }^{ S }
                \norm*{x_{ I } - x_{ J }}^2
        },
    \intertext{
    	which,
    	since 
        by \Cref{lem::diff_constant_term}
        we have
        $
            \diff_{ J \upmapsto A }^{ S }
                \norm*{x_{ I } - x_{ A } }^2
            =
            0
        $
        and
        $
            \diff_{ I \upmapsto A }^{ S }
            	\del*{
		            \diff_{ J \upmapsto A }^{ S }
			            \norm*{x_{ J } - x_{ A } }^2
			    }
            =
            0,
        $
    }
        &=
        \frac{1}{2} 
        \del*{
            \diff_{ I \upmapsto A }^{ S }
                0
            +
            0
            -
            \diff_{ I, J \upmapsto A }^{ S }
                \norm*{x_{ I } - x_{ J }}^2
        }
    	=
        - \frac{1}{2}
        \diff_{I,J \upmapsto A }^{ S }
            \norm*{ x_I - x_J}^2,
    \end{align*}
    as claimed.
\end{proof}
The obtained result suggests the following definition for 
$ \polyset[A]{S} $ when the underlying space is a~metric space.
\begin{definition}
\label{def::PAS_for_metric_spaces}
    Let $(X, \metricAlone)$ be a metric space,
    $ k \in \bN $, 
    $ S \subseteq [k] $ be nonempty, 
    and $A \subseteq [k] \setminus S$.
    Let 
    $
        \bx = \set*{ x_I }_{ I \subseteq [k] }
        \subseteq 
        X.
    $
    We define
    \begin{equation*}
        \smash{
        \polyset[A]{S}( \bx )
        \coloneq
        \sqrt{
            \frac{1}{2}
            \abs*{
                \diff_{ I, J \upmapsto A }^{ S }
                    \metric{ x_I, x_J }^2
            }
        }.}
    \end{equation*}
\end{definition}
Note that the absolute value in the above definition is necessary since, 
for~an~arbitrary metric space,
we~cannot guarantee that 
$
	\diff_{ I, J \upmapsto A }^{ S }
         \metric{ x_I, x_J }^2
$ 
will have the same sign for all tuples 
$
	\set*{ x_I }_{ I \subseteq [k]}
	\subseteq 
	X.	
$

Using the above definition of $\polyset[A]{S}$ and assuming that the other functions present in \Cref{def::poly_normed_space} 
are defined as they were for normed spaces, 
we obtain a definition for 
$\polygen{ \gen{s} }$
when the underlying spaces is a metric space. 
(Also, we again put $\polygen{\gen{0}} \equiv 1$.)
Moreover, the resulting definition is equivalent to the previously used one if the tuple $\bx$ consists of elements of an inner product space.
With this, we can now define analogs for~%
$ \bD^s_{\measureAlone}(f) $ and $ \bD^s(f)$,
where $f$ is a function defined on a metric space.
\begin{definition}[Families $\bD_{ \measureAlone }^{ s }(f)$ and $\mF_{\measureAlone}^{s}(f,G)$]
\label{def::s-hypergradient(metric)}
    Let 
    $\del{X, \metricAlone}$
    be a metric space. Let $s \ge 0$ and $k \in \bN_0$ be such that $s \in \intoc{k-1, k}$.
    Let~$\measureAlone$ be~a~%
    measure on $X$. 
    For a given measurable function 
    $f \colon X \to \bR$,
	we will denote by 
    $ \bD_{\measureAlone}^s(f)$ the family of~measurable functions
    $G \colon X \to \intcc{0,\infty}$ 
    such that
    \begin{equation}
	\label{eqdef::def::s-hypergradient(metric)::m-upper} 
        \measureAlone
        \forall 
            \bx = \set*{x_{I}}_{ I \subseteq [k] }
            \subseteq X
    \qquad
        \abs*{
        	\diff_{ I = \emptyset}^{ [k] }
        		f\del*{ x_I }
        }
        \le 
        \polygen{ \gen{s} }(\bx)
        \sum_{ I = \emptyset }^{ [k] }
            G(x_I).
    \end{equation}
    For a given $G \in \bD_{ \measureAlone }^{ s }(f)$, 
	we will denote the family of all subsets 
	of $X$ of
	full measure on~which inequality \eqref{eqdef::def::s-hypergradient(metric)::m-upper} is satisfied by 
	$\mF_{\measureAlone}^{s}(f,G)$.
    Also, as before, if $\measureAlone$ is the counting measure $\#$, we~will~write
    $ \bD^s(f)$ instead of $\bD_{\#}^s(f)$.   
\end{definition}
Using the above definition, we can finally define higher-order Sobolev spaces on~metric spaces with~a~measure.
\begin{definition}[$\multipoint{s,p}\del*{ X, \metricAlone, \measureAlone }$ space]
\label{def::mulitpointwise_on_metric}
	Let $\del{X, \metricAlone}$ be a metric space,
	$ \measureAlone $ be a measure on $X$,
	$ s \ge 0$, 
	and $p \in \intcc{1, \infty}$.
	We define the~family
	\begin{equation*}
		\multipoint{s,p}\del*{ X, \metricAlone, \measureAlone }
		\coloneq 
		\setc*{
			f \in \integrable{p}\del*{ X, \measureAlone }
		}{
			\text{%
				there exists 
				$ G \in \bD_{\measureAlone}^{s}(f) $
				such that
				$ G \in \integrable{p}\del*{ X, \measureAlone } $%
			}
		}
	\end{equation*}
	and endow it with the following seminorm and norm:
	\begin{equation*}
		\forall f \in \multipoint{s,p}\del*{ X }
	\qquad
		\seminorm{ f }_{ 
			\multipoint{s,p}\del*{ X }
		}
		\coloneq 
		\inf_{ G \in \bD_{\measureAlone}^{s}(f) }
			\norm{ G }_{ \integrable{p}\del*{ X} }
	\quad \text{and} \quad 
		\norm{ f }_{ 
			\multipoint{s,p}\del*{ X}
		}
		\coloneq 
		\norm{ f }_{ \integrable{p}\del*{ X}}
		+
		\seminorm{ f }_{ 
			\multipoint{s,p}\del*{ X} 
		}.
	\end{equation*}
	Also, much like we already did above,
	we will also write 
	$ \multipoint{s,p}(X)$
	instead of
	$ 
		\multipoint{s,p}\del*{ 
			X, \metricAlone, \measureAlone 
		}  
	$
	if the metric $\metricAlone$
	and~measure
	$ \measureAlone $
	should be inferable from the context.
\end{definition}
Let us note that if $X$ is a subset of an inner product space and the metric $\metricAlone$ is induced by the inner product, then 
$
	\multipoint{s,p}\del*{ X, \metricAlone, \measureAlone }
	\cong
	\multipoint{s}
	\integrable{p}
	\del*{ X, \measureAlone },
$
where the latter space is defined as in 
\Cref{def::multipointwise_Banach_function_space}.
Also, much like in~the mentioned definition, 
we could have defined the space
$ \multipoint{s}\mc{F}\del*{X, \measureAlone} $
for an arbitrary Banach function space
$ \mc{F}\del*{X, \measureAlone} $
and obtained the same equality of spaces.

Suppose now that we consider an arbitrary metric space
$ \del*{ X, \metricAlone }$
and a measure $\measureAlone$ on it. 
First, let us notice that 
$
	\multipoint{0,p}(X)
	\cong 
	\integrable{p}(X).
$
Indeed, since
$ 
	\polygen{\gen{0}} \equiv 1,
$
for all 
$ f \in \integrable{p}(X)$
we have
$ \abs{f} \in \bD_{\measureAlone}^{0}(f) $.
Also, if 
$ G \in \bD_{\measureAlone}^{0}(f) $,
then
$
	\abs{f} \le G
$
$ \measureAlone $-almost everywhere,
so
$
	\norm{G}_{ \integrable{p}(X) }
	\ge 
	\norm{f}_{ \integrable{p}(X) }.
$
Thus, for all $f \in \integrable{p}(X)$,
\begin{equation*}
	\norm{ f }_{ \multipoint{0,p}(X) }
	=
	\norm{ f }_{ \integrable{p}(X) }
	+
	\inf_{ G \in \bD_{\measureAlone}^{0}(f) }
		\norm{G}_{ \integrable{p}(X) }
	=
	\norm{ f }_{ \integrable{p}(X) }
	+
	\norm{ f }_{ \integrable{p}(X) }
	=
	2 \norm{f}_{ \integrable{p}(X) }.
\end{equation*}

Secondly, note that, 
for all
$ s \in \intoc{0,1}$
and  
$ p \in \intcc{1, \infty}$,
we have 
$
	\multipoint{s,p}\del*{X, \metricAlone, \measureAlone }
	\cong
	\hajlasz{s,p}\del*{X, \metricAlone, \measureAlone },
$
where 
$
	\hajlasz{s,p}\del*{X, \metricAlone, \measureAlone }
	\coloneq
	\hajlasz{1,p}\del*{X, \metricAlone^s, \measureAlone }.
$
(%
	Here,
	$ \metricAlone^s $
	is a metric on $X$
	defined by the formula
	$
		\metricAlone^s(x,y)
		\coloneq 
		\metric{x,y}^s
	$
	for all $x, y \in X$.%
)

Indeed, let us observe that, for every
$ 
	\bx = \set*{ x_I }_{ I \subseteq [1] }
	\subseteq 
	X,
$
we have
\begin{align*}
	2\polyset[\emptyset]{ \set*{1} }(\bx)^2
	=
	\smash[t]{
	\abs*{
		\diff_{I, J \upmapsto \emptyset }^{ \set*{1} }
			\metric{x_I, x_J}^2
	}}
&=
	\abs*{
		\diff_{I = \emptyset }^{ \set*{1} }
			\diff_{J = \emptyset }^{ \set*{1} }
				\metric{x_I, x_J}^2
	}
\\
&=
	\abs*{
		\metric{
			x_{\set*{1}},
			x_{\set*{1}}
		}^2
		-
		\metric{
			x_{\set*{1}},
			x_{\emptyset}
		}^2
		-
		\metric{
			x_{\emptyset},
			x_{\set*{1}}
		}^2
		+
		\metric{
			x_{\emptyset},
			x_{\emptyset}
		}^2
	}
\\
&=
	\abs*{
		- 2 \metric{ x_{ \set*{1}}, x_{\emptyset } }^2
	}
\\
&=
	2 \metric{  x_{ \set*{1}}, x_{\emptyset } }^2,
\end{align*}
so, when
$ s \in \intoc{0,1}$,
\begin{equation*}
	\polygen{\gen{s}}(\bx)
	=
	\polynum{s}(\bx)
	=
	\polyset{ \set*{1} }(\bx)^{s}
	=
	\polyset[\emptyset]{ \set*{1} }(\bx)^{s}
	=
	\metric{  x_{ \set*{1}}, x_{\emptyset } }^s.
\end{equation*}
Furthermore, denoting
\begin{equation*}
	\holderBounded{0,s}\del*{ X }
	\coloneq 
	\setc*{
		f \in \continuousAndBounded\del*{ X }	
	}{
		\exists C \ge 0
	\quad 
		\forall x, y \in X
	\qquad
		\abs*{f(x)-f(y)}
		\le 
		C \metric{x,y}^s
	},
\end{equation*}
and endowing the above space with the following seminorm and norm:
\begin{equation*}
	\forall f \in \holderBounded{0,s}\del*{ X }
	%
\qquad
	\seminorm{ f }_{ \holderBounded{0,s}\del*{ X } }
	\coloneq 
	\sup_{
		\substack{ 
			x, y \in X 
		\\
			x \ne y
		}
	}
		\frac{ \abs{ f(x) - f(y) } }{ \metric{x,y}^s }
\quad \text{and} \quad 
	\norm{ f }_{ \holderBounded{0,s}\del*{ X } }
	\coloneq 
	\norm{ f }_{ \infty }
	+
	\seminorm{ f }_{ \holderBounded{0,s}\del*{ X } },
\end{equation*}
we have
$
	\multipoint{s,\infty}\del*{X, \# }
	\cong
	\holderBounded{0,s}\del*{ X }
$
with 
$
	2\seminormAlone[{\multipoint{s,\infty}\del*{X, \# }}]
	=
	\seminormAlone[{\holderBounded{0,s}\del*{ X }}].
$
Indeed, this fact can be~proved analogously as a similar fact for spaces
$ \multipoint{s}\integrable{\infty}\del*{\Omega, \#} $
and 
$ \holderBounded{0,s}\del*{\cl{\Omega}}$,
where $\Omega$ is~an~open subset of a Euclidean space
(see \Cref{rem::up_to_first-order_characterization_holder}).
Because of this
and 
\Cref{thm::higher-order_characterization_of_holder_spaces},\footnote{Note that when $X$ is a subset of an Euclidean space, then 
$
	\multipoint{s,\infty}\del*{X, \#}
	=
	\multipoint{s}
	\integrable{\infty}\del*{X, \#},
$
where the latter family is defined as in 
\Cref{def::multipointwise_Banach_function_space}.
}
we~can treat the
$
	\multipoint{s,\infty}\del*{X, \# }
$
spaces, where
$ s \in \intoo{0,\infty}$,
as a possible definition of higher-order H\"{o}lder spaces on~metric spaces.

As such, we can hope that the family of spaces
$
	\multipoint{s,p}\del*{ X },
$
where
$ s \in \intco{0, \infty}$
and 
$ p \in \intcc{1,\infty}$,
might exhibit many properties
similar to the properties of Hajłasz--Sobolev spaces
$ \hajlasz{1,p}(X)$
or H\"{o}lder spaces
$ \holderBounded{0,\alpha}(X)$,
where $\alpha \in \intoc{0,1}$.
In the penultimate result of the paper, 
we will show that $\multipoint{s,p}(X)$ spaces are Banach.
As we will see, the~reasoning will be 
quite similar to~the~one present in 
Hajłasz's article \cite[Theorem 3]{Hajlasz}, 
further justifying our hopes.
\begin{theorem}
	Let $\del*{ X, \metricAlone } $
	be a metric space and 
	$\measureAlone$ be a measure on $X$.
	Then, 
	for all 
	$s \in \intco{0, \infty}$
	and 
	$p \in \intcc{1, \infty}$,
	$\multipoint{s,p}(X)$
	is a Banach space.
\end{theorem}
\begin{proof}
	First, let us note that for 
	$ s = 0 $
	this fact is clear since
	$
		\multipoint{s,p}(X)
		\cong 
		\integrable{p}(X)
	$
	for all $p \in \intcc{1, \infty}$.
	
	Moving forwards, suppose that $s > 0$ and 
	$p \in \intcc{1, \infty}$.
	It is clear that 
	$\multipoint{s,p}(X)$ 
	is a normed space, 
	so~it~remains to show that it is complete.
	Let
	$ \del{ f_n }_{n \in \bN} $
	be a Cauchy sequence in
	$\multipoint{s,p}(X)$.
	Then it is also a~Cauchy sequence 
	in 
	$ \integrable{p}(X)$,
	so it is convergent in 
	$ \integrable{p}(X)$
	to some
	$ f \in \integrable{p}(X)$.

	Then
	$ \del{ f_n }_{n \in \bN} $
	has a subsequence
	(denoted using the same indices as the original one)
	such that
	\begin{equation*}
		f_n \to f 
	\quad
		\text{$\measureAlone$-almost everywhere, and, for all $n \in \bN$,}
	\quad
		\norm{
			f_{n+1} - f_n
		}_{ \multipoint{s,p}(X) } 
		\le 
		2^{-n}.
	\end{equation*}
	Let $F'$ be a set of full measure such that
	$ f_n(x) \to f(x)  $
	for all $x \in F'$.
	Also, for all $n \in \bN$, 
	let 
	$
		G_n' \in \bD^s_{\measureAlone}\del*{ f_{n+1} - f_n}
	$
	be such that
	$
		\norm{ G_n' }_{
			\integrable{p}(X)
		}
		\le 
		\norm{
			f_{n+1} - f_n
		}_{ \multipoint{s,p}(X) }
		+
		2^{-n}
		\le
		2^{-n+1}.
	$
	Then, for all such $n$, 
	fix
	$ 
		F_n \in \mF^s_{\measureAlone}\del*{
			f_{n+1} - f_n, G_n'
		}
	$
	and denote
	$
		F
		\coloneq 
		F' 
		\cap 
		\bigcap_{ n = 1}^{ \infty }
			F_n.
	$
	Clearly, 
	$ F $ is of full measure.
	Next, for all $n \in \bN$, 
	define
	$
		G_n
		\coloneq 
		\sum_{ \ell = n }^{ \infty }
			G_{\ell}'
	$
	and note that
	\begin{equation*}
		\norm{ G_n }_{ \integrable{p}(X) }
		=
		\norm*{ 
			\sum_{ \ell = n }^{ \infty }
				G_{\ell}'
		}_{ \integrable{p}(X) }
		\le 
		\sum_{ \ell = n }^{ \infty }
			\norm*{ 
				G_{\ell}'
			}_{ \integrable{p}(X) }
		\le 
		\sum_{ \ell = n }^{ \infty }
			2^{-\ell+1}
		=
		2^{-n+2}
		\xrightarrow{n \to \infty}
		0.
	\end{equation*}

	Let $k \in \bN$ be such that
	$ s \in \intoc{k-1,k}$
	and fix
	$
		\bx = \set*{ x_I }_{ I \subseteq [k] }
		\subseteq
		F.
	$
	Then for all 
	$ I \subseteq [k] $
	we have
	$ f_n\del*{ x_I } \to f \del*{ x_I } $
	as $n \to  \infty $.
	Fix $m \in \bN$ and notice that
	\begin{align*}
		\abs*{
			\diff_{I = \emptyset }^{ [k] }
				\del*{ 
					f\del*{ x_I } 
					-
					f_m\del*{ x_I }
				}
		}
	\xleftarrow{n \to \infty}
	&
		\abs*{
			\diff_{I = \emptyset }^{ [k] }
				\del*{ 
					f_{n+1}\del*{ x_I } 
					-
					f_{m}\del*{ x_I }
				}
		}
	\\
	&=
		\abs*{
			\diff_{I = \emptyset }^{ [k] }
				\sum_{ \ell = m }^{ n }	
					\del*{ 
						f_{\ell+1}\del*{ x_I } 
						-
						f_{\ell}\del*{ x_I }
					}
		}
	\\
	&\le
		\sum_{ \ell = m }^{ n }	
			\abs*{
				\diff_{I = \emptyset }^{ [k] }
					\del*{ 
						f_{\ell+1}\del*{ x_I } 
						-
						f_{\ell}\del*{ x_I }
					}
			}
	\\
	&\quad \quad \le 
		\sum_{ \ell = m }^{ n }	
			\polygen{ \gen{s} }(\bx)
			\sum_{I = \emptyset }^{ [k] }
				G_{\ell}'\del*{x_I}
		\le 
		\polygen{ \gen{s} }(\bx)
			\sum_{I = \emptyset }^{ [k] }
				G_m\del*{x_I}.
	\end{align*}
	Thus,
	$
		\abs*{
			\diff_{I = \emptyset }^{ [k] }
				\del*{ 
					f\del*{ x_I } 
					-
					f_m\del*{ x_I }
				}
		}
		\le 
		\polygen{ \gen{s} }(\bx)
			\sum_{I = \emptyset }^{ [k] }
				G_m\del*{x_I}
	$
	for all
	$
		\bx = \set*{ x_I }_{ I \subseteq [k] }
		\subseteq
		F.
	$	
	Since $F$ is of full measure, we~have
	$
		G_m \in 
		\bD_{\measureAlone}^{s}(f - f_m).
	$
	Thus,
	since
	$ G_m \in \integrable{p}(X) $,
	we have
	$
		f - f_m 
		\in 
		\multipoint{s,p}(X)
	$
	and 
	$
		f \in \multipoint{s,p}(X).
	$
	Finally, note that
	\begin{align*}
		\norm*{
			f - f_m 
		}_{ \multipoint{s,p}(X) }
		=
		\norm*{
			f - f_m 
		}_{ \integrable{p}(X) }	
		+
		\inf_{ G \in \bD_{\measureAlone}^s\del*{ f - f_m} }
			\norm*{
				G
			}_{ \integrable{p}(X) }	
	&
	\\
		\le 
		\norm*{
			f - f_m 
		}_{ \integrable{p}(X) }	
		+
		\norm*{
			G_m
		}_{ \integrable{p}(X) }	
	&\xrightarrow{m \to \infty}
		0.		
	\end{align*}
	Thus, 
	$ f_n \to f $
	in 
	$  \multipoint{s,p}(X) $.
	As
	$ \del*{ f_n }_{ n \in \bN}$
	is a subsequence 
	of the original Cauchy sequence, 
	the~original sequence also converges to $f$, proving the completeness of
	$  \multipoint{s,p}(X) $.	 
\end{proof}
In the last theorem of the paper, we will prove that, if 
$ s \ge 1 $ and
$p \in \intcc{1,\infty}$,
then
$  
	\multipoint{s,p}\del*{X}
	\hookrightarrow
	\hajlasz{1,p}\del*{X}.
$
In~consequence,
we also have that
$  
	\multipoint{s,\infty}\del*{X, \#}
	\hookrightarrow
	\holderBounded{0,1}\del*{X};
$
indeed, this follows from the fact that
$
	\hajlasz{1,\infty}\del*{ X, \# }
	\cong 
	\holderBounded{0,1}\del*{X}.
$
However, before we state the theorem, 
we will first prove the following lemma.
\begin{lemma}
\label{lem::folding_hypercubes_two_elements}
    Let $(X, \metricAlone)$
    be a metric space and
    $
        x, y \in X.
    $
    Let $s > 1$
    and $ k \in \bN$
    be such that
    $ s \in \intoc{k-1,k}$.
    Define
    $
        \bx = \set{x_I}_{I\subseteq[k]}
        \subseteq 
        X
    $
    by the formula
    \begin{equation*}
        \forall I \subseteq [k]
    \qquad 
        x_I
        \coloneq 
        \begin{cases}
            x
        &
            \text{if $\card{I}$ is even},
        \\
            y
        &
            \text{otherwise}.
        \end{cases} 
    \end{equation*}
    Then 
    \begin{equation*}
        \polygen{ \gen{s} }\del*{ \bx }
        \le 
        2^{k^2}
        \partitions{ k }
        \del*{
            \metric{x,y}
            +
            \metric{x,y}^s
        }.
    \end{equation*} 
\end{lemma}
\begin{proof}
    First, let us notice that, 
    for all $I, J \subseteq [k]$,
    we have that
    $
        \metric{x_I, x_J}
        \le 
        \metric{x, y}.
    $
    In consequence, 
    for all 
    $ S \subseteq [k]$
    that are nonempty and 
    $ A \subseteq [k] \setminus S$,
    we have that
    \begin{align*}
        2
        \polyset[A]{S}\del*{ \bx }^2
        =
        \abs*{
            \diff_{ I \upmapsto A}^{S}
            	\diff_{ J \upmapsto A}^{S}
	                \metric{x_I, x_J}^2
        }
    	=
        \abs*{
            \diff_{ I = \emptyset }^{S}
            	\diff_{ J = \emptyset }^{S}
	                \metric{x_{A \cup I}, x_{A \cup J} }^2
        }  
    \hspace{-1cm}&
    \\    
    	\le 
        \sum_{ I = \emptyset }^{S}
           	\sum_{ J = \emptyset }^{S}
             \metric{x_{A \cup I}, x_{A \cup J} }^2
    &\le 
        \sum_{ I = \emptyset }^{S}
           	\sum_{ J = \emptyset }^{S}
	      \metric{x,y}^2
	\\
    &=
        \del*{ 2^{ \card{S} } }^2 \,
        \metric{x, y}^2
        \le 
        2
        \del*{ 2^{ \card{S} } }^2 \,
        \metric{x, y}^2.
    \end{align*}
    Thus, 
    $
        \polyset[A]{S}\del*{ \bx }
        \le 
        2^{ \card{S} }
        \metric{x,y}
    $
    for all such $S$ and $A$.
    Hence, for all 
    nonempty 
    $ S \subseteq [k]$,
    we~have
    \begin{equation*}
        \polyset{S}\del*{ \bx }
        =
        \sum_{ A = \emptyset }^{ [k] \setminus S }
           \polyset[A]{S}\del*{ \bx }
        \le 
        \sum_{ A = \emptyset }^{ [k] \setminus S }
            2^{ \card{S} }
            \metric{x,y} 
        =
        2^k
        \metric{x, y}.
    \end{equation*}
    In consequence, for all 
    $ \cP \in \partitions{ [k] }$
    such that
    $ \card*{ \cP } \le k-1$,
    we have
    \begin{equation*}
        \polyfam{ \cP }\del*{ \bx }
        =
        \prod_{ S \in \cP }
            \polyset{ S}\del*{ \bx }
        \le 
        \prod_{ S \in \cP }
            2^k \metric{x, y}
        =
        \del*{
            2^k \metric{x, y}
        }^{ \card{ \cP } }
        \le 
        2^{k^2}
        \del*{
            \metric{x, y}
            + 
            \metric{x, y}^s
        }
        .
    \end{equation*}
    Therefore,
    \begin{align*}
        \sum_{ j = 1 }^{ k-1 }
            \polynum{ j }\del*{ \bx }
    &=
        \sum_{ j = 1 }^{ k-1 }
            \sum_{ \cP \in \partitions[j]{ [k] } }
            \polyfam{ \cP }\del*{ \bx }
    \\    
    &\le 
        \sum_{ j = 1 }^{ k-1 }
            \sum_{ \cP \in \partitions[j]{ [k] } }     
                    2^{k^2}
        \del*{
            \metric{x, y}
            + 
            \metric{x, y}^s
        }
        =
        \del*{
            \sum_{ j = 1}^{ k-1}
                \card*{
                    \partitions[j]{ [k] }
                }
        }
        2^{k^2}
        \del*{
            \metric{x, y}
            + 
            \metric{x, y}^s
        }.
    \end{align*}
    Also, 
    \begin{align*}
        \polynum{s}\del*{ \bx }
    &=
    	\smash[t]{
        \polyset{ \set{k} }\del*{ \bx }^{ s-k+1 }
        \prod_{ j =1 }^{ k-1}
            \polyset{ \set{j}}\del*{ \bx } }
    \\
    &\le 
        \del*{ 2^k \metric{x, y} }^{ s-k+1 }
        \prod_{ j =1 }^{ k-1}
            2^k \metric{x, y}    
        =
        \del*{ 2^k \metric{x, y} }^{ s }
    	\le 
        2^{k^2}
        \del*{
            \metric{x, y}
            + 
            \metric{x, y}^s
        }.        
    \end{align*} 

    Finally,
    \begin{align*}
        \polygen{ \gen{s} }\del*{ \bx }
    &=
        \polynum{ s }\del*{ \bx }
        +
        \sum_{ j = 1}^{ k-1 }
            \polynum{ j }\del*{ \bx }
    \\
    &\le 
        2^{k^2}
        \del*{
            \metric{x, y}
            + 
            \metric{x, y}^s
        }
        +
        \del*{
            \sum_{ j = 1}^{ k-1}
                \card*{
                    \partitions[j]{ [k] }
                }
        }
        2^{k^2}
        \del*{
            \metric{x, y}
            + 
            \metric{x, y}^s
        }
    \\
    &=
        2^{k^2}
        \partitions{ k }
        \del*{
            \metric{x, y}
            + 
            \metric{x, y}^s
        },        
    \end{align*}
    where in the final equality we use the fact that
    \begin{equation*}
        1
        +
        \sum_{ j = 1}^{ k-1}
            \card*{
                \partitions[j]{ [k] }
            }
        =
        \card*{ 
            \partitions[k]{ [k] }
        }
        +
        \sum_{ j = 1}^{ k-1}
            \card*{
                \partitions[j]{ [k] }
            }
        =
        \sum_{ j = 1}^{ k}
            \card*{
                \partitions[j]{ [k] }
            }
        =
        \card*{ 
            \partitions{ [k] }
        }
        =
        \partitions{ k }.
    \end{equation*}
    Thus, the proof is complete.
\end{proof}
\begin{theorem}
    Let 
    $ \del*{ X, \metricAlone }$
    be a metric space, 
    $ \measureAlone $ be
    a measure on $X$, 
    $ s \ge 1$, 
    and 
    $ p \in \intcc{1, \infty}$.
    Then
    $
        \multipoint{s,p}\del*{X}
        \hookrightarrow
        \hajlasz{1,p}\del*{X}.
    $
\end{theorem}
\begin{proof}
    First, let us recall that 
    $
        \multipoint{1,p}\del*{X}
        \cong
        \hajlasz{1,p}\del*{X},
    $
    so we only need to prove the theorem for~%
    $ {s > 1.} $
    
    Fix 
    $ s > 1 $,
    $ f \in \multipoint{s,p}\del*{X}$,
    and $\eps > 0$.
    Then
    there exists 
    $ G \in \bD_{ \measureAlone }^{s}\del*{ f }$ 
    for which we have
    $
        \norm{ G }_{ 
            \integrable{p}(X)
        }
        \le 
        \norm{
            f
        }_{ \multipoint{s,p}\del*{X} }
        +
        \eps.
    $
    Let
    $ 
        F \in \mF_{ \measureAlone}^{s}\del*{ f, G}
    $   
    and 
    $ x, y \in F$.
    If 
    $ \metric{x, y} \ge 1$, 
    then 
    \begin{align*}
        \abs*{f(x) - f(y)}
        \le 
            \abs{f(x)} + \abs{f(y)}
    &\le 
        \metric{x, y}
         \del*{
            \abs{f(x)} + \abs{f(y)}
        } 
    \\
    &\le 
        \metric{x, y}
         \del*{
            \abs{f(x)} 
            +
            2^{k^2+1}
            \partitions{k}
            G(x)
            + 
            \abs{f(y)}
            +
            2^{k^2+1}
            \partitions{k}
            G(y)
        }.  
    \end{align*}
    Next, suppose that
    $ \metric{x, y} \le 1$
    and define 
    $
        \bx = \set{x_I}_{ I \subseteq[k]}
    $
    by the formula
    \begin{equation*}
        \forall I \subseteq [k]
    \qquad 
        x_I
        \coloneq 
        \begin{cases}
            x
        &
            \text{if $\card{I}$ is even},
        \\
            y
        &
            \text{otherwise}.
        \end{cases} 
    \end{equation*}
    Then, since by
    \Cref{lem::between_A_and_B_cardinalities}
    there are
    $ 2^{k-1}$ subsets of 
    $ [k] $ that have even cardinality
    and the same number of~subsets with odd cardinality,
    \begin{align*}
        2^{k-1}
        \abs*{
            f(x) - f(y)
        }
        =
        \abs*{
            \sum_{\substack{
                I = \emptyset
            \\
                \card*{I} \text{ is even}
            }}^{ [k] }
                f\del*{x_I}
            -
            \sum_{\substack{
                I = \emptyset
            \\
                \card*{I} \text{ is odd}
            }}^{ [k] }
                f\del*{x_I}
        }
    &=
        \abs*{
            \sum_{ I = \emptyset}^{ [k] }
                (-1)^{ \card{I} }
                f\del*{x_I}
        }
    \\
    &=
    \smash{
        \abs*{
            \sum_{ I = \emptyset}^{ [k] }
                (-1)^{ \card{[k] \setminus I} }
                f\del*{x_I}
        } 
    }     
    	=
        \abs*{
            \diff_{ I = \emptyset }^{ [k] }
                f\del*{x_I}        
        }.
    \end{align*}
    Similarly, 
    \begin{equation*}
        \sum_{ I = \emptyset}^{ [k] }
            G\del*{ x_I }
        =
        \sum_{\substack{
            I = \emptyset
        \\
            \card*{I} \text{ is even}
        }}^{ [k] }
            G\del*{x_I}
        +
        \sum_{\substack{
            I = \emptyset
        \\
            \card*{I} \text{ is odd}
        }}^{ [k] }
            G\del*{x_I}     
        =
        2^{k-1}
        \del*{
            G(x)
            +
            G(y)
        }.
    \end{equation*}
    Therefore, since
    $ G \in \bD_{\measureAlone}^{s}(f) $,
    $ 
        F \in \mF_{ \measureAlone}^{s}\del*{ f, G}
    $,
    and 
    $ \bx \subseteq F$,
    \begin{equation*}
        2^{k-1}
        \abs*{
            f(x) - f(y)
        }
        =
        \abs*{
            \diff_{ I = \emptyset }^{ [k] }
                f\del*{x_I}        
        }
        \le 
        \polygen{ \gen{s} }\del*{ \bx }
        \sum_{ I = \emptyset}^{ [k] }
            G\del*{ x_I }
        \le 
        \polygen{ \gen{s} }\del*{ \bx }
        \t 
        2^{k-1}
        \del*{
            G(x)
            +
            G(y)
        }.
    \end{equation*}
    Dividing the resulting inequality by 
    $2^{k-1}$
    and using \Cref{lem::folding_hypercubes_two_elements},
    we get that
    \begin{align*}
        \abs*{
            f(x) - f(y)
        }
    &\le 
        \polygen{ \gen{s} }\del*{ \bx }
        \del*{
            G(x)
            +
            G(y)
        }
    \\
    &\le 
        2^{k^2}
        \partitions{ k }
        \del*{
            \metric{x,y}
            +
            \metric{x,y}^s
        } 
        \del*{
            G(x)
            +
            G(y)
        }
    \\
    &\le 
        2\t 2^{k^2} \,
        \partitions{ k } 
        \metric{x,y}
        \del*{
            G(x)
            +
            G(y)
        }
    \\
    &\le 
        \metric{x, y}
         \del*{
            \abs{f(x)} 
            +
            2^{k^2+1}
            \partitions{k}
            G(x)
            + 
            \abs{f(y)}
            +
            2^{k^2+1}
            \partitions{k}
            G(y)
        },
    \end{align*}
    where the penultimate inequality follows from the fact that since
    $ \metric{x, y} \le 1$ and $s \ge 1$,
    we~have~that
    $ \metric{x, y}^s \le \metric{x, y}$.

    In consequence, since $x, y \in F$ are arbitrary, we found that
    \begin{equation*}
        \forall x, y \in F
    \qquad 
        \abs*{
            f(x) - f(y)
        }       
        \le 
        \metric{x, y}
         \del*{
            \abs{f(x)} 
            +
            2^{k^2+1}
            \partitions{k}
            G(x)
            + 
            \abs{f(y)}
            +
            2^{k^2+1}
            \partitions{k}
            G(y)
        }.  
    \end{equation*} 
    Therefore, since $F$ is of full measure,
    $
        \abs{f} 
        +
        2^{k^2+1}
        \partitions{k}
        G
    $
    is a Hajłasz gradient of $f$.
    Since $f, G \in \integrable{p}(X)$,
    it follows that 
    $ f \in \hajlasz{1,p}(X)$.
    Moreover, 
    \begin{align*}
        \norm{ f }_{
            \hajlasz{1,p}(X)        
        }
    &\le 
        \norm{ f }_{
            \integrable{p}(X)
        }
        +
        \smash{
        \norm*{
            \abs{f} 
            +
            2^{k^2+1} \,
            \partitions{k}
            G       
        }_{
            \integrable{p}(X)
        }
        }
    \\
    &\le 
        2\norm{ f }_{ \integrable{p}(X) }
        +
        2^{k^2+1} \, 
        \partitions{k}
        \norm{ G }_{ \integrable{p}(X) }
    \\
    &\le 
        2\norm{ f }_{ 
            \multipoint{s,p}\del*{X}
        }
        +
        2^{k^2+1} \,
        \partitions{k}
        \del*{
            \norm{ f }_{ 
                \multipoint{s,p}\del*{X}
            }
            +
            \eps
        }
        \le 
        2^{k^2+2} \,
        \partitions{k}
        \del*{
            \norm{ f }_{ 
                \multipoint{s,p}\del*{X}
            }
            +
            \eps
        }
        .
    \end{align*}
    Since $\eps > 0$ is arbitrary, 
    it follows that 
    $
        \norm{ f }_{
            \hajlasz{1,p}(X)        
        }
        \le 
        2^{k^2+2} \partitions{k}
        \norm{ f }_{ 
            \multipoint{s,p}\del*{X}
        }.   
    $
    Finally, because
    $ {f \in \multipoint{s,p}\del*{X} }$
    is arbitrary, 
    $
        \multipoint{s,p}\del*{X}
        \hookrightarrow
        \hajlasz{1,p}(X),
    $
    as claimed.
\end{proof}
%


\section*{Acknowledgements}

The author would like to thank Przemysław Górka for his help and useful discussions throughout the process of writing the paper.

\printbibliography

\end{document}